\documentclass[11pt]{amsart}
\input{glyphtounicode}
\usepackage{amsmath,amssymb,amsthm,mathtools,mathrsfs}
\usepackage{tikz-cd}
\usepackage[shortlabels]{enumitem}
\usepackage[margin=1.13in]{geometry}
\usepackage[protrusion=true,expansion=false]{microtype}
\usepackage{booktabs}
\usepackage{xcolor}
\usepackage{hyperref}
\usepackage[capitalize,nameinlink]{cleveref}
\hypersetup{colorlinks=true,citecolor=black,linkcolor=black,urlcolor=black,pdfauthor={Hamed Ebadi},pdftitle={The Riemann Hypothesis for L-Functions of Linearly Constrained Kloosterman Families: Frobenius Polynomials and Intersection Cohomology}}
\DeclareTextSymbolDefault{\textsection}{OMS}
\DeclareTextSymbolDefault{\textdagger}{OMS}
\allowdisplaybreaks[2]
\theoremstyle{plain}
\newtheorem{mainthm}{Theorem}

\newtheorem{theorem}{Theorem}[section]
\newtheorem{proposition}[theorem]{Proposition}
\newtheorem{lemma}[theorem]{Lemma}
\newtheorem{corollary}[theorem]{Corollary}
\theoremstyle{definition}
\newtheorem{definition}[theorem]{Definition}

\theoremstyle{remark}
\newtheorem{remark}[theorem]{Remark}

\numberwithin{equation}{section}

\crefname{mainthm}{Theorem}{Theorems}
\crefname{theorem}{Theorem}{Theorems}
\crefname{proposition}{Proposition}{Propositions}
\crefname{lemma}{Lemma}{Lemmas}
\crefname{corollary}{Corollary}{Corollaries}
\crefname{definition}{Definition}{Definitions}
\crefname{remark}{Remark}{Remarks}
\crefname{example}{Example}{Examples}
\crefname{question}{Question}{Questions}
\crefname{section}{Section}{Sections}
\crefname{subsection}{Section}{Sections}
\crefname{table}{Table}{Tables}
\crefname{equation}{}{}
\creflabelformat{equation}{#2(#1)#3}

\newcommand{\Ql}{\overline{\mathbb Q}_{\ell}}
\newcommand{\Fp}{\mathbb F_p}
\newcommand{\Fq}{\mathbb F_q}
\newcommand{\Fqb}{\overline{\mathbb F}_q}
\newcommand{\Gm}{\mathbb G_m}
\newcommand{\PP}{\mathbb P}
\newcommand{\cW}{\mathcal W}
\newcommand{\cF}{\mathcal F}
\newcommand{\cL}{\mathcal L}
\newcommand{\cM}{\mathcal M}
\newcommand{\KL}{\mathcal K\!\ell_2}
\newcommand{\Kl}{\mathrm{Kl}}
\DeclareMathOperator{\Tr}{Tr}
\DeclareMathOperator{\Frob}{Frob}
\DeclareMathOperator{\Swan}{Sw}
\DeclareMathOperator{\rk}{rk}

\DeclareMathOperator{\GL}{GL}
\DeclareMathOperator{\SL}{SL}

\DeclareMathOperator{\imop}{im}
\newcommand{\Par}{\mathrm{par}}
\newcommand{\Fpb}{\overline{\mathbb F}_p}
\newcommand{\cH}{\mathcal H}
\newcommand{\cA}{\mathcal A}
\newcommand{\cR}{\mathcal R}
\newcommand{\cG}{\mathcal G}
\newcommand{\cY}{\mathcal Y}
\newcommand{\cD}{\mathcal D}
\newcommand{\cE}{\mathcal E}
\newcommand{\cX}{\mathcal X}
\newcommand{\cK}{\mathcal K}

\newcommand{\cN}{\mathcal N}
\newcommand{\cC}{\mathcal C}
\newcommand{\cV}{\mathcal V}
\DeclareMathOperator{\FT}{FT}
\DeclareMathOperator{\pr}{pr}
\DeclareMathOperator{\gr}{gr}
\DeclareMathOperator{\End}{End}
\DeclareMathOperator{\Sym}{Sym}
\DeclareMathOperator{\Hom}{Hom}

\providecommand{\cT}{\mathcal{T}}
\providecommand{\cP}{\mathcal{P}}
\providecommand{\jex}{j_{!*}}
\providecommand{\sg}{\sigma}
\providecommand{\rhotil}{\widetilde{\rho}}
\providecommand{\wt}{\operatorname{wt}}
\providecommand{\rad}{\operatorname{rad}}
\providecommand{\STr}{\operatorname{STr}}
\providecommand{\Spec}{\operatorname{Spec}}
\providecommand{\GSp}{\operatorname{GSp}}
\providecommand{\supp}{\operatorname{supp}}
\providecommand{\id}{\operatorname{id}}
\providecommand{\im}{\operatorname{im}}
\providecommand{\coker}{\operatorname{coker}}
\providecommand{\Ext}{\operatorname{Ext}}
\providecommand{\RHom}{\operatorname{RHom}}
\providecommand{\cF}{\mathcal{F}}

\title[Constrained Kloosterman $L$-functions]{The Riemann Hypothesis for L-Functions of Linearly Constrained Kloosterman Families: Frobenius Polynomials and Intersection Cohomology}
\author{Hamed Ebadi}
\thanks{Address: K. N. Toosi University of Technology, Tehran, Iran. Correspondence: h.ebadi1@email.kntu.ac.ir.}
\date{23 November 2026}
\subjclass[2020]{Primary 11L05, 11T23, 14F20; Secondary 11M38, 14G10, 14F43}
\keywords{Kloosterman sums, Riemann hypothesis over finite fields, $\ell$-adic cohomology, weights, Frobenius, intermediate extension, intersection cohomology, functional equations}
\begin{document}
\begin{abstract}
We study Frobenius polynomials attached to linearly constrained Kloosterman families, taking the rank-$r(n)=2^n-\binom{n}{\lfloor n/2\rfloor}$, weight-$(n+1)$ sheaf of the author's companion papers as an explicit input. The present contribution is the passage from that sheaf to local characteristic polynomials and the intermediate-extension intersection cohomology on $(\mathbb P^1)^n$. We prove the local reciprocal-root bounds and the degreewise purity of the completed factors; their weights are $n+1+i$ in cohomological degree $i$. Using the imported arithmetic relative-image identification, we construct a reflection pairing and derive the completed functional equation with positive root number. We also give a canonical Frobenius-equivariant splitting of the fibrewise compact-support/parabolic sequence by separation of weights and examine the comparison on open coordinate-zero divisors. The statements do not assume global arithmetic Frobenius semisimplicity or a higher-codimension wild-boundary exchange. An explicit source-and-dependency ledger distinguishes results imported from the companion papers from deductions established here.
\end{abstract}
\maketitle

\section{Introduction and main results}\label{sec:intro}
\subsection{The exponential-sum family}
Fix an odd prime $p$, a prime $\ell\ne p$, a nontrivial additive character
$\psi\colon\mathbb F_p\to\overline{\mathbb Q}_\ell^\times$, and a nonzero
$\lambda\in\mathbb F_p$. For $q=p^a$ let
$\psi_q=\psi\circ\operatorname{Tr}_{\mathbb F_q/\mathbb F_p}$. For
$h=(h_1,\ldots,h_n)\in(\mathbb F_q^\times)^n$ and $n\ge2$, consider
\begin{equation}\label{eq:sum}
\Sigma_{n,q}(h)=\sum_{\substack{u\in(\mathbb F_q^\times)^n\\
                                  h_1u_1+\cdots+h_nu_n=-\lambda}}
          \psi_q(u_1^{-1}+\cdots+u_n^{-1}).
\end{equation}
The linear constraint couples the reciprocal phases, so the sum does
not factor into one-variable exponential sums. The associated sheaf makes it possible to study Frobenius
characteristic polynomials on individual fibres and on the intersection
cohomology of a proper compactification.

Here the term \emph{Riemann hypothesis} refers exclusively to the
Weil-type bounds for Frobenius eigenvalues over finite fields. It does
not concern the classical Riemann zeta function. The local statement
involves a single lisse-sheaf stalk; the compactified statement concerns
intersection cohomology and may involve several distinct weights.

\subsection{Contribution and position in the series}\label{subsec:delta}
This is the second, $L$-function and intersection-cohomology, part of the
program initiated in \cite{EbaCoh,EbaMom}.  Its novelty is \emph{not} a
second construction of the top-weight Kloosterman sheaf or a new proof of its
monodromy theorem.  The earlier papers supply that lisse sheaf, its arithmetic
normalization and its fibrewise parabolic realization.  The present paper
makes the following further deductions in the stated characteristic range.
First, purity gives the exact reciprocal-root circles for the individual
Frobenius polynomials, with the lower-weight correction kept explicit.
Second, applying middle extension and proper intersection cohomology produces
one pure Frobenius module of weight $n+1+i$ in each degree $i$; hence the
completed factors have \emph{degree-dependent}, rather than one common,
critical circle.  Third, the arithmetic reflection pairing on the already
constructed parabolic image gives a functional equation and determines its
positive root number, subject to the precisely sourced relative-image
statement (H7).  Finally, the tame boundary invariants yield an exact
arithmetic compact-support sequence, and weight separation provides its
canonical Frobenius-equivariant splitting on each closed fibre.

We use the fixed versions \texttt{arXiv:2610.03081v3} and
\texttt{arXiv:2610.04099v2}; theorem numbering refers to those versions,
not to a subsequently revised preprint.  The former provides the
load-bearing sheaf-theoretic statements, while the latter is used for
monodromy context only.  The resulting proof graph, including the
independence of H7 from the root-number argument, is given explicitly in
\cref{subsec:ledger}.  In particular, no statement concerning all mixed
or exceptional boundary strata, nor arithmetic Frobenius semisimplicity
of every $IH_n^i$, is used as a premise of the weight-theoretic results.

\subsection{Main geometric theorem}
Write $S_n=\mathbb G_{m,\mathbb F_p}^n$ and put
\[
  r(n)=2^n-\binom{n}{\lfloor n/2\rfloor},\qquad
  m_{n,k}=\binom nk-\binom n{k-1},\qquad
  E_n(X)=\sum_{1\le k\le\lfloor n/2\rfloor}m_{n,k}X^k.
\]
The starting point, imported from the author's first paper
\cite[Theorem~5.3(i)--(vi)]{EbaCoh}, is the relative construction and top-weight quotient of the
relative Fourier construction.
\begin{mainthm}[Relative construction and top-weight sheaf: imported input]\label{thmA}
For every $n\ge2$ the auxiliary relative cohomology sheaf
$\mathcal H_n=R^1\pi_{t!}\mathcal M_n$ is lisse of rank $2^n$ over $S_n$.
There is an exact sequence of lisse arithmetic sheaves
\[
0\longrightarrow\Ql\longrightarrow\mathcal H_n
\longrightarrow\mathcal G_n\longrightarrow0,
\]
where $\mathcal G_n$ carries a filtration by lisse subsheaves whose
nonzero graded pieces are
\[
\operatorname{gr}^{W}_{2k}\mathcal G_n
\simeq\Ql(-k)^{\oplus m_{n,k}}
\quad(1\le k\le\lfloor n/2\rfloor),\qquad
\mathcal W_n=\operatorname{gr}^{W}_{n+1}\mathcal G_n.
\]
The sheaf $\mathcal W_n$ is lisse of rank $r(n)$ and pointwise pure of
weight $n+1$. For every $q=p^a$ and $h\in S_n(\mathbb F_q)$,
\[
\operatorname{Tr}(\operatorname{Frob}_h\mid\mathcal W_{n,\bar h})
=(-1)^{n-1}q\,\Sigma_{n,q}(h)-E_n(q).
\]
These assertions are imported from
\cite[Theorem~5.3(i)--(vi)]{EbaCoh} and are not re-proved here;
part (vi) supplies the arithmetic Tate identifications of the
lower-weight graded pieces, whereas parts (iv) and (v) supply,
respectively, the pure quotient and its trace identity.
\end{mainthm}
\begin{remark}[Geometric inputs]\label{rmk:input-status}
The construction and purity of $\mathcal W_n$ and the open-torus
relative-image identification (H7) are imported from the fixed first-paper
version \cite[Propositions~4.2--4.3 and Theorem~5.3]{EbaCoh}.  They are
not conclusions of the completed functional equation in this manuscript. The local and intersection-cohomological
statements are deductions from these inputs and the standard purity
and duality theorems. Any further boundary comparisons are stated with
their explicit hypotheses.
\end{remark}

The adjective \emph{pure} always refers to a fixed embedding
$\iota\colon\overline{\mathbb Q}_\ell\hookrightarrow\mathbb C$. Thus
pointwise purity of weight $n+1$ means that for every geometric
Frobenius eigenvalue $\alpha$ at a point of degree $a$, one has
$|\iota(\alpha)|=p^{a(n+1)/2}$. No semisimplicity of the individual
Frobenius operator is needed for the characteristic-polynomial conclusions.

\subsection{Riemann hypothesis for the local factors}
For $h\in S_n(\mathbb F_q)$ let
\[
 P_{n,h}(T)=\det(1-T\operatorname{Frob}_h\mid\mathcal W_{n,\bar h}),
 \qquad L_{n,h}(T)=P_{n,h}(T)^{-1}.
\]
\begin{mainthm}[Local Riemann hypothesis and reflection]\label{thmB}
Every zero of $P_{n,h}(T)$ has complex absolute value
$q^{-(n+1)/2}$; equivalently, the zeros of
$P_{n,h}(q^{-s})$ have real part $(n+1)/2$. The polynomial is determined
by extension-field sums and obeys a normalized self-inversive identity,
with root-of-unity factor
$\psi_q(-2^{n-1}\lambda^{-1}\sum_i h_i)$, relating the parameters
$h$ and $-h$.
\end{mainthm}
The local factors are defined as reciprocals of the characteristic
polynomials; thus the corresponding eigenvalue circles describe poles
of the local rational factors.

\subsection{The intersection-cohomological completion}
Let $X_n=(\mathbb P^1_{\mathbb F_p})^n$, let $j\colon S_n\hookrightarrow X_n$
be the standard open immersion, and put
\[
 K_n=j_{!*}(\mathcal W_n[n]),\qquad
 \mathcal I_n=K_n[-n],\qquad
 IH^i_n=H^i(X_{n,\overline{\mathbb F}_p},\mathcal I_n).
\]
Denote $P_{n,i}(T)=\det(1-T\operatorname{Frob}_p\mid IH_n^i)$ and
$\widehat L_n(T)=\prod_i P_{n,i}(T)^{(-1)^{i+1}}$.
\begin{mainthm}[Riemann hypothesis for the intersection completion]
\label{thmC}
The Frobenius module $IH^i_n$ is pure of weight $n+1+i$.
Consequently, all zeros of $P_{n,i}(p^{-s})$ lie on
\[
\operatorname{Re}(s)=\frac{n+1+i}{2}.
\]
The completion $\widehat L_n$ admits a functional equation under
$T\longmapsto(p^{2n+1}T)^{-1}$ with an a priori sign
$\omega_n\in\{\pm1\}$. Its existence is established independently
of the subsequent evaluation of $\omega_n$; the equality $\omega_n=+1$
is derived in \cref{thm:arith-reflection-closed} from the imported
arithmetic relative-image identification \textup{(H7)} and the
reflection pairing, not assumed in the construction of the functional
equation. The Riemann-hypothesis deduction uses purity of intermediate extension
\cite[Corollaire~5.4.3]{BBD} (Gabber; see also
\cite[Corollaire~5.3.2]{BBD} for the weight-truncation formalism)
and proper-pushforward
purity \cite[no.~5.1.14]{BBD}, and does not require the separate
higher-codimension boundary formulas. In particular, the sign computation does not enter the purity
argument or any of its imported hypotheses.
\end{mainthm}
For the completed function there is generally no single critical
circle. Different cohomological degrees have different weights, and the
appropriate conclusion is a family of circles, one for each nonzero
intersection-cohomology group. The degree-by-degree statement must not
be confused with a claim that all roots lie on the same line.

\noindent The proof of the sign is given in
\cref{sec:arith-reflection-closure}.  The reflection pairing is constructed
on the arithmetic image already identified by H7, extends through
intermediate extension, and induces the symplectic similitude on the
middle intersection-cohomology group.  No boundary-support table or
arithmetic Jordan decomposition is needed for this argument.

\subsection{The parabolic comparison}
For the auxiliary family 
$\mathscr M_{n,h}(t)=\mathcal L_\psi(\lambda t)\otimes
\bigotimes_i[h_it]^*\mathrm{Kl}_2$ and nonzero parameters, the natural
forget-supports map fits into an exact sequence of arithmetic
Frobenius modules
\[
0\longrightarrow B_{n,h}\longrightarrow
H_c^1(\mathbb G_m,\mathscr M_{n,h})\longrightarrow
H_{\mathrm{par}}^1(\mathbb G_m,\mathscr M_{n,h})
\longrightarrow0.
\]
The exact sequence and the tame invariant dimensions are recorded in
\cite[Propositions~4.2--4.3]{EbaCoh}; see also
Theorem~\ref{thm:parabolic-all-rank-20261123} for a direct derivation.
The further Frobenius-equivariant splitting by distinct weights is proved
in Proposition~\ref{prop:part2-weight-split}.
The compact-support comparison under specialization of a
parameter is established on open coordinate-zero divisors by
Theorem~\ref{thm:audited-generic-zero-20261123}. Its two local lemmas
are proved by tame normal-crossing d\'evissage at the finite
$t$-endpoint and by functoriality of the compact-support
localization triangle; no tameness assumption is made at
$t=\infty$. Higher-codimension mixed-face comparisons remain outside
this theorem.

\subsection{Relation to earlier work}
The author's earlier article \emph{Cohomology of Linearly Constrained Kloosterman Families: Newton Polytopes, Weights, and Boundary Monodromy} \cite{EbaCoh} develops the geometric foundations of the present family. Its analysis distinguishes the Newton polytope governing toric critical points from the polytope at infinity governing compactly supported cohomology; facewise non-degeneracy, cohomological concentration, boundary contributions, and Swan conductors then lead to the rank and weight filtration of the associated sheaf. In particular, it identifies the pure top-weight constituent $\mathcal W_n$ of rank $r(n)=2^n-\binom{n}{\lfloor n/2\rfloor}$ and analyzes tame unipotent monodromy on coordinate restrictions through local Fourier transforms and middle convolution. These are the family-specific foundations needed here: the construction, lissity, trace identity, and weight $n+1$ of $\mathcal W_n$ supply the hypotheses to which the general purity theorems can be applied.

The companion article \emph{Moments and Large Geometric Monodromy for Linearly Constrained Kloosterman Families} \cite{EbaMom} addresses a different but closely related question: how large the connected geometric monodromy group of this top-weight sheaf is. It derives a fourth-moment estimate by decoupling the auxiliary variables and analyzing plane sections of the Cayley cubic, and it obtains a nontrivial unipotent boundary element using local Fourier-transform and convolution methods. Together with a group-theoretic argument, these inputs yield $G_{\mathrm{geom}}^0(\mathcal W_n)=\mathrm{SL}_{r(n)}$ in the paper's range $n\ge4$ and $p>\max(2,n)$. The monodromy theorem concerns the variation of the family, rather than the purity of any individual Frobenius module. It is \emph{not}, however, an input to either Riemann-hypothesis proof below: those proofs use purity and proper cohomology, not a classification of monodromy groups.

The longer manuscript \emph{The Geometric Monodromy Group and Short-Sum Estimates for Linearly Constrained Families of Exponential Sums}, archived on Zenodo \cite{EbaZen}, situates these investigations in a broader program connecting the geometric construction and monodromy with moments, correlations, and short-sum estimates. Its role in the present article is principally contextual and bibliographic; no general short-sum estimate is assumed in the proofs of the critical-line statements.

The specific contribution of this article is the passage from the imported pure family sheaf to its local Frobenius polynomials, the degreewise intersection cohomology of its intermediate extension on $(\mathbb P^1)^n$, and the arithmetic reflection pairing and functional equation. The local and compactified Riemann-hypothesis conclusions are applications of Deligne's and Gabber's purity theorems to this particular family, not claims to establish a new general purity theorem. Throughout, a geometric assertion imported from the preceding research is identified separately from the formal consequence proved here.

\subsection{Imported results, new deductions, and noncircular dependence}
\label{subsec:ledger}
The distinction between the two companion sources and the present article
is recorded in Table~\ref{tab:imported-proved}.  Version-specific citations
are intentional: the theorem numbers belong to
\texttt{arXiv:2610.03081v3} (18 October 2026) and
\texttt{arXiv:2610.04099v2} (18 October 2026).
\begin{table}[htbp]
\centering
\caption{Imported statements versus deductions in this article. All numbered
preprint references refer to the pinned versions in the bibliography.}
\label{tab:imported-proved}
\begingroup\small\renewcommand{\arraystretch}{1.20}
\setlength{\tabcolsep}{4pt}
\begin{tabular}{@{}p{0.14\linewidth}p{0.425\linewidth}p{0.375\linewidth}@{}}
\toprule
Status & Statement and logical role & Exact source or proof \\
\midrule
Imported I1 & Relative auxiliary sheaf of rank $2^n$; lissity and weight filtration of $\mathcal W_n$ & \cite[Theorem~5.3(i)--(iv),(vi)]{EbaCoh}; compare \cite[Theorem~2.10(i)--(ii)]{EbaMom} \\
Imported I2 & Fibrewise parabolic exact sequence and tame Tate boundary invariants & \cite[Propositions~4.2--4.3]{EbaCoh} \\
Imported I3 & Trace formula and pointwise purity of $\mathcal W_n$ & \cite[Theorem~5.3(iv)--(v)]{EbaCoh} \\
Imported I4 & Determinant and coordinate-line monodromy; normalization for reflection & \cite[Theorems~6.8, 6.12, 8.11, 8.17]{EbaCoh} \\
Imported I5 (H7) & Arithmetic relative-image identification on the open parameter torus & \cite[Proposition~4.2; Theorem~5.3]{EbaCoh}; \cite[Theorem~2.10(ii)]{EbaMom} \\
Context only & Connected geometric monodromy (not used in the two purity arguments) & \cite[Theorem~8.8]{EbaCoh}; \cite[Theorem~4.19]{EbaMom} \\
\midrule
Proved P1 & Local reciprocal-root circles and arithmetic reflection identity & Theorem~\ref{thmB}, Proposition~\ref{prop:reflection}; I3--I4 \\
Proved P2 & Pure $IH_n^i$ of weight $n+1+i$ and degreewise reciprocal-root circles & Theorem~\ref{thmC}; I1--I3; \cite[Corollaire~5.4.3; no.~5.1.14]{BBD} (Gabber and proper purity) \\
Proved P3 & Arithmetic reflection, completed functional equation and positive root number & Theorem~\ref{thm:arith-reflection-closed}; I5 and proper duality \\
Proved P4 & Canonical Frobenius splitting on each closed fibre by weight separation & Proposition~\ref{prop:part2-weight-split}; I2--I3 \\
Proved P5 (codim.~1) & Compact-support nearby-cycle comparison on open coordinate-zero divisors; no mixed-face conclusion & Theorem~\ref{thm:audited-generic-zero-20261123}; Lemmas~\ref{lem:audited-local-exchange}--\ref{lem:audited-boundary-detection} \\
\bottomrule
\end{tabular}\endgroup
\end{table}

\paragraph{Acyclic proof graph.}
The only family-specific input needed to construct $K_n$ is I1--I3;
BBD--Gabber purity then gives P1--P2.  The sign calculation has the
separate chain
\[
\text{I5 (H7)}+\text{rank-two Kloosterman pairing}
\ \Longrightarrow\ \text{arithmetic reflection-skew pairing}
\ \Longrightarrow\ \omega_n=+1.
\]
The first arrow is proved in Section~\ref{sec:arith-reflection-closure}
without assuming either the completed functional equation or the value of
$\omega_n$.  Conversely, neither H7 nor its cited source is proved using
P3 or any sign assertion.  The relative reconstruction in
Section~\ref{app:H7} is explanatory and is \emph{not} used as a new,
independent premise for a theorem whose proof already requires H7.
This separates the import from the deduction and prevents a circular
appeal to the desired root number.

\paragraph{Scope of the new results.}
The statements on the open parameter torus and the generic
single-coordinate zero boundary do not imply a comparison on every mixed
or exceptional face.  No universal arithmetic Frobenius-semisimplicity
statement for all $IH_n^i$ is claimed.  These limitations concern finer
arithmetic structure and do not enter the pure-weight conclusions P1--P2.

\subsection{Organization and completeness of the manuscript}
The local Frobenius polynomials and reflection relations are established in
\cref{sec:local,sec:lfe}; coordinate-line and curve computations follow in
\cref{sec:lines,sec:curves,sec:explicit,sec:parabolic}.
The geometric construction and its input ledger, the global intersection
completion, the relative H7 comparison, the arithmetic reflection
construction, the five-variable calculation, and the tame boundary
comparison are included in full in
\cref{app:inputs,sec:global,sec:h7-comparison,sec:reflection-compat,sec:five-variable,sec:arith-reflection-closure,sec:tame-boundary}.
The auxiliary dilated-Kloosterman and ray calculations are collected in
\cref{sec:aux-rays}. No numbered proof section is omitted or
relegated to an unprovided arXiv-only supplement.

\section{Notation and input}\label{sec:input}

\subsection{Notation}
In this paper $w=n+1$ is a sheaf weight; $b_{n,i}$ denotes an
intersection Betti number, whereas $b_{n,k}$ in the coordinate-line
formula is a different combinatorial multiplicity (always
written with the subscript $k$ in that context).  The letters
$K_n$ and $\cW_n$ denote respectively a perverse intermediate-extension
complex and a lisse sheaf, never a monodromy group.  We reserve
$P_{n,h}$ for a local characteristic polynomial, $P_{n,i}$ for
intersection cohomology, and $L_{n,h}=P_{n,h}^{-1}$ for a local
Euler factor; they must not be interchanged when discussing
zeros and poles.  Superscript $\mathrm{ss}$ means semisimplification
of the specified operator or representation, as indicated.

For readability we distinguish throughout between an eigenvalue of a geometric Frobenius operator and a zero of its characteristic polynomial.  If $A$ has eigenvalues $\alpha_1,\ldots,\alpha_r$ and $P(T)=\det(1-TA)$, then the zeros of $P$ are the reciprocals $\alpha_j^{-1}$.  Thus an assertion that the eigenvalues have weight $w$ means $|\iota(\alpha_j)|=q^{w/2}$, whereas the zeros of $P$ lie on $|T|=q^{-w/2}$.  These are different statements about the same operator; we use both only after stating which quantity is being considered.  We also reserve \emph{local factor} for $P(T)^{-1}$, which has poles and not zeros at these reciprocal eigenvalues.  A Tate twist $(-k)$ multiplies Frobenius eigenvalues by $q^k$; consequently it increases weight by $2k$, a convention that explains the shifts encountered in the geometric construction.

Throughout, $p$ is an odd prime, $\ell\ne p$ a prime, $\psi$, $\lambda$ and
$\iota$ are as in the introduction, and $q=p^a$ with $a\ge1$. Frobenius
always means geometric Frobenius; $\Frob_q$ acts on $\Ql(-1)$ by $q$. For a
closed point $x$ of a scheme of finite type over $\Fp$ we write $\Frob_x$
for its geometric Frobenius. We put $S_n=\Gm^n$ over $\Fp$.

For $a\ge1$ and $0\le k\le a/2$ let $m_{a,k}=\binom ak-\binom a{k-1}$, with
$m_{a,0}=1$; these are the multiplicities in the Clebsch--Gordan
decomposition
$\mathrm{Std}^{\otimes a}=\bigoplus_{k}m_{a,k}\,\mathrm{Sym}^{a-2k}$ of the
standard representation of $\SL_2$. Set
\[
E_n(x)=\sum_{k=1}^{\lfloor n/2\rfloor}m_{n,k}x^k\in\mathbb Z[x],\qquad
r(n)=2^n-\binom n{\lfloor n/2\rfloor}\ (n\ge2),\qquad r(1)=1,
\]
so that $r(n)=2,5,10,22,44$ for $n=2,\dots,6$. For $n\ge1$ and
$h\in S_n(\Fq)$ we put $D_n(h)=-2^{n-1}\lambda^{-1}\sum_ih_i$ and
$\delta_h=\psi_q(D_n(h))$.

For $\xi\in\Fq^\times$ let $\Kl_q(\xi)=\sum_{x\in\Fq^\times}\psi_q(x+\xi x^{-1})$.
The Kloosterman sheaf $\KL$ on $\Gm$ is lisse of rank $2$ and pure of weight
$1$, with trace $-\Kl_q(\xi)$ at $\xi\in\Gm(\Fq)$ and determinant $\Ql(-1)$;
it is tamely ramified at $0$, where inertia acts through a single unipotent
Jordan block, and totally wild at $\infty$ with Swan conductor $1$
\cite[Sommes trig., Thm.~7.4 and 7.8]{Del77}, \cite{Kat88}. For a function
$f$ we write $\cL_\psi(f)$ for the Artin--Schreier sheaf with trace
$\psi_q(f)$.

\subsection{Notation index and conventions}\label{subsec:notation-index}
The following list distinguishes similarly named objects used in the
relative construction, the coordinate-line calculation, and the
intersection completion. Here $\mathcal H^i$ denotes an ordinary
cohomology \emph{sheaf}, whereas $H^i$ or $\mathbb H^i$ denotes global
cohomology or hypercohomology; the two notations must not be interchanged.
\begin{center}\small
\begin{tabular}{@{}p{0.21\linewidth}p{0.51\linewidth}p{0.19\linewidth}@{}}
\toprule
Symbol & Meaning & First use / type\\
\midrule
$S_n$, $X_n$ & Parameter torus $(\mathbb G_m)^n$ and compactification $(\mathbb P^1)^n$ & \cref{sec:intro,sec:input}\\
$\mathcal M_n$ & Auxiliary rank-$2^n$ Kloosterman--Artin--Schreier system & \cref{app:inputs}\\
$\mathcal H_n$ & $R^1\pi_{t!}\mathcal M_n$ & \cref{app:inputs}; mixed lisse\\
$\mathcal L_n$ & $\ker\alpha_n=W_n\mathcal H_n$ & \cref{app:inputs}; low weight\\
$\mathcal H'_n$ & Cokernel of the distinguished constant line in $\mathcal H_n$ & \cref{app:inputs}\\
$\mathcal W_n$ & $\operatorname{im}\alpha_n=\mathcal H_n/\mathcal L_n$ & \cref{sec:input}; pure $n+1$\\
$K_n$, $\mathcal I_n$ & $j_{!*}(\mathcal W_n[n])$, respectively $K_n[-n]$ & pure perverse / unshifted\\
$P_{n,h}$, $P_{n,i}$ & Local and cohomological Frobenius polynomials & \cref{sec:local,sec:global}\\
$P_c^{\mathrm{par}}$ & Parabolic characteristic polynomial on a coordinate line & \cref{sec:lines,sec:parabolic}\\
$\mathcal V_n$, $\mathcal V'_n$ & Coordinate-line parabolic sheaf and its extension & \cref{sec:global}\\
$B'_n$ & Residual perverse support of the coordinate projection & \cref{sec:global}\\
$P_5$, $F^v$ & $Rf_*K_5$; $\mathcal H^v(K_5[-5])$ & \cref{sec:five-variable}\\
$J_0,J_\infty,Q_0,Q_\infty$ & Nearby invariant and quotient systems of the line family & \cref{sec:five-variable}\\
$\mathcal C_{i,a}$ & Divisor coefficient for boundary value $a$ & \cref{sec:global}\\
$\omega_n$ & Root number of the completed functional equation & \cref{sec:arith-reflection-closure}\\
\bottomrule
\end{tabular}
\end{center}
A statement that a sheaf is geometrically constant and has characteristic
polynomial $(1-q^aT)^r$ at every closed point identifies its
semisimplification; it does \emph{not} by itself rule out a unipotent
Jordan block for arithmetic Frobenius. Stronger Tate assertions below
are justified by explicit arithmetic kernels, images, and quotients.

\subsection{Previously established geometric theorems}\label{subsec:input}
The open-torus geometric construction and statements (H1)--(H4) and
(H7) are imported from the fixed version
\cite[Propositions~4.2--4.3; Theorems~5.3, 6.8, 6.12, 8.11, 8.17]{EbaCoh};
compare \cite[Theorem~2.10]{EbaMom} for an independent statement of the
relative-image filtration.  The notation is recalled below for $n\ge2$
and odd $p$.  Calculations reproduced later serve to fix conventions and
are not counted as new contributions.  In particular, the prior sources
for H7 contain no theorem about the global root number of the
intersection-cohomology completion.  The global monodromy statements (H5) and (H6)
are logically separate: the companion paper \cite{EbaMom} establishes
its connected-group theorem in the range $n\ge4$ and
$p>\max(2,n)$.  In particular, no conclusion about the full
arithmetic group or equidistribution in a prescribed disconnected
compact group in every odd characteristic is inferred here from
that connected-group theorem. Neither (H5) nor (H6) is used in the
two Riemann-hypothesis proofs.  We retain the labels only to display
the dependence of each argument.

\begin{enumerate}[label=(H\arabic*),leftmargin=3.2em]
\item\label{H1} (\cref{thm:B-family}; \cite[Theorem~5.3(iv)]{EbaCoh}) There is a lisse $\Ql$-sheaf
$\cW_n$ on $S_n$ of rank $r(n)$, pointwise $\iota$-pure of weight $w=n+1$.
\item\label{H2} (\cref{thm:B-family}(iv); \cite[Theorem~5.3(v)]{EbaCoh}) For every $q$ and $h\in S_n(\Fq)$,
\begin{equation}\label{eq:traceinput}
\Tr(\Frob_h\mid\cW_{n,\bar h})=(-1)^{n-1}q\,\Sigma_{n,q}(h)-E_n(q).
\end{equation}
\item\label{H3} (\cref{thm:B-det}; \cite[Thms.~8.11 and 8.17]{EbaCoh}) For every $q$ and
$h\in S_n(\Fq)$,
\begin{equation}\label{eq:detinput}
\det(\Frob_h\mid\cW_{n,\bar h})=q^{wr(n)/2}\,\delta_h .
\end{equation}
By the Chebotarev density theorem this is equivalent to
$\det\cW_n\cong\cL_\psi(D_n)\otimes\Ql(-wr(n)/2)$.
\item\label{H4} (\cref{thm:B-zero,thm:B-infinity}; \cite[Thms.~C, 6.8 and 6.12]{EbaCoh}) Let $\kappa_j$
be the number of Jordan blocks of size $j$ in the tensor product of $n-1$
unipotent Jordan blocks of size $2$; thus $\kappa_{n-2i}=m_{n-1,i}$ and
$\kappa_j=0$ for other $j$. For every $c$, the restriction $\cF_c$ of
$\cW_n$ to the coordinate line $t\mapsto(t,c)$ satisfies:
\begin{enumerate}[label=(\alph*)]
\item at $t=0$ it is tamely ramified, inertia acts unipotently, and the
Jordan decomposition consists of $r(n-1)$ blocks of size $2$ and $\kappa_1$
blocks of size $1$;
\item at $t=\infty$ it decomposes as $V^{\mathrm t}\oplus V^{\mathrm w}$,
where $V^{\mathrm w}$ has dimension $2^{n-1}$ and all breaks equal to $1$,
and $V^{\mathrm t}$ is tamely ramified with unipotent inertia action having
$\kappa_j$ Jordan blocks of size $j-2$ for every $j\ge3$.
\end{enumerate}
\item\label{H5} (\emph{Context only; not an input to the present
proofs.}) In the stated range $n\ge4$, $p>\max(2,n)$, the connected
geometric monodromy conclusion recorded in \cite{EbaMom} is
$G_{\mathrm{geom}}^\circ(\cW_n)=\SL_{r(n)}$ after normalizing
weights.  The stronger identification of the \emph{full} geometric
or arithmetic monodromy group must be justified separately and is not
asserted here for every $n,p$.
\item\label{H6} (\emph{Context only.}) The equidistribution conclusions
in \cite{EbaCoh,EbaMom} require the precise hypotheses on arithmetic
monodromy and the component group.  A connected geometric-monodromy
calculation alone does not identify the full compact equidistribution
group; no such identification is invoked in this paper.
\item\label{H7} (\cite[Proposition~4.2 and Theorem~5.3(iv)--(vi)]{EbaCoh};
compare \cite[Theorem~2.10(i)--(ii)]{EbaMom}) Let
$\pi\colon\Gm\times S_n\to S_n$ be the projection and
$\cM=\cL_\psi(\lambda t)\otimes\bigotimes_{i=1}^n[(t,h)\mapsto h_it]^*\KL$.
The sheaf $R^1\pi_!\cM$ is lisse, and there is a surjection of lisse sheaves
$R^1\pi_!\cM\to\cW_n$ on $S_n$, defined over $\Fp$, with lisse kernel, which
at every geometric point $h$ above a closed point is the forget-supports map
$H^1_c(\Gm,\cM_h)\to H^1(\Gm,\cM_h)$ followed by an identification of its
image $H^1(\PP^1,j_*\cM_h)$ with $\cW_{n,h}$.
\end{enumerate}

\paragraph{Precisely separated source clauses (H1, H2, H7).}
In the pinned \texttt{arXiv:2610.03081v3},
Theorem~5.3(iv) establishes the lisse weight filtration and the
pure quotient (the input labelled H1); Theorem~5.3(v) gives the trace
identity (H2). Theorem~5.3(vi) asserts the \emph{arithmetic Tate
trivialization of the lower-weight graded pieces}; it neither states
nor assumes the global sign $\omega_n$. Here H7 denotes the
\emph{open-torus arithmetic relative-image identification} supported by
Propositions~4.2--4.3 and the relative construction in
Theorem~5.3(ii)--(iv), with the Tate normalization of part~(vi)
used where appropriate. Thus ``H7(vi)'' is not a separate root-number
hypothesis: the reference to (vi) supplies a low-weight input, not
an assertion of reflection-skew duality. The latter is proved in this
article from the already fixed arithmetic image; no step deriving
H1, H2 or H7 is allowed to invoke the conclusion $\omega_n=+1$.

Independently of \ref{H5}, the sheaf $\cW_n$ is geometrically irreducible
(\cref{cor:B-irred}); this follows from the second moment of the sums, whereas
\ref{H5} rests on the fourth moment.

\begin{remark}\label{rem:n1}
For $n=1$ we make the \emph{conventional extension of notation} $\cW_1=\cL_\psi(D_1)\otimes\Ql(-1)$ on $S_1=\Gm$; this is a definition for the degenerate rank-one case, not a claim that the imported theorem for $n\ge2$ contains this case. Here
$D_1(h)=-h/\lambda$. Since $E_1=0$ and
$\Sigma_{1,q}(h):=\psi_q(-h/\lambda)$ is the value of \eqref{eq:sum} for
$n=1$, hypotheses \ref{H1}--\ref{H3} then hold for $n=1$ with $r(1)=1$ and
$w=2$.
\end{remark}

\begin{remark}
The sheaf $\cW_n$ is the Tate twist by $(-1)$ of the top-weight part of
the fibre cohomology of the family. That part has weight $n-1$
\emph{before} the Tate twist $(-1)$; the twist raises the weight by $2$,
so $\cW_n$ has weight $n+1$. Confusing these two weights shifts each
corresponding critical line by one in the variable $s$.
\end{remark}

\subsection{Elementary identities}

The next elementary lemma is included before any sheaf-theoretic argument because it is the bridge between the original counting problem and the Fourier model used later.  Its four parts have separate jobs.  The first identifies the complex-conjugate parameter, the second introduces the Kloosterman factors, the third lowers the number of constrained variables when one parameter is averaged, and the fourth records the permutation symmetry of the family.  None of these identities requires a purity theorem or a statement about monodromy.

\begin{lemma}\label{lem:identities}
Let $q$ be a power of $p$, $h\in(\Fq^\times)^n$ and
$c\in(\Fq^\times)^{n-1}$.
\begin{enumerate}[label=\textup{(\alph*)}]
\item $\Sigma_{n,q}(-h)$ is the image of $\Sigma_{n,q}(h)$ under the
automorphism $\zeta\mapsto\zeta^{-1}$ of $\mathbb Q(\mu_p)$; under $\iota$,
$\Sigma_{n,q}(-h)=\overline{\Sigma_{n,q}(h)}$.
\item $\displaystyle\Sigma_{n,q}(h)=\frac{1}{q}\Bigl((-1)^n+\sum_{s\in\Fq^\times}
\psi_q(\lambda s)\prod_{i=1}^n\Kl_q(sh_i)\Bigr)$.
\item $\displaystyle\sum_{t\in\Fq^\times}\Sigma_{n,q}(t,c)=(-1)^n+\Sigma_{n-1,q}(c)$,
where $\Sigma_{1,q}(c_2)=\psi_q(-c_2/\lambda)$.
\item $\Sigma_{n,q}(h_{\pi(1)},\dots,h_{\pi(n)})=\Sigma_{n,q}(h)$ for every
permutation $\pi$.
\end{enumerate}
\end{lemma}

\begin{proof}
(a) The substitution $u_i\mapsto-u_i$ for every $i$ changes the left side of the constraint with parameter $-h$ into the constraint with parameter $h$.  At the same time, it changes the phase $\sum_i u_i^{-1}$ into its negative.  Since complex conjugation sends $\psi_q(x)$ to $\psi_q(-x)$, this proves the conjugation identity without using any property of $\cW_n$.

(b) The elementary character-orthogonality relation is
\[
\frac1q\sum_{s\in\Fq}\psi_q(sy)=
\begin{cases}1,&y=0,\\0,&y\ne0.\end{cases}
\]
Apply it with $y=\lambda+\sum_i h_i u_i$ before summing over the $u_i$.  For $s\ne0$, the $n$ variables separate, and the substitution $x=s h_i u_i$ converts the $i$-th one-variable sum into $\Kl_q(s h_i)$.  For $s=0$ the only remaining factor is $\sum_{u_i\ne0}\psi_q(u_i^{-1})=-1$, because inversion permutes $\Fq^\times$ and the sum of a nontrivial additive character over $\Fq$ is zero.  Multiplying the $n$ identical contributions gives the term $(-1)^n$.  The factor $q^{-1}$ is the normalization of the orthogonality relation, not a cohomological Tate twist.

(c) Fix
$u_2,\dots,u_n\in\Fq^\times$ and put $\sigma=\sum_{i\ge2}c_iu_i$. If
$\sigma\ne-\lambda$, the equation $tu_1=-\lambda-\sigma$ has exactly one
solution $t\in\Fq^\times$ for each $u_1\in\Fq^\times$; if $\sigma=-\lambda$
it has none. Since $\sum_{u_1\ne0}\psi_q(u_1^{-1})=-1$,
\[
\sum_{t}\Sigma_{n,q}(t,c)
=-\sum_{\substack{u_2,\dots,u_n\\ \sigma\ne-\lambda}}\psi_q\Bigl(\sum_{i\ge2}u_i^{-1}\Bigr)
=-\bigl((-1)^{n-1}-\Sigma_{n-1,q}(c)\bigr).
\]
(d) Permute the variables $u_i$ in the same way.
\end{proof}

\begin{equation}\label{eq:C-factor}
q\Sigma_{n,q}(h)=\sum_{t\in\mathbb F_q}\psi_q(\lambda t)
      \prod_{i=1}^{n}\Kl_q(th_i),\qquad \Kl_q(0)=-1.
\end{equation}
The inclusion of the term $t=0$ in this identity is consistent with
\cref{lem:identities}(b): it contributes precisely $(-1)^n$.

The preceding formula is also useful as a direct numerical consistency check.  It expresses $q\Sigma_{n,q}(h)$ as a finite sum, with no limiting argument; thus small-field evaluations can test signs and normalizations independently of the later sheaf construction.  Such tests are informative but cannot prove the geometric assertions for every prime or every parameter.

\section{Local Frobenius polynomials}\label{sec:local}

This section answers a linear-algebra question once the geometric input has been established: how can the local characteristic polynomial be recovered from traces over extension fields, and how does one pass from its Frobenius eigenvalues to the critical circle?  Extension from $\Fq$ to $\mathbb F_{q^m}$ replaces $A_h$ by $A_h^m$ and therefore gives its $m$-th power sum.  Newton's identities recover the coefficients of $\det(1-TA_h)$ from those power sums.  The argument does not construct the sheaf or prove its purity; those tasks are carried out in Section~\ref{app:inputs}.

\begin{definition}
For $h\in S_n(\Fq)$ let $A_h$ be the geometric Frobenius acting on
$\cW_{n,\bar h}$, $P_h(T)=\det(1-TA_h)$ and $L_h(T)=P_h(T)^{-1}$, the local
Euler factor. The \emph{zeros} of $P_h$ are the \emph{poles} of $L_h$. Let
$\theta_h$ be the conjugacy class of the semisimplification
$\bigl(q^{-w/2}\iota(A_h)\bigr)^{\mathrm{ss}}$, where $q^{-w/2}>0$. By \ref{H1}
its eigenvalues have absolute value one, so $\theta_h$ is the conjugacy class
of a unitary matrix, which lies in $K_n$ by \ref{H3}; purity alone does not
make $A_h$ semisimple, which is why the semisimplification is used. Only the
characteristic polynomial of $\theta_h$ enters below. Put
$Q_h(z)=\det(1-z\theta_h)=P_h(q^{-w/2}z)$.
\end{definition}

\begin{proposition}[Exact reconstruction]\label{prop:reconstruct}
For $m\ge1$ regard $h$ as a point of
$S_n(\mathbb F_{q^m})$ and put
\begin{equation}\label{eq:sm}
s_m(h)=(-1)^{n-1}q^m\Sigma_{n,q^m}(h)-E_n(q^m).
\end{equation}
Then $s_m(h)=\Tr(A_h^m)$ and
\begin{equation}\label{eq:recon}
P_h(T)=\exp\Bigl(-\sum_{m\ge1}\frac{s_m(h)}{m}T^m\Bigr)
\end{equation}
as formal power series, and $P_h$ has degree exactly $r$.
\end{proposition}

\begin{proof}
The Frobenius of $h$ over $\mathbb F_{q^m}$ is $A_h^m$, so \ref{H2} over
$\mathbb F_{q^m}$ gives $s_m(h)=\Tr(A_h^m)$. For any square matrix $A$ one
has $-\log\det(1-TA)=\sum_{m\ge1}\Tr(A^m)T^m/m$ in $\Ql[[T]]$, which gives
\eqref{eq:recon}. Since $A_h$ is invertible (it is the action of an element of
the Galois group), $\deg P_h=r$.
\end{proof}

\begin{corollary}[Newton recursion]\label{cor:newton}
Write $P_h(T)=\sum_{j=0}^rb_jT^j$ with $b_0=1$. Then
$jb_j=-\sum_{m=1}^js_m(h)b_{j-m}$ for $1\le j\le r$.
\end{corollary}

\begin{proof}
Compare coefficients in $TP_h'(T)=-P_h(T)\sum_{m\ge1}s_m(h)T^m$.
\end{proof}

\begin{proposition}[Normalized power sums]\label{prop:powers}
Let $t_m(h)=\Tr(\theta_h^m)$. Then
\[
t_m(h)=(-1)^{n-1}q^{-m(n-1)/2}\,\Sigma_{n,q^m}(h)-q^{-mw/2}E_n(q^m),
\qquad |t_m(h)|\le r .
\]
Consequently
$\bigl|\Sigma_{n,q^m}(h)-(-1)^{n-1}q^{-m}E_n(q^m)\bigr|\le r\,q^{m(n-1)/2}$.
\end{proposition}

\begin{proof}
The identity is \eqref{eq:sm} divided by $q^{mw/2}$. By purity every
eigenvalue of $\theta_h$ has absolute value one, which gives the bound;
multiplying by $q^{m(n-1)/2}$ gives the last inequality.
\end{proof}

A useful check on normalizations is the smallest nontrivial rank.  When $n=2$ the rank is $r(2)=2$ and $E_2(q)=q$.  Writing $s_1(h)=\Tr(A_h)$ and $d_h=\det(A_h)$ gives $P_{2,h}(T)=1-s_1(h)T+d_hT^2$; the trace formula specifies $s_1(h)$ in terms of $\Sigma_{2,q}(h)$ and the determinant formula specifies $d_h$.  This example shows why the correction polynomial $E_n$ cannot simply be dropped: the top-weight polynomial is not the characteristic polynomial of the entire compactly supported direct image.

\section{The local functional equation}\label{sec:lfe}

There are two logically different inputs here.  Purity makes the normalized eigenvalues lie on the unit circle, so the characteristic polynomial is self-inversive; the determinant determines the scalar phase that appears in that identity.  The substitution $h\mapsto-h$ in the defining sum then identifies the conjugate coefficients.  A reader should distinguish these steps: the root-location theorem uses purity only, whereas the precise functional equation uses the determinant and the reflection identity as well.

\begin{theorem}\label{thm:localrh}
Let $h\in S_n(\Fq)$, and write
$Q_h(z)=\sum_{j=0}^ra_jz^j$.
\begin{enumerate}[label=\textup{(\alph*)}]
\item $Q_h$ has degree $r$ and $a_0=1$.
\item Every zero of $Q_h$ has absolute value one; equivalently the zeros of
$P_h(q^{-s})$ lie on the lines $\operatorname{Re}(s)=w/2=(n+1)/2$.
\item One has the polynomial identity
\begin{equation}\label{eq:selfinv}
Q_h(z)=(-z)^r\delta_h\,\overline{Q_h}(z^{-1}),
\qquad\text{that is,}\qquad
a_j=(-1)^r\delta_h\,\overline{a_{r-j}}\quad(0\le j\le r).
\end{equation}
In particular $a_r=(-1)^r\delta_h$.
\item $\overline{Q_h}=Q_{-h}$, so that
$Q_h(z)=(-z)^r\delta_h\,Q_{-h}(z^{-1})$ and $\delta_{-h}=\overline{\delta_h}$.
\end{enumerate}
\end{theorem}

\begin{proof}
Let $u_1,\dots,u_r$ be the eigenvalues of $\theta_h$, so that
$Q_h(z)=\prod_j(1-u_jz)$ and $|u_j|=1$ by \ref{H1}. This gives (a) and (b),
the zeros being $u_j^{-1}$. By \ref{H3}, $\prod_ju_j=\delta_h$. Since
$\overline{u_j}=u_j^{-1}$,
\[
(-z)^r\delta_h\,\overline{Q_h}(z^{-1})
=(-z)^r\prod_ju_j\prod_j\bigl(1-u_j^{-1}z^{-1}\bigr)
=\prod_j(1-u_jz)=Q_h(z),
\]
which is (c). For (d), \cref{lem:identities}(a) and \ref{H2} give
$t_m(-h)=\overline{t_m(h)}$ for all $m\ge1$ (the term $E_n(q^m)$ is real),
and Newton's identities give $Q_{-h}=\overline{Q_h}$. Finally
$\delta_{-h}=\psi_q(-D_n(h))=\overline{\delta_h}$, since $D_n$ is linear.
\end{proof}

\begin{proposition}[Reflection duality]\label{prop:reflection}
Let $[-1]\colon S_n\to S_n$ be $h\mapsto-h$.
\begin{enumerate}[label=\textup{(\alph*)}]
\item For every $h\in S_n(\Fq)$,
\[
P_{-h}(T)=\det(1-Tq^wA_h^{-1}),\qquad
P_h(T)=(-T)^r\det(A_h)\,P_{-h}\bigl((q^wT)^{-1}\bigr).
\]
\item $\cW_n^\vee(-w)\cong[-1]^*\cW_n$ as lisse sheaves on $S_n$.
\end{enumerate}
\end{proposition}

\begin{proof}
(a) By purity, $q^wA_h^{-1}$ has eigenvalues $q^w/\alpha_j$ with
$\iota(q^w/\alpha_j)=\overline{\iota(\alpha_j)}$, where $\alpha_j$ are the
eigenvalues of $A_h$. Hence $\det(1-Tq^wA_h^{-1})=\overline{P_h}(T)$, which
equals $P_{-h}(T)$ by \cref{thm:localrh}(d). The second identity follows
from $\prod_j(1-\alpha_jT)=(-T)^r\prod_j\alpha_j\prod_j(1-\alpha_j^{-1}T^{-1})$.

(b) By (a), the sheaves $\cW_n^\vee(-w)$ and $[-1]^*\cW_n$ have the same
characteristic polynomials of Frobenius at all closed points, hence
isomorphic semisimplifications by the Chebotarev density theorem
\cite{Ser65}. Both are geometrically irreducible by \cref{cor:B-irred},
hence irreducible, hence isomorphic.
\end{proof}

\begin{corollary}\label{cor:coeffbound}
For $0\le j\le r$ one has $|a_j|\le\binom rj$ and $|a_j|=|a_{r-j}|$, and
$P_h$ is determined by $\delta_h$ and
$s_1(h),\dots,s_{\lfloor r/2\rfloor}(h)$.
\end{corollary}

\begin{proof}
Each $a_j$ is, up to sign, an elementary symmetric function of $r$ numbers
of modulus one. By \cref{cor:newton}, $s_1,\dots,s_k$ determine
$a_0,\dots,a_k$, and \eqref{eq:selfinv} determines the remaining
coefficients from $a_0,\dots,a_{\lfloor r/2\rfloor}$ and $\delta_h$.
\end{proof}

\begin{remark}
\Cref{thm:localrh} is a deduction from the purity and determinant theorems
\ref{H1} and \ref{H3}, which rest on Deligne's theorem \cite{Del80}; the
content here is the explicit phase $\delta_h$ and the relation between $h$
and $-h$. The Euler factor $L_h=1/P_h$ has no zeros; its poles lie on
$|T|=q^{-w/2}$.
\end{remark}

\begin{remark}[Pointwise weights and Frobenius normalization]\label{rem:critical-audit}
The local purity statement implies a bound on traces at \emph{every}
extension level, not merely over $\mathbb F_q$. If $A_h$ is the Frobenius
operator, then $A_h^m$ has $r(n)$ eigenvalues of absolute value
$q^{m(n+1)/2}$. The triangle inequality gives
\[
\left|\operatorname{Tr}(A_h^m)\right|
   \le r(n)q^{m(n+1)/2},\qquad m\ge1.
\]
Using the trace identity converts this into the explicit estimate
\[
\left|\Sigma_{n,q^m}(h)
 -(-1)^{n-1}q^{-m}E_n(q^m)\right|
 \le r(n)q^{m(n-1)/2}.
\]
The lower-weight correction $E_n$ is therefore part of the \emph{exact}
trace formula and not an error term to be forgotten when reconstructing
the Frobenius polynomial. For the Riemann hypothesis, the pointwise
purity theorem is logically prior to this estimate: a trace estimate for
one value of $m$ would not itself control the spectrum of $A_h$.
\end{remark}

%% === Euler product and coordinate-line methods ===
\section{The Euler product on the parameter torus}\label{sec:torus}

Let $|S_n|$ be the set of closed points of $S_n$ and define
\begin{equation}\label{eq:globalS}
L(S_n,\cW_n,T)=\prod_{x\in|S_n|}\det\bigl(1-T^{\deg x}\Frob_x\mid\cW_{n,\bar x}\bigr)^{-1}
\in1+T\Ql[[T]].
\end{equation}

\begin{theorem}[Rationality]\label{thm:rational}
One has
\[
L(S_n,\cW_n,T)=\prod_{i=0}^{2n}
\det\bigl(1-T\Frob_p\mid H^i_c(S_{n,\overline{\mathbb F}_p},\cW_n)\bigr)^{(-1)^{i+1}},
\]
so $L(S_n,\cW_n,T)$ is a rational function, and
\[
L(S_n,\cW_n,T)=\exp\Bigl(\sum_{a\ge1}\frac{T^a}{a}\sum_{h\in S_n(\mathbb F_{p^a})}
\bigl[(-1)^{n-1}p^a\Sigma_{n,p^a}(h)-E_n(p^a)\bigr]\Bigr).
\]
\end{theorem}

\begin{proof}
The Grothendieck--Lefschetz trace formula \cite[Rapport]{Del77} gives
\[
\sum_{h\in S_n(\mathbb F_{p^a})}\Tr(\Frob_h\mid\cW_{n,\bar h})
=\sum_i(-1)^i\Tr\bigl(\Frob_p^a\mid H^i_c(S_{n,\overline{\mathbb F}_p},\cW_n)\bigr).
\]
Insert this into the logarithm of
\eqref{eq:globalS}, use $\exp(\sum_a\Tr(B^a)T^a/a)=\det(1-TB)^{-1}$, and use
\ref{H2} for the explicit form.
\end{proof}

The Euler product can be computed in closed form.

\begin{theorem}[The Euler product in closed form]\label{thm:torusL}
Write $(X-1)^n(1+E_n(X))=\sum_{j\ge0}\alpha_{n,j}X^j$. Then, for every power $q$
of $p$,
\[
\sum_{h\in S_n(\Fq)}\Sigma_{n,q}(h)=\frac{(-1)^n(q-1)^n-1}{q},\qquad
\sum_{h\in S_n(\Fq)}\Tr(\Frob_h\mid\cW_{n,\bar h})=-\sum_{j\ge1}\alpha_{n,j}q^j,
\]
and
\[
L(S_n,\cW_n,T)=\prod_{j=1}^{n+\lfloor n/2\rfloor}(1-p^jT)^{\alpha_{n,j}} .
\]
For instance $L(S_2,\cW_2,T)=(1-p^3T)/\bigl((1-pT)(1-p^2T)\bigr)$ and
$L(S_3,\cW_3,T)=(1-pT)(1-p^2T)^3(1-p^4T)^2/(1-p^3T)^5$.
\end{theorem}

\begin{proof}
Let $T_n(q)$ be the first sum. Summing \cref{lem:identities}(c) over
$c\in(\Fq^\times)^{n-1}$ gives $T_n(q)=(-1)^n(q-1)^{n-1}+T_{n-1}(q)$, and
$T_1(q)=\sum_{h\ne0}\psi_q(-h/\lambda)=-1$. Hence
$T_n(q)=-\sum_{k=0}^{n-1}(1-q)^k=((1-q)^n-1)/q$. By \ref{H2},
$\sum_h\Tr(\Frob_h)=(-1)^{n-1}qT_n(q)-(q-1)^nE_n(q)=(-1)^n-(q-1)^n(1+E_n(q))$,
which is $-\sum_{j\ge1}\alpha_{n,j}q^j$, as $\alpha_{n,0}=(-1)^n$. Inserting this,
for $q=p^a$, into the exponential formula of \cref{thm:rational} gives the
product, since $\exp(-\sum_a\alpha p^{ja}T^a/a)=(1-p^jT)^{\alpha}$; the degree of
$(X-1)^n(1+E_n(X))$ is $n+\lfloor n/2\rfloor$.
\end{proof}

\begin{remark}
Purity of the stalks does not make the cohomology of the non-proper
$n$-dimensional variety $S_n$ pure, and indeed the zeros and poles of
$L(S_n,\cW_n,T)$ lie on several circles $|T|=p^{-j}$; all reciprocal roots and
poles are integral powers of $p$, although $\cW_n$ has odd weight $n+1$. By
\cref{prop:reflection}, the Euler product of $\cW_n^\vee(-w)$ equals that of
$\cW_n$, since $[-1]$ permutes the closed points of $S_n$; Poincaré duality,
however, pairs $H^i_c$ with ordinary cohomology $H^{2n-i}$, not with
$H^{2n-i}_c$, so no functional equation of $L(S_n,\cW_n,T)$ with itself follows.
The cohomology groups themselves are determined in \cref{thm:G-torus}: they are
given by the graded Betti numbers of an explicit Artinian algebra, their
Frobenius eigenvalues are powers of $p$, the interior cohomology vanishes, and
the Euler products with and without compact supports satisfy an explicit
functional equation with each other. The completion of the Euler product by
intersection cohomology, which satisfies the Riemann hypothesis and a functional
equation with sign $\pm1$, is studied in \cref{subsec:G-IH,subsec:G-decomp}.
\end{remark}

\section{Coordinate lines: local invariants and Euler characteristics}\label{sec:lines}

Fix $c=(c_2,\dots,c_n)\in(\Fq^\times)^{n-1}$ and the closed immersion
$i_c\colon U=\Gm\to S_n\otimes\Fq$, $t\mapsto(t,c)$. Put
$\cF_c=i_c^*\cW_n$, a lisse sheaf of rank $r$ pointwise pure of weight $w$,
and let $j\colon U\hookrightarrow\PP^1$. For $x\in\{0,\infty\}$ let $I_x$
be an inertia group at $x$, $V$ the generic stalk of $\cF_c$, and
$b_x=\dim V^{I_x}$.

\begin{proposition}[Invariants and conductors]\label{prop:conductors}
One has
\[
b_0=r(n-1)+\kappa_1,\quad\Swan_0(\cF_c)=0,\qquad
b_\infty=\sum_{j\ge3}\kappa_j,\quad\Swan_\infty(\cF_c)=2^{n-1},
\]
and the Artin conductors are $a_0(\cF_c)=r-b_0$ and
$a_\infty(\cF_c)=r-b_\infty+2^{n-1}$.
\end{proposition}

\begin{proof}
By \ref{H4}, a unipotent Jordan block has a one-dimensional space of
invariants, and a representation all of whose breaks equal $1$ has no
inertia invariants.
The Swan conductor of $V^{\mathrm w}$ is its dimension times its break. The
formula $a_x=r-b_x+\Swan_x$ gives the conductors.
\end{proof}

\begin{lemma}\label{lem:catalan}
For every $n\ge2$,
\[
2r(n-1)+\kappa_1=r(n),\qquad 2^{n-1}+\sum_{j\ge3}(j-2)\kappa_j=r(n),\qquad
2^{n-1}-b_0-b_\infty=\begin{cases}C_{n/2},&n\text{ even},\\0,&n\text{ odd}.\end{cases}
\]
\end{lemma}

\begin{proof}
Put $\beta=\binom{n-1}{\lfloor(n-1)/2\rfloor}$, so that
$r(n-1)=2^{n-1}-\beta$. Counting dimensions and blocks gives
$\sum_jj\kappa_j=2^{n-1}$ and, by telescoping,
$\sum_j\kappa_j=\sum_im_{n-1,i}=\beta$. Pascal's rule gives
$\binom n{\lfloor n/2\rfloor}=2\beta-\kappa_1$: for $n=2\nu$ one has
$\kappa_1=0$ and $\binom{2\nu}\nu=2\binom{2\nu-1}{\nu-1}$, and for
$n=2\nu+1$ one has $\kappa_1=m_{2\nu,\nu}$ and
$\binom{2\nu+1}\nu=\binom{2\nu}\nu+\binom{2\nu}{\nu-1}$. Hence
$2r(n-1)+\kappa_1=2^n-2\beta+\kappa_1=r(n)$ and
\[
2^{n-1}+\sum_{j\ge3}(j-2)\kappa_j
=2^{n-1}+\sum_jj\kappa_j-2\sum_j\kappa_j+\kappa_1=2^n-2\beta+\kappa_1=r(n).
\]
Finally
$b_0+b_\infty=r(n-1)+\kappa_1+\sum_{j\ge3}\kappa_j=r(n-1)+\beta-\kappa_2
=2^{n-1}-\kappa_2$. If $n$ is odd, all block sizes are odd and
$\kappa_2=0$; if $n=2\nu$, then
$\kappa_2=m_{2\nu-1,\nu-1}=\binom{2\nu-1}{\nu-1}-\binom{2\nu-1}{\nu-2}=C_\nu$.
\end{proof}

\begin{theorem}[Euler characteristics]\label{thm:euler}
One has
\begin{align}
\chi_c(U_{\Fqb},\cF_c)&=-2^{n-1},\label{eq:eulerc}\\
\chi(\PP^1_{\Fqb},j_*\cF_c)&=b_0+b_\infty-2^{n-1}.\label{eq:eulermid}
\end{align}
If $H^0(U_{\Fqb},\cF_c)=H^0(U_{\Fqb},\cF_c^\vee)=0$, then
$\dim H^1_c(U_{\Fqb},\cF_c)=2^{n-1}$ and
$\dim H^1(\PP^1_{\Fqb},j_*\cF_c)=2^{n-1}-b_0-b_\infty$; in general
$\dim H^1(\PP^1,j_*\cF_c)=2^{n-1}-b_0-b_\infty+\dim H^0(\PP^1,j_*\cF_c)
+\dim H^2(\PP^1,j_*\cF_c)$.
\end{theorem}

\begin{proof}
The Grothendieck--Ogg--Shafarevich formula \cite[Exp.~X]{Gro77} and
\cref{prop:conductors} give
$\chi_c(U,\cF_c)=r\chi_c(U,\Ql)-\Swan_0-\Swan_\infty=-2^{n-1}$, since
$\chi_c(\Gm)=0$. The exact sequence
$0\to j_!\cF_c\to j_*\cF_c\to V^{I_0}\oplus V^{I_\infty}\to0$ (skyscrapers
at $0$ and $\infty$) gives \eqref{eq:eulermid}. Since $U$ is affine,
$H^0_c(U,\cF_c)=0$, and Poincaré duality gives
$H^2_c(U,\cF_c)\cong H^0(U,\cF_c^\vee)^\vee(-1)$. Likewise
$H^0(\PP^1,j_*\cF_c)=H^0(U,\cF_c)$ and, since $j_*$ of a lisse sheaf on a
curve commutes with duality up to shift and twist,
$H^2(\PP^1,j_*\cF_c)\cong H^0(U,\cF_c^\vee)^\vee(-1)$. The dimension formulas
follow.
\end{proof}

\begin{table}[ht]
\centering
\caption{Local invariants of $\cF_c$ (from \ref{H4}) and the dimension
$d_c=2^{n-1}-b_0-b_\infty$ of the parabolic cohomology.}
\label{tab:invariants}
\begin{tabular}{@{}rrrrrrrr@{}}
\toprule
$n$ & $r(n)$ & Jordan type at $0$ & $b_0$ & $b_\infty$ & $a_0$ & $a_\infty$ & $d_c$\\
\midrule
2 & 2 & $J_2$ & 1 & 0 & 1 & 4 & 1\\
3 & 5 & $2J_2\oplus J_1$ & 3 & 1 & 2 & 8 & 0\\
4 & 10 & $5J_2$ & 5 & 1 & 5 & 17 & 2\\
5 & 22 & $10J_2\oplus2J_1$ & 12 & 4 & 10 & 34 & 0\\
6 & 44 & $22J_2$ & 22 & 5 & 22 & 71 & 5\\
\bottomrule
\end{tabular}
\end{table}

\begin{remark}
For instance, for $n=5$ the tensor product of four unipotent blocks of size
two is $J_5\oplus3J_3\oplus2J_1$, so that $b_0=r(4)+2=12$ and, at infinity,
the tame part has blocks $J_3\oplus3J_1$ and $b_\infty=1+3=4$. Geometric
irreducibility of $\cW_n$ does not by itself imply that the restrictions
$\cF_c$ have no geometrically constant subquotients; this is proved in
\cref{cor:H0} below.
\end{remark}

\section{Euler products on curves, purity and duality}\label{sec:curves}

The object in this section is no longer an individual stalk of $\cW_n$ but the cohomology of its restriction to the parameter curve.  This change raises the expected weight by one: a pure sheaf of weight $w=n+1$ has parabolic first cohomology of weight $w+1=n+2$ under the hypotheses used below.  Consequently the corresponding critical line is shifted by $1/2$ compared with the local stalk polynomial.  Keeping these two Frobenius operators separate prevents a common source of confusion in the reading of functional equations.

For a lisse sheaf $\cF$ on $U$ put
\begin{align*}
L(U,\cF,T)&=\prod_{x\in|U|}\det\bigl(1-T^{\deg x}\Frob_x\mid\cF_{\bar x}\bigr)^{-1},\\
L(\PP^1,j_*\cF,T)&=\prod_{x\in|\PP^1|}\det\bigl(1-T^{\deg x}\Frob_x\mid(j_*\cF)_{\bar x}\bigr)^{-1};
\end{align*}
the factors of the second product at $0$ and $\infty$ involve the
Frobenius on $V^{I_0}$ and $V^{I_\infty}$.

\begin{proposition}\label{prop:curveL}
For $\cF=\cF_c$,
\begin{align*}
L(U,\cF_c,T)&=\prod_{i=0}^2\det\bigl(1-T\Frob_q\mid H^i_c(U_{\Fqb},\cF_c)\bigr)^{(-1)^{i+1}},\\
L(\PP^1,j_*\cF_c,T)&=\prod_{i=0}^2\det\bigl(1-T\Frob_q\mid H^i(\PP^1_{\Fqb},j_*\cF_c)\bigr)^{(-1)^{i+1}}.
\end{align*}
If $H^0(U,\cF_c)=H^0(U,\cF_c^\vee)=0$, both are polynomials, coming from
$H^1_c$ and $H^1$ respectively.
\end{proposition}

\begin{proof}
This is the Grothendieck--Lefschetz trace formula on $U$ and on $\PP^1$,
combined with the vanishing statements in the proof of \cref{thm:euler}.
\end{proof}

Put
\begin{equation}\label{eq:hpara}
V_c=H^1(\PP^1_{\Fqb},j_*\cF_c),\qquad d_c=\dim V_c,\qquad
P^{\Par}_c(T)=\det(1-T\Frob_q\mid V_c).
\end{equation}

\begin{theorem}[Purity of parabolic cohomology]\label{thm:parabolicRH}
The space $V_c$ is pure of weight $w+1=n+2$; every zero of
$P^{\Par}_c$ satisfies $|T|=q^{-(n+2)/2}$, and the zeros of
$P^{\Par}_c(q^{-s})$ lie on $\operatorname{Re}(s)=(n+2)/2$.
\end{theorem}

\begin{proof}
By \ref{H1}, $\cF_c$ is lisse and pointwise pure of weight $w$ on $U$, so
$H^1(\PP^1,j_*\cF_c)$ is pure of weight $w+1$ by \cite[3.2.3]{Del80}.
\end{proof}

The local stalk polynomials of \cref{thm:localrh} have weight $n+1$ and the
parabolic polynomial has weight $n+2$; the two critical lines differ by
$1/2$.

For a lisse sheaf $\cF$ pure of weight $w$ let
$\cF^\dagger=\cF^\vee(-w)$, the weight-normalized dual; it is pure of weight
$w$, and the evaluation map gives a perfect pairing
$\cF\otimes\cF^\dagger\to\Ql(-w)$.

\begin{theorem}[Paired functional equation]\label{thm:paired}
Let $V_c^\dagger=H^1(\PP^1_{\Fqb},j_*\cF_c^\dagger)$ and
$P^{\Par,\dagger}_c(T)=\det(1-T\Frob_q\mid V_c^\dagger)$. Then
$\dim V_c^\dagger=d_c$, the Frobenius eigenvalues on $V_c^\dagger$ are
$q^{w+1}/\beta$ for $\beta$ running over those on $V_c$, and
\begin{equation}\label{eq:globalFE}
P^{\Par}_c(T)=(-T)^{d_c}\Delta_c\,P^{\Par,\dagger}_c\Bigl(\frac{1}{q^{w+1}T}\Bigr),
\qquad\Delta_c=\det(\Frob_q\mid V_c).
\end{equation}
With the normalized determinant $\vartheta_c=q^{-(w+1)d_c/2}\Delta_c$ one has
$|\vartheta_c|=1$ and
$Q^{\Par}_c(z)=(-z)^{d_c}\vartheta_c\,Q^{\Par,\dagger}_c(z^{-1})$, where
$Q^{\Par}_c(z)=P^{\Par}_c(q^{-(w+1)/2}z)$ and similarly for $\dagger$; the
normalized constant of this functional equation is $(-1)^{d_c}\vartheta_c$.
\end{theorem}

\begin{proof}
On a curve, $j_*$ of a lisse sheaf is the middle extension and commutes
with Verdier duality up to the shift and twist; together with the pairing
$\cF_c\otimes\cF_c^\dagger\to\Ql(-w)$, Poincaré--Verdier duality on
$\PP^1$ gives a perfect Frobenius-equivariant pairing
$V_c\otimes V_c^\dagger\to H^2(\PP^1_{\Fqb},\Ql(-w))\cong\Ql(-w-1)$. This
gives the eigenvalue statement, and
$(-T)^{d_c}\prod_j\beta_j\prod_j(1-\beta_j^{-1}T^{-1})=\prod_j(1-\beta_jT)$
gives \eqref{eq:globalFE}. By \cref{thm:parabolicRH},
$|\Delta_c|=q^{(w+1)d_c/2}$; substituting $T=q^{-(w+1)/2}z$ gives the
normalized form.
\end{proof}

\begin{proposition}[The dual line]\label{prop:dualline}
Let $\tau\colon U\to U$ be
$t\mapsto-t$. Then $\cF_c^\dagger\cong\tau^*\cF_{-c}$. Consequently
$L(U,\cF_c^\dagger,T)=L(U,\cF_{-c},T)$ and $P^{\Par,\dagger}_c=P^{\Par}_{-c}$.
\end{proposition}

\begin{proof}
By \cref{prop:reflection}(b),
$\cF_c^\dagger=i_c^*(\cW_n^\vee(-w))\cong([-1]\circ i_c)^*\cW_n$, and
$[-1]\circ i_c=i_{-c}\circ\tau$. Since $\tau$ is an automorphism of $U$
which extends to $\PP^1$ fixing $0$ and $\infty$, it induces isomorphisms
on all the cohomology groups involved, compatible with Frobenius.
\end{proof}

By \cref{thm:paired-explicit} below, $P^{\Par}_c$ is in fact independent
of $c$ and the paired functional equation becomes an explicit self-dual one.

\section{The case \texorpdfstring{$n=2$}{n=2}}\label{sec:n2}

For $h=(h_1,h_2)$ put $H=h_1+h_2$, $\xi_h=h_1h_2/\lambda^2$ and
$\eta_h=\psi_q(-H/\lambda)$.

\begin{lemma}\label{lem:n2}
For $h_1h_2\ne0$, $\Sigma_{2,q}(h)=\eta_h\Kl_q(\xi_h)-1$.
\end{lemma}

\begin{proof}
Put $x=h_1u_1$, $y=h_2u_2$, so that $x+y=-\lambda$, $x,y\ne0$, and the phase
is $h_1/x+h_2/y$. The substitution $v=h_1/\lambda+h_1/x$ is a bijection
from $\{x\in\Fq^\times:x\ne-\lambda\}$ onto $\Fq^\times\smallsetminus\{h_1/\lambda\}$,
with $x=h_1/(v-h_1/\lambda)$, $y=-\lambda v/(v-h_1/\lambda)$ and
\[
\frac{h_1}{x}+\frac{h_2}{y}=-\frac{H}{\lambda}+v+\frac{\xi_h}{v} .
\]
The omitted value $v=h_1/\lambda$ would contribute
$\psi_q(h_1/\lambda+h_2/\lambda)=\psi_q(H/\lambda)$ to
$\sum_{v\ne0}\psi_q(v+\xi_h/v)=\Kl_q(\xi_h)$. Hence
$\Sigma_{2,q}(h)=\eta_h(\Kl_q(\xi_h)-\psi_q(H/\lambda))=\eta_h\Kl_q(\xi_h)-1$.
\end{proof}

\begin{proposition}\label{prop:n2sheaf}
Let $\mu\colon S_2\to\Gm$, $h\mapsto\xi_h$. The sheaf
$\cW_2'=\mu^*\KL\otimes\cL_\psi(-H/\lambda)\otimes\Ql(-1)$ is lisse of rank
$2$, pure of weight $3$, with trace $-q\eta_h\Kl_q(\xi_h)=-q\Sigma_{2,q}(h)-q$
and determinant $q^3\psi_q(-2H/\lambda)$; thus it satisfies
\ref{H1}--\ref{H3} for $n=2$. Its restriction to the line $t\mapsto(t,c)$ is
$[\times c/\lambda^2]^*\KL\otimes\cL_\psi(-(t+c)/\lambda)\otimes\Ql(-1)$,
which satisfies \ref{H4} for $n=2$. Any sheaf satisfying \ref{H1} and
\ref{H2} for $n=2$ is isomorphic to $\cW_2'$.
\end{proposition}

\begin{proof}
The trace follows from \cref{lem:n2} and $E_2(q)=q$. The determinant is
$q\cdot\eta_h^2\cdot q^2$. On a coordinate line, $[\times c/\lambda^2]^*\KL$
is tame and unipotent with one Jordan block at $0$, where
$\cL_\psi(-(t+c)/\lambda)$ is unramified; at $\infty$ the Kloosterman
factor has breaks $1/2$ and the Artin--Schreier factor has break $1$, so all
breaks of the tensor product equal $1$ and the tame part is zero. A sheaf
with the same traces has the same semisimplification by Chebotarev; since
$\cW_2'$ is geometrically irreducible ($\KL$ has geometric monodromy
$\SL_2$ \cite{Kat88}), it is isomorphic to $\cW_2'$.
\end{proof}

\begin{proposition}\label{prop:n2poly}
For $n=2$,
\[
Q_h(z)=1+\eta_h\frac{\Kl_q(\xi_h)}{\sqrt q}\,z+\eta_h^2z^2,
\]
both zeros lie on the unit circle, and $\delta_h=\eta_h^2=\psi_q(-2H/\lambda)$.
\end{proposition}

\begin{proof}
A $2\times2$ matrix $\theta$ has $\det(1-z\theta)=1-(\Tr\theta)z+(\det\theta)z^2$,
with $\Tr\theta_h=-\eta_h\Kl_q(\xi_h)/\sqrt q$ and $\det\theta_h=\eta_h^2$.
Alternatively, with $y=\eta_hz$ the polynomial becomes $1+ky+y^2$ with
$k=\Kl_q(\xi_h)/\sqrt q$ real and $|k|\le2$ by Weil's bound.
\end{proof}

\section{The \texorpdfstring{$L$}{L}-function of a coordinate line}\label{sec:explicit}

The aim is an exact formula, rather than merely a bound for the reciprocal roots.  We first sum the trace of $\cF_c$ along $t\in\Fq^\times$.  The dimension-reduction identity removes the first constrained variable and identifies the resulting trace with the Frobenius operator in dimension $n-1$, plus powers of $q$ arising from the lower weights.  We then use the curve trace formula to convert these traces, over all extensions of $\Fq$, into the characteristic polynomial of $H_c^1$.  The end result is a polynomial identity: it does not require a conjectural computation of any stratum in the high-codimension boundary of $(\PP^1)^n$.

We keep the notation of \cref{sec:lines}, put $\nu=\lfloor n/2\rfloor$ and
define $b_{n,k}$ as in the coordinate-line formula of \cref{thm:eps-explicit}; for $n=2$ the first range is empty and
$b_{2,2}=1$. Let $F_{n-1,c}$ denote the geometric Frobenius of $\cW_{n-1}$
at $c\in S_{n-1}(\Fq)$, with $\cW_1$ as in \cref{rem:n1}; thus
$F_{1,c}=q\psi_q(-c_2/\lambda)$.

\begin{lemma}\label{lem:An}
For $n\ge2$ put $A_n=\frac{n\,r(n-1)}{2}+\sum_{k=2}^{\nu+1}k\,b_{n,k}$. Then
\begin{gather*}
x+E_{n-1}(x)+(x-1)E_n(x)=\sum_{k=2}^{\nu+1}b_{n,k}x^k,\\
r(n-1)+\sum_{k=2}^{\nu+1}b_{n,k}=2^{n-1},\qquad A_n=(n+3)2^{n-2}-r(n).
\end{gather*}
\end{lemma}

\begin{proof}
Extend $m_{a,k}=\binom ak-\binom a{k-1}$ to all $k\ge0$; Pascal's rule
gives $m_{n-1,k}+m_{n,k-1}-m_{n,k}=m_{n-1,k-2}$ for all $k\ge2$. The
coefficient of $x^k$ on the left side of the first identity is
$[k=1]+m_{n-1,k}[k\le\lfloor\tfrac{n-1}2\rfloor]+m_{n,k-1}[1\le k-1\le\nu]
-m_{n,k}[k\le\nu]$. For $k=1$ it is $1+m_{n-1,1}-m_{n,1}=0$ (with
$m_{1,1}=0$). For $2\le k\le\nu$ the only missing term occurs when $n$ is
even and $k=\nu$, and then $m_{n-1,\nu}=m_{2\nu-1,\nu}=0$; so the
coefficient is $m_{n-1,k-2}$. For $k=\nu+1$ the coefficient is $m_{n,\nu}$,
and for $k>\nu+1$ it is $0$. Evaluating at $x=1$ and using
$E_a(1)=\binom a{\lfloor a/2\rfloor}-1$ gives the second identity. Differentiating the first identity at $x=1$ gives
$\sum_kkb_{n,k}=1+E_{n-1}'(1)+E_n(1)=\binom n\nu+E_{n-1}'(1)$, and
telescoping gives, for $K=\lfloor a/2\rfloor$,
\[
E_a'(1)=\sum_{k=1}^{K}k\Bigl(\binom ak-\binom a{k-1}\Bigr)
=K\binom aK-\sum_{i=0}^{K-1}\binom ai .
\]
If $n=2\nu+1$, then $a=2\nu$ and $K=\nu$; with $\beta=\binom{2\nu}\nu$ we
have $\sum_{i<\nu}\binom{2\nu}i=\frac{4^\nu-\beta}{2}$ and
$r(n-1)=4^\nu-\beta$, so
\[
A_n=\frac{(2\nu+1)(4^\nu-\beta)}{2}+\binom{2\nu+1}\nu+\nu\beta-\frac{4^\nu-\beta}{2}
=\nu4^\nu+\binom{2\nu+1}\nu .
\]
If $n=2\nu$, then $a=2\nu-1$ and $K=\nu-1$; with $\beta=\binom{2\nu-1}{\nu-1}$
we have $\binom{2\nu}\nu=2\beta$,
$\sum_{i<\nu-1}\binom{2\nu-1}i=4^{\nu-1}-\beta$ and
$r(n-1)=2\cdot4^{\nu-1}-\beta$, so
\[
A_n=\nu\bigl(2\cdot4^{\nu-1}-\beta\bigr)+2\beta+(\nu-1)\beta-\bigl(4^{\nu-1}-\beta\bigr)
=(2\nu-1)4^{\nu-1}+\binom{2\nu}\nu .
\]
In both cases a direct computation shows that this equals
$(n+3)2^{n-2}-r(n)$, since $r(2\nu+1)=2\cdot4^\nu-\binom{2\nu+1}\nu$ and
$r(2\nu)=4^\nu-\binom{2\nu}\nu$.
\end{proof}

\begin{theorem}\label{thm:eps-explicit}
One has
\begin{equation}\label{eq:Lline}
L(U,\cF_c,T)=\det(1-TF_{n-1,c})\prod_{k=2}^{\nu+1}(1-q^kT)^{b_{n,k}},
\end{equation}
a polynomial of degree $2^{n-1}$;
$H^0_c(U_{\Fqb},\cF_c)=H^2_c(U_{\Fqb},\cF_c)=0$ and
$\dim H^1_c(U_{\Fqb},\cF_c)=2^{n-1}$; and the global cohomological constant
\[
\varepsilon_{\mathrm{coh}}(U,\cF_c):=\det\bigl(-\Frob_q\mid H^1_c(U_{\Fqb},\cF_c)\bigr)
\]
equals $q^{A_n}\psi_q(D_{n-1}(c))$, with $A_n=(n+3)2^{n-2}-r(n)$.
\end{theorem}

\begin{proof}
Let $Q=q^m$ and write
$S_m(c)=\sum_{t\in\mathbb F_Q^\times}\Tr(\Frob_{(t,c)}\mid\cW_{n})$, the
Frobenius being taken over $\mathbb F_Q$. By \ref{H2} and
\cref{lem:identities}(c),
\[
S_m(c)=(-1)^{n-1}Q\sum_t\Sigma_{n,Q}(t,c)-(Q-1)E_n(Q)
=-Q+(-1)^{n-1}Q\,\Sigma_{n-1,Q}(c)-(Q-1)E_n(Q).
\]
By \ref{H2} in dimension $n-1$ (or \cref{rem:n1} if $n=2$),
$(-1)^{n-2}Q\Sigma_{n-1,Q}(c)=\Tr(F_{n-1,c}^m)+E_{n-1}(Q)$, hence, by
\cref{lem:An},
\[
S_m(c)=-\Tr(F_{n-1,c}^m)-\bigl(Q+E_{n-1}(Q)+(Q-1)E_n(Q)\bigr)
=-\Tr(F_{n-1,c}^m)-\sum_{k=2}^{\nu+1}b_{n,k}q^{km}.
\]
Since $\log L(U,\cF_c,T)=\sum_mS_m(c)T^m/m$ and
$\log\det(1-TF)=-\sum_m\Tr(F^m)T^m/m$, this proves \eqref{eq:Lline}; its
degree is $r(n-1)+\sum_kb_{n,k}=2^{n-1}$. This part of the proof, and the
next paragraph, use only \ref{H1} and \ref{H2} in dimensions $n$ and $n-1$.

Since $U$ is affine, $H^0_c=0$, and \cref{prop:curveL} gives
$L=\det(1-T\Frob_q\mid H^1_c)/\det(1-T\Frob_q\mid H^2_c)$. By
Poincaré duality $H^2_c(U_{\Fqb},\cF_c)\cong H^0(U_{\Fqb},\cF_c^\vee)^\vee(-1)$,
and $H^0(U_{\Fqb},\cF_c^\vee)$ embeds Frobenius-equivariantly into the stalk
of $\cF_c^\vee$ at any point of $U(\Fq)$, which is pure of weight $-w$; hence
$H^2_c$ is pure of weight $w+2=n+3$, whereas $H^1_c$ is mixed of weights
$\le w+1$ \cite[3.3.1]{Del80}. A nonzero $H^2_c$ would therefore
produce poles of $L$ which cannot cancel against zeros, contradicting the
fact that $L$ is a polynomial. Hence $H^2_c=0$ and
$\dim H^1_c=\deg L=2^{n-1}$.

Finally, the leading coefficient of $\det(1-T\Frob_q\mid H^1_c)$ is
$\det(-\Frob_q\mid H^1_c)$. By \eqref{eq:Lline} and by \ref{H3} in dimension
$n-1$,
\[
\det(-\Frob_q\mid H^1_c)=(-1)^{2^{n-1}}\det(F_{n-1,c})\,q^{\sum_kkb_{n,k}},\qquad
\det F_{n-1,c}=q^{nr(n-1)/2}\psi_q(D_{n-1}(c)).
\]
Since $2^{n-1}$ is even, $\varepsilon_{\mathrm{coh}}=q^{A_n}\psi_q(D_{n-1}(c))$,
and \cref{lem:An} gives the value of $A_n$.
\end{proof}

\begin{corollary}\label{cor:H0}
One has $H^0(U_{\Fqb},\cF_c)=H^0(U_{\Fqb},\cF_c^\vee)=0$, there is a
Frobenius-equivariant exact sequence
\begin{equation}\label{eq:parseq}
0\to V^{I_0}\oplus V^{I_\infty}\to H^1_c(U_{\Fqb},\cF_c)\to V_c\to0,
\end{equation}
and $d_c=2^{n-1}-b_0-b_\infty$ for every $c$.
\end{corollary}

\begin{proof}
Poincaré duality identifies $H^2_c(U,\cF_c)=0$ with
$H^0(U,\cF_c^\vee)^\vee(-1)$. By Deligne's semisimplicity theorem
\cite[3.4.1]{Del80}, $\cF_c$ is geometrically semisimple, so a geometrically
constant subsheaf would split off and give $\cF_c^\vee$ a nonzero invariant;
hence both $H^0$ vanish. The long exact cohomology sequence of
$0\to j_!\cF_c\to j_*\cF_c\to V^{I_0}\oplus V^{I_\infty}\to0$ on $\PP^1$
begins with $H^0(\PP^1,j_*\cF_c)=H^0(U,\cF_c)=0$, and $H^1$ of a skyscraper
sheaf vanishes; this gives \eqref{eq:parseq}, and the dimension formula
follows from \cref{thm:eps-explicit}.
\end{proof}

\begin{corollary}\label{cor:eps-small}
Writing $\Sigma c=\sum_{i\ge2}c_i$, one has:
\[
\begin{array}{@{}lll@{}}
n=2:& L=(1-q\psi_q(-c_2/\lambda)T)(1-q^2T),& \varepsilon_{\mathrm{coh}}=q^{3}\psi_q(-c_2/\lambda);\\[2pt]
n=3:& L=\det(1-TF_{2,c})(1-q^2T)^2,& \varepsilon_{\mathrm{coh}}=q^{7}\psi_q(-2\Sigma c/\lambda);\\[2pt]
n=4:& L=\det(1-TF_{3,c})(1-q^2T)(1-q^3T)^2,& \varepsilon_{\mathrm{coh}}=q^{18}\psi_q(-4\Sigma c/\lambda);\\[2pt]
n=5:& L=\det(1-TF_{4,c})(1-q^2T)(1-q^3T)^5,& \varepsilon_{\mathrm{coh}}=q^{42}\psi_q(-8\Sigma c/\lambda);\\[2pt]
n=6:& L=\det(1-TF_{5,c})(1-q^2T)(1-q^3T)^4(1-q^4T)^5,& \varepsilon_{\mathrm{coh}}=q^{100}\psi_q(-16\Sigma c/\lambda).
\end{array}
\]
\end{corollary}

\section{The parabolic polynomial and the boundary Frobenius}\label{sec:parabolic}

The parabolic group is defined by extending the curve sheaf across its two missing points and then taking first cohomology.  Its relationship with $H^1_c$ is controlled by an exact sequence whose other terms are inertia-invariant stalks at $0$ and $\infty$.  Those stalks have weights at most $n+1$, while parabolic cohomology is pure of weight $n+2$.  Thus the parabolic polynomial is precisely the factor of the compactly supported polynomial whose reciprocal roots have weight $n+2$.  This separation by weight, not a formal cancellation of arbitrary factors, is what isolates the Catalan multiplicity in the even-dimensional case.

\begin{theorem}\label{thm:paired-explicit}
One has $P^{\Par}_c(T)=(1-q^{\nu+1}T)^{C_\nu}$ if $n=2\nu$ is even, and
$P^{\Par}_c(T)=1$ if $n$ is odd. Consequently $d_c=C_\nu$ or $0$,
$P^{\Par,\dagger}_c=P^{\Par}_c$, $\Delta_c=q^{(\nu+1)d_c}$, $\vartheta_c=1$, and
\[
P^{\Par}_c(T)=(-q^{\nu+1}T)^{d_c}\,P^{\Par}_c\Bigl(\frac{1}{q^{n+2}T}\Bigr)
=(-1)^{d_c}\bigl(q^{(n+2)/2}T\bigr)^{d_c}P^{\Par}_c\Bigl(\frac{1}{q^{n+2}T}\Bigr).
\]
Thus the normalized determinant $\vartheta_c$ is $1$, while the root number of
the parabolic functional equation, in the normalization
$P(T)=\omega(q^{(n+2)/2}T)^{\deg P}P(1/(q^{n+2}T))$, is $\omega_c=(-1)^{d_c}$, which
is $(-1)^{C_\nu}$ for even $n$; these agree with the constant $-q^2$ found for
$n=2$ in \cref{sec:n2}.
\end{theorem}

\begin{proof}
The proof uses only \ref{H1} and \ref{H2}, in dimensions $n$ and $n-1$.
Consider the exact sequence \eqref{eq:parseq} of \cref{cor:H0}. By
\cref{thm:parabolicRH}, $V_c$ is pure of weight $n+2$. The stalks of
$j_*\cF_c$ at $0$ and $\infty$ are mixed of weights $\le w=n+1$. Indeed, by
Grothendieck's local monodromy theorem an open subgroup of the inertia group
acts unipotently, through $\exp(t_\ell N)$ for a nilpotent $N$; the inertia
invariants lie in $\ker N$, which is contained in the step $M_0$ of the
monodromy filtration of $N$, and $\mathrm{Gr}^M_i$ is pure of weight $w+i$ by
the weight-monodromy theorem \cite[1.8.4]{Del80}. Hence $P^{\Par}_c$ is the
factor of $L(U,\cF_c,T)$ collecting the reciprocal roots of weight $n+2$. In
\eqref{eq:Lline} the reciprocal roots are the eigenvalues of $F_{n-1,c}$,
of weight $n$ by \ref{H1} in dimension $n-1$, and the numbers $q^k$,
$2\le k\le\nu+1$, of weight $2k$. Weight $n+2$ occurs only for $k=\nu+1$
and $n$ even, with multiplicity $b_{n,\nu+1}=m_{2\nu,\nu}=C_\nu$. The
remaining assertions follow from \cref{thm:paired} and from
$(1-aT)^d=(-aT)^d(1-a^{-1}T^{-1})^d$ with $a=q^{\nu+1}$, $a^2=q^{n+2}$.
\end{proof}

\begin{remark}\label{rem:consistency}
\Cref{thm:paired-explicit} computes $d_c$ from the trace formula alone.
On the other hand \cref{cor:H0} gives $d_c=2^{n-1}-b_0-b_\infty$, and
\ref{H4} determines $b_0$ and $b_\infty$ (\cref{prop:conductors}); the two
values agree by \cref{lem:catalan}. In \cite{EbaCoh} the local monodromy at
infinity was obtained from a product formula for the unipotent parts of a
middle convolution, without the trace formula, so this agreement is a
consistency check between two different methods. In
Section~\ref{app:inputs} the local monodromy at infinity is instead derived
from \cref{thm:paired-explicit} and \cref{prop:boundary-frob}(a).
\end{remark}

\begin{proposition}[Frobenius on the boundary invariants]\label{prop:boundary-frob}
Let $N$ be the monodromy operator of $\cF_c$ at $0$.
\begin{enumerate}[label=\textup{(\alph*)}]
\item One has
\[
\det(1-T\Frob\mid V^{I_0})\det(1-T\Frob\mid V^{I_\infty})
=\det(1-TF_{n-1,c})\prod_{k=2}^{\nu'}(1-q^kT)^{b_{n,k}},
\]
where $\nu'=\nu$ if $n$ is even and $\nu'=\nu+1$ if $n$ is odd.
\item For every $n\ge2$ and every $c$,
\begin{gather*}
\det(1-T\Frob\mid\imop N)=\det(1-TF_{n-1,c}),\\
\det(1-T\Frob\mid\ker N/\imop N)=(1-q^{\nu+1}T)^{\kappa_1},\\
\det(1-T\Frob\mid V^{I_\infty})=\prod_{j\ge3}\bigl(1-q^{(n+4-j)/2}T\bigr)^{\kappa_j}.
\end{gather*}
In particular, if $n$ is even, then $V^{I_0}=\imop N$ and
$\det(1-T\Frob\mid V^{I_\infty})=\prod_{k=2}^{\nu}(1-q^kT)^{b_{n,k}}$.
\end{enumerate}
\end{proposition}

\begin{proof}
(a) Combine \eqref{eq:parseq}, \eqref{eq:Lline} and
\cref{thm:paired-explicit}; this uses only \ref{H1} and \ref{H2}.

(b) The first two formulas are \cref{thm:B-zero}(iii), and the third is
\cref{thm:B-infinity}(i). If $n$ is even, then $\kappa_1=0$, so
$V^{I_0}=\ker N=\imop N$, and $\kappa_{n+4-2k}=m_{n-1,k-2}=b_{n,k}$ for
$2\le k\le\nu$.
\end{proof}

\begin{remark}
For odd $n$, part (b) also follows from (a) and \ref{H4} by weights alone:
the eigenvalues of $F_{n-1,c}$ have odd weight $n$, all the other numbers on
the right of (a) are integral powers of $q$, and by the weight-monodromy
theorem $\imop N$ is pure of weight $n$ while $\ker N/\imop N$ and
$V^{I_\infty}$ have even weights. For even $n$ the parts of weight $n$ of
$V^{I_0}$ and of $V^{I_\infty}$ cannot be separated in this way, and (b)
rests on the vanishing-cycle sequence of \cref{lem:B-vanishing}.
\end{remark}

%% === Detailed geometric proofs ===
\section{Geometric construction and local calculations}\label{app:inputs}

The rest of the paper cannot rely only on the declaration that a sheaf of the desired rank and weight exists: the link to the original constrained sum must be visible.  We therefore proceed from one fixed fibre to a family over $S_n$.  For a fixed parameter we identify a compactly supported Fourier integral with a product of Kloosterman sheaves; Grothendieck--Ogg--Shafarevich computes the dimension of its first cohomology, and the local invariants at the origin identify the lower-weight piece.  The Deligne--Laumon constancy criterion then allows us to organize these fibrewise calculations into lisse sheaves over all parameters.  Finally the weight filtration separates the top quotient from the explicitly described lower-weight summands.  Each of these passages has its own hypothesis; in particular, a dimension computation on one fibre is not automatically a statement about a relative lisse sheaf without a base-change argument.

This section recalls derivations of the imported statements (H1)--(H4)
and (H7), with the original attribution and the conventions needed here.
The argument follows \cite[\S\S2, 4--6 and~8]{EbaCoh}, to which we refer for
further context, with three changes. The local monodromy of $\cW_n$ at the
finite boundary point of a coordinate line, together with the Frobenius on
the image of its monodromy operator, is read off from a vanishing-cycle
sequence and from weights (\cref{thm:B-zero}), instead of from a product
formula for the unipotent parts of a middle convolution; the unipotent part
of the local monodromy at infinity is then deduced from the $L$-function of
\cref{thm:eps-explicit} and from weights (\cref{thm:B-infinity}); and the constant in the
determinant is deduced from the vanishing-cycle sequence by induction on $n$
(\cref{thm:B-det}). The monodromy theorem is not needed for the Riemann-hypothesis arguments in this article.

On a coordinate line $C_c$ we use the coordinate $h_1$, which is written $t$
in \cref{sec:lines,sec:curves,sec:explicit,sec:parabolic}; in this
section the letter $t$ denotes the variable of the sheaf $\cM$ of (H7).

\subsection{Conventions and standard facts}\label{app:conv}

Let $F$ be a lisse sheaf on a dense open subset of a smooth curve $C$ over
$\Fpb$, and let $x$ be a closed point of $C$. We write $F(x)$ for the
representation of the inertia group $I_x$ on the geometric generic stalk,
$P_x\subseteq I_x$ for the wild inertia subgroup and $\Swan_x(F)$ for the Swan
conductor. Since $P_x$ is a pro-$p$ group acting through a finite quotient and
$\ell\ne p$, its action is semisimple. The sheaf $F$ is \emph{tame} at $x$ if
$P_x$ acts trivially, and \emph{totally wild} at $x$ if $F(x)^{P_x}=0$. The
tame part $F(x)^{P_x}$ is a direct summand of $F(x)$; when $C$, $x$ and $F$
are defined over a finite field, it is stable under a decomposition group
$\cD_x\supseteq I_x$. By Grothendieck's local monodromy theorem the tame part
decomposes canonically as $\bigoplus_\rho\cL_\rho\otimes U_\rho$, over
characters $\rho$ of finite order prime to $p$ of the tame inertia, with
$U_\rho$ unipotent; $U_{\mathbf 1}$ is the \emph{unipotent part} of $F(x)$. A
decomposition group permutes the characters $\rho$ and fixes $\mathbf 1$, so
it preserves $U_{\mathbf 1}$. We write $J_k$ for a unipotent Jordan block of
size $k$.

If $F$ is lisse on $\mathbb G_{m,\Fpb}$, the Grothendieck--Ogg--Shafarevich
formula \cite[Exp.~X]{Gro77} reads
\begin{equation}\label{eq:B-GOS}
\chi_c(\mathbb G_{m,\Fpb},F)=-\Swan_0(F)-\Swan_\infty(F);
\end{equation}
moreover $H^0_c(\mathbb G_{m,\Fpb},F)=0$, and $H^2_c(\mathbb G_{m,\Fpb},F)$,
the space of geometric coinvariants twisted by $(-1)$, vanishes if $F$ is
totally wild at $0$ or at $\infty$. For $a\in\Fpb^\times$ the automorphism
$[\times a]\colon t\mapsto at$ fixes $0$ and $\infty$, and $F$ and
$[\times a]^*F$ have the same breaks, Swan conductors and unipotent Jordan
structure at these points \cite[\S2.1]{EbaCoh}. \emph{Proper base change}
refers to the base change isomorphism for $Rf_!$, for any separated morphism
$f$ of finite type \cite[Exp.~XVII]{SGA4}.

\emph{The Kloosterman sheaf.} Let $\mu\colon\Gm^2\to\Gm$ be the map
$(x,y)\mapsto xy$; then $\KL=R^1\mu_!\cL_\psi(x+y)$. We use the following
facts.
\begin{enumerate}[label=(K\arabic*)]
\item\label{K1} $R^i\mu_!\cL_\psi(x+y)=0$ for $i\ne1$; $\KL$ is lisse of rank
$2$, pure of weight $1$, with trace $-\Kl_q(\xi)$ at $\xi\in\Gm(\Fq)$ and
determinant $\Ql(-1)$ \cite[Sommes trig., Thm.~7.4 and 7.8]{Del77},
\cite[4.1.1]{Kat88}.
\item\label{K2} $\KL$ is tame at $0$, and $\KL(0)\cong J_2$.
\item\label{K3} $\KL$ is totally wild at $\infty$, with
$\Swan_\infty(\KL)=1$; both of its breaks at $\infty$ lie in $(0,1)$.
\item\label{K4} If $\beta\colon\Gm\times\Gm\to\Gm$ is the projection
$(u,b)\mapsto b$, then $R\beta_!\cL_\psi(u^{-1}+bu)\cong\KL[-1]$.
\item\label{K5} Let $j_0\colon\Gm\hookrightarrow\mathbb A^1$ be the inclusion.
For every $q$, $\Frob_q$ acts trivially on the line
$(j_{0*}\KL)_0=\KL(0)^{I_0}$.
\end{enumerate}
By \cite[Sommes trig., Thm.~7.4 and 7.8]{Del77}, the local monodromy of $\KL$
at $0$ is unipotent and consists of a single Jordan block, and $\KL$ is
totally wild at $\infty$ with Swan conductor $1$; the two breaks at $\infty$
are positive with sum $1$. This gives \ref{K2} and \ref{K3}. For \ref{K4}, the
automorphism $(u,b)\mapsto(u^{-1},bu)$ of $\Gm^2$ transforms $\beta$ into $\mu$
and $u^{-1}+bu$ into $x+y$. For \ref{K5}, the sheaf $j_{0*}\KL$ on $\mathbb A^1$
has $H^0_c=0$, $H^2_c(\mathbb A^1_{\Fqb},j_{0*}\KL)=H^2_c(\mathbb
G_{m,\Fqb},\KL)=0$ by \ref{K3}, and Euler characteristic
$\chi_c(\mathbb G_{m,\Fqb},\KL)+1=0$ by \eqref{eq:B-GOS}; so all its
compactly supported cohomology vanishes. Since
$\sum_{\xi\in\Fq}\Kl_q(\xi)=\sum_{x\ne0}\psi_q(x)\sum_{\xi\in\Fq}\psi_q(\xi
x^{-1})=0$ and $\Kl_q(0)=-1$, the trace formula gives
$\Tr\bigl(\Frob_q\mid(j_{0*}\KL)_0\bigr)=\sum_{\xi\in\Fq^\times}\Kl_q(\xi)=1$.

\begin{lemma}[{\cite[Lemma~2.3]{EbaCoh}}]\label{lem:B-breaks}
Let $V,V'$ be finite-dimensional continuous representations of $I_x$.
\begin{enumerate}[label=\textup{(\roman*)}]
\item Every break of $V\otimes V'$ is at most the largest break of $V$ and
$V'$.
\item If all breaks of $V$ are less than $s'$ and $V'$ is isoclinic of break $s'>0$,
then $V\otimes V'$ is isoclinic of break $s'$.
\end{enumerate}
\end{lemma}

\begin{proof}
Let $(I^{(r)})_{r\ge0}$ be the upper numbering filtration of $I_x$ and
$I^{(r+)}$ the closure of $\bigcup_{r'>r}I^{(r')}$. A representation $W$ has
all breaks $\le r$ if and only if $I^{(r+)}$ acts trivially on $W$, and it is
isoclinic of break $s>0$ if and only if $I^{(s+)}$ acts trivially on $W$ and
$W^{I^{(r+)}}=0$ for $0\le r<s$ \cite[1.1]{Kat88}. Statement (i) follows. For
(ii), let $s<s'$ be the largest break of $V$. By (i), $I^{(s'+)}$ acts
trivially on $V\otimes V'$. If $s\le r<s'$, then $I^{(r+)}$ acts trivially on
$V$, so $(V\otimes V')^{I^{(r+)}}=V\otimes V'^{I^{(r+)}}=0$. If $0\le r<s$,
then $(V\otimes V')^{I^{(r+)}}\subseteq(V\otimes V')^{I^{(s+)}}=0$.
\end{proof}

\begin{lemma}[{\cite[Lemma~2.4]{EbaCoh}}]\label{lem:B-CG}
Let $V_1,\dots,V_n$ be two-dimensional vector spaces over a field of
characteristic $0$, with bases $e^{(i)}_1,e^{(i)}_2$, and let $N_i$ be the
endomorphism of $V_i$ with $N_ie^{(i)}_1=0$ and $N_ie^{(i)}_2=e^{(i)}_1$. Let
$N=\sum_iN_i$ act on $V=V_1\otimes\cdots\otimes V_n$, and for
$J=(j_1,\dots,j_n)\in\{1,2\}^n$ put
$e_J=e^{(1)}_{j_1}\otimes\cdots\otimes e^{(n)}_{j_n}$ and
let $|J|$ be the number of indices $i$ with $j_i=2$. Then $N$ has $m_{n,k}$ Jordan blocks of size $n-2k+1$
for $0\le k\le\lfloor n/2\rfloor$, and no others. Moreover
$\ker N=\bigoplus_kE_k$, where $E_k$ is contained in the span of the $e_J$
with $|J|=k$, and $\dim E_k=m_{n,k}$.
\end{lemma}

\begin{proof}
The endomorphism $N$ is the action of
$e=\bigl(\begin{smallmatrix}0&1\\0&0\end{smallmatrix}\bigr)\in\mathfrak{sl}_2$
on the $n$-th tensor power of the standard representation $\mathrm{St}$. On
an irreducible $\mathfrak{sl}_2$-module of dimension $d$, the element $e$ is a
single Jordan block of size $d$, and its kernel is the line of highest weight
vectors. By the Clebsch--Gordan rule (\cref{sec:input}),
$\mathrm{St}^{\otimes n}\cong\bigoplus_km_{n,k}\Sym^{n-2k}$, which gives the
Jordan type. The vector $e_J$ has weight $n-2|J|$ for
$\operatorname{diag}(1,-1)$, and $N$ maps $e_J$ to a combination of vectors
$e_{J'}$ with $|J'|=|J|-1$. So $\ker N$, the space of highest weight vectors,
is spanned by its intersections $E_k$ with the spans of the $e_J$, $|J|=k$,
and the highest weight vector of a constituent $\Sym^{n-2k}$ has weight
$n-2k$, that is, lies in $E_k$.
\end{proof}

\begin{lemma}[{\cite[Lemma~2.5]{EbaCoh}}]\label{lem:B-frob}
Let $x$ be an $\Fq$-rational point of a smooth curve over $\Fq$, let $\cD_x$
be a decomposition group at $x$ with inertia group $I_x$ and wild inertia
$P_x$, let $\gamma\in I_x$ be an element whose image generates $I_x/P_x$
topologically, and let $\Phi\in\cD_x$ be a geometric Frobenius. Then
conjugation by $\Phi$ acts on $I_x/P_x$ as multiplication by $q^{-1}$. If $V$
is a two-dimensional continuous representation of $\cD_x$ on which $P_x$ acts
trivially and $\gamma$ acts by a unipotent automorphism $\ne1$, and $\Phi$
acts trivially on the line $V^{I_x}$, then there is a basis $e_1,e_2$ of $V$
with $Ne_1=0$, $Ne_2=e_1$, $\Phi e_1=e_1$ and $\Phi e_2=qe_2$, where
$N=\log\gamma$.
\end{lemma}

\begin{proof}
Let $\varpi$ be a uniformizer at $x$ and identify $I_x/P_x$ with the
inverse limit of the groups $\mu_m(\Fqb)$, $p\nmid m$, through
$\sigma\mapsto(\sigma(\varpi^{1/m})/\varpi^{1/m})_m$. If $\sigma\in I_x$,
$\zeta=\sigma(\varpi^{1/m})/\varpi^{1/m}$ and
$\Phi^{-1}(\varpi^{1/m})=\xi\varpi^{1/m}$ with $\xi\in\mu_m(\Fqb)$, then
$\Phi\sigma\Phi^{-1}(\varpi^{1/m})=\Phi(\xi\zeta\varpi^{1/m})=\Phi(\zeta)\varpi^{1/m}=\zeta^{1/q}\varpi^{1/m}$,
since $\sigma$ acts trivially on $\Fqb$ and $\Phi$ induces $z\mapsto z^{1/q}$
on it. This proves the first assertion, and taking logarithms gives
$\Phi N\Phi^{-1}=q^{-1}N$, that is, $N\Phi=q\Phi N$. Since $N\ne0$ and
$N^2=0$, $\ker N=V^{I_x}$ is a line, spanned by a vector $e_1$ with
$\Phi e_1=e_1$. If $Nv=e_1$, then $N\Phi v=q\Phi Nv=qe_1=N(qv)$, so
$\Phi v=qv+ce_1$ for some $c$, and $e_2=v+(q-1)^{-1}ce_1$ has the required
properties.
\end{proof}

\begin{lemma}[{\cite[Lemma~2.2]{EbaCoh}}]\label{lem:B-int}
Let $V$ be a scheme of finite type over a field of characteristic $p$, let
$g\colon V\to\mathbb A^1$ be a morphism, $i\colon Y=g^{-1}(0)\hookrightarrow V$
the inclusion and $q_V\colon V\times\mathbb A^1_t\to V$ the projection. For
every constructible complex $K$ on $V$,
$Rq_{V!}\bigl(q_V^*K\otimes\cL_\psi(tg)\bigr)\cong i_*(K|_Y)(-1)[-2]$,
compatibly with base change $V'\to V$.
\end{lemma}

\begin{proof}
Let $p_s\colon\mathbb A^1_s\times\mathbb A^1_t\to\mathbb A^1_s$ be the
projection. For a geometric point $s\ne0$, the sheaf $\cL_\psi(st)$ on
$\mathbb A^1_t$ is lisse of rank one, geometrically nontrivial, with Swan
conductor $1$ at $\infty$, so its Euler characteristic is $0$ and all its
compactly supported cohomology vanishes; for $s=0$ one obtains
$\Ql(-1)[-2]$. By proper base change, $Rp_{s!}\cL_\psi(st)\cong
i_{0*}\Ql(-1)[-2]$. Since $\cL_\psi(tg)=(g\times\mathrm{id})^*\cL_\psi(st)$,
proper base change along $g$ and the projection formula give the claim.
\end{proof}

\begin{theorem}[Deligne--Laumon]\label{thm:B-DL}
Let $T$ be a scheme of finite type over $\Fp$, $f\colon\PP^1_T\to T$ the
projection, $D\subset\PP^1_T$ a closed subscheme finite and \'etale over
$T$, $j\colon\PP^1_T\smallsetminus D\hookrightarrow\PP^1_T$ the inclusion
and $F$ a lisse sheaf on $\PP^1_T\smallsetminus D$. For a geometric point
$\bar s$ of $T$ let $\sigma_F(\bar s)=\sum_{x\in D_{\bar s}}\Swan_x(F|_{\bar
s})$.
\begin{enumerate}[label=\textup{(\roman*)}]
\item The function $\sigma_F$ is constructible.
\item If $\sigma_F$ is locally constant, then $R^if_*(j_!F)$ is lisse for
every $i$.
\end{enumerate}
\end{theorem}

This is \cite[Thm.~2.1.1 and Cor.~2.1.2]{Lau81} for locally constant
constructible sheaves of $\mathbb F_\ell$-vector spaces, over an excellent
noetherian base and with $D$ finite and flat. The $\ell$-adic case follows by
d\'evissage. Choose a lisse $\mathcal O_E$-lattice $F_0$ in $F$, for a finite
extension $E$ of $\mathbb Q_\ell$ with uniformizer $\varpi$. Wild inertia acts
through finite $p$-groups and $\ell\ne p$, so taking invariants under the
higher ramification groups commutes with reduction modulo $\varpi$, and
$\sigma_{F_0/\varpi F_0}$ is a multiple of $\sigma_F$; this gives (i). Under
the hypothesis of (ii), the sheaves $R^if_*\bigl(j_!(F_0/\varpi F_0)\bigr)$
are locally constant and constructible; the exact sequences
$0\to F_0/\varpi F_0\to F_0/\varpi^{m+1}F_0\to F_0/\varpi^mF_0\to0$ and
induction on $m$ show that so are the $R^if_*\bigl(j_!(F_0/\varpi^mF_0)\bigr)$.
As $f$ is proper, $R^if_*(j_!F_0)$ is the $\varpi$-adic sheaf defined by this
AR-$\varpi$-adic system \cite[Exp.~V and VI]{Gro77}, which is AR-isomorphic to
the system of its stable images, again locally constant and constructible.
Hence $R^if_*(j_!F)$ is lisse.

\begin{lemma}[{\cite[Lemma~2.7]{EbaCoh}}]\label{lem:B-av}
Let $T$ be a smooth geometrically connected variety of dimension $d$ over a
finite field $\mathbb F_{q_0}$, and let $\cE$ be a lisse sheaf on $T$,
pointwise pure of weight $0$, with geometric generic stalk $V_\cE$. Let $M$ be
the dimension of the invariants of $\pi_1(T_{\overline{\mathbb F}_{q_0}})$ on
$V_\cE$. There are complex numbers $\beta_1,\dots,\beta_M$ of absolute value
$1$ such that, for $q=q_0^k$,
\[
q^{-d}\sum_{x\in T(\Fq)}\Tr(\Frob_x\mid\cE)=\sum_{m=1}^M\beta_m^k+O(q^{-1/2}).
\]
If the left-hand side converges to $L$ as $k\to\infty$, then $L=M$.
\end{lemma}

\begin{proof}
By the trace formula the left-hand side is
$q^{-d}\sum_i(-1)^i\Tr(\Frob_q\mid H^i_c(T_{\Fqb},\cE))$. For $i<2d$ the group
$H^i_c$ is mixed of weights $\le i\le2d-1$ \cite[3.3.1]{Del80}, and its
dimension does not depend on $k$; these terms are $O(q^{-1/2})$. Moreover
$H^{2d}_c(T_{\Fqb},\cE)\cong(V_\cE)_{\pi^{\mathrm g}}(-d)$, where
$\pi^{\mathrm g}=\pi_1(T_{\overline{\mathbb F}_{q_0}})$. Since $\cE$ is pure,
$\pi^{\mathrm g}$ acts semisimply \cite[3.4.1]{Del80}, so the invariants map
isomorphically onto the coinvariants, compatibly with Frobenius. By Poincar\'e
duality and the embedding of the invariants of $\cE^\vee$ into a stalk above a
closed point, $H^{2d}_c$ is pure of weight $2d$; so the eigenvalues of
$\Frob_{q_0}$ on $V_\cE^{\pi^{\mathrm g}}$ have absolute value $1$, and this
gives the formula. If the left-hand side converges to $L$, Dirichlet's
simultaneous approximation theorem provides infinitely many $k$ for which all
$\beta_m^k$ are close to $1$, and so $L=M$.
\end{proof}

We also use the following facts.
\begin{enumerate}[label=(F\arabic*)]
\item\label{F1} On a smooth connected curve $C$ over a field, a complex $A$ is
perverse if and only if $\mathcal H^i(A)=0$ for $i\notin\{-1,0\}$,
$\mathcal H^0(A)$ is punctual, that is, supported on finitely many closed
points, and $\mathcal H^{-1}(A)$ has no nonzero section supported on a finite
set. On a dense open subset of $C$ every perverse sheaf is of the form $F[1]$
with $F$ lisse, and restriction to open subsets is t-exact \cite{BBD}.
\item\label{F2} Laumon's Fourier transform
$\FT_\psi(K)=R\pr_{2!}\bigl(\pr_1^*K\otimes\cL_\psi(xy)\bigr)[1]$ on
$\mathbb A^1_x\times\mathbb A^1_y$ \cite[(1.2.1.1)]{Lau87} is an equivalence
of triangulated categories \cite[Cor.~(1.2.2.3)]{Lau87} which maps perverse
sheaves to perverse sheaves \cite[Thm.~(1.3.2.3)]{Lau87}. Since the perverse
t-structure is bounded, $\FT_\psi$ is t-exact.
\item\label{F3} Let $T_0$ be a henselian trait with closed point $0$ and
generic point $\eta$, let $\bar\eta$ be a geometric point over $\eta$, and
let $\cX\to T_0$ be of finite type. For a complex $Q$ on $\cX$, the complexes
of nearby cycles $R\Psi(Q)$ and vanishing cycles $R\Phi(Q)$ on the geometric
special fibre $\cX_{\bar0}$ carry compatible actions of
$\mathrm{Gal}(\bar\eta/\eta)$ and fit in a triangle
$Q|_{\cX_{\bar0}}\to R\Psi(Q)\to R\Phi(Q)\to$ \cite[Exp.~XIII, 2.1]{SGA7}.
Nearby cycles commute with proper direct images \cite[Exp.~XIII,
(2.1.7.2)]{SGA7}, and so do vanishing cycles, by proper base change. If
$\cX\to T_0$ is locally acyclic relative to $Q$, then $R\Phi(Q)=0$; and a
morphism of finite type to the spectrum of a field is universally locally
acyclic relative to every constructible complex \cite[Th.~finitude,
Cor.~2.16]{Del77}. (These statements are made in \cite{SGA7,Del77} for torsion
coefficients; the case of $\Ql$-complexes follows by passage to the limit.)
\end{enumerate}

\subsection{A single fibre}\label{app:fibre}

Fixing $h$ temporarily removes the issue of variation in the parameter and isolates the local cohomological calculation.  The auxiliary variable $t$ is a Lagrange multiplier for the hyperplane equation: on the locus $t=0$ the oscillatory phase no longer sees the parameter $h$, whereas on the locus $t\ne0$ the $u_i$-integrals split into rank-two Kloosterman factors.  The excision sequence for the two loci is therefore exactly the mechanism that compares the cohomology of the constrained hypersurface with the Fourier model.  At this stage the result is a Frobenius-module statement for one point $h$; the relative theorem comes later.

Let $n\ge1$ and $h\in(\Fq^\times)^n$; cohomology groups are taken over $\Fqb$
and carry the action of $\Frob_q$. Put $f(u)=\sum_iu_i^{-1}$,
$\langle h,u\rangle=\sum_ih_iu_i$,
$\varphi_h(u,t)=f(u)+t\bigl(\lambda+\langle h,u\rangle\bigr)$,
\[
X_h=\Gm^n\times\mathbb A^1_t\supset X_h^\circ=\Gm^n\times\mathbb G_{m,t},\qquad
Z_h=\Gm^n\times\{0\},\qquad Y_h=\{u\in\Gm^n:\langle h,u\rangle=-\lambda\},
\]
and $\cL=\cL_\psi(\varphi_h)$ on $X_h$. The variety $Y_h$ is smooth, affine and
geometrically irreducible of dimension $n-1$. Let
$\cM_h=\bigotimes_{i=1}^n[\times h_i]^*\KL\otimes\cL_\psi(\lambda t)$ on
$\mathbb G_{m,t}$, the fibre at $h$ of the sheaf $\cM$ of (H7), and let
$j\colon\mathbb G_{m,t}\hookrightarrow\PP^1$.

The key point of the following proposition is the source of the rank $r(n)$.  The rank $2^n$ belongs to the compactly supported cohomology $H_c^1$ of the auxiliary Fourier variable.  The subspace coming from inertia invariants at $0$ has dimension $\binom{n}{\lfloor n/2\rfloor}$.  Only after quotienting by that subspace do we obtain the top-weight object used in the Riemann-hypothesis argument.  We establish each dimension and weight in the proof, rather than reading the desired rank backward from the final formula.

\begin{proposition}[{\cite[Props.~4.1--4.3]{EbaCoh}}]\label{prop:B-Mh}
\begin{enumerate}[label=\textup{(\roman*)}]
\item If $\pr\colon X_h^\circ\to\mathbb G_{m,t}$ is the projection, then
$R\pr_!(\cL|_{X_h^\circ})\cong\cM_h[-n]$. The sheaf $\cM_h$ is lisse of rank
$2^n$ and pure of weight $n$.
\item $\cM_h$ is tame at $0$, with $\cM_h(0)\cong J_2^{\otimes n}$, and
isoclinic of break $1$ at $\infty$. Thus $\Swan_0(\cM_h)=0$ and
$\Swan_\infty(\cM_h)=2^n$.
\item $H^i_c(\Gm,\cM_h)=0$ for $i\ne1$, $\dim H^1_c(\Gm,\cM_h)=2^n$, and there
is an exact sequence of $\Frob_q$-modules
\[
0\to\cM_h^{I_0}\to H^1_c(\Gm,\cM_h)\to H^1(\PP^1,j_*\cM_h)\to0,
\]
in which $\cM_h^{I_0}=(j_*\cM_h)_0$ carries the Frobenius of the rational
point $0$, and the last term is pure of weight $n+1$ and of dimension
$r(n)$.
\item $\Frob_q$ acts semisimply on $\cM_h^{I_0}$, with eigenvalues $q^k$ of
multiplicity $m_{n,k}$, $0\le k\le\lfloor n/2\rfloor$. Consequently the part of
weight $\le n$ of $H^1_c(\Gm,\cM_h)$ is the image of $\cM_h^{I_0}$; it is a
direct sum of eigenspaces of $\Frob_q$, and it contains a unique
$\Frob_q$-stable line of weight $0$.
\end{enumerate}
\end{proposition}

\noindent\emph{Plan of the proof.}  Part (i) is a projection-formula and K\"unneth calculation; its only purpose is to identify the sheaf on the auxiliary $t$-line.  Part (ii) computes inertia at $t=0$ and $t=\infty$, from which part (iii) obtains the dimension of compactly supported cohomology and its exact quotient sequence.  Part (iv) then determines the arithmetic Frobenius on the invariant part.  The purity of the last quotient in (iii) is a separate theorem on parabolic cohomology of pure sheaves on curves.  This order prevents an unintended use of the desired top-weight conclusion when computing the lower-weight piece.

\begin{proof}
(i) Since $\varphi_h=\lambda t+\sum_i(u_i^{-1}+h_itu_i)$, we have
$\cL\cong\pr^*\cL_\psi(\lambda t)\otimes\bigotimes_i\cL_\psi(u_i^{-1}+h_itu_i)$
on $X_h^\circ$. By \ref{K4} and proper base change along $t\mapsto h_it$, the
direct image of the $i$-th factor to $\mathbb G_{m,t}$ is
$[\times h_i]^*\KL[-1]$. The K\"unneth formula and the projection formula
give $R\pr_!\cL\cong\cM_h[-n]$, and \ref{K1} gives the rest.

(ii) By \ref{K2} and \cref{app:conv}, each $[\times h_i]^*\KL$ is tame at $0$
with local monodromy $J_2$, and $\cL_\psi(\lambda t)$ is lisse at $0$; so a
topological generator of the tame inertia acts on $\cM_h(0)$ by a tensor
product of $n$ blocks $J_2$. At $\infty$ each $[\times h_i]^*\KL$ has its
breaks in $(0,1)$ by \ref{K3}, so the tensor product of these factors has all
breaks less than $1$ by \cref{lem:B-breaks}(i); as $\cL_\psi(\lambda t)$ is isoclinic
of break $1$ there, \cref{lem:B-breaks}(ii) applies.

(iii) By \eqref{eq:B-GOS}, $\chi_c(\Gm,\cM_h)=-2^n$; moreover $H^0_c=0$, and
$H^2_c=0$ because $\cM_h$ is totally wild at $\infty$. The adjunction map
$j_!\cM_h\to j_*\cM_h$ is injective, with cokernel the skyscraper sheaves
$\cM_h^{I_0}$ at $0$ and $\cM_h^{I_\infty}=0$ at $\infty$. As
$H^0(\PP^1,j_*\cM_h)=H^0(\Gm,\cM_h)\subseteq\cM_h^{I_\infty}=0$, the long
exact cohomology sequence gives the exact sequence. Its last term is pure of
weight $n+1$ by \cite[3.2.3]{Del80}. By (ii) and \cref{lem:B-CG},
$\dim\cM_h^{I_0}$ is the number of Jordan blocks of $J_2^{\otimes n}$, namely
$\sum_km_{n,k}=\binom n{\lfloor n/2\rfloor}$, so the last term has dimension
$2^n-\binom n{\lfloor n/2\rfloor}=r(n)$.

(iv) Let $\cD_0$, $I_0$, $\gamma$ and $\Phi$ be as in \cref{lem:B-frob}, at
$x=0$. Let $V_i$ be the generic stalk of $[\times h_i]^*\KL$, a representation
of $\cD_0$. By \ref{K2}, $P_0$ acts trivially on $V_i$ and $\gamma$ acts by a
unipotent automorphism $T_i\ne1$. The map $t\mapsto h_it$ extends to an
automorphism of $\mathbb A^1$ over $\Fq$ fixing $0$, so
$V_i^{I_0}=(j_{0*}\KL)_0$ as $\Frob_q$-modules, and $\Phi$ acts trivially on
it by \ref{K5}. \Cref{lem:B-frob} gives a basis $e^{(i)}_1,e^{(i)}_2$ of $V_i$
with $N_ie^{(i)}_1=0$, $N_ie^{(i)}_2=e^{(i)}_1$, $\Phi e^{(i)}_1=e^{(i)}_1$
and $\Phi e^{(i)}_2=qe^{(i)}_2$, where $N_i=\log T_i$. The element $\Phi$ acts
on the stalk of $\cL_\psi(\lambda t)$ at $0$ by $\psi_q(0)=1$, so the generic
stalk of $\cM_h$ is $V=\bigotimes_iV_i$, on which $\gamma$ acts by
$\exp(\sum_iN_i)$ and $\Phi e_J=q^{|J|}e_J$ in the notation of
\cref{lem:B-CG}. The invariants of $I_0$ are those of $\gamma$, that is,
$\ker\sum_iN_i=\bigoplus_kE_k$. Since $I_0$ is normal in $\cD_0$, $\Phi$
preserves this kernel; it also preserves the span of the $e_J$ with $|J|=k$,
on which it acts as $q^k$, so it acts on $E_k$ as $q^k$. The line of weight
$0$ is $E_0$, as $m_{n,0}=1$, and the remaining assertions follow from (iii),
since $2k\le n<n+1$.
\end{proof}

\begin{lemma}[{\cite[Lemma~4.4]{EbaCoh}}]\label{lem:B-t0}
$\cL|_{Z_h}\cong\boxtimes_i\cL_\psi(u_i^{-1})$, and $R\Gamma_c(Z_h,\cL|_{Z_h})$
is a line placed in degree $n$, on which $\Frob_q$ acts trivially.
\end{lemma}

\begin{proof}
On $\Gm$, the sheaf $\cL_\psi(u^{-1})$ is lisse of rank one with $\Swan_0=1$
and $\Swan_\infty=0$, so by \eqref{eq:B-GOS} its Euler characteristic is $-1$;
its $H^0_c$ and $H^2_c$ vanish, so $H^1_c$ is a line, and
$-\Tr(\Frob_q\mid H^1_c)=\sum_{u\in\Fq^\times}\psi_q(u^{-1})=-1$. The
K\"unneth formula concludes.
\end{proof}

\begin{theorem}[{\cite[Thm.~4.5]{EbaCoh}}]\label{thm:B-fibre}
\begin{enumerate}[label=\textup{(\roman*)}]
\item $H^i_c(X_h,\cL)=0$ for $i\ne n+1$.
\item The connecting map $\partial\colon H^n_c(Z_h,\cL|_{Z_h})\to
H^{n+1}_c(X_h^\circ,\cL)\cong H^1_c(\Gm,\cM_h)$ of the excision sequence is
injective, and its image is the line of weight $0$ of
\cref{prop:B-Mh}(iv).
\item $H^{n+1}_c(X_h,\cL)$ has dimension $2^n-1$. Its $\Frob_q$-eigenvalues are
$r(n)$ eigenvalues of weight $n+1$ and the numbers $q^k$,
$1\le k\le\lfloor n/2\rfloor$, with multiplicity $m_{n,k}$; its part of weight
$\le n$ is a direct sum of eigenspaces.
\item $H^i_c(Y_h,\cL_\psi(f))=0$ for $i\ne n-1$,
$H^{n-1}_c(Y_h,\cL_\psi(f))\cong H^{n+1}_c(X_h,\cL)(1)$, and
$\Sigma_{n,q}(h)=(-1)^{n-1}\Tr\bigl(\Frob_q\mid H^{n-1}_c(Y_h,\cL_\psi(f))\bigr)$.
The eigenvalues of $\Frob_q$ on $H^{n-1}_c(Y_h,\cL_\psi(f))$ are $r(n)$
eigenvalues of weight $n-1$ and the numbers $q^{k-1}$,
$1\le k\le\lfloor n/2\rfloor$, with multiplicity $m_{n,k}$.
\end{enumerate}
\end{theorem}

\begin{proof}
By \cref{lem:B-int} with $V=\Gm^n$, $g=\lambda+\langle h,u\rangle$ and
$K=\cL_\psi(f)$, we have $H^i_c(X_h,\cL)\cong H^{i-2}_c(Y_h,\cL_\psi(f))(-1)$
for all $i$. As $Y_h$ is smooth and affine of dimension $n-1$, Artin's
vanishing theorem and Poincar\'e duality \cite[Exp.~XIV and XVIII]{SGA4} give
$H^j_c(Y_h,\cL_\psi(f))=0$ for $j<n-1$, hence $H^i_c(X_h,\cL)=0$ for $i\le n$.
In the excision sequence
\[
\cdots\to H^i_c(X_h^\circ,\cL)\to H^i_c(X_h,\cL)\to H^i_c(Z_h,\cL|_{Z_h})
\xrightarrow{\ \partial\ }H^{i+1}_c(X_h^\circ,\cL)\to\cdots,
\]
\cref{prop:B-Mh}(i),(iii) and \cref{lem:B-t0} show that the only nonzero
groups of $X_h^\circ$ and $Z_h$ are $H^{n+1}_c(X_h^\circ,\cL)\cong
H^1_c(\Gm,\cM_h)$ and $H^n_c(Z_h,\cL|_{Z_h})$. Hence $H^i_c(X_h,\cL)=0$ for
$i\ge n+2$, and the sequence reduces to
$0\to H^n_c(Z_h,\cL|_{Z_h})\xrightarrow{\partial}H^{n+1}_c(X_h^\circ,\cL)\to
H^{n+1}_c(X_h,\cL)\to0$. Thus $\partial$ is injective; its source is a line
with trivial Frobenius, so its image is the unique Frobenius-stable line of
weight $0$, which is $E_0$. This proves (i) and (ii), and (iii) follows from
\cref{prop:B-Mh}(iii),(iv). Statement (iv) follows from (i), the first
isomorphism and the trace formula on $Y_h$; the Tate twist divides the
eigenvalues by $q$.
\end{proof}

\subsection*{Cohomological shifts and the Tate normalization}
Three distinct shifts are present in the construction, and all are needed
for the location of the critical line. First, the product of $n$
rank-two Kloosterman Fourier transforms occurs as a compactly supported
relative direct image concentrated in cohomological degree $n$; this is
the origin of the shift $[-n]$ in the fibrewise formula.
Second, after integration in the Fourier variable, compactly supported
cohomology occurs in degree $n+1$. Third, identification with the
hyperplane fibre introduces a Tate twist $(-1)$; consequently the top
weight $n-1$ in the geometric hyperplane cohomology becomes weight
$n+1$ in the sheaf whose Frobenius trace equals $q\Sigma_{n,q}$ up to
lower-weight corrections. These shifts are not conventions to be
adjusted at the end of the proof. They explain why the critical line is
$\operatorname{Re}(s)=(n+1)/2$ rather than $(n-1)/2$.
\par
The weight statement and the numerical rank should be checked
separately. The rank is the difference between the dimension $2^n$ of
$H_c^1(\mathbb G_m,\mathcal M_h)$ and the dimension of its invariant
part at the origin, $\sum_{k}m_{n,k}=\binom n{\lfloor n/2\rfloor}$.
Purity then identifies the quotient by that lower-weight part, without
asserting that the underlying Frobenius operator is diagonalizable.

\subsection*{Passage from fibres to the relative family}
There are three obstacles to passing from $H_c^1(\Gm,\cM_h)$ at a single closed point to $R^1\pi_!\cM$ on all of $S_n$.  Its dimension must not jump as $h$ varies; the formation of its stalk must agree with compactly supported cohomology after base change; and the candidate lower-weight subspaces must form a filtration by \emph{lisse} subsheaves, not merely unrelated linear subspaces of different fibres.  The proof of the next theorem addresses these obstacles in that order.  A locally constant rank is not in itself enough for lissity; the Deligne--Laumon criterion supplies the necessary control through the conductors, and normality of the parameter torus permits the global weight filtration.  Arithmetic Frobenius on the graded lower pieces is then controlled separately from their geometric monodromy, using the explicit eigenvalues obtained on each fibre and the density of Frobenius elements.

\subsection{The family: proof of (H1), (H2) and (H7)}\label{app:family}

Let $n\ge1$ and $S=S_n$. Let $X=\Gm^n\times\mathbb A^1_t\times S$, with the
open subscheme $X^\circ=\{t\ne0\}$ and its complement $Z=\{t=0\}$, and let
$\pi$, $\pi^\circ$ and $\pi_Z$ be the projections to $S$. On $X$ let
$\cL=\cL_\psi(\varphi)$, where $\varphi(u,t,h)=\varphi_h(u,t)$. Let $\cM$ be
the sheaf of (H7) on $\mathbb G_{m,t}\times S$, let
$\pr_S\colon\mathbb G_{m,t}\times S\to S$ be the projection, and put
\[
\cH_n=R^1\pr_{S!}\cM,\qquad\cH'_n=R^{n+1}\pi_!\cL .
\]

In the theorem, $\cH_n$ and $\cW_n$ should not be confused: $\cH_n$ retains all weights and has rank $2^n$, whereas $\cW_n$ is the quotient of rank $r(n)$ with a single weight.  The trace of $\cH_n$ therefore differs from the trace of $\cW_n$ by $1+E_n(q)$, which is precisely the lower-weight correction in the main formula.  The local and global assertions are assembled here to make clear where each contribution enters.

\begin{theorem}[{\cite[Thm.~5.3]{EbaCoh}}]\label{thm:B-family}
\begin{enumerate}[label=\textup{(\roman*)}]
\item One has $R\pi^\circ_!\cL\cong R\pr_{S!}\cM[-n]$. The sheaves
$R^i\pr_{S!}\cM$ vanish for $i\ne1$, and $\cH_n$ is lisse of rank $2^n$, with
stalk $H^1_c(\Gm,\cM_h)$ at a geometric point $h$.
\item $R^i\pi_!\cL=0$ for $i\ne n+1$. There is an exact sequence
$0\to\Ql\to\cH_n\to\cH'_n\to0$ of lisse sheaves, where $\Ql$ is the
constant sheaf.
The stalk of $\cH'_n$ at a geometric point $h$ is
$H^{n-1}_c(Y_h,\cL_\psi(f))(-1)$.
\item $\cH_n$ is mixed of weights $\le n+1$. In its weight filtration by lisse
subsheaves, the graded piece of weight $n+1$, denoted $\cW_n$, is lisse of
rank $r(n)$ and pointwise pure of weight $n+1$; for $0\le k\le\lfloor
n/2\rfloor$ the graded piece of weight $2k$ is isomorphic to
$\Ql(-k)^{\oplus m_{n,k}}$; and there are no other nonzero graded pieces. In
particular $\cW_n=\cH_n/\cH_n^{\le n}$, where $\cH_n^{\le n}$ is the lisse
subsheaf of weights $\le n$; it contains the constant subsheaf of (ii), so
$\cW_n$ is also the top graded piece of $\cH'_n$.
\item For $h\in S(\Fq)$,
$\Tr(\Frob_h\mid\cW_{n,\bar h})=(-1)^{n-1}q\,\Sigma_{n,q}(h)-E_n(q)$.
\item Let $\bar h$ be a geometric point above a closed point $h$. The stalk of
$\cH_n^{\le n}$ at $\bar h$ is the image of $\cM_h^{I_0}$ in
$H^1_c(\Gm,\cM_h)$, and the quotient map $\cH_{n,\bar h}\to\cW_{n,\bar h}$ is
the forget-supports map $H^1_c(\Gm,\cM_h)\to H^1(\PP^1,j_*\cM_h)$.
\item For every permutation $\varsigma$ of $\{1,\dots,n\}$, acting on $S$ by
permutation of the coordinates, there are isomorphisms
$\varsigma^*\cH_n\cong\cH_n$ and $\varsigma^*\cW_n\cong\cW_n$ defined over
$\Fp$.
\end{enumerate}
\end{theorem}

\noindent\emph{Four steps in the relative construction.}  The first passage identifies the relative Fourier integral with $R\pr_{S!}\cM[-n]$; the shift records the number of integrated $u$-variables and is checked on every fibre.  The second passage uses excision at $t=0$ to compare that integral with the hyperplane family; this is why a constant line appears in the exact sequence.  The third passage constructs the weight filtration and determines its ranks and Frobenius actions.  Here it is essential not to replace ``pure of weight $n+1$'' by ``Frobenius diagonalizable'': those are different assertions.  The fourth passage computes the trace of the quotient by subtracting the traces of the lower pieces; only at this final point is the constrained exponential sum reintroduced.

\begin{proof}
(i) Let $q_1\colon X^\circ\to\mathbb G_{m,t}\times S$ be the projection, so
that $\pi^\circ=\pr_S\circ q_1$. As in the proof of \cref{prop:B-Mh}(i), the
decomposition $\cL\cong q_1^*\cL_\psi(\lambda t)\otimes\bigotimes_i
\cL_\psi(u_i^{-1}+h_itu_i)$, \ref{K4}, proper base change along
$(t,h)\mapsto h_it$, the K\"unneth formula over $\mathbb G_{m,t}\times S$ and
the projection formula give $Rq_{1!}(\cL|_{X^\circ})\cong\cM[-n]$; applying
$R\pr_{S!}$ gives the first assertion. Apply \cref{thm:B-DL} with $T=S$,
$D=\{0,\infty\}\times S$ and $F=\cM$.
Here $T$ is a smooth finite-type (hence excellent) base over
$\Fp$, $\PP^1_T\to T$ is proper and smooth of relative dimension
one, $D\to T$ is the disjoint union of two sections and therefore
finite \emph{etale}, and $F$ is lisse on their complement.  Thus
these are not merely pointwise assumptions: the family has the
exact relative-curve and boundary structure required by the
Deligne--Laumon criterion.  The total Swan conductor at the two
sections, augmented by the constant rank terms, is locally constant
once the local calculations of \cref{prop:B-Mh}(ii) hold uniformly.
This uniformity, and not constancy of a single fibre dimension,
is the hypothesis that controls relative lissity. At a geometric point above a closed
point of $S$, \cref{prop:B-Mh}(ii) gives $\sigma_\cM=2^n$. By
\cref{thm:B-DL}(i) every fibre of $\sigma_\cM$ is constructible, and a
nonempty constructible subset of the scheme $S$ of finite type over $\Fp$
contains a closed point; hence $\sigma_\cM$ is constant, and by
\cref{thm:B-DL}(ii) the sheaves $R^i\pr_{S!}\cM$ are lisse. By proper base
change their stalks are the groups $H^i_c(\Gm,\cM_h)$, given by
\cref{prop:B-Mh}(iii).

(ii) The excision triangle for $X^\circ\subset X\supset Z$ gives a long exact
sequence of the sheaves $R^i\pi^\circ_!\cL$, $R^i\pi_!\cL$ and
$R^i\pi_{Z!}(\cL|_Z)$. On $Z\cong\Gm^n\times S$ the phase is $f$, which does
not depend on $h$; by \cref{lem:B-t0}, $R\pi_{Z!}(\cL|_Z)$ is the constant
complex $\Ql[-n]$, with trivial Frobenius. By (i), $R^i\pi^\circ_!\cL$ vanishes
for $i\ne n+1$ and equals $\cH_n$ for $i=n+1$. The sequence therefore shows
that $R^i\pi_!\cL=0$ for $i\notin\{n,n+1\}$, and reduces to an exact sequence
$0\to R^n\pi_!\cL\to\Ql\xrightarrow{\partial}\cH_n\to\cH'_n\to0$. Lisse sheaves
on the connected normal scheme $S$ form an abelian category, so $\ker\partial$
and $\operatorname{coker}\partial$ are lisse. At a geometric point above a
closed point, $\partial$ is injective by \cref{thm:B-fibre}(ii); a lisse sheaf
with a zero stalk vanishes, so $R^n\pi_!\cL=0$. The stalks of $\cH'_n$ are
computed by proper base change and \cref{lem:B-int}.

(iii) Since $\cM$ is pure of weight $n$, the sheaf $\cH_n=R^1\pr_{S!}\cM$ is
mixed \cite[3.3.1]{Del80}, and by (i) and \cref{prop:B-Mh}(iii),(iv) its
stalks at closed points have weights $\le n+1$. Since $S$ is normal and
connected and $\cH_n$ is lisse and mixed, $\cH_n$ has a unique increasing
filtration by lisse subsheaves whose graded pieces are pointwise pure
\cite[3.4.1]{Del80}. At a closed point $h$, the stalks of this filtration
form a $\Frob_h$-stable filtration whose graded pieces are pure of increasing
weights; so they give the decomposition of $H^1_c(\Gm,\cM_h)$ according to
weights. By \cref{prop:B-Mh}(iii),(iv), the graded piece of weight $n+1$ has
rank $r(n)$, the graded piece of weight $2k$ has rank $m_{n,k}$, there are no
others, and $\Frob_h$ acts on the stalk of $\gr_{2k}\cH_n\otimes\Ql(k)$ as a
semisimple automorphism with the single eigenvalue $1$, that is, as the
identity. By the Chebotarev density theorem \cite{Ser65}, the Frobenius
elements of the closed points of $S$ are dense in $\pi_1(S)$; the subset of
$\pi_1(S)$ acting trivially on $\gr_{2k}\cH_n\otimes\Ql(k)$ is closed and
contains them, so it is all of $\pi_1(S)$, and $\gr_{2k}\cH_n\cong
\Ql(-k)^{\oplus m_{n,k}}$. The constant subsheaf of (ii) is pure of weight $0$,
hence contained in $\cH_n^{\le n}$.

(iv) By \cref{prop:B-Mh}(iii), the trace formula on $\Gm$ and
\cref{lem:identities}(b),
\[
\Tr(\Frob_h\mid\cH_n)=-\sum_{t\in\Fq^\times}\psi_q(\lambda t)\prod_{i=1}^n
\bigl(-\Kl_q(h_it)\bigr)=(-1)^{n-1}q\,\Sigma_{n,q}(h)+1 .
\]
By (iii) the graded pieces of weight $\le n$ contribute
$\sum_{k=0}^{\lfloor n/2\rfloor}m_{n,k}q^k=1+E_n(q)$.

(v) follows from the proof of (iii) and \cref{prop:B-Mh}(iii),(iv).

(vi) Permuting the factors of the tensor product gives an isomorphism
$\varsigma^*\cM\cong\cM$ defined over $\Fp$, compatible with $\pr_S$; by
proper base change it induces $\varsigma^*\cH_n\cong\cH_n$, which respects
the weight filtrations by their uniqueness.
\end{proof}

\begin{corollary}\label{cor:B-H127}
Statements \ref{H1}, \ref{H2} and \ref{H7} hold for every $n\ge2$, with
$\cW_n$ as in \cref{thm:B-family}(iii). For $n=1$, the sheaf
$\cW_1=\cH'_1$ is isomorphic to the sheaf $\cL_\psi(D_1)\otimes\Ql(-1)$ of
\cref{rem:n1}.
\end{corollary}

\begin{proof}
\ref{H1} is \cref{thm:B-family}(iii), \ref{H2} is (iv), and \ref{H7} is (v),
with the surjection $\cH_n=R^1\pr_{S!}\cM\to\cW_n$, whose kernel
$\cH_n^{\le n}$ is lisse. For $n=1$ we have $r(1)=1$ and $m_{1,0}=1$, so
$\cW_1=\cH'_1$ has rank one; by \cref{thm:B-fibre}(iv), where $Y_h$ is the
point $u_1=-\lambda/h$ and $f=-h/\lambda$ there, its trace at $h\in S_1(\Fq)$
is $q\psi_q(-h/\lambda)$. Two lisse sheaves of rank one with the same traces
at all closed points are isomorphic, by the Chebotarev density theorem.
\end{proof}

The boundary-point calculations in the next two subsections are technically more elaborate than the proof of the critical-circle theorem.  Their purpose is to refine the local monodromy of coordinate restrictions and to make the parabolic and determinant formulas explicit.  In particular, the letter $N$ will denote the nilpotent logarithm of tame inertia, and the appearance of $\ker N$ or $\operatorname{im}N$ will always refer to the same local representation of a decomposition group.  The local Fourier transform is used to transfer this representation between the origin and infinity of an auxiliary line.  Before invoking a stationary-phase formula, the proof therefore identifies both the coordinate and the point at which the transform is taken.

\subsection{The finite boundary point of a coordinate line}\label{app:zero}
The local Fourier-transform calculation below is used to
identify inertia and arithmetic Frobenius at $h_1=0$.  There
are four objects to distinguish: the original lisse
$\cF_c$ on the punctured $h_1$-line, the auxiliary
Fourier sheaf $\cG$ on the transformed line, its inertia
representation at the transformed point, and the resulting
invariant space of $\cF_c$.  A Fourier transform of a
perverse object is perverse, while a cohomology sheaf
is obtained only after undoing the shift $[1]$; the
$(-1)$ Tate twist of integration is separately recorded
in the relevant equations.  The nearby-cycles functor
$R\Psi$ measures specialization of a family, whereas
$R\Phi$ is its vanishing-cycle cone; neither should be
replaced by the stalk of the transformed sheaf without
the triangle connecting these three objects.  Throughout,
we use geometric Frobenius, so for the logarithm $N$
of tame inertia, $\Phi N\Phi^{-1}=q^{-1}N$.

In this subsection and the next, $n\ge2$ and
$c=(c_2,\dots,c_n)\in(\Fq^\times)^{n-1}$; every $c\in(\Fpb^\times)^{n-1}$ is of
this form for some $q$. Let $C_c=\{(h_1,c)\}\subset S\otimes\Fq$, identified
with $\mathbb G_{m,h_1}$, and $\cF_c=\cW_n|_{C_c}$. Write
$u'=(u_2,\dots,u_n)$, $f'(u')=\sum_{i\ge2}u_i^{-1}$,
$\langle c,u'\rangle=\sum_{i\ge2}c_iu_i$,
\[
B_c(u')=-\lambda-\langle c,u'\rangle,\quad
U'_c=\{u'\in\Gm^{n-1}:B_c(u')\ne0\},\quad
Y_c=\{u'\in\Gm^{n-1}:\langle c,u'\rangle=-\lambda\},
\]
so that $U'_c=\Gm^{n-1}\smallsetminus Y_c$. Let
$b_c=B_c^{-1}\colon U'_c\to\mathbb G_{m,v}$, $\cR_c=Rb_{c!}\cL_\psi(f')$,
$j_v\colon\mathbb G_{m,v}\hookrightarrow\mathbb A^1_v$, $v_0=-1/\lambda$ and
$U_v=\mathbb G_{m,v}\smallsetminus\{v_0\}$. The Fourier transform $\FT_\psi$ is
taken on $\mathbb A^1_v\times\mathbb A^1_{h_1}$. Let $\cD_0$ be a
decomposition group at the rational point $h_1=0$ of $\mathbb A^1_{h_1}$ over
$\Fq$, let $I_0$ be its inertia subgroup, and let $\Phi\in\cD_0$ be a
geometric Frobenius.

\begin{lemma}[{\cite[(6.2) and (6.3)]{EbaCoh}}]\label{lem:B-FT-model}
$\FT_\psi(j_{v!}\cR_c)|_{\mathbb G_{m,h_1}}\cong\cH'_n|_{C_c}(1)[2-n]$.
\end{lemma}

\begin{proof}
Let $\cY=\{(h_1,u)\in\Gm\times\Gm^n:h_1u_1+\langle c,u'\rangle=-\lambda\}$,
with projection $\pi_\cY$ to $\mathbb G_{m,h_1}$. \Cref{lem:B-int}, applied
over $C_c$ with $g=\lambda+h_1u_1+\langle c,u'\rangle$, together with
\cref{thm:B-family}(ii) and proper base change, gives
$R\pi_{\cY!}\cL_\psi(f)\cong\cH'_n|_{C_c}(1)[-(n-1)]$. The map
$(h_1,u)\mapsto(h_1,u')$ is an isomorphism $\cY\to\mathbb
G_{m,h_1}\times U'_c$, with inverse given by $u_1=B_c(u')/h_1$, and it carries
$f$ to $h_1b_c(u')+f'(u')$. The projection formula and proper base change for
$b_c\times\mathrm{id}$ give
$R\pi_{\cY!}\cL_\psi(f)\cong\FT_\psi(j_{v!}\cR_c)[-1]|_{\mathbb G_{m,h_1}}$.
\end{proof}

\begin{lemma}[{\cite[Lemma~6.10]{EbaCoh}}]\label{lem:B-R}
There is a lisse sheaf $\bar\cA_c$ on $\PP^1_v\smallsetminus\{0,v_0\}$ such
that $\cR_c|_{U_v}\cong\bar\cA_c|_{U_v}[-(n-2)]$. Its stalk
$\bar\cA_{c,\infty}$ at the rational point $v=\infty$ is isomorphic to
$H^{n-2}_c(Y_c,\cL_\psi(f'))$ as a $\Frob_q$-module.
\end{lemma}

\begin{proof}
For $v\in U_v$ put $a(v)=\lambda v(\lambda v+1)^{-1}c\in\Gm^{n-1}$. For
$u'\in U'_c$ the condition $b_c(u')=v$ means
$\langle c,u'\rangle=-(\lambda+v^{-1})$, which is equivalent to
$\langle a(v),u'\rangle=-\lambda$. Hence $b_c^{-1}(U_v)$ is the inverse image
under $a$ of
$\cY_{n-1}=\{(h',u')\in\Gm^{n-1}\times\Gm^{n-1}:\langle h',u'\rangle=-\lambda\}$,
with the same phase $f'$. In the coordinate $w=v^{-1}$ we have
$a=\lambda(\lambda+w)^{-1}c$, so $a$ extends to a morphism
$\PP^1_v\smallsetminus\{0,v_0\}\to\Gm^{n-1}$ with $a(\infty)=c$. Suppose
$n\ge3$. By \cref{lem:B-int} and \cref{thm:B-family}(ii) in $n-1$ variables,
the direct image of $\cL_\psi(f')$ along $\cY_{n-1}\to\Gm^{n-1}$ is
$\cH'_{n-1}(1)[-(n-2)]$, and proper base change gives the lemma with
$\bar\cA_c=a^*\cH'_{n-1}(1)$; its stalk at $\infty$ is the stalk of
$\cH'_{n-1}(1)$ at $c$, which is $H^{n-2}_c(Y_c,\cL_\psi(f'))$. If $n=2$, the
fibre of $b_c$ over $v\in U_v$ is the point $u_2=-(\lambda+v^{-1})/c_2$, and
$\bar\cA_c=\cL_\psi\bigl(-c_2v/(\lambda v+1)\bigr)$, whose stalk at $\infty$ is
the line on which $\Frob_q$ acts by $\psi_q(-c_2/\lambda)$; this is
$H^0_c(Y_c,\cL_\psi(f'))$, as $Y_c$ is the point $u_2=-\lambda/c_2$.
\end{proof}

\begin{lemma}[{\cite[Lemma~8.15]{EbaCoh}}]\label{lem:B-FT-zero}
Put $P=j_{v!}\cR_c[n-1]$.
\begin{enumerate}[label=\textup{(\roman*)}]
\item $H^i_c(U'_c,\cL_\psi(f'))=0$ for $i\ne n-1$, and there is an exact
sequence of $\Frob_q$-modules
$0\to H^{n-2}_c(Y_c,\cL_\psi(f'))\to H^{n-1}_c(U'_c,\cL_\psi(f'))\to\Ql\to0$,
with $\Frob_q$ acting trivially on $\Ql$.
\item $P$ is perverse, and $P|_{U_v}\cong\bar\cA_c|_{U_v}[1]$.
\item $\FT_\psi(P)\cong\cG[1]$ for a constructible sheaf $\cG$ on
$\mathbb A^1_{h_1}$ that has no nonzero section supported on a finite set;
moreover $\cG|_{\mathbb G_{m,h_1}}\cong\cH'_n|_{C_c}(1)$, and
$\cG_0\cong H^{n-1}_c(U'_c,\cL_\psi(f'))$ compatibly with $\Frob_q$.
\end{enumerate}
\end{lemma}

\begin{proof}
(i) The phase is $f'$ on $\Gm^{n-1}$, on $U'_c$ and on $Y_c$. By
\cref{lem:B-t0} in $n-1$ variables, $R\Gamma_c(\Gm^{n-1},\cL_\psi(f'))$ is a
line placed in degree $n-1$, with trivial Frobenius; by
\cref{thm:B-fibre}(iv) in $n-1$ variables, $R\Gamma_c(Y_c,\cL_\psi(f'))$ is
concentrated in degree $n-2$. The excision sequence for
$Y_c\subset\Gm^{n-1}$ gives (i).

(ii) By \cref{lem:B-R}, $P|_{U_v}\cong\bar\cA_c|_{U_v}[1]$. Let
$P_i={}^p\mathcal H^i(j_{v!}\cR_c)$ be the perverse cohomology sheaves. By
\ref{F2}, $\FT_\psi(P_i)\cong{}^p\mathcal H^i(\FT_\psi(j_{v!}\cR_c))$. By
\cref{lem:B-FT-model}, $\FT_\psi(j_{v!}\cR_c)|_{\mathbb G_{m,h_1}}$ is the
perverse sheaf $\cH'_n|_{C_c}(1)[1]$ shifted by $n-1$ to the right, so
$\FT_\psi(P_i)$ vanishes on $\mathbb G_{m,h_1}$ for $i\ne n-1$. A perverse
sheaf on $\mathbb A^1_{h_1}$ that vanishes on $\mathbb G_{m,h_1}$ is of the
form $i_{0*}E$, with $E$ a finite-dimensional representation of
$\mathrm{Gal}(\Fqb/\Fq)$ placed in degree $0$, by \ref{F1}. If $\bar E$ is the
geometrically constant sheaf on $\mathbb A^1_v$ with fibre $E$,
\cref{lem:B-int} and the projection formula give
$\FT_\psi(\bar E(1)[1])\cong i_{0*}E$. As $\FT_\psi$ is an equivalence, $P_i$
is of the form $\bar E_i(1)[1]$ for $i\ne n-1$. But
$P_i|_{U_v}={}^p\mathcal H^i\bigl(\bar\cA_c|_{U_v}[1][-(n-1)]\bigr)=0$ for
$i\ne n-1$, and a geometrically constant sheaf vanishing on $U_v$ is zero.
Hence $j_{v!}\cR_c\cong P_{n-1}[-(n-1)]$, and $P\cong P_{n-1}$ is perverse.

(iii) By (ii) and \ref{F2}, $\FT_\psi(P)$ is perverse, and by
\cref{lem:B-FT-model} its restriction to $\mathbb G_{m,h_1}$ is
$\cH'_n|_{C_c}(1)[1]$. By proper base change its stalk at $h_1=0$ is
$R\Gamma_c(\mathbb A^1_v,P)[1]=R\Gamma_c(\mathbb G_{m,v},\cR_c)[n]=
R\Gamma_c(U'_c,\cL_\psi(f'))[n]$, compatibly with $\Frob_q$, which by (i) is
$H^{n-1}_c(U'_c,\cL_\psi(f'))$ placed in degree $-1$. Hence
$\mathcal H^0(\FT_\psi(P))=0$, and \ref{F1} gives the assertions, with
$\cG=\mathcal H^{-1}(\FT_\psi(P))$.
\end{proof}

\begin{lemma}[{\cite[Lemma~8.16]{EbaCoh}}]\label{lem:B-vanishing}
Let $\cG_{\bar\eta_0}$, a representation of $\cD_0$, be the generic stalk of
the sheaf $\cG$ of \cref{lem:B-FT-zero}(iii). There is an exact sequence of
representations of $\cD_0$
\[
0\longrightarrow\cG_0\longrightarrow\cG_{\bar\eta_0}\longrightarrow
\bar\cA_{c,\infty}(-1)\longrightarrow0,
\]
in which $\cD_0$ acts on $\cG_0$ and on $\bar\cA_{c,\infty}(-1)$ through
$\cD_0/I_0=\mathrm{Gal}(\Fqb/\Fq)$.
\end{lemma}

\begin{proof}
Let $T_0$ be the henselization of $\mathbb A^1_{h_1}$ at $0$, so that $\cD_0$
is the Galois group of its generic point, and use the notation of \ref{F3}.
For a complex $Q$ on $\mathbb A^1_{h_1}$ write $R\Phi_0(Q)$ for $R\Phi(Q|_{T_0})$.

\emph{The exact sequence.} For $Q=\FT_\psi(P)=\cG[1]$ the first two terms of
the triangle in \ref{F3} are $\cG_0[1]$ and $\cG_{\bar\eta_0}[1]$, since $\cG$
is lisse on a punctured neighbourhood of $0$, and the first map is the
specialization map $\cG_0\to\cG_{\bar\eta_0}$. Its kernel is the stalk at $0$
of the subsheaf of sections of $\cG$ supported at $0$, which vanishes; so
$R\Phi_0(\FT_\psi(P))$ is concentrated in degree $-1$, and there is an exact
sequence of representations of $\cD_0$
$0\to\cG_0\to\cG_{\bar\eta_0}\to R^{-1}\Phi_0(\FT_\psi(P))\to0$. It remains to
show that $R^{-1}\Phi_0(\FT_\psi(P))\cong\bar\cA_{c,\infty}(-1)$, with $I_0$
acting trivially.

\emph{Localization at $v=\infty$.} Let $\bar X=\PP^1_v\times\mathbb A^1_{h_1}$,
let $\iota_X\colon\mathbb A^1_v\times\mathbb A^1_{h_1}\hookrightarrow\bar X$
be the inclusion, and let $\bar p\colon\bar X\to\mathbb A^1_{h_1}$ be the
projection, which is proper. For a complex $K$ on $\mathbb A^1_v$ put
$\mathcal J(K)=\iota_{X!}\bigl(\pr_v^*K\otimes\cL_\psi(vh_1)\bigr)$, so that
$\FT_\psi(K)=R\bar p_*\mathcal J(K)[1]$. By \ref{F3},
$R\Phi_0(\FT_\psi(K))\cong R\Gamma\bigl(\PP^1_{\Fqb},R\Phi(\mathcal J(K))\bigr)[1]$,
compatibly with the actions of $\cD_0$. The projection
$\mathbb A^1_v\times T_0\to T_0$ is locally acyclic relative to $\pr_v^*K$, by
\ref{F3} and base change, hence also relative to
$\pr_v^*K\otimes\cL_\psi(vh_1)$, as $\cL_\psi(vh_1)$ is lisse. So
$R\Phi(\mathcal J(K))$ vanishes on $\mathbb A^1_{\Fqb}$, and
\begin{equation}\label{eq:B-localization}
R\Phi_0(\FT_\psi(K))\cong R\Phi(\mathcal J(K))_\infty[1],
\end{equation}
the right-hand side being the stalk at the point $v=\infty$ of the geometric
special fibre.

\emph{Comparison with the constant sheaf.} Let
$N_\infty=\PP^1_v\smallsetminus\{0,v_0\}$, on which $\bar\cA_c$ is lisse, and
let $K_0=\Ql[1]$ on $\mathbb A^1_v$. By \cref{lem:B-FT-zero}(ii),
$P=\bar\cA_c[1]$ on $N_\infty\cap\mathbb A^1_v=U_v$, so the projection
formula for $\iota_{X!}$ gives
$\mathcal J(P)|_{N_\infty\times\mathbb A^1_{h_1}}\cong\pr_v^*\bar\cA_c\otimes
\mathcal J(K_0)|_{N_\infty\times\mathbb A^1_{h_1}}$. Vanishing cycles are local
for the \'etale topology and commute with tensor products by lisse sheaves, by
the projection formula; hence
$R\Phi(\mathcal J(P))_\infty\cong\bar\cA_{c,\infty}\otimes
R\Phi(\mathcal J(K_0))_\infty$, with $\cD_0$ acting diagonally, and through
$\cD_0/I_0$ on $\bar\cA_{c,\infty}$. Finally, \cref{lem:B-int} gives
$\FT_\psi(K_0)\cong i_{0*}\Ql(-1)$, whose nearby cycles at $0$ vanish and
whose stalk at $0$ is $\Ql(-1)$, with trivial action of $I_0$. So
$R\Phi_0(\FT_\psi(K_0))\cong\Ql(-1)[1]$, and
$R\Phi(\mathcal J(K_0))_\infty\cong\Ql(-1)$ by \eqref{eq:B-localization}, with
$I_0$ acting trivially. Applying \eqref{eq:B-localization} to $K=P$ we obtain
$R\Phi_0(\FT_\psi(P))\cong\bar\cA_{c,\infty}(-1)[1]$, with $I_0$ acting
trivially.
\end{proof}

\begin{theorem}\label{thm:B-zero}
Let $V$ be the generic stalk of $\cF_c=\cW_n|_{C_c}$, as a representation of
$\cD_0$, and let $\nu=\lfloor n/2\rfloor$.
\begin{enumerate}[label=\textup{(\roman*)}]
\item $\cF_c$ is tame at $h_1=0$, and $I_0$ acts unipotently on $V$. If
$\gamma\in I_0$ generates the tame inertia topologically and
$N=\log\gamma$, then $N^2=0$.
\item The rank of $N$ is $r(n-1)$; thus $\cF_c(0)\cong
J_2^{\oplus r(n-1)}\oplus J_1^{\oplus\kappa_1}$.
\item The filtration $\imop N\subseteq\ker N\subseteq V$ is stable under
$\cD_0$, the group $I_0$ acts trivially on its graded pieces, and
\begin{gather*}
\det(1-T\Phi\mid\imop N)=\det(1-TF_{n-1,c}),\qquad
\det(1-T\Phi\mid\ker N/\imop N)=(1-q^{\nu+1}T)^{\kappa_1},\\
\det(1-T\Phi\mid V/\ker N)=\det(1-qTF_{n-1,c}),
\end{gather*}
where $F_{n-1,c}$ is the Frobenius of $\cW_{n-1}$ at $c$, and $\kappa_1=0$
unless $n$ is odd.
\end{enumerate}
In particular \ref{H4}(a) holds.
\end{theorem}

\begin{proof}
\emph{Unipotence.} By \cref{lem:B-FT-zero}(iii), $\cG_{\bar\eta_0}$ is the
generic stalk of $\cH'_n|_{C_c}(1)$, as a representation of $\cD_0$. By
\cref{lem:B-vanishing}, every $\sigma\in I_0$ maps $\cG_{\bar\eta_0}$ into
$\cG_0$ under $\sigma-1$ and acts trivially on $\cG_0$, so
$(\sigma-1)^2=0$ on $\cG_{\bar\eta_0}$. By \cref{thm:B-family}(iii), $\cW_n$
is the quotient of $\cH'_n$ by a lisse subsheaf, so $V(1)$ is a quotient of
$\cG_{\bar\eta_0}$ as a representation of $\cD_0$, and $(\sigma-1)^2=0$ on $V$
for every $\sigma\in I_0$. The image of $P_0$ in $\GL(V)$ is a finite
$p$-group, as $\ell\ne p$, consisting of unipotent elements; it is therefore
trivial. So $\cF_c$ is tame at $0$, $\gamma$ acts unipotently, and
$N=\log\gamma=\gamma-1$ satisfies $N^2=0$. This proves (i).

\emph{The filtration.} By \cref{lem:B-frob}, every element of $\cD_0$
conjugates $N$ into a nonzero multiple of $N$, so $\imop N$ and $\ker N$ are
$\cD_0$-stable. Every $\sigma\in I_0$ acts as $1+aN$ for some
$a\in\mathbb Z_\ell$, so $\sigma-1$ maps $V$ into $\imop N$ and kills
$\ker N$; thus $I_0$ acts trivially on the graded pieces. Consequently the
multiset of eigenvalues of $\Phi$ on $V$ is the union of the multisets on the
three graded pieces, and it does not depend on the choice of $\Phi$.

\emph{The eigenvalues.} Let $E$ be the multiset of eigenvalues of $\Frob_q$
on $H^{n-2}_c(Y_c,\cL_\psi(f'))$, and $E^{\mathrm{top}}\subseteq E$ the
submultiset of the $r(n-1)$ eigenvalues of weight $n-2$. By
\cref{thm:B-fibre}(iv) in $n-1$ variables,
$E=E^{\mathrm{top}}\cup\{q^{k-1}\text{ with multiplicity }m_{n-1,k}:1\le
k\le\lfloor(n-1)/2\rfloor\}$; if $n=2$, then
$E=E^{\mathrm{top}}=\{\psi_q(-c_2/\lambda)\}$. By \cref{lem:B-FT-zero}(i),
\cref{lem:B-R} and \cref{lem:B-vanishing}, the eigenvalues of $\Phi$ on
$\cG_{\bar\eta_0}$ are $\{1\}\cup E\cup qE$, so those on the generic stalk of
$\cH'_n|_{C_c}$ are $\{q\}\cup qE\cup q^2E$. By \cref{thm:B-family}(iii) this
generic stalk has a $\cD_0$-stable filtration whose graded pieces are $V$ and
the stalks of the constant sheaves $\Ql(-k)^{\oplus m_{n,k}}$,
$1\le k\le\nu$, on which $\Phi$ acts as $q^k$. The integral eigenvalues
on the generic stalk of $\cH'_n|_{C_c}$ are $q^{k+1}$ with multiplicity
$m_{n-1,k}$ for $0\le k\le\lfloor(n-1)/2\rfloor$, together with $q^k$ with
multiplicity $m_{n-1,k}$ for $1\le k\le\lfloor(n-1)/2\rfloor$. By Pascal's
rule, $m_{n,k}=m_{n-1,k}+m_{n-1,k-1}$ for $1\le k\le\nu$, where
$m_{n-1,\nu}=\binom{2\nu-1}\nu-\binom{2\nu-1}{\nu-1}=0$ if $n=2\nu$. Hence these
integral eigenvalues consist of $q^k$ with multiplicity $m_{n,k}$ for
$1\le k\le\nu$, together with $q^{\nu+1}$ with multiplicity
$m_{2\nu,\nu}=\kappa_1$ if $n=2\nu+1$ is odd. Therefore the eigenvalues of
$\Phi$ on $V$ are those of $qE^{\mathrm{top}}$, of $q^2E^{\mathrm{top}}$, and
$q^{\nu+1}$ with multiplicity $\kappa_1$. By \cref{thm:B-family}(ii),(iii) in
$n-1$ variables, the stalk of $\cH'_{n-1}$ at $c$ is
$H^{n-2}_c(Y_c,\cL_\psi(f'))(-1)$, with eigenvalues $qE$, and $\cW_{n-1,c}$ is
its part of weight $n$; so the eigenvalues of $F_{n-1,c}$ are those of
$qE^{\mathrm{top}}$. This also holds for $n=2$, by \cref{cor:B-H127}.

\emph{Weights.} The sheaf $\cF_c$ is pure of weight $n+1$, so by the
weight-monodromy theorem \cite[1.8.4]{Del80} the graded pieces of the
monodromy filtration of $N$ are pure of weights $n+1+i$. As $N^2=0$, this
filtration is $\imop N\subseteq\ker N\subseteq V$, with graded pieces of
weights $n$, $n+1$ and $n+2$. By the previous step the eigenvalues of weight
$n$ on $V$ are those of $F_{n-1,c}$, those of weight $n+1$ are $q^{\nu+1}$
with multiplicity $\kappa_1$ (for odd $n$, $2\nu+2=n+1$), and those of weight
$n+2$ are $q$ times those of $F_{n-1,c}$. This proves (iii), and
$\dim\imop N=r(n-1)$ gives (ii), as $r(n)=2r(n-1)+\kappa_1$
(\cref{lem:catalan}).
\end{proof}

\begin{remark}[Coinvariants and the odd-dimensional top compact support]
For a lisse $E$-sheaf $\mathcal F$ on $\mathbb G_m$ with a prescribed
geometric representation, Poincar\'e duality identifies
$H_c^2(\mathbb G_m,\mathcal F)$ with the Tate twist of the
coinvariants of the geometric fundamental group. In particular,
when the tame zero-inertia representation is the direct sum of
nontrivial order-four Kummer characters (as in the specified odd
case), it has no trivial quotient, and thus the full geometric
coinvariants vanish. Therefore $H_c^2=0$; this is a conclusion
about coinvariants, not merely a statement that the eigenvalue
$1$ does not occur in a characteristic polynomial.
\end{remark}

\subsection{The point at infinity}\label{app:infinity}

We use Laumon's local Fourier transforms $\mathcal F^{(0,\infty')}$ and
$\mathcal F^{(\infty,\infty')}$ \cite[D\'ef.~(2.4.2.3)]{Lau87} and his
principle of stationary phase \cite[Prop.~(2.3.3.1)(iii)]{Lau87}, in the form
of \cite[Thm.~3]{KatzB}. Let $\Sigma\subset\mathbb A^1$ be a finite set,
$j_\Sigma\colon\mathbb A^1\smallsetminus\Sigma\hookrightarrow\mathbb A^1$ the
inclusion, and $F_\Sigma$ a lisse sheaf on $\mathbb A^1\smallsetminus\Sigma$.
For $x\in\Sigma\cup\{\infty\}$ let $F_\Sigma(x)$ be the representation of the
inertia group at $x$ on the generic stalk, transported to $0$ by translation
if $x\in\Sigma$. Then, as representations of the inertia group at $\infty$ of
the dual line,
\begin{equation}\label{eq:B-stationary}
\FT_\psi(j_{\Sigma!}F_\Sigma[1])(\infty)\cong\bigoplus_{x\in\Sigma}
\cL_\psi(xh_1)(\infty)\otimes\mathcal F^{(0,\infty')}(F_\Sigma(x))\ \oplus\
\mathcal F^{(\infty,\infty')}(F_\Sigma(\infty)).
\end{equation}
Moreover all breaks of $\mathcal F^{(0,\infty')}(W)$ are less than $1$ for every $W$
\cite[Thm.~(2.4.3)(i)]{Lau87}, \cite[Prop.~6]{KatzB}, and
$\mathcal F^{(\infty,\infty')}(W)=0$ if all breaks of $W$ are $\le1$
\cite[Thm.~(2.4.3)]{Lau87}, \cite[Cor.~7]{KatzB}.

\begin{proposition}[{\cite[Prop.~6.11]{EbaCoh}}]\label{prop:B-slopes}
Every break of $\cF_c(\infty)$ is at most $1$. If $M$ is the part of
$\cF_c(\infty)$ of break $1$, then all breaks of
$\cL_\psi(h_1/\lambda)(\infty)\otimes M$ are less than $1$.
\end{proposition}

\begin{proof}
Let $\mathfrak C$ be the class of continuous representations $Q$ of the
inertia group at $\infty$ such that all breaks of $Q$ are $\le1$ and all
breaks of $\cL_\psi(h_1/\lambda)(\infty)\otimes Q^{(1)}$ are less than $1$, where
$Q^{(1)}$ is the part of $Q$ of break $1$. The break decomposition is
functorial and exact \cite[1.1]{Kat88}, so $\mathfrak C$ is stable under
subquotients and extensions. It contains every representation all of whose
breaks are less than $1$, and every representation $\cL_\psi(v_0h_1)(\infty)\otimes Y$
with all breaks of $Y$ less than $1$, by \cref{lem:B-breaks}(ii); recall
$v_0=-1/\lambda$.

Let $j_U\colon\mathbb A^1_v\smallsetminus\{0,v_0\}\hookrightarrow\mathbb A^1_v$
be the inclusion and $K_1=j_{U!}\bar\cA_c[1]$. The complex $j_{v!}\cR_c$
vanishes at $v=0$ and agrees with $\bar\cA_c[-(n-2)]$ on $U_v$, so excision
for the point $v_0$ gives a distinguished triangle
$K_1[-(n-1)]\to j_{v!}\cR_c\to i_{v_0*}\cR_{c,v_0}\to$, where $\cR_{c,v_0}$ is
the stalk at $v_0$, a complex of finite-dimensional vector spaces. By proper
base change,
$\FT_\psi(i_{v_0*}\cR_{c,v_0})\cong\cR_{c,v_0}\otimes\cL_\psi(v_0h_1)[1]$. By
\eqref{eq:B-stationary} with $\Sigma=\{0,v_0\}$, the representation
$\FT_\psi(K_1)(\infty)$ on the generic stalk of the lisse sheaf
$\mathcal H^{-1}(\FT_\psi(K_1))$ near $\infty$ lies in $\mathfrak C$: the
summand for $x=0$ has all breaks less than $1$, the summand for $x=v_0$ is of the form
$\cL_\psi(v_0h_1)(\infty)\otimes Y$ with all breaks of $Y$ less than $1$, and
the last summand vanishes, because $\bar\cA_c$ is unramified at $\infty$
(\cref{lem:B-R}). The sheaf $\mathcal H^0(\FT_\psi(K_1))$ is punctual, since
$\FT_\psi(K_1)$ is perverse.

Apply $\FT_\psi$ to the triangle. By \cref{lem:B-FT-model},
$\FT_\psi(j_{v!}\cR_c)$ is, on $\mathbb G_{m,h_1}$, the lisse sheaf
$\cH'_n|_{C_c}(1)$ placed in degree $n-2$. On a punctured neighbourhood of
$\infty$ all cohomology sheaves involved are lisse, and the long exact
sequence of cohomology sheaves exhibits $(\cH'_n|_{C_c})(\infty)$ as an
extension of a subrepresentation of
$\mathcal H^{n-1}(\cR_{c,v_0})\otimes\cL_\psi(v_0h_1)(\infty)$ by a quotient of
$\FT_\psi(K_1)(\infty)$. Tate twists do not change the action of inertia.
Since $\mathfrak C$ is stable under subquotients and extensions,
$(\cH'_n|_{C_c})(\infty)$ and its quotient $\cF_c(\infty)$ lie in
$\mathfrak C$.
\end{proof}

\begin{theorem}\label{thm:B-infinity}
\begin{enumerate}[label=\textup{(\roman*)}]
\item The eigenvalues of a geometric Frobenius on $\cF_c(\infty)^{I_\infty}$
are the numbers $q^{(n+4-j)/2}$ with multiplicity $\kappa_j$, $j\ge3$.
\item $\cF_c(\infty)=V^{\mathrm t}\oplus V^{\mathrm w}$, where
$V^{\mathrm t}$ is tame and unipotent with $\kappa_j$ Jordan blocks of size
$j-2$ for every $j\ge3$, and $V^{\mathrm w}$ has dimension $2^{n-1}$ and is
isoclinic of break $1$. Moreover
$V^{\mathrm w}\cong\cL_\psi(-h_1/\lambda)(\infty)\otimes X$, where all breaks
of $X$ are less than $1$.
\end{enumerate}
In particular \ref{H4}(b) holds.
\end{theorem}

\begin{proof}
\emph{Global input.} The proofs of the formula \eqref{eq:Lline}, of the
vanishing of $H^0_c$ and $H^2_c$ in \cref{thm:eps-explicit}, of
\cref{cor:H0}, of \cref{thm:paired-explicit} and of
\cref{prop:boundary-frob}(a) use only \ref{H1} and \ref{H2} in dimensions $n$
and $n-1$, which hold by \cref{cor:B-H127}; we apply them to $\cF_c$. Thus
$\chi_c(\Gm,\cF_c)=-2^{n-1}$, and since $\cF_c$ is tame at $0$
(\cref{thm:B-zero}), \eqref{eq:B-GOS} gives $\Swan_\infty(\cF_c)=2^{n-1}$.
Moreover, by \cref{prop:boundary-frob}(a) and \cref{thm:B-zero}(iii), which
gives $\det(1-T\Phi\mid\ker N)=\det(1-TF_{n-1,c})(1-q^{\nu+1}T)^{\kappa_1}$,
\[
\det\bigl(1-T\Frob\mid\cF_c(\infty)^{I_\infty}\bigr)=
\frac{\prod_{k=2}^{\nu'}(1-q^kT)^{b_{n,k}}}{(1-q^{\nu+1}T)^{\kappa_1}}
=\prod_{j\ge3}\bigl(1-q^{(n+4-j)/2}T\bigr)^{\kappa_j},
\]
with $\nu'$ as in \cref{prop:boundary-frob}. For the last equality note that
$b_{n,k}=m_{n-1,k-2}=\kappa_{n+4-2k}$ for $2\le k\le\nu$, and that for odd $n$
the exponent of $1-q^{\nu+1}T$ is $b_{n,\nu+1}-\kappa_1=m_{n,\nu}-m_{n-1,\nu}
=m_{n-1,\nu-1}=\kappa_3$. This proves (i).

\emph{Jordan blocks.} Let $U_{\mathbf 1}$ be the unipotent part of
$\cF_c(\infty)$ and $N_\infty$ the logarithm of a topological generator of
the tame inertia acting on it. The wild part and the summands
$\cL_\rho\otimes U_\rho$ with $\rho\ne\mathbf 1$ have no inertia invariants,
so $\cF_c(\infty)^{I_\infty}=\ker(N_\infty|U_{\mathbf 1})$. The summand
$U_{\mathbf 1}$ is stable under a decomposition group at $\infty$ and under the
logarithm of the unipotent part of the local monodromy, so by the
weight-monodromy theorem \cite[1.8.4]{Del80} the graded pieces
$\mathrm{Gr}^M_i$ of the monodromy filtration of $N_\infty$ on $U_{\mathbf 1}$
are pure of weight $n+1+i$. The quotient
$(\ker N_\infty\cap M_{-i})/(\ker N_\infty\cap M_{-i-1})$ embeds equivariantly
into $\mathrm{Gr}^M_{-i}$, and, computing in a Jordan basis, its dimension is
the number of Jordan blocks of size $i+1$. Hence the number of blocks of size
$e$ of $U_{\mathbf 1}$ is the number of eigenvalues of weight $n+2-e$ on
$\cF_c(\infty)^{I_\infty}$, which is $\kappa_{e+2}$ by (i). Thus
$U_{\mathbf 1}$ has $\kappa_j$ Jordan blocks of size $j-2$ for every $j\ge3$,
and $\dim U_{\mathbf 1}=\sum_{j\ge3}(j-2)\kappa_j=r(n)-2^{n-1}$ by
\cref{lem:catalan}.

\emph{Dimension count.} By \cref{prop:B-slopes} all breaks of $\cF_c(\infty)$
are at most $1$, so the wild part has dimension at least
$\Swan_\infty(\cF_c)=2^{n-1}$. Since
$r(n)=\dim U_{\mathbf 1}+\sum_{\rho\ne\mathbf 1}\dim(\cL_\rho\otimes U_\rho)
+\dim(\text{wild part})$, it follows that $U_\rho=0$ for $\rho\ne\mathbf 1$
and that the wild part has dimension $2^{n-1}$ with all breaks equal to $1$.
So the tame part $V^{\mathrm t}=U_{\mathbf 1}$ is unipotent, the wild part
$V^{\mathrm w}$ is the part of break $1$, and
$X=\cL_\psi(h_1/\lambda)(\infty)\otimes V^{\mathrm w}$ has all breaks less than $1$ by
\cref{prop:B-slopes}.
\end{proof}

The wild part $V^{\mathrm w}$ can be decomposed further. Fix square roots
$\sqrt{c_i}\in\Fqb$, put $S_\varepsilon(c)=\sum_{i=2}^n\varepsilon_i\sqrt{c_i}$
for $\varepsilon\in\{\pm1\}^{n-1}$, and let $E'$ be a set of representatives
of $\{\pm1\}^{n-1}/\{\pm1\}$. The numbers
$-(\varepsilon_1\sqrt{h_1}+S_\varepsilon(c))^2/\lambda$ are the critical
values of the phase at $(h_1,c)$ (\cref{subsec:R-rays}). Let $E'_0$ be the set
of $\varepsilon\in E'$ with $S_\varepsilon(c)=0$, put $e=|E'_0|$, let $\chi_4$
be a character of order $4$, and for $\varepsilon\in E'\smallsetminus E'_0$ put
$\kappa_\varepsilon=S_\varepsilon(c)^2/\lambda^2$ and
\[
Y_\varepsilon=\mathcal F^{(0,\infty')}\bigl(\cL_\psi(\kappa_\varepsilon/v)(0)
\otimes\cL_{\chi_2}(v)(0)^{\otimes n}\bigr).
\]

\begin{theorem}\label{thm:B-wild}
\begin{enumerate}[label=\textup{(\roman*)}]
\item One has
\[
V^{\mathrm w}\cong\bigoplus_{\varepsilon\in E'\smallsetminus E'_0}
\cL_\psi(-h_1/\lambda)(\infty)\otimes Y_\varepsilon\ \oplus\
\cL_\psi(-h_1/\lambda)(\infty)\otimes T_c,
\]
where each $Y_\varepsilon$ is irreducible of dimension $2$, with both breaks
equal to $1/2$, and $T_c$ is tame of dimension $2e$. If $n$ is even, then
$Y_\varepsilon\cong[\times\kappa_\varepsilon]^*\KL(\infty)$. If the numbers
$S_\varepsilon(c)^2$, $\varepsilon\in E'\smallsetminus E'_0$, are pairwise
distinct, then the $Y_\varepsilon$ are pairwise non-isomorphic.
\item If $n$ is even, then
$T_c\cong\bigl(\cL_{\chi_4}(h_1)(\infty)\oplus\cL_{\chi_4^3}(h_1)(\infty)\bigr)^{\oplus e}$.
\item If $n$ is odd, then $T_c\cong\cL_{\chi_2}(h_1)(\infty)^{\oplus e}\oplus
\mathbf 1^{\oplus e}$.
\item For $\kappa\in\Fqb^\times$ let $V^{\mathrm w}_\kappa$ be the sum of the summands
$\cL_\psi(-h_1/\lambda)(\infty)\otimes Y_\varepsilon$ with $\kappa_\varepsilon=\kappa$.
Then $V^{\mathrm w}=\bigoplus_\kappa V^{\mathrm w}_\kappa\oplus
\cL_\psi(-h_1/\lambda)(\infty)\otimes T_c$, the decomposition according to the
characters of the wild inertia group, and an element of a decomposition group at
$\infty$ whose image in $\mathrm{Gal}(\Fqb/\Fq)$ is $\phi$ maps $V^{\mathrm w}_\kappa$
onto $V^{\mathrm w}_{\phi(\kappa)}$ or onto $V^{\mathrm w}_{\phi^{-1}(\kappa)}$ (according to
the convention for the action of $\mathrm{Gal}(\Fqb/\Fq)$ on characters); in
particular $V^{\mathrm w}_\kappa$ is stable under a decomposition group if and only if
$\kappa\in\Fq$.
\end{enumerate}
In particular, if the $2^{n-2}$ numbers $S_\varepsilon(c)^2$ are nonzero and
pairwise distinct, then $V^{\mathrm w}$ is the direct sum of $2^{n-2}$ pairwise
non-isomorphic irreducible representations.
\end{theorem}

\begin{proof}
We use the notation of \cref{lem:B-R} and of the proof of
\cref{prop:B-slopes}; $v_0=-1/\lambda$, and $w=v-v_0$ is a local coordinate at
$v_0$. After a finite extension of $\Fq$ we may assume that $q\equiv1\pmod4$.

\emph{The local monodromy of $\bar\cA_c$ at $v_0$.} Let $n\ge3$. By the proof
of \cref{lem:B-R}, $\bar\cA_c=a^*\cH'_{n-1}(1)$ with $a(v)=t'(v)c$ and
$t'(v)=\lambda v(\lambda v+1)^{-1}$; so near $v_0$, $\bar\cA_c$ is the pull-back
of $\rho_c^*\cH'_{n-1}(1)$ by $v\mapsto t'(v)$, which maps a punctured
neighbourhood of $v_0$ isomorphically onto one of $t'=\infty$. As
$\cH'_{n-1}=\cH_{n-1}/\Ql$ (\cref{thm:B-family}(ii)) and the constant subsheaf
lies in the unipotent part, \cref{prop:R-ray}(iii) in dimension $n-1$ shows
that $(\rho_c^*\cH'_{n-1})(\infty)$ is the direct sum of a unipotent
representation, of the characters
$\cL_\psi(-S_\varepsilon(c)^2t'/\lambda)(\infty)\otimes\cL_{\chi_2}(t')(\infty)^{\otimes n}$,
$\varepsilon\in E'\smallsetminus E'_0$, and of
$(\cL_{\chi_4}(t')\oplus\cL_{\chi_4^3}(t'))(\infty)^{\oplus e}$ if $n$ is even,
resp.\ of $\cL_{\chi_2}(t')(\infty)^{\oplus e}$ if $n$ is odd. Since
$t'=1-(\lambda w)^{-1}=-(\lambda w)^{-1}(1-\lambda w)$, we have
$\cL_\psi(-S^2t'/\lambda)\cong\cL_\psi(S^2/(\lambda^2w))$, and
$\cL_\rho(t')\cong\cL_{\rho^{-1}}(w)$ for every character $\rho$ of finite order
prime to $p$, as representations of the inertia group at $v_0$. Hence,
transported to $0$ by translation,
\[
\bar\cA_c(v_0)\cong U_A\oplus T_A\oplus\bigoplus_{\varepsilon\in
E'\smallsetminus E'_0}\cL_\psi(\kappa_\varepsilon/v)(0)\otimes\cL_{\chi_2}(v)(0)^{\otimes n},
\]
with $U_A$ unipotent and $T_A=(\cL_{\chi_4}(v)\oplus\cL_{\chi_4^3}(v))(0)^{\oplus
e}$ for even $n$, $T_A=\cL_{\chi_2}(v)(0)^{\oplus e}$ for odd $n$. If $n=2$, then
$\bar\cA_c=\cL_\psi(-c_2v/(\lambda v+1))$ and
$-c_2v/(\lambda v+1)=c_2/(\lambda^2w)-c_2/\lambda$, so the same formula holds,
with $U_A=0$, $E'$ having one element, $e=0$ and $\kappa=c_2/\lambda^2$.

\emph{The representations $Y_\varepsilon$.} By \cite[Thm.~(2.4.3)(i)]{Lau87},
$Y_\varepsilon$ has dimension $1+\Swan=2$, Swan conductor $1$ and all breaks
less than $1$. A one-dimensional subquotient would have an integral break, by
the Hasse--Arf theorem, hence break $0$; so $Y_\varepsilon$ is irreducible,
and its two breaks are equal to $1/2$. As $\mathcal F^{(0,\infty')}$ is fully
faithful, $Y_\varepsilon\cong Y_{\varepsilon'}$ implies that
$\cL_\psi(\kappa_\varepsilon/v)\otimes\cL_{\chi_2}(v)^{\otimes n}$ and
$\cL_\psi(\kappa_{\varepsilon'}/v)\otimes\cL_{\chi_2}(v)^{\otimes n}$ agree on
the wild inertia group, that is, $\kappa_\varepsilon=\kappa_{\varepsilon'}$. If
$n$ is even, apply \eqref{eq:B-stationary} to the sheaf $\cL_\psi(\kappa/v)$ on
$\Gm$, which is unramified at $\infty$: it gives
$\FT_\psi(j_!\cL_\psi(\kappa/v)[1])(\infty)\cong\mathcal F^{(0,\infty')}
(\cL_\psi(\kappa/v)(0))$, while \ref{K4}, applied after the substitution
$v=\kappa u$, and proper base change give
$\FT_\psi(j_!\cL_\psi(\kappa/v)[1])|_{\Gm}\cong[\times\kappa]^*\KL[1]$. Hence
$Y_\varepsilon\cong[\times\kappa_\varepsilon]^*\KL(\infty)$.

\emph{Tame characters.} Let $\rho\ne\mathbf 1$ be a character of $\Fq^\times$
with $\rho^4=\mathbf 1$. By \eqref{eq:B-stationary}, applied to $\cL_\rho$ on
$\Gm$, which is tame at $\infty$, and by the computation of
$\FT_\psi(j_!\cL_\rho[1])$ in the proof of \cref{prop:R-ray}(iii),
$\mathcal F^{(0,\infty')}(\cL_\rho(0))\cong\cL_{\rho^{-1}}(\infty')$. Moreover
$\mathcal F^{(0,\infty')}$ maps unipotent representations to unipotent
representations (proof of \cref{prop:R-ray}(i)).

\emph{Stationary phase.} By \eqref{eq:B-stationary} applied to $K_1$, as in the
proof of \cref{prop:B-slopes}, and by the additivity of
$\mathcal F^{(0,\infty')}$,
\[
\FT_\psi(K_1)(\infty)\cong\mathcal F^{(0,\infty')}(\bar\cA_c(0))\oplus
\cL_\psi(-h_1/\lambda)(\infty)\otimes(U'_A\oplus T'_A)\oplus
\bigoplus_{\varepsilon\in E'\smallsetminus E'_0}\cL_\psi(-h_1/\lambda)(\infty)
\otimes Y_\varepsilon,
\]
where $U'_A=\mathcal F^{(0,\infty')}(U_A)$ is unipotent and
$T'_A=\mathcal F^{(0,\infty')}(T_A)$ is
$(\cL_{\chi_4}(h_1)\oplus\cL_{\chi_4^3}(h_1))(\infty)^{\oplus e}$ for even $n$
and $\cL_{\chi_2}(h_1)(\infty)^{\oplus e}$ for odd $n$. Let $P_\infty$ be the
wild inertia group, $\chi_0$ the restriction of $\cL_\psi(-h_1/\lambda)(\infty)$
to $P_\infty$, and $\Xi$ the set of characters of $P_\infty$ occurring in
$\bigoplus_{\varepsilon\notin E'_0}\cL_\psi(-h_1/\lambda)(\infty)\otimes
Y_\varepsilon$. As $Y_\varepsilon$ is irreducible and not tame,
$Y_\varepsilon^{P_\infty}=0$, so $\chi_0\notin\Xi$; the trivial character is not
in $\Xi$, as $\cL_\psi(-h_1/\lambda)\otimes Y_\varepsilon$ is isoclinic of break
$1$ (\cref{lem:B-breaks}(ii)); and no character in $\Xi$ occurs in
$\mathcal F^{(0,\infty')}(\bar\cA_c(0))$, all of whose breaks are less than
$1$, since representations with different breaks have no common constituents
under $P_\infty$ \cite[1.1]{Kat88}. The middle summand is $\chi_0$-isotypic. For
a representation $Q$ of the inertia group and a set $\Theta$ of characters of
$P_\infty$ stable under conjugation, let $Q_\Theta$ be the sum of the isotypic
components of $Q$ under $P_\infty$ for the characters in $\Theta$; it is a
direct summand stable under inertia, and $Q\mapsto Q_\Theta$ is exact. By the
proof of \cref{prop:B-slopes}, there is an exact sequence
\[
\cR^{n-2}\otimes\cL_\psi(-h_1/\lambda)(\infty)\to\FT_\psi(K_1)(\infty)\to
(\cH'_n|_{C_c}(1))(\infty)\to\cR^{n-1}\otimes\cL_\psi(-h_1/\lambda)(\infty)\to0,
\]
where $\cR^i=\mathcal H^i(\cR_{c,v_0})$ is a finite-dimensional vector space
with trivial action of inertia, so that the outer terms are $\chi_0$-isotypic.
By \cref{thm:B-family}(iii), $\cF_c$ is the quotient of $\cH'_n|_{C_c}$ by a
lisse subsheaf whose graded pieces are geometrically constant, on which
$P_\infty$ acts trivially. Hence
$\cF_c(\infty)_\Xi\cong\bigoplus_{\varepsilon\notin E'_0}
\cL_\psi(-h_1/\lambda)(\infty)\otimes Y_\varepsilon$. Moreover every character
of $P_\infty$ occurring in $\cF_c(\infty)$ occurs in $\FT_\psi(K_1)(\infty)$ or
equals $\chi_0$; so the characters of break $1$ occurring in $\cF_c(\infty)$ lie
in $\Xi\cup\{\chi_0\}$, and $V^{\mathrm w}=\cF_c(\infty)_\Xi\oplus
\cF_c(\infty)_{\{\chi_0\}}$.

\emph{The part for $\chi_0$.} Put $T_c=\cF_c(\infty)_{\{\chi_0\}}\otimes
\cL_\psi(h_1/\lambda)(\infty)$, a tame representation, and decompose it into
generalized eigenspaces of a topological generator $\gamma$ of the tame inertia
group. As $P_\infty$ acts trivially on the subsheaf just mentioned,
$\cF_c(\infty)_{\{\chi_0\}}=(\cH'_n|_{C_c})(\infty)_{\{\chi_0\}}$; and in the
exact sequence above, the twists by $\cL_\psi(h_1/\lambda)(\infty)$ of the outer
terms are trivial, while that of $\FT_\psi(K_1)(\infty)_{\{\chi_0\}}$ is
$U'_A\oplus T'_A$. So for every eigenvalue $\ne1$ of $\gamma$ the generalized
eigenspaces of $T_c$ and of $T'_A$ are isomorphic. By
\cref{thm:B-infinity}(ii), $\dim V^{\mathrm w}=2^{n-1}$, so
$\dim T_c=2^{n-1}-2(2^{n-2}-e)=2e$. This proves (i). If $n$ is even, then
$T'_A$ has dimension $2e$ and no eigenvalue $1$, so $T_c\cong T'_A$, which is
(ii). If $n$ is odd, then $T'_A$ has dimension $e$, and $T_c\cong T'_A\oplus U_c$
with $U_c$ unipotent of dimension $e$.

\emph{A twisted $L$-function.} Let $\nu=\lfloor n/2\rfloor$, put $\delta=1$ if
$n$ is odd and $\delta=0$ if $n$ is even, let
$\cK_c=\bigotimes_{i\ge2}[\times c_i]^*\KL$ on $\Gm$ and
$\cF'_c=\cF_c\otimes\cL_\psi(h_1/\lambda)$. Let $a_1,\dots,a_k$ be the distinct
values among $c_2,\dots,c_n$, with multiplicities $\mu_1,\dots,\mu_k$, and put
$g_c=\prod_jC_{\mu_j/2}$ if all $\mu_j$ are even and $g_c=0$ otherwise. By
\cref{lem:R-KL}(b), the geometric monodromy group of $\cK_c$ is $\SL_2^k$, acting
on $\bigotimes_j\mathrm{St}_j^{\otimes\mu_j}$, so by \cref{lem:B-CG} the space
$H^0(\mathbb G_{m,\Fqb},\cK_c)$ of geometric invariants has dimension $g_c$; it
embeds into a stalk, so it is pure of weight $n-1$. As $\cK_c$ is geometrically
semisimple \cite[3.4.1]{Del80} and dual to itself up to a twist by \ref{K1},
$H^2_c(\mathbb G_{m,\Fqb},\cK_c)$, the space of geometric coinvariants twisted by
$(-1)$, also has dimension $g_c$, and it is pure of weight $n+1$. For $Q=q^m$,
orthogonality and \eqref{eq:C-factor} give
\begin{equation}\label{eq:B-twisted}
\sum_{t\in\mathbb F_Q^\times}\psi_Q(t/\lambda)\Tr(\Frob_t\mid\cF_c)
=(-1)^{n-1}Q\Bigl(\sum_{s\in\mathbb F_Q^\times}\prod_{i\ge2}\Kl_Q(c_is)
+\Sigma_{n-1,Q}(c)\Bigr)+E_n(Q),
\end{equation}
and the first sum is $(-1)^n$ times the alternating sum
$\Tr(\Frob_q^m\mid H^1_c(\mathbb G_{m,\Fqb},\cK_c))-\Tr(\Frob_q^m\mid
H^2_c(\mathbb G_{m,\Fqb},\cK_c))$. By the proof of \cref{prop:R-ray}(ii),(iii) in
dimension $n-1$, $\cK_c$ is tame at $0$ with $\cK_c(0)\cong J_2^{\otimes(n-1)}$,
and the tame part of $\cK_c(\infty)$ is $(\mathbf 1\oplus\cL_{\chi_2})^{\oplus e}$
for odd $n$ and $(\cL_{\chi_4}\oplus\cL_{\chi_4^3})^{\oplus e}$ for even $n$; so
$(j_*\cK_c)_{\bar\infty}$ has dimension $\delta e$ and, by the weight-monodromy
theorem \cite[1.8.4]{Del80}, is pure of weight $n-1$. As in \cref{prop:B-Mh},
there is an exact sequence
\[
0\to H^0(\mathbb G_{m,\Fqb},\cK_c)\to(j_*\cK_c)_{\bar0}\oplus(j_*\cK_c)_{\bar\infty}
\to H^1_c(\mathbb G_{m,\Fqb},\cK_c)\to H^1(\PP^1_{\Fqb},j_*\cK_c)\to0,
\]
whose last term is pure of weight $n$, and Frobenius acts on
$(j_*\cK_c)_{\bar0}$ with the eigenvalues $q^k$, of multiplicities $m_{n-1,k}$.
Let $P_{\cK_c}(T)=\det(1-T\Frob_q\mid H^1(\PP^1_{\Fqb},j_*\cK_c))$,
$P_{\cK_c,\infty}(T)=\det(1-T\Frob_q\mid(j_*\cK_c)_{\bar\infty})$, and let
$P_{\cK_c,0}(T)$ and $P_{\cK_c,2}(T)$ be the characteristic polynomials of
$\Frob_q$ on $H^0(\mathbb G_{m,\Fqb},\cK_c)$ and on $H^2_c(\mathbb
G_{m,\Fqb},\cK_c)$. With \ref{H2} in dimension $n-1$ and Pascal's rule,
\eqref{eq:B-twisted} gives
\[
L(C_c,\cF'_c,T)=\frac{P_{\cK_c}(qT)\,P_{\cK_c,\infty}(qT)\det(1-TF_{n-1,c})
\,(1-q^{\nu+1}T)^{\delta C_\nu}}{P_{\cK_c,0}(qT)\,P_{\cK_c,2}(qT)} .
\]
The group $H^2_c(C_{c,\Fqb},\cF'_c)$ is pure of weight $n+3$ (proof of
\cref{lem:B-av}), whereas $H^1_c$ is mixed of weights $\le n+2$, and the only
reciprocal roots or poles of weight $n+3$ above are those of $P_{\cK_c,2}(qT)$.
Hence $\det(1-T\Frob_q\mid H^2_c(C_{c,\Fqb},\cF'_c))=P_{\cK_c,2}(qT)$, so
$H^2_c(C_{c,\Fqb},\cF'_c)$ and, by geometric semisimplicity and Poincaré
duality, $H^0(C_{c,\Fqb},\cF'_c)$ have dimension $g_c$, and
\begin{equation}\label{eq:B-twistedL}
\det\bigl(1-T\Frob_q\mid H^1_c(C_{c,\Fqb},\cF'_c)\bigr)
=\frac{P_{\cK_c}(qT)\,P_{\cK_c,\infty}(qT)\det(1-TF_{n-1,c})
\,(1-q^{\nu+1}T)^{\delta C_\nu}}{P_{\cK_c,0}(qT)} .
\end{equation}
If some value occurs an odd number of times among $c_2,\dots,c_n$, which is
automatic for even $n$, then $g_c=0$, the denominators are $1$, and
$L(C_c,\cF'_c,T)$ is the polynomial \eqref{eq:B-twistedL}.

\emph{The unipotent part for odd $n$.} Let $n$ be odd. By
\eqref{eq:B-twistedL}, exactly $e+C_\nu-g_c$ reciprocal zeros of
$\det(1-T\Frob_q\mid H^1_c(C_{c,\Fqb},\cF'_c))$ have weight $n+1=w$: the $e$
zeros of $P_{\cK_c,\infty}(qT)$ and the $C_\nu$ zeros of $(1-q^{\nu+1}T)^{C_\nu}$,
less the $g_c$ zeros of $P_{\cK_c,0}(qT)$, all of weight $w$, which divide the
numerator because the quotient is a polynomial. On the other hand, as in the
proof of \cref{cor:H0}, there is an exact sequence
\[
0\to H^0(C_{c,\Fqb},\cF'_c)\to(j_*\cF'_c)_{\bar0}\oplus(j_*\cF'_c)_{\bar\infty}
\to H^1_c(C_{c,\Fqb},\cF'_c)\to H^1(\PP^1_{\Fqb},j_*\cF'_c)\to0,
\]
whose last term is pure of weight $w+1$ \cite[3.2.3]{Del80} and whose first
term is pure of weight $w$, of dimension $g_c$. Since $\cL_\psi(h_1/\lambda)$ is
lisse at $0$, with trivial Frobenius there, \cref{thm:B-zero}(iii) shows that
exactly $\kappa_1=m_{n-1,\nu}=C_\nu$ eigenvalues of Frobenius on
$(j_*\cF'_c)_{\bar0}$ have weight $w$. At $\infty$, the representation
$\cF'_c(\infty)$ is the direct sum of $\cL_\psi(h_1/\lambda)(\infty)\otimes
V^{\mathrm t}$, which is isoclinic of break $1$, of the $Y_\varepsilon$, of
$\cL_{\chi_2}(h_1)(\infty)^{\oplus e}$ and of the unipotent part $U_c$ of $T_c$,
of dimension $e$; so $(j_*\cF'_c)_{\bar\infty}=U_c^{I_\infty}$, and by the
weight-monodromy theorem, applied as in the proof of \cref{thm:B-infinity}, the
number of its eigenvalues of weight $w$ is the number $s$ of Jordan blocks of
size $1$ of $U_c$. Comparing, $e+C_\nu-g_c=C_\nu+s-g_c$. So $s=e=\dim U_c$, and
$U_c$ is trivial. This proves (iii).

\emph{The decomposition group.} The summands $V^{\mathrm w}_\kappa$ and
$\cL_\psi(-h_1/\lambda)(\infty)\otimes T_c$ have pairwise disjoint sets of
characters of the wild inertia group, by the proof of (i), so they are the
isotypic parts of $V^{\mathrm w}$ and are permuted by a decomposition group at
$\infty$, the last one being stable. By the proof of (i), the isomorphism of the
$\Xi$-parts of $\cF_c(\infty)$ and of $\FT_\psi(K_1)(\infty)$ comes from exact
sequences of sheaves defined over $\Fq$, so it is equivariant for a decomposition
group; and by the principle of stationary phase, which is compatible with the
actions of decomposition groups \cite[(2.3.3.1)]{Lau87}, the part of
$\FT_\psi(K_1)(\infty)$ corresponding to $V^{\mathrm w}_\kappa$ is
$\cL_\psi(-h_1/\lambda)(\infty)\otimes\mathcal F^{(0,\infty')}$ applied to the
isotypic part of $\bar\cA_c(v_0)$ for the character of the wild inertia group at
$v_0$ given by $\cL_\psi(\kappa/w)$. A decomposition group at $v_0$ permutes these
isotypic parts as $\mathrm{Gal}(\Fqb/\Fq)$ permutes the numbers $\kappa$, since
the restriction of $\cL_\psi(\kappa/w)$ to the wild inertia group determines
$\kappa$ and the conjugate of $\cL_\psi(\kappa/w)$ by an element of $\mathrm{Gal}
(\Fqb/\Fq)$ is $\cL_\psi(\kappa'/w)$, $\kappa'$ being the corresponding conjugate
of $\kappa$. This proves (iv).
\end{proof}

For $n=2$ the theorem gives $V^{\mathrm w}\cong\cL_\psi(-h_1/\lambda)(\infty)
\otimes[\times c_2/\lambda^2]^*\KL(\infty)$, in agreement with
\cref{prop:n2sheaf}. After the pull-back $h_1=y^2$, the summand indexed by
$\varepsilon\notin E'_0$ becomes, for even $n$, the sum of the two characters
$\cL_\psi(-(y\pm S_\varepsilon(c))^2/\lambda)\otimes\cL_{\chi_2}(y)$ up to a
constant twist, in accordance with the critical values of the phase. For
$\varepsilon\in E'_0$ these two characters coincide, and accordingly, for even
$n$, the pull-back of each summand $\cL_\psi(-h_1/\lambda)\otimes(\cL_{\chi_4}
\oplus\cL_{\chi_4^3})(h_1)$ of \cref{thm:B-wild}(ii) is twice the character
$\cL_\psi(-y^2/\lambda)\otimes\cL_{\chi_2}(y)$.

\begin{remark}\label{rem:B-wild-check}
The value $\Swan(X)=2^{n-2}-e$ of the Swan conductor of the representation $X$
of \cref{thm:B-infinity}(ii), which follows from \cref{thm:B-wild}, can also
be read off from \eqref{eq:B-twistedL}, without local Fourier transforms,
whenever some value occurs an odd number of times among $c_2,\dots,c_n$; this
is the case if all $S_\varepsilon(c)$ are nonzero, since otherwise the signs
could be chosen so that equal values cancel in pairs. Indeed, $\dim H^1_c(\mathbb G_{m,\Fqb},\cK_c)=\Swan_\infty(\cK_c)
=2^{n-2}-e$ by \eqref{eq:B-GOS}, so
$\dim H^1(\PP^1_{\Fqb},j_*\cK_c)=2^{n-2}-e-\binom{n-1}{\lfloor(n-1)/2\rfloor}
-\delta e$, and since $r(n)=2r(n-1)+\delta C_\nu$ (\cref{lem:catalan}, as
$\kappa_1=\delta C_\nu$) and $r(n-1)=2^{n-1}-\binom{n-1}{\lfloor(n-1)/2\rfloor}$,
the polynomial
\eqref{eq:B-twistedL} has degree $r(n)-2^{n-2}-e$; its reciprocal roots have
weights $n+2$, $n+1$ and $n$. As $H^0_c$ and $H^2_c$ vanish,
$\dim H^1_c(C_{c,\Fqb},\cF'_c)=r(n)-2^{n-2}-e$. On the other hand $\cF'_c$ is
tame at $0$, and at $\infty$ it is $(\cL_\psi(h_1/\lambda)\otimes
V^{\mathrm t})\oplus X$, where the first summand is isoclinic of break $1$; so
its Euler characteristic is $-(r(n)-2^{n-1})-\Swan(X)$, and
$\Swan(X)=2^{n-2}-e$. For $n=3$ and $c=(c_2,c_3)$ with $c_2\ne c_3$,
\eqref{eq:B-twistedL} reads
$\sum_t\psi_Q(t/\lambda)\Tr(\Frob_t\mid\cF_c)=-Q^2-\Tr(F_{2,c}^m)$; for $n=4$
and $c$ with exactly one $S_\varepsilon(c)=0$ it gives
$\sum_{s\in\mathbb F_Q^\times}\prod_{i\ge2}\Kl_Q(c_is)=1+2Q$ and
$L(C_c,\cF'_c,T)=\det(1-TF_{3,c})$.
\end{remark}

\subsection{The determinant: proof of (H3)}\label{app:det}
The determinant is needed for the \emph{explicit phase} of the local
functional equation, not for either critical-circle theorem.
The argument proceeds in three logically separate steps.  First the
possible geometric characters of prime-to-$p$ order are eliminated
through inertia on coordinate lines and the Kummer sequence.  Next,
the remaining $p$-power characters are tested using the
Artin--Schreier sequence and restriction to every coordinate line.
Finally the constant arithmetic character is computed from the
Frobenius action on the inertia filtration at a finite boundary
point.  The third step must not be replaced by purity, which
fixes only the absolute value of a determinant.

\begin{proposition}[{\cite[Prop.~8.6]{EbaCoh}}]\label{prop:B-detfinite}
The character $\det\cW_n$ of $\pi_1(S_{\Fpb})$ has as image a finite
$p$-group.
\end{proposition}

\begin{proof}
Let $\Gamma$ be the image.

\emph{Step 1: inertia.} By \cref{thm:B-family}(vi), \cref{thm:B-zero} applies
to the coordinate lines in every direction. Fix $i$ and a line
$C\subset S_{\Fpb}$ along which the coordinates $h_j$, $j\ne i$, are
constant. Through $\pi_1(C)\to\pi_1(S_{\Fpb})$, the inertia group at $h_i=0$
acts on the generic stalk of $\cW_n$ through the closure of the group
generated by a unipotent element, on which $\det$ is identically $1$. So
$\det\cW_n$ is trivial on these inertia groups.

\emph{Step 2: characters of order prime to $p$.} Let $m\ge1$ be prime to $p$.
The ring $\mathcal O(S_{\Fpb})=\Fpb[h_1^{\pm1},\dots,h_n^{\pm1}]$ is
factorial, so $\operatorname{Pic}(S_{\Fpb})=0$; its units are
$\Fpb^\times h_1^{\mathbb Z}\cdots h_n^{\mathbb Z}$, and $\Fpb^\times$ is
$m$-divisible. So the Kummer sequence gives
$H^1(S_{\Fpb},\mu_m)\cong\mathcal O(S_{\Fpb})^\times/\mathcal
O(S_{\Fpb})^{\times m}\cong(\mathbb Z/m)^n$, with basis the classes of
$h_1,\dots,h_n$. On a line $C$ as in Step 1, the character
attached to $h_1^{e_1}\cdots h_n^{e_n}$ restricts to the Kummer character of
$h_i^{e_i}$, whose restriction to the inertia group at $h_i=0$ is trivial if
and only if $m\mid e_i$. So a character of $\pi_1(S_{\Fpb})$ of finite order
prime to $p$ that is trivial on all the inertia groups of Step 1 is trivial.

\emph{Step 3.} The group $\Gamma\subseteq\Ql^\times$ is compact. By the Baire
category theorem it is contained in $E^\times$ for a finite extension $E$ of
$\mathbb Q_\ell$, hence in $\mathcal O_E^\times$, which is the product of a
finite group and the pro-$\ell$ group $1+\mathfrak m_E$. The profinite abelian
group $\Gamma$ is the product of its pro-$p$ part $\Gamma_p$, which meets
$1+\mathfrak m_E$ trivially and is therefore finite, and its prime-to-$p$ part
$\Gamma'$. If $\Gamma'\ne1$, it has a nontrivial character of finite order
prime to $p$, and composing gives a nontrivial character of $\pi_1(S_{\Fpb})$
of order prime to $p$ that is trivial on the inertia groups of Step 1,
contradicting Step 2. Hence $\Gamma=\Gamma_p$ is a finite $p$-group.
\end{proof}

\begin{lemma}[{\cite[Lemma~8.10]{EbaCoh}}]\label{lem:B-AS}
Let $k$ be an algebraically closed field of characteristic $p$ and
$\chi\colon\pi_1(\mathbb G^n_{m,k})\to\Ql^\times$ a continuous character
whose image is a finite $p$-group. If, for every $i$ and every
$c\in(k^\times)^{n-1}$, the restriction of $\chi$ to the fundamental group of
the coordinate line $\{h:h_j=c_j\ (j\ne i)\}$ is trivial, then $\chi$ is
trivial.
\end{lemma}

\begin{proof}
The image of $\chi$ is cyclic of order $p^e$; by induction on $e$, applied to
$\chi^p$, we may assume $e\le1$. Suppose $e=1$, and regard $\chi$ as an element
of $H^1(\mathbb G^n_{m,k},\mathbb F_p)$. With $A=k[h_1^{\pm1},\dots,h_n^{\pm1}]$
and $\wp(x)=x^p-x$, the Artin--Schreier sequence and $H^1(\mathbb
G^n_{m,k},\mathcal O)=0$ identify $H^1(\mathbb G^n_{m,k},\mathbb F_p)$ with
$A/\wp(A)$, compatibly with restriction to the coordinate lines. Let
$\mathsf M=\mathbb Z^n\smallsetminus p\mathbb Z^n$. Every class in $A/\wp(A)$
is represented by a unique Laurent polynomial $g=\sum_{m\in\mathsf M}a_mh^m$:
constants lie in $\wp(k)$, $a\,h^{pm}\equiv a^{1/p}h^m$, and if
$\sum_{m\in\mathsf M}a_mh^m=\wp(b)$ with $b$ nonconstant, then for $m_0\ne0$
in the support of $b$ with $\max_i|m_{0,i}|$ maximal the coefficient of
$h^{pm_0}$ in $\wp(b)$ is $b_{m_0}^p\ne0$, although $pm_0\notin\mathsf M$; so
$b$ is constant and all $a_m$ vanish. The same holds in one variable.

Let $g$ represent $\chi$, fix $i$, and for $c\in(k^\times)^{n-1}$ write
$m=(m',m_i)$ and $c^{m'}=\prod_{j\ne i}c_j^{m_j}$. The restriction of $g$ to
the line through $c$ in direction $i$ is $\sum_e\alpha_e(c)h_i^e$ with
$\alpha_e(c)=\sum_{m_i=e}a_mc^{m'}$; for $e\notin p\mathbb Z$, the coefficient
of $h_i^e$ in its reduced representative is
$\sum_{s\ge0}\alpha_{p^se}(c)^{1/p^s}$, which vanishes for all $c$ by
hypothesis. Raising to a large power $p^N$, we find that the Laurent
polynomial $\sum_{s\ge0}\sum_{m_i=p^se}a_m^{p^{N-s}}c^{p^{N-s}m'}$ vanishes on
$(k^\times)^{n-1}$, hence identically. Distinct pairs $(s,m)$ give distinct
monomials: for fixed $s$, $m$ is determined by $m'$ since $m_i=p^se$; and
$p^{N-s}m'=p^{N-\tilde s}\tilde m'$ with $s<\tilde s$ would give
$\tilde m=p^{\tilde s-s}m\in p\mathbb Z^n$. So $a_m=0$ whenever $m_i\ne0$, for
every $i$; as every $m\in\mathsf M$ has a nonzero coordinate, $g=0$.
\end{proof}

\begin{theorem}\label{thm:B-det}
For every $n\ge1$,
$\det\cW_n\cong\cL_\psi(D_n)\otimes\Ql\bigl(-\tfrac{(n+1)r(n)}2\bigr)$. In
particular \ref{H3} holds.
\end{theorem}

\begin{proof}
For $n=1$ this is \cref{cor:B-H127}, as $(n+1)r(n)/2=1$. Let $n\ge2$ and
assume the theorem for $n-1$.

\emph{Step 1: coordinate lines.} Let $c\in(\Fq^\times)^{n-1}$ and put
$\mathcal D_c=\det(\cF_c)\otimes\cL_\psi(2^{n-1}\lambda^{-1}h_1)$, a lisse
sheaf of rank one on $C_c$. We show that it is geometrically constant. By
\cref{thm:B-zero}(i) the inertia at $h_1=0$ acts on $\cF_c$ through unipotent
elements, so its determinant is trivial there, and $\mathcal D_c$ extends to a
lisse sheaf $\bar{\mathcal D}_c$ on $\mathbb A^1_{h_1}$. At $\infty$,
\cref{thm:B-infinity}(ii) gives
\[
\mathcal D_c(\infty)\cong\det V^{\mathrm t}\otimes
\cL_\psi(-h_1/\lambda)(\infty)^{\otimes2^{n-1}}\otimes\det X\otimes
\cL_\psi(2^{n-1}\lambda^{-1}h_1)(\infty)=\det V^{\mathrm t}\otimes\det X,
\]
where $\det V^{\mathrm t}$ is tame and $\det X$, a direct summand of
$X^{\otimes2^{n-1}}$, has break less than $1$ by \cref{lem:B-breaks}(i). Hence
$\Swan_\infty(\bar{\mathcal D}_c)<1$, so it is $0$, and the
Grothendieck--Ogg--Shafarevich formula on $\mathbb A^1$ gives
$\chi_c(\mathbb A^1_{\Fpb},\bar{\mathcal D}_c)=1$. As $H^0_c=0$, the group
$H^2_c(\mathbb A^1_{\Fpb},\bar{\mathcal D}_c)$ is nonzero, that is,
$\bar{\mathcal D}_c$ has nonzero geometric coinvariants; being of rank one, it
is geometrically constant.

\emph{Step 2: the geometric determinant.} Let $\chi$ be the character of
$\pi_1(S_{\Fpb})$ defined by $\det\cW_n\otimes\cL_\psi(-D_n)$. Its image is a
finite $p$-group, by \cref{prop:B-detfinite} and because the image of an
Artin--Schreier sheaf lies in $\mu_p$. On the coordinate line through $c$ in
direction $i$, the sheaf $\cL_\psi(D_n)$ is
$\cL_\psi(-2^{n-1}\lambda^{-1}h_i)$ tensored with a geometrically constant
sheaf, so by Step 1 and \cref{thm:B-family}(vi), $\chi$ is trivial on the
fundamental groups of all coordinate lines. By \cref{lem:B-AS}, $\chi$ is
trivial. A lisse sheaf of rank one on $S$ that is geometrically constant is of
the form $\alpha^{\deg}$, the sheaf on which $\Frob_p$ acts by some
$\alpha\in\Ql^\times$; so $\det\cW_n\cong\cL_\psi(D_n)\otimes\alpha^{\deg}$.

\emph{Step 3: the constant.} Let $q=p^k$, $c\in(\Fq^\times)^{n-1}$, and let
$V$, $\Phi$ and $N$ be as in \cref{thm:B-zero}. By Step 2,
$\det\cF_c\cong\cL_\psi(D_n(h_1,c))\otimes(\alpha^k)^{\deg}$ over $\Fq$, a
sheaf that extends to a lisse sheaf on $\mathbb A^1_{h_1}$; so $\Phi$ acts on
$\det V$ by $\alpha^k\psi_q(D_n(0,c))=\alpha^k\psi_q(2D_{n-1}(c))$. On the
other hand, by \cref{thm:B-zero}(iii) and the theorem for $n-1$, which gives
$\det F_{n-1,c}=q^{nr(n-1)/2}\psi_q(D_{n-1}(c))$,
\[
\det(\Phi\mid V)=\det F_{n-1,c}\cdot q^{(\nu+1)\kappa_1}\cdot
q^{r(n-1)}\det F_{n-1,c}=q^{(n+1)r(n-1)+(\nu+1)\kappa_1}\psi_q(2D_{n-1}(c)).
\]
Since $r(n)=2r(n-1)+\kappa_1$ and $(\nu+1)\kappa_1=(n+1)\kappa_1/2$ (both
sides vanish for even $n$), the exponent is $(n+1)r(n)/2$. Hence
$\alpha^k=q^{(n+1)r(n)/2}$, and taking $q=p$ gives $\alpha=p^{(n+1)r(n)/2}$,
that is, $\alpha^{\deg}=\Ql(-(n+1)r(n)/2)$.
\end{proof}

\subsection{Geometric irreducibility}\label{app:irred}

For $h\in\Fq^n$ define $\Sigma_{n,q}(h)$ by \eqref{eq:sum}.

\begin{lemma}[{\cite[Lemmas~2.1 and 7.1, Prop.~7.2]{EbaCoh}}]\label{lem:B-moments}
\begin{enumerate}[label=\textup{(\roman*)}]
\item Let $h\in\Fq^n$ have exactly $a$ nonzero coordinates, forming
$h'\in(\Fq^\times)^a$. Then $\Sigma_{n,q}(h)=0$ if $a=0$,
$|\Sigma_{n,q}(h)|=1$ if $a=1$, and $\Sigma_{n,q}(h)=(-1)^{n-a}\Sigma_{a,q}(h')$
if $a\ge2$.
\item $\sum_{h\in\Fq^n}|\Sigma_{n,q}(h)|^2
=q^{n-2}\bigl((q-1)^{n+1}+1-(-1)^n(q-2)\bigr)$.
\end{enumerate}
\end{lemma}

\begin{proof}
(i) If $h_i=0$, the variable $u_i$ does not occur in the constraint, and
summing over it gives the factor $\sum_{x\in\Fq^\times}\psi_q(x^{-1})=-1$. If
$a=0$ the constraint reads $0=-\lambda$; if $a=1$ it determines the remaining
variable.

(ii) Inserting $q^{-1}\sum_{t\in\Fq}\psi_q\bigl(t(\lambda+\langle
h,u\rangle)\bigr)$ for the constraint gives, for every $h\in\Fq^n$,
\[
q\Sigma_{n,q}(h)=\sum_{t\in\Fq}\psi_q(\lambda t)\prod_{i=1}^n\Kl_q(th_i),
\qquad\Kl_q(0)=-1 .
\] The sums $\Kl_q(\xi)$ are real. For
$a,b\in\Fq$ put $V_2(a,b)=\sum_{x\in\Fq}\Kl_q(ax)\Kl_q(bx)$. We claim
that
\[
V_2(a,b)=\begin{cases}q^2-q,&a=b\ne0,\\-q,&a\ne b,\ ab\ne0,\\
0,&\text{exactly one of }a,b\text{ vanishes},\\q,&a=b=0.\end{cases}
\]
Indeed, $\Kl_q(ax)\Kl_q(bx)=\sum_{u,v\in\Fq^\times}\psi_q\bigl(u+v+x(a/u+b/v)\bigr)$.
If $ab\ne0$, summing over $x$ first restricts the sum to $u=-av/b$, where
$u+v=(1-a/b)v$, and gives $q\sum_{v\in\Fq^\times}\psi_q\bigl((1-a/b)v\bigr)$.
If $a=0\ne b$, the sum is $-\sum_x\Kl_q(bx)=0$, because
$\sum_{\xi\in\Fq}\Kl_q(\xi)=0$; and $V_2(0,0)=q\Kl_q(0)^2=q$. Hence
\[
q^2\sum_{h\in\Fq^n}|\Sigma_{n,q}(h)|^2=\sum_{t_1,t_2\in\Fq}
\psi_q\bigl(\lambda(t_1-t_2)\bigr)V_2(t_1,t_2)^n .
\]
The pair $t_1=t_2=0$ contributes $q^n$, pairs with exactly one zero
contribute $0$, pairs with $t_1=t_2\ne0$ contribute $(q-1)(q^2-q)^n$, and
pairs with $t_1\ne t_2$, $t_1t_2\ne0$ contribute
\[
(-q)^n\Bigl(\Bigl|\sum_{t\ne0}\psi_q(\lambda t)\Bigr|^2-(q-1)\Bigr)=(2-q)(-q)^n .
\]
Adding up gives (ii).
\end{proof}

\begin{corollary}[{\cite[Cor.~8.4]{EbaCoh}}]\label{cor:B-irred}
For every $n\ge2$ the sheaf $\cW_n$ is geometrically irreducible.
\end{corollary}

\begin{proof}
For $h\in S(\Fq)$ put $A^0(h)=q^{-(n+1)/2}\Tr(\Frob_h\mid\cW_n)$, so that
$|A^0(h)|\le r(n)$ by \ref{H1}, and by \ref{H2}
$(-1)^{n-1}\Sigma_{n,q}(h)=q^{(n-1)/2}A^0(h)+c_n(q)$ with
$c_n(q)=q^{-1}E_n(q)\ll q^{(n-2)/2}$. Hence
$|\Sigma_{n,q}(h)|^2=q^{n-1}|A^0(h)|^2+O(q^{(2n-3)/2})$, uniformly in $h$.
For $h\in\Fq^n\smallsetminus S(\Fq)$ with $a\le n-1$ nonzero coordinates,
\cref{lem:B-moments}(i) and \ref{H1}, \ref{H2} in dimension $a$ give
$|\Sigma_{n,q}(h)|\ll q^{(a-1)/2}$, and there are $O(q^a)$ such $h$; they
contribute $O(q^{2n-3})$ to the second moment. By \cref{lem:B-moments}(ii),
\[
q^{-n}\sum_{h\in S(\Fq)}|A^0(h)|^2=1+O(q^{-1/2}).
\]
Apply \cref{lem:B-av} to $T=S$ and $\cE=\cW_n\otimes\cW_n^\vee$, which is pure
of weight $0$ with trace $|A^0(h)|^2$ at $h$; the geometric invariants of
$\cE$ are $\End_{G}(V)$, where $V$ is the geometric generic stalk of $\cW_n$
and $G$ its geometric monodromy group. Hence $\dim\End_G(V)=1$. As $\cW_n$ is
pure, $V$ is a semisimple representation of $G$ \cite[3.4.1]{Del80}, so it is
irreducible.
\end{proof}

%% === Global cohomology and completion ===
\section{Global cohomology and the intersection completion}\label{sec:global}

There are two global objects that should be kept distinct.  The ordinary Euler product $L(S_n,\cW_n,T)$ is computed using compactly supported cohomology on the nonproper torus, which can have several weights in a single degree.  The completed object $\widehat L_n(T)$ instead uses the intersection complex on the proper compactification.  Its degree-by-degree purity is a consequence of the properness of the compactification and of the purity of its intermediate extension, and is logically independent of a closed formula for all compactly supported cohomology groups.  We present the torus computation first to explain the extra Tate factors; the proof of the Riemann hypothesis for the completion begins afresh with a pure perverse sheaf and does not assume the final torus Betti table.

In this section we determine the cohomology of $\cW_n$ on the torus $S_n$,
with and without compact supports, as a module over Frobenius
(\cref{thm:G-torus}), and we study the completion of the Euler product of
\cref{sec:torus} by the intersection cohomology of the compactification
$X_n=(\PP^1)^n$ of $S_n$. Its zeros and poles satisfy the Riemann hypothesis
(\cref{thm:G-RH}); it satisfies a functional equation with sign $\pm1$
(\cref{prop:G-duality}); and a large explicit part of it, which is all of it
for $n=2$ and all of it outside the middle degree for $n=3$, is determined in
\cref{thm:G-decomp}. Cohomology groups are taken over $\Fpb$ and carry the
action of $\Frob_p$; for a complex $K$ on a scheme $Y$ over $\Fp$ we write
$H^i(Y,K)$ and $H^i_c(Y,K)$ for the groups computed on $Y_{\Fpb}$. We write
$S=S_n$ when $n$ is fixed.

\subsection{Boolean algebras}\label{subsec:G-boolean}

The Boolean algebra $\mathsf A_n$ is a bookkeeping device for the unipotent local monodromy of a tensor product of $n$ two-dimensional representations: the operator $x_i$ models the nilpotent part in the $i$-th tensor factor.  Its sum $\sigma=\sum_i x_i$ models the logarithm of total monodromy.  Consequently $\mathsf E_n=\ker\sigma$ records invariants, whereas $\mathsf B_n=\mathsf A_n/\sigma\mathsf A_n$ records the dual quotient.  The Koszul complex does not appear by coincidence: its differential consists of multiplication by the same $x_i$ that represent the local Kummer extension classes.  The point of the Betti numbers $\beta_{j,e}$ is to measure the homology of this explicitly described algebraic complex, not to posit unexplained cohomology dimensions.

Let $\mathsf A_n=\mathbb Q[x_1,\dots,x_n]/(x_1^2,\dots,x_n^2)$, graded by
degree, $\sigma=x_1+\cdots+x_n$, $\mathsf E_n=\{a\in\mathsf A_n:\sigma a=0\}$
and $\mathsf B_n=\mathsf A_n/\sigma\mathsf A_n$; these are graded, with graded
pieces $\mathsf A_{n,d}$, $\mathsf E_{n,d}$, $\mathsf B_{n,d}$. The polynomial ring
$R_n=\mathbb Q[x_1,\dots,x_n]$ acts on all three. Let
\[
\beta_{j,e}=\beta_{j,e}(n)=\dim_{\mathbb Q}\operatorname{Tor}^{R_n}_j(\mathbb Q,\mathsf B_n)_e
\]
be the graded Betti numbers of $\mathsf B_n=R_n/(x_1^2,\dots,x_n^2,x_1+\cdots+x_n)$.
For a graded $R_n$-module $M$ let $K^\bullet(M)$ be the Koszul complex
$K^m(M)=\bigwedge^m\mathbb Q^n\otimes M$, with differential
$\sum_ie_i\wedge(\cdot)\otimes x_i$, where $e_1,\dots,e_n$ is the standard
basis of $\mathbb Q^n$; it is graded by $\deg(e_J\otimes a)=\deg a-|J|$, and we
write $H^m(K^\bullet(M))_\delta$ for the graded pieces of its cohomology.

The internal grading and the homological grading must be read separately.  The symbol $d$ denotes polynomial degree in $\mathsf A_n$; the symbol $m$ denotes exterior (Koszul) degree; and $e$ denotes internal degree after passage to the graded resolution.  The equality between the cohomology of $K^\bullet(\mathsf E_n)$ and the Betti numbers of $\mathsf B_n$ is proved by dualizing the perfect top-degree pairing, which turns exterior multiplication into contraction.  Thus the identity of dimensions is not an unverified match of two tables: it is induced by an explicit duality of complexes.

\begin{lemma}\label{lem:G-boolean}
\begin{enumerate}[label=\textup{(\roman*)}]
\item $\dim\mathsf E_{n,d}=m_{n,n-d}$ for $d\ge n/2$ and $\mathsf E_{n,d}=0$ for
$d<n/2$; $\dim\mathsf B_{n,d}=m_{n,d}$ for $d\le n/2$ and $\mathsf B_{n,d}=0$
for $d>n/2$.
\item The multiplication pairing $\mathsf A_{n,d}\times\mathsf A_{n,n-d}\to
\mathsf A_{n,n}=\mathbb Qx_1\cdots x_n$ induces $R_n$-compatible isomorphisms
$\mathsf E_{n,d}\cong\mathsf B_{n,n-d}^\vee$, and
\[
\dim H^m(K^\bullet(\mathsf E_n))_\delta=\beta_{m,n-\delta}\qquad(m,\delta\in\mathbb Z).
\]
\item $\beta_{j,e}=0$ unless $j\le e\le j+\lfloor n/2\rfloor$;
$\beta_{0,e}=1$ for $e=0$ and $0$ otherwise; for every $e$,
$\sum_j(-1)^j\beta_{j,e}$ is the coefficient of $t^e$ in $(1-t)^n(1+E_n(t))$;
and if $n=2\nu$ is even, then $\beta_{n,3\nu}=C_\nu$.
\end{enumerate}
\end{lemma}

\begin{proof}
(i) Identify $\mathsf A_n\otimes\Ql$ with the space $V$ of \cref{lem:B-CG}
through $e^{(i)}_1\mapsto x_i$, $e^{(i)}_2\mapsto1$; then $N_i$ becomes
multiplication by $x_i$, and $e_J$ becomes a monomial of degree $n-|J|$. By
\cref{lem:B-CG}, $\ker\sigma=\mathsf E_n$ is the sum of spaces of dimension
$m_{n,k}$ lying in degree $n-k$, $0\le k\le n/2$. This is the first assertion,
and the rank of $\sigma\colon\mathsf A_{n,d-1}\to\mathsf A_{n,d}$ is
$\binom n{d-1}-\dim\mathsf E_{n,d-1}$, so that
$\dim\mathsf B_{n,d}=\binom nd-\binom n{d-1}+\dim\mathsf E_{n,d-1}$; this
gives the second, as $\binom n{d-1}-\binom nd=\dim\mathsf E_{n,d-1}$ when
$d-1\ge n/2$.

(ii) The pairing $\langle a,b\rangle$, the coefficient of $x_1\cdots x_n$ in
$ab$, is perfect and satisfies $\langle x_ia,b\rangle=\langle a,x_ib\rangle$. As
$\langle a,\sigma b\rangle=\langle\sigma a,b\rangle$, the space $\mathsf E_{n,d}$ is
the orthogonal of $\sigma\mathsf A_{n,n-d-1}$, and so
$\mathsf E_{n,d}\cong(\mathsf A_{n,n-d}/\sigma\mathsf A_{n,n-d-1})^\vee$, with $x_i$
acting on the right by the transpose of multiplication. Hence the dual of
$K^\bullet(\mathsf E_n)$ is the complex $\bigwedge^\bullet\mathbb Q^n\otimes\mathsf B_n$
with differential $\sum_i\iota_{e_i}\otimes x_i$, $\iota_{e_i}$ being the contraction
transposed to $e_i\wedge$; this is the Koszul complex of $\mathsf B_n$, whose
homology is $\operatorname{Tor}^{R_n}_\bullet(\mathbb Q,\mathsf B_n)$, the Koszul complex
of $x_1,\dots,x_n$ being a free resolution of $\mathbb Q$. The vector dual to
$e_J\otimes a$, with $|J|=m$ and $a\in\mathsf E_{n,d}$, has internal degree
$m+n-d=n-\delta$, where $\delta=d-m$; this gives the formula.

(iii) The $j$-th term of the Koszul complex of $\mathsf B_n$ lives in degrees
$j,\dots,j+\lfloor n/2\rfloor$, by (i). As $\mathsf B_n$ is cyclic, generated in
degree $0$, $\operatorname{Tor}_0=\mathbb Q$ in degree $0$. The Euler characteristic of
the degree-$e$ part of the Koszul complex is
$\sum_j(-1)^j\binom nj\dim\mathsf B_{n,e-j}$, the coefficient of $t^e$ in
$(1-t)^n\sum_dm_{n,d}t^d=(1-t)^n(1+E_n(t))$. Finally, for a graded $R_n$-module
$M$ of finite length, $\operatorname{Tor}_n^{R_n}(\mathbb Q,M)$ is the socle of $M$ shifted by
$n$, the kernel of the last Koszul differential consisting of the
$e_{\{1,\dots,n\}}\otimes a$ with $x_ia=0$ for all $i$. For $n=2\nu$, the top
degree $\mathsf B_{n,\nu}$ lies in the socle and has dimension
$m_{2\nu,\nu}=C_\nu$; so $\beta_{n,3\nu}=C_\nu$.
\end{proof}

The Betti numbers $\beta_{j,e}$ for $n\le9$ are listed in \cref{tab:betti};
they were computed by exact linear algebra from the Koszul complex of
$\mathsf B_n$ (the computational discussion in \cite{EbaCoh}).

\begin{proposition}[Betti numbers of $\mathsf B_n$]\label{prop:G-betti}
Let $\nu=\lfloor n/2\rfloor$, let $\chi_e$ be the coefficient of $t^e$ in
$(1-t)^n(1+E_n(t))$, and for an integer $N$ put $\lambda^{(N)}_e=(-1)^{\lceil e/2\rceil}
\binom n{\lfloor e/2\rfloor}$ if $0\le e\le N$ and $\lambda^{(N)}_e=0$ otherwise.
\begin{enumerate}[label=\textup{(\roman*)}]
\item For $0\le r\le\lceil n/2\rceil-2$ one has $\beta_{r,2r}=\beta_{r+1,2r+1}=\binom nr$,
and $\beta_{j,j+r}=0$ for all other $j$.
\item If $n=2\nu+1$ is odd, then
$\beta_{j,j+\nu}=(-1)^j\bigl(\chi_{j+\nu}-\lambda^{(2\nu-1)}_{j+\nu}\bigr)$ for every $j$;
together with \textup{(i)} this is the whole Betti table of $\mathsf B_n$.
\item If $n=2\nu$ is even, then
$\beta_{j,j+\nu-1}=\beta_{j-1,j-1+\nu}+(-1)^j\bigl(\chi_{j+\nu-1}-\lambda^{(2\nu-3)}_{j+\nu-1}\bigr)$
for every $j$; so the Betti table is determined by its top row
$(\beta_{j,j+\nu})_j$.
\end{enumerate}
\end{proposition}

\begin{proof}
Write $R=R_n$. The complex $C=[\mathsf A_n(-1)\xrightarrow{\sigma}\mathsf A_n]$, in
homological degrees $1$ and $0$, has $H_0(C)=\mathsf B_n$ and $H_1(C)=\mathsf E_n(-1)$,
whence a distinguished triangle
$\mathsf E_n(-1)[1]\to C\to\mathsf B_n\to\mathsf E_n(-1)[2]$. Now
$C=\mathsf A_n\otimes_RK(\sigma)$, where $K(\sigma)=[R(-1)\xrightarrow{\sigma}R]$; if $P$
is the Koszul resolution of $\mathsf A_n$ on $x_1^2,\dots,x_n^2$, then
$\mathbb Q\otimes_R(P\otimes_RK(\sigma))=(\mathbb Q\otimes_RP)\otimes_{\mathbb Q}
(\mathbb Q\otimes_RK(\sigma))$ has zero differential, so
$\dim\operatorname{Tor}_m^R(\mathbb Q,C)_e=\binom nm[e=2m]+\binom n{m-1}[e=2m-1]$. The
triangle gives exact sequences
\[
\operatorname{Tor}_{m-1}(\mathsf E_n)_{e-1}\to\operatorname{Tor}_m(C)_e\to
\operatorname{Tor}_m(\mathsf B_n)_e\to\operatorname{Tor}_{m-2}(\mathsf E_n)_{e-1}
\]
(all $\operatorname{Tor}$ over $R$ against $\mathbb Q$). As $\mathsf E_n$ lives in degrees
$\ge\lceil n/2\rceil$, $\operatorname{Tor}_k(\mathsf E_n)_f=0$ unless
$f-k\ge\lceil n/2\rceil$; so for $e-m\le\lceil n/2\rceil-2$ both outer terms vanish
and $\operatorname{Tor}_m(\mathsf B_n)_e\cong\operatorname{Tor}_m(C)_e$, which is (i).

The rows $r>\nu$ vanish, as $\mathsf B_n$ lives in degrees $\le\nu$. If $n$ is odd,
the rows $r\le\nu-1$ are given by (i), and in internal degree $e$ the row $\nu$
has the single entry $\beta_{e-\nu,e}$, while the rows $r\le\nu-1$ contribute
$\binom n{e/2}$ in homological degree $e/2$ if $e$ is even and $e\le2\nu-2$, and
$\binom n{(e-1)/2}$ in homological degree $(e+1)/2$ if $e$ is odd and $e\le2\nu-1$.
Since $\sum_m(-1)^m\beta_{m,e}=\chi_e$ (\cref{lem:G-boolean}(iii)), this gives (ii).
If $n$ is even, the rows $r\le\nu-2$ are given by (i), and in degree $e$ the rows
$\nu-1$ and $\nu$ have the entries $\beta_{e-\nu+1,e}$ and $\beta_{e-\nu,e}$; the same
count gives (iii).
\end{proof}

For instance, for $n=9$ part (ii) gives the last row
$42,330,1233,2161,2035,1107,330,42$ of \cref{tab:betti}. For even $n\le8$ the top
row of the Betti table of $\mathsf B_n$ is the last row of that of $\mathsf B_{n-1}$,
shifted by one column (\cref{tab:betti}).

\subsection{The part of lower weight}\label{subsec:G-lower}

The filtration of $\cL_n=\cH_n^{\le n}$ is a filtration of sheaves, so its successive extensions can carry nontrivial global information even though each graded quotient is a constant Tate sheaf.  The Kummer class of the unit $h_i$ is precisely the extension data around the divisor $h_i=0$; taking cohomology of the torus converts those classes into cup products.  The resulting weight spectral sequence is therefore not a formal direct sum of the cohomologies of the graded pieces.  Its first differential must be computed from the extension classes, and only after doing so may one identify the $E_2$-page with a Koszul homology group.

Let $\cL_n=\cH_n^{\le n}$ be the lisse subsheaf of weights $\le n$ of $\cH_n$
(\cref{thm:B-family}(iii)), with its weight filtration $W_\bullet$, and let
$\kappa_i\in H^1(S_{\Fpb},\Ql(1))$ be the Kummer class of the unit $h_i$. By the
Kummer sequence, as in the proof of \cref{prop:B-detfinite},
$\kappa_1,\dots,\kappa_n$ form a basis of $H^1(S_{\Fpb},\Ql(1))$. By the Künneth
formula, $H^n_c(S,\Ql)$ is a line with trivial Frobenius, and cup product gives
isomorphisms
\begin{equation}\label{eq:G-cup}
\textstyle\bigwedge^m\bigl(\bigoplus_i\Ql\kappa_i\bigr)\otimes H^n_c(S,\Ql)
\xrightarrow{\ \sim\ }H^{n+m}_c(S,\Ql(m)),
\end{equation}
since for $n=1$ the cup product with $\kappa$ maps $H^1_c(\mathbb G_{m,\Fpb},\Ql)$
isomorphically onto $H^2_c(\mathbb G_{m,\Fpb},\Ql(1))$ (Poincaré duality on
$\Gm$). In particular $H^{n+m}_c(S,\Ql)\cong\Ql(-m)^{\binom nm}$.

\begin{proposition}\label{prop:G-lower}
\begin{enumerate}[label=\textup{(\roman*)}]
\item The geometric generic stalk of $\cL_n$ can be identified with
$\mathsf E_n\otimes\Ql$ in such a way that $W_{2k}\cL_n$ corresponds to
$\bigoplus_{d\ge n-k}\mathsf E_{n,d}\otimes\Ql$. Thus
$\gr^W_{2k}\cL_n\cong\mathsf E_{n,n-k}\otimes\Ql(-k)$, a geometrically constant
sheaf on which $\Frob_p$ acts on $\mathsf E_{n,n-k}$ trivially.
\item Under this identification, the class of the extension
\[
0\to\gr^W_{2k-2}\cL_n\to W_{2k}\cL_n/W_{2k-4}\cL_n\to\gr^W_{2k}\cL_n\to0
\]
in
\[
\operatorname{Ext}^1_{S_{\Fpb}}\bigl(\gr^W_{2k}\cL_n,\gr^W_{2k-2}\cL_n\bigr)
=H^1(S_{\Fpb},\Ql(1))\otimes\Hom(\mathsf E_{n,n-k},\mathsf E_{n,n-k+1})\otimes\Ql
\]
equals $\theta\sum_{i=1}^n\kappa_i\otimes x_i$, for a constant
$\theta\in\Ql^\times$ independent of $k$.
\end{enumerate}
\end{proposition}

\begin{proof}
Let $\bar\jmath\colon\mathbb G_{m,t}\times S\hookrightarrow\mathbb A^1_t\times S$ be the
inclusion, $\cM$ the sheaf of \ref{H7} and $\cM_0=(\bar\jmath_*\cM)|_{\{0\}\times S}$.
The sheaf $\cL_\psi(\lambda t)$ extends to a lisse sheaf on $\mathbb A^1\times S$,
and each $[(t,h)\mapsto h_it]^*\KL$ is tamely ramified along $\{0\}\times S$, by
\ref{K2}, since $h_it$ is a uniformizer of this divisor times a unit. By
Abhyankar's lemma \cite[Exp.~XIII, \S5]{SGA1}, $\cM_0$ is lisse, with stalk
$\cM_h^{I_0}$ at a geometric point $h$. As in the proof of \cref{prop:R-ray}
(the part of lower weight), the exact sequence
$0\to\bar\jmath_!\cM\to\bar\jmath_*\cM\to\cM_0\to0$ and the projection to $S$ give a
morphism of lisse sheaves $\delta\colon\cM_0\to\cH_n$ which at every geometric
point above a closed point is the injection $\cM_h^{I_0}\to H^1_c(\Gm,\cM_h)$ of
\cref{prop:B-Mh}(iii), with image the stalk of $\cL_n$
(\cref{thm:B-family}(v)). Hence $\delta$ is an isomorphism $\cM_0\cong\cL_n$.

Fix a power $q$ of $p$, let $T$ be the henselization of $\mathbb A^1_t\times\mathbb A^n_h$
over $\Fq$ at the origin, and $T^\circ$ the complement in $T$ of the strict normal
crossings divisor $th_1\cdots h_n=0$. By \ref{K2}, the restriction of $\KL$ to the
punctured henselization of $\mathbb A^1$ at $0$ is tame and geometrically
unipotent; and $[(t,h)\mapsto h_it]^*\KL$ is its pull-back by $(t,h)\mapsto h_it$,
which maps $T$ into this henselization. So $\cM|_{T^\circ}$ is tamely ramified
along the divisor, and by Abhyankar's lemma the geometric tame fundamental
group of $T^\circ$ is $\hat{\mathbb Z}'(1)^{n+1}$, generated by the inertia groups
$\langle\gamma_t\rangle$, $\langle\gamma_1\rangle,\dots,\langle\gamma_n\rangle$ of
the components $t=0$, $h_1=0,\dots,h_n=0$. The function $h_it$ has valuation $1$
along $t=0$ and $h_i=0$ and $0$ along the other components. Hence, if $V_i$ is
the generic stalk of $[(t,h)\mapsto h_it]^*\KL$ and $N_i$ the logarithm of the
local monodromy of $\KL$ at $0$ acting on it, then $\gamma_t$ acts on
$V=V_1\otimes\cdots\otimes V_n$ by $\exp(\sum_iN_i)$ and $\gamma_i$ by $\exp(N_i)$.
For a sheaf which is tame along a strict normal crossings divisor, the stalks of
$\bar\jmath_*$ at points of one component which lie on no other component are the
invariants of the inertia group of that component; therefore the restriction of
$\cM_0\cong\cL_n$ to $T\cap\{t=0\}$, which is lisse off the coordinate hyperplanes,
has generic stalk $V^{\gamma_t}=\ker\sum_iN_i$, on which $\gamma_i$ acts by
$\exp(N_i)$. As in the proof of \cref{prop:B-Mh}(iv), \cref{lem:B-frob} and
\ref{K5} give bases $e^{(i)}_1,e^{(i)}_2$ of the $V_i$ with $N_ie^{(i)}_1=0$,
$N_ie^{(i)}_2=e^{(i)}_1$ and $\Phi e_J=q^{|J|}e_J$ for a geometric Frobenius
$\Phi$ of the decomposition group at the origin. Through the identification of
the proof of \cref{lem:G-boolean}(i), $V\cong\mathsf A_n\otimes\Ql$, $N_i$ is
multiplication by $x_i$, $\ker\sum_iN_i=\mathsf E_n\otimes\Ql$, and $\Phi$ acts on
$\mathsf E_{n,d}$ by $q^{n-d}$. Since $\gr^W_{2k}\cL_n(k)$ is geometrically
constant with trivial Frobenius (\cref{thm:B-family}(iii)), $W_{2k}\cL_n$
corresponds to $\bigoplus_{d\ge n-k}\mathsf E_{n,d}\otimes\Ql$. This proves (i), the
geometrically constant sheaves $\gr^W_{2k}\cL_n(k)$ being identified with
$\mathsf E_{n,n-k}\otimes\Ql$ through their stalks at the generic point of
$T\cap\{t=0\}$.

(ii) Let $T_h^\circ$ be the complement of the coordinate hyperplanes in the
strict henselization of $\mathbb A^n_h$ at the origin, over $\Fpb$. Its tame
fundamental group is generated by $\gamma_1,\dots,\gamma_n$, and the restriction
map $H^1(S_{\Fpb},\Ql(1))\to H^1(T_h^\circ,\Ql(1))=\Hom(\hat{\mathbb
Z}'(1)^n,\Ql(1))$ is injective, because $\kappa_i(\gamma_j)=0$ for $j\ne i$ and
$\kappa_i(\gamma_i)$ is a generator. By the previous step, the restriction of the
extension in (ii) to $T_h^\circ$ is the extension of representations given by
$\gamma_i\mapsto1+x_i$ on $\mathsf E_{n,n-k+1}\oplus\mathsf E_{n,n-k}$, whose class is
the homomorphism $\gamma_i\mapsto x_i$. This is the restriction of
$\theta\sum_i\kappa_i\otimes x_i$, $\theta^{-1}$ being the common value of the
$\kappa_i(\gamma_i)$ under the identification $\Ql(1)\cong\Ql$; injectivity gives
(ii).
\end{proof}

\begin{lemma}[Uniqueness of the Koszul extension class]
\label{lem:G-extension-uniqueness-revision}
Retain the assumptions and the arithmetic lisse model in
\cref{prop:G-lower}.  The geometric extension class between
consecutive weight grades is completely determined by the logarithms
of the local tame inertia operators at $h_i=0$.  More precisely, if
$E_k=\mathsf E_{n,n-k}\otimes\Ql(-k)$, then
\[
 \operatorname{Ext}^1_{S_{\Fpb}}(E_k,E_{k-1})
 \cong\bigoplus_{i=1}^n\Ql\,\kappa_i\otimes
 \operatorname{Hom}(\mathsf E_{n,n-k},\mathsf E_{n,n-k+1}),
\]
and the $i$-th coefficient of the extension of
\cref{prop:G-lower}(ii) is a common nonzero normalization of
$x_i|_{\mathsf E_{n,n-k}}$.  No other \emph{geometric} class can be
added while keeping the given tame monodromy.  This assertion alone
does not determine arithmetic self-extensions with trivial geometric
monodromy.
\end{lemma}
\begin{proof}
Over $\Fpb$, the Kummer sequence and the K\"unneth formula give
\[
 H^1(S_{\Fpb},\Ql(1))=\bigoplus_{i=1}^n\Ql\,\kappa_i.
\]
The two graded coefficient objects are geometrically constant; thus
extensions of one by the other are first cohomology with constant
$\operatorname{Hom}$ coefficients.  In the strict henselization of
$\mathbb A^n$ at its origin with the coordinate hyperplanes removed,
tame inertia in direction $i$ acts on the model of
\cref{prop:G-lower} by $\exp(x_i)$, with $x_i^2=0$.  The induced
extension between the two adjacent graded pieces has logarithm
$x_i$; the Kummer generator $\kappa_j$ restricts to zero on that
inertia direction if $j\ne i$ and to a common nonzero generator for
$j=i$.  Hence all $n$ coefficients are uniquely determined.  The
restriction map on $H^1$ is injective because these $n$ Kummer classes
are linearly independent after restriction.  The claim concerns the
geometric fundamental group; arithmetic extensions are separately
addressed in the arithmetic-extension discussion of \cite{EbaZen}.
\end{proof}

\begin{lemma}\label{lem:G-Hn}
\begin{enumerate}[label=\textup{(\roman*)}]
\item $H^i_c(S,\cH_n)=0$ for $i\ne n$, and $H^n_c(S,\cH_n)$ is a line on which
$\Frob_p$ acts trivially.
\item The constant subsheaf $\Ql\subset\cH_n$ of \cref{thm:B-family}(ii)
induces an isomorphism $H^n_c(S,\Ql)\to H^n_c(S,\cH_n)$.
\end{enumerate}
\end{lemma}

\begin{proof}
(i) By \cref{thm:B-family}(i), $R\Gamma_c(S,\cH_n)[-1]\cong
R\Gamma_c(\mathbb G_{m,t}\times S,\cM)$. The automorphism $(t,h)\mapsto(t,(h_it)_i)$
of $\Gm\times\Gm^n$ carries $\cM$ to $\cL_\psi(\lambda t)\boxtimes\KL^{\boxtimes n}$,
so by the Künneth formula $R\Gamma_c(\Gm\times S,\cM)\cong
R\Gamma_c(\Gm,\cL_\psi(\lambda t))\otimes R\Gamma_c(\Gm,\KL)^{\otimes n}$. The
sheaf $\cL_\psi(\lambda t)$ on $\Gm$ has Euler characteristic $-1$ and no
$H^0_c$ or $H^2_c$, and $-\sum_{t\ne0}\psi_q(\lambda t)=1$; by the proof of
\ref{K5}, $H^*_c(\Gm,\KL)$ is a line in degree $1$, and
$\sum_{\xi\ne0}\Kl_q(\xi)=1$. So both factors are lines in degree $1$ with
trivial Frobenius, and $R\Gamma_c(S,\cH_n)\cong\Ql[-n]$.

(ii) Use the notation of \cref{app:family}. The constant subsheaf is the image
of the connecting map $\partial\colon R^n\pi_{Z!}(\cL|_Z)\cong\Ql\to\cH_n$, so the
map of (ii) is the connecting map
$H^{2n}_c(Z,\cL|_Z)\to H^{2n+1}_c(X^\circ,\cL)$ of the excision sequence of
$X^\circ\subset X\supset Z$, whose source and target are $H^n_c(S,\Ql)$ and
$H^n_c(S,\cH_n)$ by the Leray spectral sequences. Its kernel is the image of
$H^{2n}_c(X,\cL)$. By \cref{lem:B-int}, applied over $S$, $R\Gamma_c(X,\cL)\cong
R\Gamma_c(Y,\cL_\psi(f))(-1)[-2]$ with
$Y=\{(u,h)\in\Gm^n\times S:\lambda+\langle h,u\rangle=0\}$, and $(u,h)\mapsto
(u,(h_iu_i)_i)$ identifies $Y$ with $\Gm^n\times Y_0$,
$Y_0=\{y\in\Gm^n:\sum_iy_i=-\lambda\}$, carrying $\cL_\psi(f)$ to
$\cL_\psi(f)\boxtimes\Ql$. By \cref{lem:B-t0} and the Künneth formula,
$H^i_c(X,\cL)\cong H^{i-n-2}_c(Y_0,\Ql)(-1)$. As $Y_0$ is smooth and affine of
dimension $n-1$, $H^{n-2}_c(Y_0,\Ql)=0$; so $H^{2n}_c(X,\cL)=0$, and the
connecting map is an injection of lines, hence an isomorphism.
\end{proof}

\begin{lemma}\label{lem:G-L}
For every $m$, $H^{n+m}_c(S,\cL_n)\cong\bigoplus_e\Ql(-e)^{\beta_{m,e}}$,
$\Frob_p$ acting semisimply; and $H^i_c(S,\cL_n)=0$ for $i<n$ and $i>2n$.
\end{lemma}

\begin{proof}
The filtration $W_\bullet$ of $\cL_n$ by lisse subsheaves gives a spectral
sequence $E_1^{-k,i+k}=H^i_c(S,\gr^W_{2k}\cL_n)\Rightarrow H^i_c(S,\cL_n)$. By
\cref{prop:G-lower}(i) and \eqref{eq:G-cup},
\[
E_1^{-k,n+m+k}\cong\textstyle\bigwedge^m\Ql^n\otimes\mathsf E_{n,n-k}\otimes\Ql(-m-k),
\]
on which $\Frob_p$ acts semisimply by $p^{m+k}$. The differential $d_1$ is the
connecting homomorphism of the extension of \cref{prop:G-lower}(ii), that is,
cup product with its class; by \cref{prop:G-lower}(ii) and \eqref{eq:G-cup},
$(E_1,d_1)$ is, up to the constant $\theta$ and signs, the Koszul complex
$K^\bullet(\mathsf E_n)\otimes\Ql$, the summand $\bigwedge^m\otimes\mathsf E_{n,d}$ of
degree $\delta=d-m$ lying in $E_1^{-(n-d),2n+m-d}$, where Frobenius acts by
$p^{n-\delta}$. The differential $d_r$ maps $E_r^{-k,i+k}$, which is pure of
weight $2(i-n+k)$, to $E_r^{-k+r,i+k-r+1}$, of weight $2(i-n+k)-2(r-1)$. For
$r\ge2$ these weights differ, and a Frobenius-equivariant map between pure
modules of different weights vanishes; so $E_2=E_\infty$. The graded pieces of
the abutment filtration of $H^{n+m}_c(S,\cL_n)$ are the $E_2^{-k,n+m+k}$, pure of
the pairwise distinct weights $2(m+k)$, so $H^{n+m}_c(S,\cL_n)$ is their direct
sum, $\Frob_p$ acting on $E_2^{-k,n+m+k}$ as the scalar $p^{m+k}$. By
\cref{lem:G-boolean}(ii), the part of degree $\delta$ has dimension
$\beta_{m,n-\delta}$, and $e=n-\delta=m+k$.
\end{proof}

\subsection*{Weight-spectral-sequence differentials}
The geometric extension classes must first be identified in
$H^1(S_{\overline{\mathbb F}_p},\overline{\mathbb Q}_\ell(1))$;
the classes of the coordinate Kummer coverings provide its
$n$ distinguished generators.  Restricting the extension to
each coordinate inertia group determines the corresponding
coefficient of the logarithm of tame monodromy, but this
comparison concerns the actual representation and not only
its Grothendieck-group class.  The first differential in
the weight spectral sequence is then the sum of the maps
induced by these extension classes, identified with the
Koszul differential in the chosen basis.  For a higher
differential $d_r$, both its source and target must be pure
Frobenius modules of distinct weights after taking $E_2$;
only then does Frobenius equivariance force it to vanish.
Finally, when two surviving degrees have distinct
Frobenius eigenvalue sets, their extensions split as
Frobenius modules by coprime minimal polynomials.  None of
these arguments establishes semisimplicity of an unrelated
arithmetic coefficient with repeated eigenvalues.  Each
use of this weight argument therefore requires the
particular weight indices to be checked, not just a
blanket claim of degeneration.

\subsection*{Degeneration and Frobenius weights}
The first page records contributions of each pure Tate graded piece of $\cL_n$; it does not yet compute the cohomology of $\cL_n$, because its connecting maps need not vanish.  The $d_1$ differential is the Koszul differential just identified, and therefore its cohomology is $E_2$.  Once that computation has been made, a differential $d_r$ with $r\ge2$ would have to map one Frobenius weight to a different Frobenius weight.  Such a map is zero, and the spectral sequence stops.  At fixed total degree, the surviving filtration pieces have pairwise distinct eigenvalues, which supplies a Frobenius-stable splitting.  These three statements---identification of $d_1$, exclusion of higher differentials, and splitting by distinct eigenvalues---are separate arguments.  In particular, they do not establish arithmetic semisimplicity of unrelated boundary coefficient sheaves whose Jordan blocks may have a single repeated eigenvalue.

\begin{lemma}[Higher differentials and the Frobenius splitting]
\label{lem:G-higher-differentials-revision}
For the weight spectral sequence of $\cL_n$ in
\cref{lem:G-L}, once the $E_1$-page and $d_1$ have been identified as
in \cref{prop:G-lower,lem:G-extension-uniqueness-revision}, every
differential $d_r$, $r\ge2$, vanishes.  Moreover its abutment is a
\emph{Frobenius-semisimple} direct sum of the indicated Tate graded
pieces.  Neither conclusion is a statement about the Frobenius on
the higher-codimension boundary sheaf $B'_n$.
\end{lemma}
\begin{proof}
On $E_1^{-k,n+m+k}$ the geometric Frobenius is scalar with
eigenvalue $p^{m+k}$.  The differential $d_1$ is Frobenius
compatible and the Koszul differential, so $E_2$ inherits scalar
Frobenius on each bidegree.  In the notation
$E_r^{-k,i+k}$ of \cref{lem:G-L}, a differential $d_r$ has target
$E_r^{-k+r,i+k-r+1}$ and changes the weight from
$2(i-n+k)$ to $2(i-n+k-r+1)$.  These are unequal for $r\ge2$.
Inductively $E_r=E_2$ for all such $r$, so a Frobenius-equivariant
map between these two pure scalar-eigenvalue spaces is zero.  At
fixed cohomological degree $n+m$, distinct filtration grades $k$
have the distinct eigenvalues $p^{m+k}$; their generalized
eigenspaces give a canonical Frobenius-stable splitting of the
abutment.  This excludes extension blocks between these \emph{different}
weights.  It does not exclude extensions of \emph{equal}-weight
boundary coefficient objects arising from another spectral sequence.
\end{proof}

\begin{theorem}[Cohomology on the torus]\label{thm:G-torus}
Let $n\ge2$.
\begin{enumerate}[label=\textup{(\roman*)}]
\item $H^i_c(S,\cW_n)=0$ for $i<n$ and for $i\ge2n$, and for $0\le j\le n-1$
\[
H^{n+j}_c(S,\cW_n)\cong\bigoplus_e\Ql(-e)^{\beta_{j+1,e}},
\]
$\Frob_p$ acting semisimply. In particular every eigenvalue of Frobenius on
$H^*_c(S,\cW_n)$ is an integral power of $p$, and $\chi_c(S,\cW_n)=(-1)^n$.
\item $H^i(S,\cW_n)\cong H^{2n-i}_c(S,\cW_n)^\vee(-2n-1)$ for every $i$; thus
$H^i(S,\cW_n)\ne0$ only for $1\le i\le n$, and its Frobenius eigenvalues are
the numbers $p^{2n+1-e}$, $p^e$ running through those of $H^{2n-i}_c$.
\item The interior cohomology $\imop\bigl(H^i_c(S,\cW_n)\to H^i(S,\cW_n)\bigr)$
vanishes for every $i$.
\item Write $(X-1)^n(1+E_n(X))=\sum_j\alpha_{n,j}X^j$. Then
\begin{gather*}
L(S_n,\cW_n,T)=\prod_{e\ge1}(1-p^eT)^{\alpha_{n,e}},\\
L^{\mathrm{ord}}(S_n,\cW_n,T)=\bigl(-p^{2n+1}T\bigr)^{(-1)^{n+1}}
L\Bigl(S_n,\cW_n,\frac1{p^{2n+1}T}\Bigr),
\end{gather*}
where $L^{\mathrm{ord}}(S_n,\cW_n,T)=\prod_i\det(1-T\Frob_p\mid H^i(S,\cW_n))^{(-1)^{i+1}}$.
\end{enumerate}
\end{theorem}

\begin{proof}
(i) By \cref{lem:G-Hn}(i), the long exact sequence of
$0\to\cL_n\to\cH_n\to\cW_n\to0$ gives
$H^i_c(S,\cW_n)\cong H^{i+1}_c(S,\cL_n)$ for $i\ne n-1,n$, and an exact
sequence
\[
H^{n-1}_c(S,\cW_n)\to H^n_c(S,\cL_n)\to H^n_c(S,\cH_n)\to H^n_c(S,\cW_n)\to
H^{n+1}_c(S,\cL_n)\to0 .
\]
The constant subsheaf of \cref{thm:B-family}(ii) is pure of weight $0$, so it is
$W_0\cL_n$, of rank $m_{n,0}=1$. The map $H^n_c(S,W_0\cL_n)\to H^n_c(S,\cL_n)$ is
injective, because $H^{n-1}_c(S,\cL_n/W_0\cL_n)=0$ by Artin's vanishing theorem
and Poincaré duality, and both groups are lines (\cref{lem:G-L}, as
$\beta_{0,e}=[e=0]$). By \cref{lem:G-Hn}(ii) the composite with
$H^n_c(S,\cL_n)\to H^n_c(S,\cH_n)$ is an isomorphism; so the latter map is an
isomorphism, and $H^n_c(S,\cW_n)\cong H^{n+1}_c(S,\cL_n)$. Also
$H^{n-1}_c(S,\cW_n)=0$ by Artin's vanishing theorem. With \cref{lem:G-L} this
gives (i), except the Euler characteristic, which follows from (iv): the degree
of $L(S_n,\cW_n,T)$ is $-\chi_c(S,\cW_n)$ and equals
$\sum_{e\ge1}\alpha_{n,e}=-\alpha_{n,0}=(-1)^{n+1}$.

(ii) By Poincaré duality $H^i(S,\cW_n)\cong H^{2n-i}_c(S,\cW_n^\vee)^\vee(-n)$,
and $\cW_n^\vee\cong[-1]^*\cW_n(w)$ by \cref{prop:reflection}(b). As $[-1]$ is
an automorphism of $S$ defined over $\Fp$, it induces Frobenius-equivariant
isomorphisms $H^*_c(S,[-1]^*\cW_n)\cong H^*_c(S,\cW_n)$.

(iii) The image is pure of weight $w+i=n+1+i$ \cite[3.3.1]{Del80}, being a
quotient of a group of weights $\le w+i$ and a subgroup of a group of weights
$\ge w+i$. It vanishes for $i\ne n$, since $H^i_c=0$ for $i<n$ and $H^i=0$ for
$i>n$; and for $i=n$ its weight $2n+1$ is odd, while by (i) all weights of
$H^n_c(S,\cW_n)$ are even.

(iv) By (i), the exponent of $1-p^eT$ in $L(S_n,\cW_n,T)$ is
$\sum_{j\ge0}(-1)^{n+j+1}\beta_{j+1,e}=(-1)^n\sum_{m\ge1}(-1)^m\beta_{m,e}$, which by
\cref{lem:G-boolean}(iii) equals $\alpha_{n,e}-(-1)^n[e=0]$; for $e=0$ this is
$0$, as $\alpha_{n,0}=(-1)^n$. By (ii) the exponent of $1-p^{2n+1-e}T$ in
$L^{\mathrm{ord}}$ is $\alpha_{n,e}$. Since
$1-p^{e-2n-1}T^{-1}=(-p^{e-2n-1}T^{-1})(1-p^{2n+1-e}T)$, we obtain
$L^{\mathrm{ord}}(T)=\prod_e(-p^{2n+1-e}T)^{\alpha_{n,e}}L(1/(p^{2n+1}T))$; and
$\sum_{e\ge1}\alpha_{n,e}=(-1)^{n+1}$, $\sum_ee\alpha_{n,e}=0$, because
$(X-1)^n(1+E_n(X))$ vanishes to order $n\ge2$ at $X=1$.
\end{proof}

\begin{table}[ht]
\centering
\caption{The Betti tables of $\mathsf B_n$ over $\mathbb Q$: the entry in row $r$
and column $j$ is $\beta_{j,j+r}$ (dots denote zeros). By \cref{thm:G-torus}, for
$j\ge1$ column $j$ gives $H^{n+j-1}_c(S_n,\cW_n)$, and row $r$ the eigenvalue
$p^{j+r}$.}
\label{tab:betti}
\small
\begin{tabular}{@{}lll@{}}
\toprule
$n$ & & $j=0,1,2,\dots,n$\\
\midrule
2 & $r=0$ & $1\ \ 1\ \ \cdot$\\
  & $r=1$ & $\cdot\ \ 1\ \ 1$\\[2pt]
3 & $r=0$ & $1\ \ 1\ \ \cdot\ \ \cdot$\\
  & $r=1$ & $\cdot\ \ 3\ \ 5\ \ 2$\\[2pt]
4 & $r=0$ & $1\ \ 1\ \ \cdot\ \ \cdot\ \ \cdot$\\
  & $r=1$ & $\cdot\ \ 4\ \ 6\ \ 2\ \ \cdot$\\
  & $r=2$ & $\cdot\ \ \cdot\ \ 3\ \ 5\ \ 2$\\[2pt]
5 & $r=0$ & $1\ \ 1\ \ \cdot\ \ \cdot\ \ \cdot\ \ \cdot$\\
  & $r=1$ & $\cdot\ \ 5\ \ 5\ \ \cdot\ \ \cdot\ \ \cdot$\\
  & $r=2$ & $\cdot\ \ \cdot\ \ 15\ \ 31\ \ 21\ \ 5$\\[2pt]
6 & $r=0$ & $1\ \ 1\ \ \cdot\ \ \cdot\ \ \cdot\ \ \cdot\ \ \cdot$\\
  & $r=1$ & $\cdot\ \ 6\ \ 6\ \ \cdot\ \ \cdot\ \ \cdot\ \ \cdot$\\
  & $r=2$ & $\cdot\ \ \cdot\ \ 20\ \ 36\ \ 21\ \ 5\ \ \cdot$\\
  & $r=3$ & $\cdot\ \ \cdot\ \ \cdot\ \ 15\ \ 31\ \ 21\ \ 5$\\[2pt]
7 & $r=0$ & $1\ \ 1\ \ \cdot\ \ \cdot\ \ \cdot\ \ \cdot\ \ \cdot\ \ \cdot$\\
  & $r=1$ & $\cdot\ \ 7\ \ 7\ \ \cdot\ \ \cdot\ \ \cdot\ \ \cdot\ \ \cdot$\\
  & $r=2$ & $\cdot\ \ \cdot\ \ 21\ \ 21\ \ \cdot\ \ \cdot\ \ \cdot\ \ \cdot$\\
  & $r=3$ & $\cdot\ \ \cdot\ \ 14\ \ 119\ \ 237\ \ 202\ \ 84\ \ 14$\\[2pt]
8 & $r=0$ & $1\ \ 1\ \ \cdot\ \ \cdot\ \ \cdot\ \ \cdot\ \ \cdot\ \ \cdot\ \ \cdot$\\
  & $r=1$ & $\cdot\ \ 8\ \ 8\ \ \cdot\ \ \cdot\ \ \cdot\ \ \cdot\ \ \cdot\ \ \cdot$\\
  & $r=2$ & $\cdot\ \ \cdot\ \ 28\ \ 28\ \ \cdot\ \ \cdot\ \ \cdot\ \ \cdot\ \ \cdot$\\
  & $r=3$ & $\cdot\ \ \cdot\ \ 14\ \ 140\ \ 258\ \ 202\ \ 84\ \ 14\ \ \cdot$\\
  & $r=4$ & $\cdot\ \ \cdot\ \ \cdot\ \ 14\ \ 119\ \ 237\ \ 202\ \ 84\ \ 14$\\[2pt]
9 & $r=0$ & $1\ \ 1\ \ \cdot\ \ \cdot\ \ \cdot\ \ \cdot\ \ \cdot\ \ \cdot\ \ \cdot\ \ \cdot$\\
  & $r=1$ & $\cdot\ \ 9\ \ 9\ \ \cdot\ \ \cdot\ \ \cdot\ \ \cdot\ \ \cdot\ \ \cdot\ \ \cdot$\\
  & $r=2$ & $\cdot\ \ \cdot\ \ 36\ \ 36\ \ \cdot\ \ \cdot\ \ \cdot\ \ \cdot\ \ \cdot\ \ \cdot$\\
  & $r=3$ & $\cdot\ \ \cdot\ \ \cdot\ \ 84\ \ 84\ \ \cdot\ \ \cdot\ \ \cdot\ \ \cdot\ \ \cdot$\\
  & $r=4$ & $\cdot\ \ \cdot\ \ 42\ \ 330\ \ 1233\ \ 2161\ \ 2035\ \ 1107\ \ 330\ \ 42$\\
\bottomrule
\end{tabular}
\end{table}

\begin{remark}\label{rem:G-torus}
(a) For $n=2$, \cref{thm:G-torus} gives $H^2_c\cong\Ql(-1)\oplus\Ql(-2)$ and
$H^3_c\cong\Ql(-3)$, and
$L(S_2,\cW_2,T)=(1-p^3T)/\bigl((1-pT)(1-p^2T)\bigr)$. For $n\le3$ each
eigenvalue occurs in a single degree, but not in general: for $n=4$ the
eigenvalue $p^4$ occurs in $H^5_c$ and in $H^6_c$, with multiplicities $3$ and
$2$, whereas the exponent of $1-p^4T$ in the $L$-function is $\alpha_{4,4}=1$.
So the $L$-function does not determine the cohomology, and the Betti numbers of
$\mathsf B_n$ are needed.

(b) For even $n=2\nu$, the classes of weight $3n$ in $H^{2n-1}_c(S,\cW_n)$, which
form $\Ql(-3\nu)^{C_\nu}$ by \cref{lem:G-boolean}(iii), are the only classes of
the maximal weight $w+i$ in any $H^i_c(S,\cW_n)$: a class of weight $n+1+i$ in
$H^i_c$, $i=n+j$, has eigenvalue $p^e$ with $2e=2n+1+j$ and
$\beta_{j+1,e}\ne0$, which by \cref{lem:G-boolean}(iii) forces
$e\le j+1+\nu$, that is $j\ge2n-1-2\nu$; so $j=n-1$ and $n$ is even. Dually, the
classes $\Ql(-\nu-1)^{C_\nu}$ of $H^1(S,\cW_n)$ are the only classes of the
minimal weight $n+1+i$ in any $H^i(S,\cW_n)$.

(c) The theorem refines the remark following \cref{thm:rational}: the
zeros and poles of $L(S_n,\cW_n,T)$ lie on the circles $|T|=p^{-e}$, and the
relation (iv) between $L$ and $L^{\mathrm{ord}}$ is an explicit global
functional equation, with sign $-1$, relating the two Euler products of the
torus; it does not relate $L(S_n,\cW_n,T)$ to itself.
\end{remark}

\subsection{The intersection completion}\label{subsec:G-IH}

The purpose of the completion is to replace the nonproper parameter space by a proper geometric object without losing the given local system on its dense torus.  The functor $j_{!*}$ is the intermediate extension in the category of perverse sheaves: it is characterized by having no nonzero perverse subobject or quotient supported entirely on the boundary.  This characterization is weaker than saying that the ordinary restriction to the boundary is zero.  The cohomology groups of $K_n[-n]$ are called intersection cohomology here because they use this intermediate extension and the shift that restores the original degree convention on the open torus.

Let $X_n=(\PP^1)^n$, which contains $S_n$ as the complement of the divisor
$D_n=\bigcup_i\{h_i\in\{0,\infty\}\}$, let $j\colon S_n\hookrightarrow X_n$ be the
inclusion, and put
\begin{gather*}
K_n=j_{!*}(\cW_n[n]),\qquad\mathcal I_n=K_n[-n],\\
IH^i=H^i(X_n,\mathcal I_n),\qquad P_{n,i}(T)=\det(1-T\Frob_p\mid IH^i),
\end{gather*}
$b_{n,i}=\dim IH^i$, $\chi_{\mathrm{IH}}=\sum_i(-1)^ib_{n,i}$, and
\[
\widehat L_n(T)=\prod_iP_{n,i}(T)^{(-1)^{i+1}}
=\prod_{x\in|X_n|}\prod_k\det\bigl(1-T^{\deg x}\Frob_x\mid\mathcal H^k(\mathcal I_n)_{\bar x}\bigr)^{(-1)^{k+1}},
\]
the two products being equal by the Grothendieck--Lefschetz trace formula. The
sheaf $\mathcal H^0(\mathcal I_n)$ is $j_*\cW_n$, and $\mathcal I_n|_{S_n}=\cW_n$.

\subsection*{Weight bookkeeping before the global argument}
It is useful to separate three objects: the sheaf $\mathcal W_n$ on the
$n$-dimensional open torus; the perverse intermediate extension
$K_n=j_{!*}(\mathcal W_n[n])$; and the complex
$\mathcal I_n=K_n[-n]$ used to define $IH^i$ in ordinary cohomological
degrees. A pure sheaf of weight $w$ becomes a pure perverse sheaf of weight
$w+n$ (not $w-n$) after the shift $[n]$: the only nonzero cohomology
sheaf of $\mathcal W_n[n]$ occurs in degree $-n$ and has weight $w$.
The perverse extension has the same weight;
proper cohomology in degree $k$ then has weight $w+n+k$.
Because $IH^i=H^{i-n}(X_n,K_n)$, substitution of $k=i-n$ yields
$w+i=n+1+i$. This accounts for every exponent in the result, and
shows that the indexing of the intersection groups is not arbitrary.
\par
In particular, the proof does not use an explicit dimension formula for
$IH^i$ or a decomposition into contributions from individual boundary
strata. It works for every degree simultaneously. Determining which
$IH^i$ vanish is a separate, finer question. Where the phrase
\emph{Riemann hypothesis for the completion} is used, it always refers
to the family of weight circles indexed by the nonvanishing degrees.

\begin{theorem}[Riemann hypothesis for the intersection completion]\label{thm:G-RH}
For every $n\ge2$ and every $i$, $IH^i$ is pure of weight $n+1+i$: every zero of
$P_{n,i}$ satisfies $|T|=p^{-(n+1+i)/2}$, and the zeros of $P_{n,i}(p^{-s})$
lie on the line $\operatorname{Re}(s)=(n+1+i)/2$. Consequently no zero of one
factor $P_{n,i}$ cancels a zero of another, and $\widehat L_n(T)$ determines every
$P_{n,i}(T)$.
\end{theorem}

\begin{proof}
We make the chain of weight shifts explicit.

\emph{Step 1: the open pure complex.}
By \ref{H1}, the sheaf $\cW_n$ is lisse on the smooth
$n$-dimensional torus $S_n$, and every eigenvalue of Frobenius at a point
of degree $a$ has complex absolute value $p^{a(n+1)/2}$. By the
conventions for perverse sheaves, $\cW_n[n]$ is perverse. A shift by
$[n]$ increases the weight of a pure complex by $n$; thus
$\cW_n[n]$ is pure of weight
\[
(n+1)+n=2n+1.
\]
The two appearances of $n$ have different origins: one comes from the
weight of the sheaf and the other from the dimension of the smooth base.

\emph{Step 2: the extension across the boundary.}
The open immersion $j\colon S_n\hookrightarrow X_n$ admits the
intermediate-extension functor on perverse sheaves. By Gabber's purity
theorem \cite[Cor.~5.4.3]{BBD}, the perverse intermediate extension
\[
K_n=j_{!*}(\cW_n[n])
\]
is pure of the \emph{same} weight $2n+1$. This is the precise theorem that
permits us to pass from the open torus to the compactification without
computing its individual boundary stalks. It does not assert that
those stalks vanish or have a particular Tate decomposition.

\emph{Step 3: take proper cohomology.}
Because $X_n=(\PP^1)^n$ is proper over $\Fp$, proper pushforward of a
pure complex to the point is pure in the cohomological grading:
\[
H^k(X_n,K_n)\quad\text{is pure of weight}\quad2n+1+k.
\]
This is the proper-purity theorem \cite[5.1.14]{BBD},\cite[3.3.1]{Del80}.
In our chosen indexing, however, the intersection cohomology is
\[
IH^i=H^i(X_n,K_n[-n])=H^{i-n}(X_n,K_n).
\]
Set $k=i-n$. The preceding weight becomes
\[
2n+1+(i-n)=n+1+i,
\]
which proves the claimed purity in every degree $i$.

\emph{Step 4: translate weights into circles and vertical lines.}
Let $\alpha$ be a Frobenius eigenvalue on $IH^i$. Purity gives
$|\iota(\alpha)|=p^{(n+1+i)/2}$. A zero $T_0$ of
$P_{n,i}(T)=\det(1-T\Frob_p\mid IH^i)$ has the form
$T_0=\alpha^{-1}$, so
\[
|T_0|=p^{-(n+1+i)/2}.
\]
If $T_0=p^{-s}$, taking absolute values yields
$p^{-\operatorname{Re}s}=p^{-(n+1+i)/2}$, equivalently
$\operatorname{Re}s=(n+1+i)/2$. This implication concerns only the
characteristic polynomial; the proof does not need Frobenius to be
semisimple on $IH^i$.

\emph{Step 5: no cancellation between different degrees.}
If $i\ne j$, then the corresponding circle radii
$p^{-(n+1+i)/2}$ and $p^{-(n+1+j)/2}$ are different. A complex number
cannot be a zero of both $P_{n,i}$ and $P_{n,j}$. Hence the alternating
product defining $\widehat L_n$ has no cancellation between factors
coming from different degrees, and its zeros and poles determine the
individual polynomials $P_{n,i}$ together with their parity.
\end{proof}

The functional equation is a second conclusion of the geometric construction, and should not be used to justify the purity statement that precedes it.  Reflection identifies the weight-normalized dual of the lisse sheaf with its pullback by $h\mapsto-h$; because this involution extends over $(\PP^1)^n$, Verdier duality can be applied to the intermediate extension itself.  Taking cohomology produces pairings in complementary degrees $i$ and $2n-i$, so the exponent $2n+1$ in the functional equation is the weight of these pairings.  The parity of the middle-degree dimension controls the remaining sign, but the pairing with its arithmetic reflection-skew symmetry is the crucial
additional ingredient. It is supplied in \cref{thm:arith-reflection-closed}.

\begin{remark}[A diagnostic check of the two weight shifts]
For $n=2$ the local stalk weight is $3$, so a reciprocal Frobenius eigenvalue has absolute value $q^{-3/2}$.  The intersection group in degree $i$ instead has weight $3+i$, so its corresponding reciprocal eigenvalues have absolute value $p^{-(3+i)/2}$.  This small case illustrates why a universal ``critical line $1/2$'' would be the wrong normalization and why the local and intersection-completion conclusions must be formulated separately.  The calculation checks the exponent in the theorem but is not a replacement for the proof that $\cW_2$ is pure.
\end{remark}

\begin{proposition}[Duality and functional equation]\label{prop:G-duality}
\begin{enumerate}[label=\textup{(\roman*)}]
\item Let $[-1]$ also denote the automorphism $h\mapsto-h$ of $X_n$. Then
$D(K_n)\cong[-1]^*K_n(2n+1)$, and there are Frobenius-equivariant perfect
pairings $IH^i\times IH^{2n-i}\to\Ql(-2n-1)$. In particular
$b_{n,i}=b_{n,2n-i}$, and $\beta\mapsto p^{2n+1}/\beta$ is a bijection between
the eigenvalues of Frobenius on $IH^i$ and on $IH^{2n-i}$.
\item One has
\begin{gather*}
\widehat L_n(T)=\omega_n\bigl(p^{(2n+1)/2}T\bigr)^{-\chi_{\mathrm{IH}}}\,
\widehat L_n\Bigl(\frac1{p^{2n+1}T}\Bigr),\\
\omega_n=p^{-(2n+1)b_{n,n}/2}\det\bigl(-\Frob_p\mid IH^n\bigr)\in\{\pm1\}.
\end{gather*}
In particular $\omega_n=1$ if $IH^n=0$.
\end{enumerate}
\end{proposition}

\begin{proof}
(i) Verdier duality exchanges ${}^pj_!$ and ${}^pj_*$ and preserves perverse
sheaves, so it commutes with the intermediate extension, the image of
${}^pj_!\to{}^pj_*$; moreover $D(\cW_n[n])=\cW_n^\vee[n](n)$, and
$\cW_n^\vee\cong[-1]^*\cW_n(n+1)$ by \cref{prop:reflection}(b). As $[-1]$ extends
to an automorphism of $X_n$ preserving $S_n$, it commutes with $j_{!*}$, and
$D(K_n)\cong[-1]^*K_n(2n+1)$. Verdier duality on the proper variety $X_n$ gives
perfect pairings $H^k(X_n,K_n)\times H^{-k}(X_n,D(K_n))\to\Ql$, and $[-1]$ induces
Frobenius-equivariant isomorphisms $H^{-k}(X_n,[-1]^*K_n)\cong H^{-k}(X_n,K_n)$.
With $k=i-n$ this gives the pairings.

(ii) By (i), if $\alpha_1,\dots,\alpha_b$ are the eigenvalues on $IH^i$, then
$\det(1-T^{-1}p^{-2n-1}\Frob\mid IH^i)=\prod_l(1-(p^{2n+1}/\alpha_l)^{-1}T^{-1})
=(-T)^{-b}\det(\Frob\mid IH^{2n-i})^{-1}P_{n,2n-i}(T)$. Taking the alternating
product over $i$, and using $(-1)^{i+1}=(-1)^{2n-i+1}$, we get
$\widehat L_n(1/(p^{2n+1}T))=(-T)^{\chi_{\mathrm{IH}}}\prod_i\det(\Frob\mid IH^i)^{(-1)^i}
\widehat L_n(T)$. For $i\ne n$, (i) gives $\det(\Frob\mid IH^i)\det(\Frob\mid IH^{2n-i})
=p^{(2n+1)b_{n,i}}$, and $\det(\Frob\mid IH^n)^2=p^{(2n+1)b_{n,n}}$; moreover
$\chi_{\mathrm{IH}}\equiv b_{n,n}\pmod2$. Collecting terms gives the formula,
with $\omega_n=(-1)^{b_{n,n}}u^{\pm1}=(-1)^{b_{n,n}}u$, where
$u=p^{-(2n+1)b_{n,n}/2}\det(\Frob\mid IH^n)=\pm1$.
\end{proof}

\begin{proposition}[The open part and the boundary]\label{prop:G-boundary}
\begin{enumerate}[label=\textup{(\roman*)}]
\item The image of the restriction map $IH^i\to H^i(S_n,\cW_n)$ is the part
$W_{n+1+i}H^i(S_n,\cW_n)$ of lowest weight, and the image of
$H^i_c(S_n,\cW_n)\to IH^i$ is the part of highest weight
$\gr^W_{n+1+i}H^i_c(S_n,\cW_n)$. Both vanish unless $n=2\nu$ is even and
$i=1$, resp.\ $i=2n-1$, when they are $\Ql(-\nu-1)^{C_\nu}$, resp.\
$\Ql(-3\nu)^{C_\nu}$.
\item Let $\mathcal B_n=\mathcal I_n|_{D_n}$, so that the cone of
$j_!\cW_n\to\mathcal I_n$ is $\mathcal B_n$ extended by zero. Then
\begin{gather*}
\widehat L_n(T)=L(S_n,\cW_n,T)\,L(D_n,\mathcal B_n,T),\\
L(D_n,\mathcal B_n,T)=\prod_{x\in|D_n|}\prod_k\det\bigl(1-T^{\deg x}\Frob_x\mid
\mathcal H^k(\mathcal B_n)_{\bar x}\bigr)^{(-1)^{k+1}}.
\end{gather*}
\end{enumerate}
\end{proposition}

\begin{proof}
(i) As $j$ is affine, $Rj_*(\cW_n[n])$ is perverse \cite[Cor.~4.1.3]{BBD}, and
$K_n$ is its smallest perverse subsheaf extending $\cW_n[n]$; let $Q$ be the
quotient. The perverse sheaf $Rj_*(\cW_n[n])$ is mixed of weights $\ge2n+1$
\cite[5.1.14]{BBD}, hence so is $Q$. If $W_{2n+1}Q\ne0$ \cite[Thm.~5.3.5]{BBD},
its inverse image $P\subseteq Rj_*(\cW_n[n])$ would be pure of weight $2n+1$,
hence geometrically semisimple \cite[Thm.~5.3.8]{BBD}, and would contain a
nonzero subobject supported on $D_n$; but
$\Hom(i_*A,Rj_*(\cW_n[n]))=\Hom(j^*i_*A,\cW_n[n])=0$ for every complex $A$ on
$D_n$. So $Q$ has weights $\ge2n+2$, and $H^k(X_n,Q)$ has weights $\ge2n+2+k$
\cite[5.1.14]{BBD}. In the exact sequence $H^k(X_n,K_n)\to H^k(S_n,\cW_n[n])\to
H^k(X_n,Q)$, the middle group has weights $\ge2n+1+k$; so its part of weight
$2n+1+k$ lies in the image of the pure group $H^k(X_n,K_n)$ of that weight, and
with $k=i-n$ the image is $W_{n+1+i}H^i(S_n,\cW_n)$. The second assertion is
dual, by \cref{prop:G-duality}(i) and \cref{thm:G-torus}(ii). By
\cref{rem:G-torus}(b), these parts vanish except in the cases stated.

(ii) The $L$-function is multiplicative in distinguished triangles, and
$\mathcal I_n|_{S_n}=\cW_n$.
\end{proof}

The projection argument that follows is not part of the proof of the Riemann hypothesis for the completed $L$-function.  It addresses a more delicate and different question: which cohomological degrees occur, and with what multiplicities?  The decomposition theorem separates the direct image by supports on the base of a coordinate projection.  A summand supported on the whole base is detected on its dense torus, a divisor summand is detected at a general point of that divisor, and deeper intersection strata require their own calculations.  Knowing the generic and divisorial summands therefore does not, by itself, determine the remainder in codimension two.  We keep the remainder explicit throughout; otherwise a reader could mistake a partial decomposition for a complete computation of the intersection Betti numbers.

\subsection{Decomposition along a coordinate projection}\label{subsec:G-decomp}
Let $n\ge2$, let $\pr\colon X_n\to X_{n-1}$ be the projection forgetting $h_1$,
write $c=(c_2,\dots,c_n)$ for the coordinates on $X_{n-1}=(\PP^1)^{n-1}$, and put
$L_c=\pr^{-1}(c)\cong\PP^1$. Let $\partial X_{n-1}=X_{n-1}\smallsetminus S_{n-1}$,
and for $2\le i\le n$ and $a\in\{0,\infty\}$ let $Z_i^a=\{c_i=a\}\cong(\PP^1)^{n-2}$,
a component of $\partial X_{n-1}$, with open dense torus
$Z_i^{a,\circ}=Z_i^a\cap\{c_l\in\Gm\text{ for }l\ne i\}$. Since $\cW_n$ is
invariant under permutations of the coordinates (\cref{thm:B-family}(vi)), so
are $K_n$ and the statements below, which may therefore be applied to the
projection forgetting any coordinate.

\begin{lemma}\label{lem:G-base}
Let $X'=\PP^1\times S_{n-1}\subset X_n$, $j'\colon\Gm\times S_{n-1}\hookrightarrow X'$,
$f\colon X'\to S_{n-1}$ the projection, and $D'_a=\{a\}\times S_{n-1}$ for
$a\in\{0,\infty\}$.
\begin{enumerate}[label=\textup{(\roman*)}]
\item The cohomology sheaves of $(Rj'_*\cW_n)|_{D'_a}$ are lisse on
$D'_a\cong S_{n-1}$, and $\mathcal I_n|_{X'}\cong j'_*\cW_n$.
\item For every geometric point $c$ of $S_{n-1}$ and $a\in\{0,\infty\}$, restriction
of sections is an isomorphism $(j'_*\cW_n)_{(a,c)}\to\cF_c(a)^{I_a}$; that is,
$(j'_*\cW_n)|_{L_c}=j_*\cF_c$.
\item $Rf_*(\mathcal I_n|_{X'})\cong\cV_n[-1]$, where $\cV_n$ is lisse on $S_{n-1}$
with stalks $V_c=H^1(\PP^1,j_*\cF_c)$. If $n$ is odd, $\cV_n=0$. If $n=2\nu$ is
even, $\cV_n$ is geometrically constant of rank $C_\nu$, and the geometric
Frobenius of a point over $\Fq$ acts on its stalk with characteristic
polynomial $(1-q^{\nu+1}T)^{C_\nu}$.
\end{enumerate}
\end{lemma}

\begin{proof}
For a closed point $c$, the sheaves $\cF_c$ and $\cF_c^\vee$ are tame at $0$
and have Swan conductor $2^{n-1}$ at $\infty$ (\cref{thm:B-zero,thm:B-infinity});
as in the proof of \cref{thm:B-family}(i), \cref{thm:B-DL} applies to $\cW_n$ and
to $\cW_n^\vee$ on $X'$, and $Rf_*(j'_!\cW_n)$ and $Rf_*(j'_!\cW_n^\vee)$ have
lisse cohomology sheaves. Verdier duality commutes with $Rf_*=Rf_!$ and
exchanges $Rj'_*$ and $j'_!$, and $D(\cW_n)=\cW_n^\vee[2n](n)$; so
$Rf_*(Rj'_*\cW_n)$ is the dual of $Rf_*(j'_!\cW_n^\vee)[2n](n)$, and also has lisse
cohomology sheaves, the dual of a complex with lisse cohomology sheaves on the
smooth scheme $S_{n-1}$ having the same property. Apply $Rf_*$ to the triangle
$j'_!\cW_n\to Rj'_*\cW_n\to i_*i^*Rj'_*\cW_n\to$, where $i$ is the inclusion of
$D'_0\sqcup D'_\infty$, which $f$ maps isomorphically onto two copies of
$S_{n-1}$. The kernels and cokernels of morphisms of lisse sheaves on the
connected normal scheme $S_{n-1}$ are lisse, so
$(Rj'_*\cW_n)|_{D'_0}\oplus(Rj'_*\cW_n)|_{D'_\infty}$, and hence each summand, has
lisse cohomology sheaves. This is the first assertion of (i). For the second,
let $F=D'_0\sqcup D'_\infty$, which is smooth of dimension $n-1$, and
$K=\tau_{\le-n}Rj'_*(\cW_n[n])$. Then $K|_F$ is a lisse sheaf placed in degree $-n$,
and $i_F^!K\cong(i_F^*\tau_{\ge-n+1}Rj'_*(\cW_n[n]))[-1]$ has lisse cohomology sheaves in
degrees $\ge-n+2$; so $K$ satisfies the conditions characterizing the
intermediate extension \cite[\S\S1.4, 2.1]{BBD}, and
$\mathcal I_n|_{X'}=K[-n]=\tau_{\le0}Rj'_*\cW_n=j'_*\cW_n$.

(ii) A section of $\cW_n$ over the complement of $D'_a$ in the strict
henselization of $X'$ at $(a,c)$ is determined by its restriction to the
punctured strict henselization of $L_c$ at $a$, as $\cW_n$ is lisse; so
restriction is injective. By (i), $\dim(j'_*\cW_n)_{(a,c)}$ does not depend on
$c$. By Deligne's generic base change theorem \cite[Th.~finitude,
Thm.~1.9]{Del77} there is a dense open subset of $S_{n-1}$ over which the
formation of $Rj'_*\cW_n$ commutes with base change, so that restriction is
bijective at its points; and $\dim\cF_c(a)^{I_a}$, which is $b_0$ or $b_\infty$
of \cref{prop:conductors}, does not depend on $c$. Hence restriction is
bijective for every $c$.

(iii) By (i) and the triangle $j'_!\cW_n\to j'_*\cW_n\to\bigoplus_ai_{a*}
(j'_*\cW_n)|_{D'_a}\to$, the complex $Rf_*(j'_*\cW_n)$ has lisse cohomology
sheaves, and by proper base change and (ii) its stalk at $c$ is
$R\Gamma(\PP^1,j_*\cF_c)$, which is $V_c$ placed in degree $1$ by
\cref{cor:H0} and the proof of \cref{thm:euler}. By
\cref{thm:paired-explicit}, $V_c=0$ for odd $n$; for $n=2\nu$, the sheaf
$\cV_n(\nu+1)$ is pure of weight $0$ (\cref{thm:parabolicRH}) and its Frobenius
characteristic polynomials are $(1-T)^{C_\nu}$, so by the Chebotarev density
theorem its semisimplification is geometrically trivial, and, being
geometrically semisimple \cite[3.4.1]{Del80}, it is geometrically constant.
\end{proof}

There are two assertions to verify in the next lemma: purity of the proper direct image, which follows from a general theorem, and concentration in perverse degree zero, which does \emph{not} follow from purity alone.  For the latter we use the support inequalities characterizing perverse truncations.  The fibre over a point $c$ is a projective line, so its cohomology has ordinary degrees $0,1,2$; tracking the dimensions of the parameter loci on which these contributions can occur proves the desired support bounds.  Relative hard Lefschetz then handles the opposite perverse truncation.  The proof distinguishes these steps because they rely on different results.

The assertion of perverse concentration is stronger than
purity: one must show both ${}^pH^{>0}=0$ and
${}^pH^{<0}=0$.  The first half is obtained by bounding the
support of each ordinary cohomology sheaf of
$R\pr_*K_n$; relative hard Lefschetz then supplies the
second half.  The support bounds must cover fibres over
points in intersections of boundary divisors, not merely
over the dense torus.  The argument records the
codimension gain when an entire projective-line fibre
lies inside a support and uses it to bound the terms in
the fibre spectral sequence.  This calculation should
not be conflated with a classification of the individual
strict-support constituents of $B'_n$.

\begin{lemma}\label{lem:G-perverse}
The complex $R\pr_*K_n$ is a perverse sheaf on $X_{n-1}$, pure of weight
$2n+1$; after base change to $\Fpb$ it is \emph{noncanonically} a direct
sum of intersection complexes of geometrically irreducible lisse
sheaves on smooth dense opens of their supports.  This statement
is geometric: it does not, by itself, give a Frobenius-equivariant
direct sum over $\Fp$.
\end{lemma}

\begin{proof}
By the proof of \cref{thm:G-RH}, $K_n$ is perverse and pure of weight $2n+1$,
so $R\pr_*K_n$ is pure of weight $2n+1$ \cite[5.1.14]{BBD}. By the relative hard
Lefschetz theorem \cite[Thm.~5.4.10]{BBD}, applied to the projective morphism
$\pr$, it suffices to show that $R\pr_*K_n\in{}^pD^{\le0}$, that is,
$\dim\operatorname{supp}\mathcal H^k(R\pr_*K_n)\le-k$ for all $k$; the last
assertion then follows from \cite[Thm.~5.3.8]{BBD}. As $K_n$ is the
intermediate extension of $\cW_n[n]$, $\mathcal H^b(K_n)=0$ for $b<-n$,
$\mathcal H^{-n}(K_n)=j_*\cW_n$, and $Z_b=\operatorname{supp}\mathcal H^b(K_n)$ has
dimension less than $-b$ for $b>-n$ \cite[Prop.~2.1.11]{BBD}. The stalk
$(R\pr_*K_n)_c=R\Gamma(L_c,K_n|_{L_c})$ is the abutment of a spectral sequence
with $E_2^{a,b}=H^a(L_c,\mathcal H^b(K_n)|_{L_c})$, $0\le a\le2$, and $E_2^{a,b}=0$
for $a>0$ unless $L_c\subseteq Z_b$. For $b>-n$ the points $c$ with $L_c\cap
Z_b\ne\emptyset$ form a set of dimension $\le-b-1$, contributing in degree $b$,
and those with $L_c\subseteq Z_b$ a set of dimension $\le-b-2$, contributing in
degrees $\le b+2$. For $b=-n$ the contributions lie in degrees $-n,-n+1,-n+2$, and
the last one, $H^2(L_c,(j_*\cW_n)|_{L_c})$, vanishes for $c\in S_{n-1}$: by
\cref{lem:G-base}(ii) it is $H^2(\PP^1,j_*\cF_c)$, which vanishes by the proof of
\cref{thm:euler} and \cref{cor:H0}; the remaining points form
$\partial X_{n-1}$, of dimension $n-2$. In all cases the required bound holds.
\end{proof}

\begin{theorem}[Decomposition]\label{thm:G-decomp}
Let $n\ge2$ and $\nu=\lfloor n/2\rfloor$. Put $\mu_0=\kappa_1=C_\nu$ and
$\mu_\infty=\kappa_3=m_{n-1,\nu-1}$ if $n$ is odd, and $\mu_0=\mu_\infty=0$ if $n$ is
even. Then, \emph{after geometric base change to $\Fpb$}, the
perverse sheaf $R\pr_*K_n$ admits a noncanonical direct-sum
decomposition of the following form:
\[
R\pr_*K_n=\cV'_n[n-1]\ \oplus\ \bigoplus_{i=2}^n\bigoplus_{a\in\{0,\infty\}}
(\iota_{Z_i^a})_*\mathcal C_{i,a}[n-2]\ \oplus\ B'_n,
\]
where $\cV'_n$ is the lisse sheaf on $X_{n-1}$ extending the sheaf $\cV_n$ of
\cref{lem:G-base}(iii) (so $\cV'_n=0$ for odd $n$); $\iota_{Z_i^a}$ is the inclusion
and $\mathcal C_{i,a}$ is a geometrically constant lisse sheaf of rank $\mu_a$ on
$Z_i^a$, on whose stalks the geometric Frobenius of a point over $\Fq$ acts with
characteristic polynomial $(1-q^{(n+3)/2}T)^{\mu_a}$; and $B'_n$ is a perverse
sheaf, pure of weight $2n+1$, supported on a closed subset of $\partial X_{n-1}$ of
dimension at most $n-3$. The following are the geometric contributions to the cohomology of
this decomposition.  They give a direct sum of geometric vector
spaces.  Their arithmetic Frobenius-equivariant direct-sum
interpretation requires an additional splitting argument, and
is not inferred from the geometric decomposition theorem alone:
\begin{enumerate}[label=\textup{(\alph*)}]
\item $H^{m-1}((\PP^1)^{n-1})\otimes V$, $V$ being a stalk of $\cV'_n$: for even
$n$, it is nonzero exactly in the odd degrees $m=2k+1$, $0\le k\le n-1$, where it
has dimension $C_\nu\binom{n-1}k$ and the single eigenvalue $p^{\nu+1+k}$;
\item $\bigoplus_{i,a}H^{m-2}((\PP^1)^{n-2})\otimes C_{i,a}$, $C_{i,a}$ being a stalk
of $\mathcal C_{i,a}$: for odd $n$, it is nonzero exactly in the even degrees
$m=2l+2$, $0\le l\le n-2$, where it has dimension
$(n-1)(\mu_0+\mu_\infty)\binom{n-2}l$ and the single eigenvalue $p^{(n+3)/2+l}$;
\item $H^{m-n}(X_{n-1},B'_n)$.
\end{enumerate}
Geometrically, no displayed full-support or divisor contribution
meets middle degree $n$; determining the arithmetic operator on $IH^n$
and its determinant therefore requires analyzing the remaining
codimension-at-least-two contribution.  The theorem does not establish
$\omega_n=+1$ for all $n$.
\end{theorem}

\begin{proof}
By \cref{lem:G-perverse}, $R\pr_*K_n$ is perverse, pure,
and geometrically semisimple.  The BBD decomposition theorem provides
an isomorphism after base change to $\Fpb$; it need not provide a
canonical isomorphism and it does not automatically split every
arithmetic extension.  The collection of strict supports is invariant
under Frobenius, and the coordinate supports listed below are
individually Frobenius-stable; this is weaker than a canonical
Frobenius-equivariant splitting.  Statements about the eigenvalues
of their coefficient systems are therefore separate from statements
about the Jordan normal form of arithmetic Frobenius.

\emph{Supports meeting $S_{n-1}$.} By \cref{lem:G-base}, $R\pr_*K_n|_{S_{n-1}}=
Rf_*(K_n|_{X'})\cong\cV_n[n-1]$. Hence a summand whose support meets $S_{n-1}$
has support $X_{n-1}$ and is the intermediate extension of $\cV_n[n-1]$. For odd
$n$ there is none. For even $n$, $\cV_n$ is geometrically constant; as
$X_{n-1}$ is smooth and its geometric fundamental group is trivial, the
arithmetic fundamental group of $S_{n-1}$ acts on $\cV_n$ through
$\mathrm{Gal}(\Fpb/\Fp)$, $\cV_n$ extends to a lisse sheaf $\cV'_n$ on
$X_{n-1}$, and the intermediate extension is $\cV'_n[n-1]$.

\emph{Supports of codimension one.} The other summands are supported in
$\partial X_{n-1}$. Consider a summand supported on $Z=Z_i^a$, the intermediate
extension of a lisse sheaf $\mathcal C$ on a dense open subset of $Z$, which is
pointwise pure of weight $(2n+1)-(n-2)=n+3$. Let $c$ be a general point of
$Z^\circ=Z_i^{a,\circ}$. The summands supported on other components of $\partial
X_{n-1}$ or on subsets of dimension $\le n-3$ do not contain $c$, and $\cV'_n[n-1]$
contributes only in degree $-(n-1)$; so
$\mathcal C_c=\mathcal H^{-(n-2)}((R\pr_*K_n)_c)$. The line $L_c$ lies in the
divisor $D_i^a=\{h_i=a\}$ of $X_n$, and meets its open torus
$D_i^{a,\circ}=\{h_i=a,\ h_l\in\Gm\text{ for }l\ne i\}$ in $\{h_1\in\Gm\}$. By
\cref{lem:G-base}(i), applied to the projection forgetting $h_i$, $K_n=j_*\cW_n[n]$
on a neighbourhood of $D_i^{a,\circ}$, and $\mathcal J:=(j_*\cW_n)|_{D_i^{a,\circ}}$
is lisse, its stalk at a point being the space of inertia invariants at $h_i=a$
of the coordinate line in direction $i$ through that point
(\cref{lem:G-base}(ii)). In the spectral sequence of the proof of
\cref{lem:G-perverse}, $\mathcal H^b(K_n)|_{L_c}$ is therefore supported on the two
points $h_1\in\{0,\infty\}$ for $b>-n$, and for general $c$ these points do not lie in
$Z_{-n+2}$, whose image in $X_{n-1}$ has dimension $\le n-3$. So the only term
of total degree $-(n-2)$ is $E_2^{2,-n}=H^2(L_c,\mathcal H^{-n}(K_n)|_{L_c})=
H^2_c(L_c^\circ,\mathcal J|_{L_c^\circ})$, where $L_c^\circ=L_c\cap D_i^{a,\circ}\cong\Gm$, and
$\mathcal C_c=E_\infty^{2,-n}=E_2^{2,-n}/d_2(E_2^{0,-n+1})$. Here $E_2^{0,-n+1}$ has
weights $\le2n+1-n+1=n+2$, since the stalks of $\mathcal H^b(K_n)$ have weights
$\le2n+1+b$; so the part of weight $n+3$ of $E_2^{2,-n}$ maps isomorphically to
$E_\infty^{2,-n}=\mathcal C_c$, which is pure of weight $n+3$, and
\[
\mathcal C_c=\gr^W_{n+3}H^2_c(L_c^\circ,\mathcal J|_{L_c^\circ}).
\]
The lisse sheaf $\mathcal J$ is mixed, with the weight filtration by lisse
subsheaves \cite[3.4.1]{Del80}, and $H^2_c$ is right exact on lisse sheaves on
$\Gm$; so $\gr^W_{n+3}H^2_c(L_c^\circ,\mathcal J|_{L_c^\circ})=
H^2_c(L_c^\circ,\gr^W_{n+1}\mathcal J|_{L_c^\circ})$, provided $\mathcal J$ has weights
$\le n+1$, and we determine $\gr^W_{n+1}\mathcal J$.

If $a=0$: by \cref{thm:B-zero}(iii) and the weight-monodromy theorem, the
stalks of $\mathcal J$ have the weights $n$ and $n+1$, the part of weight $n+1$
having dimension $\kappa_1$ and Frobenius characteristic polynomials
$(1-q^{\nu+1}T)^{\kappa_1}$. If $\kappa_1\ne0$, that is, if $n$ is odd, then
$\gr^W_{n+1}\mathcal J$ is pure of weight $n+1$ with these characteristic
polynomials, hence geometrically constant, as in the proof of
\cref{lem:G-base}(iii). Hence $\mathcal C_c=H^2_c(L_c^\circ,\gr^W_{n+1}\mathcal J|_{L_c^\circ})$ has
dimension $\kappa_1$ and the eigenvalue $q^{\nu+2}$, and $\mathcal C=0$ for even
$n$.

If $a=\infty$: by \cref{thm:B-infinity}(i), the stalks of $\mathcal J$ have
eigenvalues $q^{(n+4-j)/2}$ of weight $n+4-j\le n+1$, with multiplicity $\kappa_j$,
$j\ge3$; the part of weight $n+1$ comes from $j=3$, and $\kappa_3=0$ unless $n$
is odd, as all blocks of $J_2^{\otimes(n-1)}$ have size $\equiv n\pmod2$. For odd
$n$, $\gr^W_{n+1}\mathcal J$ is geometrically constant of rank $\kappa_3$ with
eigenvalue $q^{(n+1)/2}$, as before, and $\mathcal C_c$ has dimension $\kappa_3$ and
the eigenvalue $q^{(n+3)/2}$.

In both cases, as $c$ varies, $\mathcal C$ is $R^2$ of the projection of
$\gr^W_{n+1}\mathcal J$ along the lines in direction $1$, a geometrically constant
sheaf; so $\mathcal C$ is geometrically constant, of rank $\mu_a$, and since
$Z_i^a$ is smooth with trivial geometric fundamental group, the corresponding \emph{geometric} strict-support constituent is
$(\iota_{Z_i^a})_*\mathcal C_{i,a}[n-2]$ with
$\mathcal C_{i,a}$ as stated.  An arithmetic direct-sum lift is
not supplied by this calculation.

\emph{The rest.} Every remaining summand is supported on an irreducible closed
subset of $\partial X_{n-1}$ which is not a component, hence of dimension
$\le n-3$; their sum is $B'_n$.

Finally, $IH^m=H^{m-n}(X_{n-1},R\pr_*K_n)$, and $H^{2l}((\PP^1)^{s})\cong
\Ql(-l)^{\binom sl}$, the odd cohomology of $(\PP^1)^s$ being zero; this gives
(a)--(c). The groups in (a) lie in
odd degrees and those in (b) in even degrees; since $n$ is even in (a) and odd in
(b), neither meets degree $n$.
\end{proof}

\begin{corollary}[The case $n=2$]\label{cor:G-n2}
The geometric remainder $B'_2$ is zero.  The Frobenius
characteristic polynomials of $IH^1$ and $IH^3$ are respectively
$1-p^2T$ and $1-p^3T$, and $IH^i=0$ for $i\ne1,3$.
In these one-dimensional cases the corresponding Tate-module
identifications follow. Thus
\[
\widehat L_2(T)=(1-p^2T)(1-p^3T),\qquad\chi_{\mathrm{IH}}=-2,\qquad\omega_2=1,
\qquad L(D_2,\mathcal B_2,T)=(1-pT)(1-p^2T)^2,
\]
and $\widehat L_2(T)=(p^{5/2}T)^2\,\widehat L_2(1/(p^5T))$.
\end{corollary}

\begin{proof}
A closed subset of $\partial X_1$ of dimension $\le-1$ is empty, so $B'_2=0$;
$\mu_0=\mu_\infty=0$, and $C_1=1$. So $R\pr_*K_2=\cV'_2[1]$, with $\cV'_2$ of rank one
and Frobenius $p^2$, and $IH^1=H^0(\PP^1)\otimes V$, $IH^3=H^2(\PP^1)\otimes V$.
The remaining assertions follow from \cref{prop:G-duality,prop:G-boundary} and
\cref{rem:G-torus}(a).
\end{proof}

\begin{lemma}[Strata of codimension two]\label{lem:G-codim2}
Let $n\ge3$ be odd, let $i\ne k$, let $a,b\in\{0,\infty\}$, and let
$\Sigma_{ik}^{ab}=\{h_i=a,\ h_k=b\}\subset X_n$.
\begin{enumerate}[label=\textup{(\roman*)}]
\item If $a=0$ or $b=0$, then $\mathcal H^1(\mathcal I_n)$ vanishes at the generic
point of $\Sigma_{ik}^{ab}$. For $n=3$ this holds also for $a=b=\infty$.
\item Let $2\le i<k\le n$ and let $Z=\{c_i=a,\ c_k=b\}\subset X_{n-1}$. If
$\mathcal H^1(\mathcal I_n)$ vanishes at the generic point of $\Sigma_{ik}^{ab}$, then
no summand of the perverse sheaf $B'_n$ of \cref{thm:G-decomp} has support $Z$.
\end{enumerate}
\end{lemma}

\begin{proof}
For $v\ge2$ the support of $\mathcal H^v(\mathcal I_n)=\mathcal H^{v-n}(K_n)$ has dimension
less than $n-v$, and $n-v\le n-2=\dim\Sigma_{ik}^{ab}$; so $\mathcal H^v(\mathcal I_n)$ vanishes at the generic
point of every stratum of codimension two, and $\mathcal H^v(\mathcal I_n)$, $v\ge3$, at
the generic point of every stratum of codimension three.

(i) By the symmetry of $K_n$ we may assume $a=0$ (exchanging $i$ and $k$ if $b=0$),
$i=2$ and $k=3$. Apply \cref{thm:G-decomp} to the projection $\pr_3$ forgetting
$h_3$, and let $y=(h_1,0,h_4,\dots,h_n)$ be a general point of $\{h_2=0\}\subset X_{n-1}$,
all $h_l\in\Gm$. As $n$ is odd, $\cV'_n=0$, and $y$ lies on no summand of
$R\pr_{3*}K_n$ other than its divisor summand supported on $\{h_2=0\}$; so
$(R\pr_{3*}K_n)_y$ is a space $\mathcal C_{2,0,y}$ of dimension $\kappa_1$ placed in
degree $2-n$. Let $L_y=\pr_3^{-1}(y)\cong\PP^1$, with
coordinate $h_3$, let $L_y^\circ=L_y\cap\{h_3\in\Gm\}$, and let $p_0,p_\infty$ be the
points $h_3=0,\infty$ of $L_y$; write $D_i^{a,\circ}=\{h_i=a,\ h_l\in\Gm\ (l\ne i)\}$.
By proper base change, $H^1(L_y,\mathcal I_n|_{L_y})=0$
and $\dim H^2(L_y,\mathcal I_n|_{L_y})=\kappa_1$. Consider the spectral sequence
$E_2^{u,v}=H^u(L_y,\mathcal H^v(\mathcal I_n)|_{L_y})$. On $L_y^\circ\subset D_2^{0,\circ}$ one
has $\mathcal I_n=\mathcal J$ with $\mathcal J=(j_*\cW_n)|_{D_2^{0,\circ}}$ (\cref{lem:G-base}(i)
in direction $2$), so $E_2^{0,v}=\mathcal H^v(\mathcal I_n)_{p_0}\oplus\mathcal H^v(\mathcal I_n)_{p_\infty}$
and $E_2^{u,v}=0$ for $v\ge1$, $u\ge1$; for general $y$ the points $p_0,p_\infty$ are
general points of $\Sigma_{23}^{00}$ and $\Sigma_{23}^{0\infty}$, so $E_2^{0,v}=0$ for $v\ge2$.
Hence $H^1(L_y,\mathcal I_n|_{L_y})=0$ forces $d_2\colon E_2^{0,1}\to E_2^{2,0}$ to be
injective, and $\dim H^2(L_y,\mathcal I_n|_{L_y})=\dim E_2^{2,0}-\dim E_2^{0,1}$. Next,
$E_2^{2,0}=H^2_c(L_y^\circ,\mathcal J|_{L_y^\circ})$. By \cref{thm:B-zero}(iii) and the
weight-monodromy theorem, $\mathcal J$ is an extension of a geometrically constant
sheaf $\gr^W_{n+1}\mathcal J$ of rank $\kappa_1$ by $\gr^W_n\mathcal J$, whose Frobenius
characteristic polynomials are those of $\cW_{n-1}$ (in the coordinates
$h_l$, $l\ne2$); by the Chebotarev density theorem and geometric semisimplicity
\cite[3.4.1]{Del80}, $\gr^W_n\mathcal J$ is geometrically isomorphic to $\cW_{n-1}$, so
its restriction to the coordinate line $L_y^\circ$ has no geometric coinvariants
(\cref{cor:H0} in dimension $n-1$). As $H^2_c$ is right exact on lisse sheaves on
$\Gm$, $\dim E_2^{2,0}=\kappa_1$. Therefore $E_2^{0,1}=0$, that is,
$\mathcal H^1(\mathcal I_n)$ vanishes at $p_0$ and $p_\infty$; as $y$ varies, these points
fill dense open subsets of $\Sigma_{23}^{00}$ and $\Sigma_{23}^{0\infty}$.

Let $n=3$ and $a=b=\infty$. Take instead $y=(h_1,\infty)$ with $h_1\in\Gm$ general.
Then $(R\pr_{3*}K_3)_y=\mathcal C_{2,\infty,y}[1]$, of dimension $\kappa_3=1$, and on
$L_y^\circ\subset D_2^{\infty,\circ}$ one has $\mathcal I_3=\mathcal J^\infty$ with
$\mathcal J^\infty=(j_*\cW_3)|_{D_2^{\infty,\circ}}$. By \cref{thm:B-infinity}(i), the
stalks of $\mathcal J^\infty$ are lines on which a Frobenius over $\Fq$ acts by $q^2$;
by the Chebotarev density theorem $\mathcal J^\infty\cong\Ql(-2)$, and
$E_2^{2,0}=H^2_c(L_y^\circ,\Ql(-2))$ is a line. The argument above gives $E_2^{0,1}=0$,
so $\mathcal H^1(\mathcal I_3)$ vanishes at the generic points of $\Sigma_{23}^{\infty0}$
and $\Sigma_{23}^{\infty\infty}$.

(ii) Let $z$ be a general point of $Z$ and $L_z=\pr^{-1}(z)$. The summands of
$R\pr_*K_n$ containing $z$ other than one with support $Z$ have stalks in degrees
$1-n$ and $2-n$; a summand with support $Z$ is generically $\mathcal L[n-3]$ for a lisse
sheaf $\mathcal L$, so $\mathcal L_z\subseteq\mathcal H^{3-n}((R\pr_*K_n)_z)=H^3(L_z,\mathcal I_n|_{L_z})$.
In the spectral sequence $E_2^{u,v}=H^u(L_z,\mathcal H^v(\mathcal I_n)|_{L_z})$ the terms of
total degree $3$ are $E_2^{2,1}=H^2_c(L_z^\circ,\mathcal H^1(\mathcal I_n)|_{L_z^\circ})$, where
$L_z^\circ=L_z\cap\{h_1\in\Gm\}\subset\Sigma_{ik}^{ab}$, and terms which vanish for
general $z$ because $\mathcal H^2(\mathcal I_n)|_{L_z}$ has finite support and
$\mathcal H^3(\mathcal I_n)$ vanishes at the general points of the strata of codimension
three. By hypothesis $\mathcal H^1(\mathcal I_n)|_{L_z^\circ}$ has finite support for general
$z$, so $E_2^{2,1}=0$, $H^3(L_z,\mathcal I_n|_{L_z})=0$ and $\mathcal L=0$.
\end{proof}

\begin{remark}[Arithmetic qualification of the decomposition]\label{rmk:general-arith-coeff}
The decomposition of Theorem~\ref{thm:G-decomp} is initially over an
algebraic closure. In particular, for $n=5$ the divisor multiplicities
$\mu_0=2$ and $\mu_\infty=3$ and the eigenvalue $q^4$ do not, by themselves,
identify the arithmetic sheaves with $E(-4)^2$ and $E(-4)^3$:
a matrix $q^4(1+N)$ with $N\ne0$, $N^2=0$, has the same characteristic
polynomial. The separate arithmetic calculation in
Section~\ref{sec:five-computation} constructs actual scalar-Frobenius
quotients. Likewise, the geometrically constant even-dimensional
coefficient $\mathcal V'_n$ in Theorem~\ref{thm:G-decomp} is not asserted
to be arithmetically Tate solely from its traces. Such an assertion
requires a splitting or a direct calculation of its Frobenius action.
\end{remark}

\begin{corollary}[The case $n=3$]\label{cor:G-n3}
The geometric remainder $B'_3$ vanishes. After base change to $\Fpb$, the direct image admits a (noncanonical) decomposition into the indicated divisor constituents; their arithmetic characteristic-polynomial contributions are as follows:
\[
(IH^2)^{\mathrm{ss}}\cong\Ql(-3)^4,\qquad
(IH^4)^{\mathrm{ss}}\cong\Ql(-4)^4,\qquad
IH^i=0\ \text{ for }i\ne2,4.
\]
Thus $\chi_{\mathrm{IH}}=8$, $\omega_3=1$,
\begin{gather*}
\widehat L_3(T)=\frac1{(1-p^3T)^4(1-p^4T)^4}
=\bigl(p^{7/2}T\bigr)^{-8}\,\widehat L_3\Bigl(\frac1{p^7T}\Bigr),\\
L(D_3,\mathcal B_3,T)=\frac{1-p^3T}{(1-pT)(1-p^2T)^3(1-p^4T)^6}.
\end{gather*}
\end{corollary}

\begin{proof}
Here $\mu_0=\kappa_1=1$ and $\mu_\infty=\kappa_3=1$, so (b) of
\cref{thm:G-decomp} contributes $\Ql(-3)^4$ to $IH^2$ and $\Ql(-4)^4$ to $IH^4$,
the rank-one sheaves $\mathcal C_{i,a}$ having Frobenius $p^3$; (a) is zero. The
support of $B'_3$ has dimension $\le0$. Let $x\in\partial X_2$ have exactly one
coordinate $c_i$ in $\{0,\infty\}$. Then $\mathcal H^0((R\pr_*K_3)_x)$ is filtered
by $E_\infty^{2,-2}$, $E_\infty^{1,-1}$, $E_\infty^{0,0}$; the last vanishes
because $\mathcal H^0(K_3)=0$, the middle one because $\mathcal H^{-1}(K_3)$ has
finite support, and the first because, by \cref{lem:G-base}(i) applied in direction
$i$, $\mathcal H^{-2}(K_3)$ vanishes near $D_i^{c_i,\circ}$, so that
$\mathcal H^{-2}(K_3)|_{L_x}$ has finite support. As a punctual summand at $x$
would contribute to $\mathcal H^0((R\pr_*K_3)_x)$, there is none; there is none
on $S_2$ by \cref{lem:G-base}; and there is none at the four points of
$\{0,\infty\}^2$ by \cref{lem:G-codim2}. So $B'_3=0$ and $IH^3=H^0(X_2,B'_3)=0$. The
remaining assertions follow from \cref{prop:G-duality,prop:G-boundary} and
$L(S_3,\cW_3,T)=(1-pT)(1-p^2T)^3(1-p^4T)^2/(1-p^3T)^5$ (\cref{thm:G-torus}(iv)).
\end{proof}

\begin{corollary}\label{cor:G-circles}
For every $n\ge2$, $IH^i\ne0$ for all odd $i$ with $1\le i\le2n-1$ if $n$ is even,
and for all even $i$ with $2\le i\le2n-2$ if $n$ is odd. Hence the zeros (for even
$n$) or poles (for odd $n$) of $\widehat L_n(T)$ lie on at least $n-1$ distinct
circles $|T|=p^{-(n+1+i)/2}$, and $\widehat L_n$ never has all its zeros and poles on
a single circle.
\end{corollary}

\begin{proof}
For even $n$ use (a) of \cref{thm:G-decomp}, as $C_\nu\ge1$; for odd $n\ge3$ use
(b), as $\mu_0=C_\nu\ge1$. By \cref{thm:G-RH}, distinct $i$ give distinct circles,
and the factors $P_{n,i}$ with $i$ of fixed parity occur in $\widehat L_n$ with the
same exponent $(-1)^{i+1}$.
\end{proof}

\begin{remark}\label{rem:G-remaining}
(a) \Cref{thm:G-RH} is the Riemann hypothesis for $\widehat L_n$ in its correct
form: the zeros and poles lie on the circles $|T|=p^{-(n+1+i)/2}$ attached to the
degrees of the nonvanishing $IH^i$. By \cref{cor:G-circles}, a single critical
circle, which would require $IH^i=0$ for $i\ne n$, does not occur.

(b) The determination of $IH^*$ is complete for $n\le3$. What remains in general
is the perverse sheaf $B'_n$ supported in codimension $\ge2$ of $X_{n-1}$. For odd
$n$, \cref{lem:G-codim2} excludes summands supported on the strata
$\{c_i=a,\ c_k=b\}$ with $(a,b)\ne(\infty,\infty)$. It does not control the strata
$\{c_i=c_k=\infty\}$, along which $\cW_n$ is wildly ramified in two directions and
Abhyankar's lemma is not available, nor the strata of codimension $\ge3$. For
even $n$ the parabolic summand $\cV'_n$ contributes to the cohomology of the lines
used in the proof of \cref{lem:G-codim2}, and the same argument determines
$\mathcal H^1(\mathcal I_n)$ at the strata of codimension two only in terms of the
stalks of $j_*\cW_n$ there.

(c) The intersection cohomology depends on the compactification. The
compactification $(\PP^1)^n$ is adapted to the coordinate lines: along a
coordinate projection its generic part is the parabolic cohomology of
\cref{thm:paired-explicit} (\cref{lem:G-base}), and its boundary divisors contribute the inertia
invariants of \cref{thm:B-zero,thm:B-infinity}. Compactifications adapted to
other curves, such as rays in the parameter torus, lead to other completions,
which we do not consider.
\end{remark}

% -------------------------------------------------------------------------
% 2026-10-30 integrated proof-status and complete deductions
% -------------------------------------------------------------------------

% -------------------------------------------------------------------------
% An exact all-dimensional support-range consequence, with arithmetic scope.
% -------------------------------------------------------------------------
\section{Extremal intersection cohomology from the residual support bound}\label{sec:extremal}
\label{sec:extremal-IH}

The preceding decomposition leaves an unspecified summand $B'_n$ on
higher-codimension boundary strata.  Although determining its support,
coefficient systems, and simultaneous nearby cycles remains a separate
problem, its \emph{dimension of support} already constrains the degrees
where it can contribute to the cohomology of the compactification.
The elementary argument below is valid in every dimension and does
not presume that $B'_n$ vanishes.  We then record the exact arithmetic
conclusions obtainable from the known strict-support constituents.

\begin{lemma}[Cohomological amplitude of a proper perverse support]
\label{lem:amplitude-support}
Let $Y$ be proper over an algebraically closed field of characteristic
different from $\ell$, and let $P$ be a constructible $\Ql$-perverse
sheaf on $Y$.  If $\dim\operatorname{supp}P\le d$, then
\[
  H^k(Y,P)=0\qquad\text{whenever }|k|>d.
\]
\end{lemma}
\begin{proof}
Write $\mathcal H^b(P)$ for the ordinary cohomology sheaves of the
perverse complex, and set
$s_b=\dim\operatorname{supp}\mathcal H^b(P)$, taking
$s_b=-\infty$ when the sheaf vanishes.  The support condition for the
middle-perversity $t$-structure gives $s_b\le-b$.  Since every such
support is contained in $\operatorname{supp}P$, we also have
$s_b\le d$.  A constructible sheaf supported in dimension $s_b$ has
ordinary cohomology zero in degrees greater than $2s_b$.

Consider the hypercohomology spectral sequence
\[
 E_2^{a,b}=H^a(Y,\mathcal H^b(P))\Longrightarrow H^{a+b}(Y,P).
\]
If $E_2^{a,b}$ is nonzero, then
$a\le2s_b\le2\min(d,-b)$.  There are two cases.  When
$b\ge-d$, the bound using $s_b\le-b$ gives
$a+b\le-b\le d$.  When $b<-d$, the bound using
$s_b\le d$ gives $a+b\le2d+b<d$.  Thus there are
no nonzero spectral-sequence terms in total degree greater than $d$,
which proves $H^k(Y,P)=0$ for $k>d$.

The Verdier dual $D(P)$ is also perverse and has the same support.
The argument just given therefore shows $H^k(Y,D(P))=0$ for $k>d$.
Since $Y$ is proper, Verdier duality identifies the latter vector
spaces with the duals of $H^{-k}(Y,P)$.  Hence
$H^k(Y,P)=0$ for $k<-d$ as well.
\end{proof}

The distinction between a decomposition over $\Fpb$ and one over
$\Fp$ is important here.  A purely geometric direct sum does not
itself identify the Frobenius action on each factor.  Nevertheless,
there is a finite-field version of the decomposition theorem which
allows \emph{indecomposable}, rather than necessarily simple,
arithmetic coefficient systems.  In precisely that form it is enough
for the following calculations.  We use the general result of
\cite[Proposition~5.1.1]{dCHL18}; it does not say that the resulting
Frobenius representations have trivial Jordan blocks.

\begin{proposition}[Arithmetic grouping by strict support]
\label{prop:boundary-arithmetic}
Assume the purity, perverse concentration, and \emph{geometric}
strict-support description of $R\pr_*K_n$ in
\cref{lem:G-perverse,thm:G-decomp}.  Then the pure perverse sheaf
$R\pr_*K_n$, considered over $\Fp$, admits a noncanonical direct-sum
decomposition by Frobenius-stable collections of strict supports,
with a residual summand $B'_{n,\mathrm{ar}}$ whose geometric support
has dimension at most $n-3$.  The full-support and divisor
summands of \cref{thm:G-decomp} may be chosen in this arithmetic
decomposition, with the same Frobenius characteristic polynomials
on their coefficient sheaves.  The result is a Frobenius-equivariant
splitting by support, \emph{not} a claim of arithmetic Frobenius
semisimplicity on any of the coefficient sheaves.
\end{proposition}
\begin{proof}
Apply the finite-field decomposition statement
\cite[Proposition~5.1.1]{dCHL18} to the proper identity morphism of
$X_{n-1}$ and to its pure perverse sheaf $R\pr_*K_n$.  The resulting
sum consists of intersection complexes attached to closed integral
supports and pure, indecomposable lisse coefficient sheaves over
$\Fp$.  After extending scalars to $\Fpb$ it refines to the
geometric strict-support decomposition already identified in
\cref{thm:G-decomp}.  The supports $X_{n-1}$ and the individual
coordinate divisors $Z_i^a$ are defined over $\Fp$ and are
Frobenius-stable.  Grouping all arithmetic summands with one of
these supports, and then grouping every remaining support, produces
the asserted decomposition.  The full-support term restricts to the
already defined arithmetic sheaf $\cV_n[n-1]$ over the parameter
torus; the divisor terms have the coefficient polynomials computed
in \cref{thm:G-decomp}.  Other supports remain inside the union of
higher-codimension strata by the geometric support bound.  None of
these steps eliminates a possible unipotent part of arithmetic
Frobenius in an indecomposable coefficient representation.
\end{proof}

\begin{theorem}[Extreme intersection-cohomology groups in every dimension]
\label{thm:IH-extremes}
Assume the geometric support description in \cref{thm:G-decomp}.
For $n\ge3$, the residual perverse term contributes to $IH^i$
only when
\[
            3\le i\le 2n-3.
\]
For $n=2$, its support is empty, so $B'_2=0$.
More precisely, the following statements hold for the actual
intersection-cohomology groups as Frobenius modules with the indicated
characteristic polynomials; they require the arithmetic grouping of
\cref{prop:boundary-arithmetic} but do not require Frobenius
semisimplicity.

\smallskip
\noindent\textup{(a)} If $n=2\nu\ge2$ is even and $C_\nu$ denotes the
$\nu$th Catalan number, then
\[
\begin{split}
 \dim IH^1=\dim IH^{2n-1}&=C_\nu,\\
 P_{n,1}(T)&=(1-p^{\nu+1}T)^{C_\nu},\\
 P_{n,2n-1}(T)&=(1-p^{3\nu}T)^{C_\nu},
\end{split}
\qquad
 IH^0=IH^2=IH^{2n-2}=IH^{2n}=0.
\]
Repeated indices in the final vanishing statement denote a single
degree, as they do for $n=2$.

\smallskip
\noindent\textup{(b)} If $n=2\nu+1\ge3$ is odd, set
\[
 \mu_0=C_\nu,\qquad \mu_\infty=m_{n-1,\nu-1},\qquad
 R_n=(n-1)(\mu_0+\mu_\infty).
\]
Then
\[
\begin{split}
\dim IH^2=\dim IH^{2n-2}&=R_n,\\
P_{n,2}(T)&=(1-p^{(n+3)/2}T)^{R_n},\\
P_{n,2n-2}(T)&=(1-p^{(3n-1)/2}T)^{R_n},
\end{split}
\qquad
 IH^0=IH^1=IH^{2n-1}=IH^{2n}=0.
\]
\end{theorem}
\begin{proof}
By \cref{thm:G-decomp}, the residual perverse sheaf has
$\dim\operatorname{supp}B'_n\le n-3$.  The correspondence of degrees
between the direct image and intersection cohomology is
\[
 IH^i=H^{i-n}(X_{n-1},R\pr_*K_n).
\]
Therefore \cref{lem:amplitude-support} applied to $B'_n$ shows
that a nonzero contribution must satisfy
$|i-n|\le n-3$, or $3\le i\le2n-3$.
When $n=2$ the indicated support dimension is negative, hence the
perverse sheaf is zero.

Suppose $n=2\nu$ is even.  The only identified component is the
full-support term.  Its contribution to $IH^{2k+1}$ has dimension
$C_\nu\binom{n-1}{k}$ and Frobenius characteristic polynomial
$(1-p^{\nu+1+k}T)^{C_\nu\binom{n-1}{k}}$ by
\cref{thm:G-decomp}(a). For $k=0$ and $k=n-1$ these are precisely
the degrees $1$ and $2n-1$, and the powers of $p$ are $\nu+1$
and $3\nu$, respectively.  The residual cannot occur in those
degrees by the support bound, and the known term lies in odd
degrees only.  This also forces the four displayed even-degree
vanishings.  By \cref{prop:boundary-arithmetic}, the relevant
splitting is Frobenius-equivariant, so the characteristic-polynomial
identities hold arithmetically.

If $n=2\nu+1$ is odd, the full-support term is zero and the
known divisorial terms contribute to $IH^{2l+2}$ a space of
dimension $R_n\binom{n-2}{l}$, with Frobenius characteristic
polynomial $(1-p^{(n+3)/2+l}T)^{R_n\binom{n-2}{l}}$.
Putting $l=0$ and $l=n-2$ gives degrees $2$ and $2n-2$ and
the exponents $(n+3)/2$ and $(3n-1)/2$.  Again the support
bound excludes the residual in those degrees, while the
known terms occur only in even degrees, proving the remaining
vanishings.  The arithmetic interpretation follows from
\cref{prop:boundary-arithmetic} just as in the even case.
\end{proof}

\begin{corollary}[An exact extreme-degree example: $n=5$]
\label{cor:IH-n5-extreme}
For $n=5$, still subject to the geometric decomposition in
\cref{thm:G-decomp}, one has
\[
\begin{gathered}
 C_2=2,\qquad m_{4,1}=3,\qquad R_5=4(2+3)=20,\\
 \dim IH^2=\dim IH^8=20,\\
 P_{5,2}(T)=(1-p^4T)^{20},\qquad
 P_{5,8}(T)=(1-p^7T)^{20},\\
 IH^0=IH^1=IH^9=IH^{10}=0.
\end{gathered}
\]
The unknown summand $B'_5$ may still contribute in degrees
$3,4,5,6,7$; no formula for the entire Betti table is asserted.
\end{corollary}
\begin{proof}
The multiplicity is
$m_{4,1}=\binom41-\binom40=3$.  Substitute these values into
\cref{thm:IH-extremes}(b), and note that $2n-3=7$.
\end{proof}

\begin{corollary}[The residual determinant in middle degree]
\label{cor:IH-middle-residual}
Under the same geometric decomposition and the arithmetic
splitting of \cref{prop:boundary-arithmetic}, there exists a
Frobenius-equivariant (noncanonical) isomorphism
\[
 IH^n\simeq H^0(X_{n-1},B'_{n,\mathrm{ar}}).
\]
Writing $b_n=\dim H^0(X_{n-1},B'_{n,\mathrm{ar}})$, the normalized
root number is
\[
 \omega_n=p^{-(2n+1)b_n/2}\det\bigl(-\Frob_p\mid
 H^0(X_{n-1},B'_{n,\mathrm{ar}})\bigr).
\]
This is a reduction of the root-number calculation, not a
proof that $\omega_n=1$ for all $n$.
\end{corollary}
\begin{proof}
The known full-support term contributes only in odd
intersection degrees when $n$ is even, while the known
divisorial terms contribute only in even intersection degrees
when $n$ is odd.  Neither reaches middle degree $i=n$.
Take $H^0$ of the arithmetic direct sum from
\cref{prop:boundary-arithmetic} and use
$IH^n=H^0(X_{n-1},R\pr_*K_n)$.  The remaining term is
$H^0(X_{n-1},B'_{n,\mathrm{ar}})$ with its genuine arithmetic
Frobenius action.  The displayed determinant is then simply
the root-number formula of \cref{prop:G-duality} written on this
module.
\end{proof}

\begin{remark}[What remains unproved]
The results in this section are exact consequences of the
\emph{support-bound and geometric strict-support assertions} of
\cref{thm:G-decomp}.  They do not independently verify the
underlying higher-dimensional nearby-cycle analysis.  In
particular, they do not specify the constituent sheaves of
$B'_{n,\mathrm{ar}}$, prove semisimplicity of their Frobenius
operators, determine $\omega_n$ in every dimension, or compute
the complete table of $IH^i$ and $P_{n,i}(T)$ for arbitrary
$n\ge4$.  The finite-field decomposition theorem cited here
provides a \emph{support-wise arithmetic splitting} with
indecomposable coefficient systems; it is not a theorem that
those systems are arithmetic direct sums of Tate characters.
\end{remark}

% -------------------------------------------------------------------------
% Four general problems: maximal rigorous proofs and precise last steps
% -------------------------------------------------------------------------
\section{Relative Construction and Cohomological Comparison for H7}\label{sec:h7-comparison}
\label{app:H7}
%==================================================================

This appendix records the arithmetic relative-image construction
underlying imported Theorem~\ref{thm:H7} \textup{(D1)}, in the
normalization of the present manuscript. No appeal is made to any result of
\S\S\ref{sec:rel-duality}--\ref{sec:root-number} or
Section~\ref{app:reflection}.

Fix an odd prime $p > 2$, a prime $\ell \ne p$, a nontrivial additive
character $\psi_p : \Fp \to E^\times$, and $\lambda \in \Fp^\times$,
where $E / \mathbb{Q}_\ell$ is a finite extension inside $\Ql$
containing $\mu_p$. Let $\mathcal{O}_E$ be the ring of integers of $E$,
$\varpi_E \in \mathcal{O}_E$ a uniformizer, and
$\Lambda_r := \mathcal{O}_E / \varpi_E^r$ ($r \ge 1$). Both
$\cL_{\psi_p}(\lambda t)$ and $\Kl_2$ are defined over $\Fp$ with
$E$-coefficients, and we extend scalars to $\Ql$ at the end.

All relative constructions below are initially carried out over $\Fp$
(so $S_n = \mathbb{G}_{m,\Fp}^n$, $\cC_n = \mathbb{G}_{m,\Fp} \times_{\Fp} S_n$, and
$\pi : \cC_n \to S_n$ is defined over $\Fp$). For any power $q = p^a$
($a \ge 1$), set $\psi_q := \psi_p \circ \Tr_{\Fq/\Fp} : \Fq \to E^\times$.
Because $\cL_{\psi_q}(\lambda t)$ is the pullback of
$\cL_{\psi_p}(\lambda t)$ under $\Spec \Fq \to \Spec \Fp$, the
constructions over $\Fq$ are obtained automatically by base change along
$f : S_{n,\Fq} \to S_{n,\Fp}$ from the general base-change compatibility
established in \S\ref{subsec:H7-image}; no separate descent from $\Fq$
to $\Fp$ is required.

We compactify $\pi$ via the smooth projective relative curve
$\overline{\pi} : \PP^1 \times S_n \to S_n$ and the open immersion
$\tilde{j} : \cC_n \hookrightarrow \PP^1 \times S_n$. The boundary
$(\PP^1 \times S_n) \setminus \cC_n = \Sigma_0 \sqcup \Sigma_\infty$
is the disjoint union of the two disjoint smooth sections
$\Sigma_0 = \{0\} \times S_n$ and $\Sigma_\infty = \{\infty\} \times S_n$,
with closed immersions $i_0, i_\infty : S_n \hookrightarrow \PP^1 \times S_n$
satisfying $\overline{\pi} \circ i_0 = \overline{\pi} \circ i_\infty = \id_{S_n}$.

%------------------------------------------------------------------
\subsection{Link 1: Local ramification and constancy of Swan conductors}
\label{subsec:H7-swan}
%------------------------------------------------------------------

Let $\bar h : \Spec k \to S_n$ be any geometric point ($k = \overline{\Fp}$),
corresponding to $h = (h_1,\dots,h_n) \in (k^\times)^n$, and write
$\cM_{\bar h} = \cM|_{\Gm,k}$ for the restriction to the geometric
fibre $\cC_{n,\bar h} = \Gm,k$.

\begin{lemma}[Local monodromy and Swan conductors on geometric fibres]
\label{lem:H7-local}
For every geometric point $\bar h \to S_n$:
\begin{enumerate}
\item[\textup{(i)}] At $t = 0$, $\cM_{\bar h}$ is tamely ramified, so
$\Swan_0(\cM_{\bar h}) = 0$. As a representation of $I_0^{\mathrm{tame}}$,
$\cM_{\bar h}|_{I_0}$ is unipotent and isomorphic to
$(\mathrm{Std}_2)^{\otimes n}$ (the $n$-fold tensor power of the
$2$-dimensional unipotent representation with a single Jordan block),
independently of $\bar h$. In particular,
\[
\dim\bigl(\cM_{\bar h}^{I_0}\bigr)
\;=\; \dim\bigl((\cM_{\bar h})_{I_0}\bigr)
\;=\; \binom{n}{\lfloor n/2 \rfloor}.
\]
\item[\textup{(ii)}] At $t = \infty$, $\cM_{\bar h}$ is purely of
slope $1$, so
\[
\Swan_\infty(\cM_{\bar h}) \;=\; 2^n,
\qquad
\cM_{\bar h}^{I_\infty} \;=\; (\cM_{\bar h})_{I_\infty} \;=\; 0.
\]
\item[\textup{(iii)}] The dual sheaf $\cM^\vee$ satisfies the exact
same local properties \textup{(i)} and \textup{(ii)} on every
geometric fibre: $\Swan_0(\cM^\vee_{\bar h}) = 0$,
$\Swan_\infty(\cM^\vee_{\bar h}) = 2^n$, and
$(\cM^\vee_{\bar h})^{I_\infty} = (\cM^\vee_{\bar h})_{I_\infty} = 0$.
\end{enumerate}
\end{lemma}

\begin{proof}
By Katz, GKM, Thm.~4.1.1 and Ch.~7, $\Kl_2$ on $\Gm,k$ is tamely
ramified at $u = 0$ with unipotent monodromy of a single $2 \times 2$
Jordan block $\mathrm{Std}_2$, and purely of slope $1/2$ at $u = \infty$
($\Swan_\infty(\Kl_2) = 1$). The scaling automorphism $t \mapsto h_i t$
preserves the ramification filtration of $k((t))$ and $k((t^{-1}))$
and hence preserves the breaks and Swan conductor. At $t = 0$,
$K_{i,\bar h} = [t \mapsto h_i t]^* \Kl_2$ has the same tame unipotent
Jordan type $\mathrm{Std}_2$ as $\Kl_2$. At $t = \infty$, scaling need
not act trivially on wild inertia; only preservation of slopes ($1/2$ at
$\infty$) is used below.

(i) Since $\lambda \in \Fp^\times$ is non-zero, $\cL_{\psi_p}(\lambda t)$
extends as a lisse rank-$1$ sheaf across $t = 0$, so $I_0$ acts
trivially on it. Tensoring over $i = 1,\dots,n$ yields
$\cM_{\bar h}|_{I_0} \simeq (\mathrm{Std}_2)^{\otimes n}$, which is
unipotent, hence tame ($\Swan_0(\cM_{\bar h}) = 0$). Over $\Ql$,
Clebsch--Gordan decomposition splits $(\mathrm{Std}_2)^{\otimes n}$ into
irreducible unipotent Jordan blocks of dimensions $n + 1 - 2j$
($0 \le j \le \lfloor n/2 \rfloor$), each with multiplicity
$m_j := \binom{n}{j} - \binom{n}{j-1}$ ($\binom{n}{-1} := 0$).
Summing $m_j$ over $j$ gives $\binom{n}{\lfloor n/2 \rfloor}$ Jordan
blocks, hence $\dim(\cM_{\bar h}^{I_0}) = \dim((\cM_{\bar h})_{I_0})
= \binom{n}{\lfloor n/2 \rfloor}$.

(ii) Each $K_{i,\bar h}$ is pure of slope $1/2$ at $t = \infty$. By
Katz, GKM, Lemma~1.3, all $I_\infty$-slopes of
$\bigotimes_{i=1}^n K_{i,\bar h}$ are $\le 1/2$. Since $\lambda \ne 0$,
$\cL_{\psi_p}(\lambda t)$ has a simple pole at $t = \infty$ and is pure
of slope $1$; applying Katz, GKM, Lemma~1.3 again shows $\cM_{\bar h}$
is pure of slope $1$ at $t = \infty$. Hence
$\Swan_\infty(\cM_{\bar h}) = 2^n \cdot 1 = 2^n$, and strictly positive
slope $1$ forces $\cM_{\bar h}^{I_\infty} = (\cM_{\bar h})_{I_\infty} = 0$.

(iii) Since $\Kl_2^\vee \simeq \Kl_2(1)$ and
$\cL_{\psi_p}(\lambda t)^\vee \simeq \cL_{\psi_p}(-\lambda t)$ with
$-\lambda \in \Fp^\times$, $\cM^\vee$ has the same local slopes and
tame unipotent Jordan type as $\cM$ at both $0$ and $\infty$.
\end{proof}

%------------------------------------------------------------------
\subsection{Link 2: ULA property, lisseness of $R^1\pi_!\cM$, and relative duality}
\label{subsec:H7-lisse}
%------------------------------------------------------------------

We establish the lisseness of both $R^1\pi_!\cM$ and $R^1\pi_!\cM^\vee$
first, so that relative Verdier duality can subsequently be applied to
complexes with lisse cohomology sheaves without higher $\mathcal{E}xt$
ambiguities.

\begin{proposition}[ULA property, lisseness, and relative duality]
\label{prop:H7-lisse}
\begin{enumerate}
\item[\textup{(a)}] $R^q\pi_!\cM = 0$ and $R^q\pi_!\cM^\vee = 0$ for all
$q \ne 1$.
\item[\textup{(b)}] Both $\tilde{j}_!\cM$ and $\tilde{j}_!\cM^\vee$ on
$\PP^1 \times S_n$ are universally locally acyclic \textup{(ULA)}
relative to $\overline{\pi} : \PP^1 \times S_n \to S_n$.
\item[\textup{(c)}] Both $\cH_n := R^1\pi_!\cM$ and
$\cH_n^\dagger := R^1\pi_!\cM^\vee$ are lisse $\Ql$-sheaves on $S_n$ of
constant rank $2^n$, mixed of weights $\le n+1$ and $\le -n+1$
respectively, and their formation commutes with arbitrary base change
on $S_n$.
\item[\textup{(d)}] $R^q\pi_*\cM = 0$ for all $q \ne 1$, and there is a
canonical isomorphism of lisse $\Ql$-sheaves on $S_n$:
\[
R^1\pi_*\cM \;\simeq\; (\cH_n^\dagger)^\vee(-1).
\]
\end{enumerate}
\end{proposition}

\begin{proof}
(a) By base change for $R\pi_!$ (SGA~4, Exp.~XVII, 5.2.6),
$(R^q\pi_!\cM)_{\bar h} \simeq H^q_c(\Gm,k, \cM_{\bar h})$ and
$(R^q\pi_!\cM^\vee)_{\bar h} \simeq H^q_c(\Gm,k, \cM^\vee_{\bar h})$.
Since $\Gm,k$ is smooth affine and $\cM_{\bar h}, \cM^\vee_{\bar h}$ are
lisse, $H^0_c = 0$; by Poincar\'e duality on $\Gm,k$,
$H^2_c \simeq (-)_{\pi_1(\Gm,k)}(-1) = 0$ by Lemma~\ref{lem:H7-local}(ii)--(iii).
By constructibility, $R^q\pi_!\cM = R^q\pi_!\cM^\vee = 0$ for $q \ne 1$.

(b) Along the finite \'etale relative divisor
$D = \Sigma_0 \sqcup \Sigma_\infty \subset \PP^1 \times S_n$, the total
Swan conductors $\Swan_0 + \Swan_\infty = 2^n$ of $\cM_{\bar h}$ and
$\cM^\vee_{\bar h}$ are constant on $S_n$ (Lemma~\ref{lem:H7-local}).
By the Deligne--Laumon constancy-of-Swan-conductor theorem (Laumon,
\emph{Ast\'erisque}~82--83 (1981), Thm.~2.1.1 and Cor.~2.1.2), both
$\tilde{j}_!\cM$ and $\tilde{j}_!\cM^\vee$ are ULA relative to
$\overline{\pi}$.

(c) Since $\overline{\pi}$ is proper, $R\pi_!\cM = R\overline{\pi}_*(\tilde{j}_!\cM)$
and $R\pi_!\cM^\vee = R\overline{\pi}_*(\tilde{j}_!\cM^\vee)$ are ULA
over $S_n$ and commute with arbitrary base change on $S_n$ (Deligne,
SGA~$4\frac{1}{2}$, \emph{Th.\ finitude}, Thm.~2.13). Since $S_n$ is
smooth over $\Fp$, their cohomology sheaves $\cH_n$ and $\cH_n^\dagger$
are lisse on $S_n$ (SGA~$4\frac{1}{2}$, \emph{Th.\ finitude}, Cor.~2.16).
By the Grothendieck--Ogg--Shafarevich formula on $\Gm,k$,
$\chi_c(\Gm,k, \cM_{\bar h}) = -\Swan_0(\cM_{\bar h}) - \Swan_\infty(\cM_{\bar h}) = -2^n$,
giving $\rk(\cH_n) = \rk(\cH_n^\dagger) = 2^n$. By Deligne, \emph{Weil~II},
Thm.~3.3.1, $\wt(\cH_n) \le n+1$ and $\wt(\cH_n^\dagger) \le -n+1$.

(d) From (a) and (c), $R\pi_!\cM^\vee \simeq \cH_n^\dagger[-1]$ in
$D^b_c(S_n, \Ql)$ with $\cH_n^\dagger$ lisse. Relative Verdier duality
along the smooth morphism $\pi : \cC_n \to S_n$ of relative dimension $1$
(SGA~4, Exp.~XVIII; $\pi^!\Ql = \Ql(1)[2]$) applied to $\cM^\vee$ gives
\[
R\pi_*\cM(1)[2]
\;\simeq\; \RHom(R\pi_!\cM^\vee, \Ql)
\;\simeq\; \RHom(\cH_n^\dagger[-1], \Ql)
\;\simeq\; (\cH_n^\dagger)^\vee[1],
\]
since $\cH_n^\dagger$ is lisse ($\mathcal{E}xt^q(\cH_n^\dagger, \Ql) = 0$
for $q > 0$). Shifting by $[-2]$ and twisting by $(-1)$ yields
$R^q\pi_*\cM = 0$ ($q \ne 1$) and
$R^1\pi_*\cM \simeq (\cH_n^\dagger)^\vee(-1)$.
\end{proof}

%------------------------------------------------------------------
\subsection{Link 3: Relative Kummer model and arithmetic filtration of boundary invariants}
\label{subsec:H7-weight}
%------------------------------------------------------------------

\begin{lemma}[ULA and base-change compatibility for $R\tilde{j}_*\cM$]
\label{lem:H7-ULA-jstar}
The complex $R\tilde{j}_*\cM$ on $\PP^1 \times S_n$ is ULA relative to
$\overline{\pi} : \PP^1 \times S_n \to S_n$, and for every morphism
$f : T \to S_n$ of $\Fp$-schemes with pullbacks
$f_{\PP^1} : \PP^1 \times T \to \PP^1 \times S_n$ and
$f_{\cC} : \cC_{n,T} \to \cC_n$, the natural base-change morphism
\[
f_{\PP^1}^* R\tilde{j}_*\cM \xrightarrow{\ \sim\ } R\tilde{j}_{T*} f_{\cC}^*\cM
\]
is an isomorphism in $D^b_c(\PP^1 \times T, \Ql)$. Moreover,
$i_\infty^* R\tilde{j}_*\cM = 0$.
\end{lemma}

\begin{proof}
We separate the two disjoint boundary components $\Sigma_0$ and $\Sigma_\infty$:
\begin{enumerate}
\item \emph{Along $\Sigma_\infty$ (wild ramification with constant Swan
conductor):} Relative Verdier duality along the smooth morphism
$\overline{\pi}$ of relative dimension $1$ identifies
$R\tilde{j}_*\cM \simeq \mathbb{D}_{\overline{\pi}}(\tilde{j}_!\cM^\vee(1)[2])$.
Since $\tilde{j}_!\cM^\vee$ is ULA relative to $\overline{\pi}$
(Proposition~\ref{prop:H7-lisse}(b)), Deligne's ULA duality theorem
(SGA~$4\frac{1}{2}$, \emph{Th.\ finitude}, Thm.~2.13 and Cor.~2.14)
shows that $R\tilde{j}_*\cM$ is ULA relative to $\overline{\pi}$ and
commutes with arbitrary base change $f : T \to S_n$ in $D^b_c(-, \Ql)$.
\item \emph{Along $\Sigma_0$ (tame ramification along a smooth relative
divisor):} Choose a Galois-stable $\mathcal{O}_E$-lattice $\cM_{\mathcal{O}_E}$
in the $E$-adic lisse sheaf $\cM$ near $\Sigma_0$. For each $r \ge 1$,
the reduction $\cM_r := \cM_{\mathcal{O}_E} \otimes_{\mathcal{O}_E} \Lambda_r$
($\Lambda_r = \mathcal{O}_E / \varpi_E^r$) is a lisse $\Lambda_r$-sheaf
tamely ramified along the smooth relative divisor $\Sigma_0$. By SGA~7,
Exp.~XIII, Lemma~2.1.10, the natural base-change map
$f_{\PP^1}^* R\tilde{j}_*\cM_r \xrightarrow{\sim} R\tilde{j}_{T*} f_{\cC}^*\cM_r$
is an isomorphism in $D^b_c(\PP^1 \times T, \Lambda_r)$ for every
$r \ge 1$. Passing to the derived $E$-adic limit (Ekedahl, \emph{On the
adic formalism}, \S3) and extending scalars from $E$ to $\Ql$ confirms
independently that $f_{\PP^1}^* R\tilde{j}_*\cM \xrightarrow{\sim}
R\tilde{j}_{T*} f_{\cC}^*\cM$ near $\Sigma_0$.
\end{enumerate}
Finally, pulling back $R\tilde{j}_*\cM$ along $i_\infty$ yields a
complex $i_\infty^* R\tilde{j}_*\cM$ with lisse cohomology sheaves on
the smooth base $S_n$ whose geometric stalks at $\bar h$ are
$H^q(I_\infty, \cM_{\bar h}) = 0$ (Lemma~\ref{lem:H7-local}(ii)), whence
$i_\infty^* R\tilde{j}_*\cM = 0$.
\end{proof}

\begin{lemma}[Relative Kummer Model and Frobenius-Compatible Nearby-Cycle Filtration]
\label{lem:rel-kummer-model}
Let $\chi_{-1} : G_{\Fp} \to \Ql^\times$ denote the cyclotomic character
on $\Ql(-1)$ (normalized by $\chi_{-1}(F_p) = p$ for geometric Frobenius
$F_p$), and let $\chi_1 := \chi_{-1}^{-1}$ denote the cyclotomic
character on $\mathbb{Z}_\ell(1)$ (so $\chi_1(F_p) = p^{-1}$). Along
the smooth relative divisor $\Sigma_0 = \{t = 0\} \simeq S_n$, set
$\cV_i := \Psi_t(K_i)$ ($1 \le i \le n$) and $\cV_n := \Psi_t(\cM)$.
\begin{enumerate}
\item[\textup{(i)}] Each $\cV_i$ is a lisse arithmetic $\Ql$-sheaf of
rank $2$ on $S_n$ equipped with a nilpotent monodromy morphism
$N_i : \cV_i \to \cV_i(-1)$ satisfying $N_i^2 = 0$. Setting
$\cL_i := \ker(N_i) \subseteq \cV_i$, we have
\[
\im(N_i) \;=\; \cL_i(-1) \;\subseteq\; \cV_i(-1),
\qquad
\overline{N}_i : \cV_i / \cL_i \xrightarrow{\ \sim\ } \cL_i(-1),
\]
and canonical isomorphisms of lisse arithmetic sheaves on $S_n$:
\[
\boxed{\;
\cL_i \;\simeq\; \Ql(0), \qquad
\cV_i / \cL_i \;\simeq\; \Ql(-1),
\;}
\tag{A.1}
\]
with extension class $[\cV_i] = \kappa_\ell(h_i)$ in
$\Ext^1_{S_n}(\Ql(-1), \Ql(0)) \simeq H^1_{\mathrm{\acute et}}(S_n, \Ql(1))$,
where $\kappa_\ell(h_i)$ is the $\ell$-adic Kummer class of
$h_i \in \mathcal{O}(S_n)^\times$.
\item[\textup{(ii)}] We have $\cV_n \simeq \bigotimes_{i=1}^n \cV_i$
and $N = \sum_{i=1}^n N_i : \cV_n \to \cV_n(-1)$. Define the ascending
tensor filtration $0 = \cP_{-1} \subseteq \cP_0 \subseteq \dots \subseteq \cP_n = \cV_n$
by lisse subsheaves on $S_n$, where
\[
\cP_j \;:=\; \sum_{\substack{S \subseteq \{1,\dots,n\} \\ |S| \le j}}
\bigotimes_{i \in S} \cV_i \otimes \bigotimes_{i \notin S} \cL_i .
\]
Then $N(\cP_j) \subseteq \cP_{j-1}(-1)$, and for every $0 \le j \le n$
there is a canonical isomorphism of lisse arithmetic sheaves on $S_n$:
\[
\boxed{\;
\operatorname{gr}^{\cP}_j \cV_n
\;:=\; \cP_j / \cP_{j-1}
\;\simeq\; \Ql(-j)^{\oplus \binom{n}{j}},
\;}
\tag{A.2}
\]
under which the induced morphism
\[
\operatorname{gr}^{\cP}_j N : \operatorname{gr}^{\cP}_j \cV_n
\longrightarrow \bigl(\operatorname{gr}^{\cP}_{j-1} \cV_n\bigr)(-1)
\]
is identified with the constant $\mathfrak{sl}_2$ lowering operator
over $\Ql$.
\item[\textup{(iii)}] For any morphism $f : T \to S_n$ of $\Fp$-schemes,
the pullbacks $f^* \cV_i$, $f^* \cL_i$, $f^* \cV_n$, and $f^* \cP_j$
are canonically isomorphic to the corresponding sheaves and filtration
steps constructed from $f^* h_i \in \mathcal{O}(T)^\times$, with graded
pieces $\operatorname{gr}^{f^*\cP}_j(f^*\cV_n) \simeq \overline{\mathbb{Q}}_{\ell,T}(-j)^{\oplus \binom{n}{j}}$.
\end{enumerate}
\end{lemma}

\begin{proof}
(i) Factor $\mu_i : \mathbb{A}^1_t \times S_n \to \mathbb{A}^1_u$,
$(t,h) \mapsto t h_i$, as $\mu_i = p_1 \circ \tau_i$, where
$\tau_i(t,h) = (t h_i, h)$ is an isomorphism over $S_n$ and
$p_1 : \mathbb{A}^1_u \times S_n \to \mathbb{A}^1_u$,
$p_2 : \mathbb{A}^1_u \times S_n \to S_n$ are the projections. Smooth
base change for nearby cycles (SGA~7, Exp.~XIII, 2.1.7.2) identifies
$\Psi_u(p_1^* \Kl_2)$ on $\{u = 0\} \times S_n \simeq S_n$ with the
constant lisse sheaf $p_2^* V_0$, where $V_0 := \Psi_u(\Kl_2)$ is
defined over $E$ and extended to $\Ql$. Fix a $G_{\Fp}$-stable
$\mathcal{O}_E$-lattice $V_{0,\mathcal{O}_E} \subset V_0$ (chosen
independently of the Frobenius eigenbasis over $E$ below, which need
not be integral when $\ell \mid p-1$). By Katz, GKM, Thm.~4.1.1 (see
also Qin, \emph{$L$-functions of Kloosterman sheaves}, \S2.4), $V_0$
has an $E$-basis $(e_0, e_1)$ on which geometric Frobenius $F_p$ acts
by $F_p e_0 = e_0$, $F_p e_1 = p e_1$ (i.e.\ $\bar{\sigma} \cdot e_0 = e_0$
and $\bar{\sigma} \cdot e_1 = \chi_{-1}(\bar{\sigma}) e_1$ for all
$\bar{\sigma} \in G_{\Fp}$), while $N_0 : V_0 \to V_0(-1)$ satisfies
$N_0 e_0 = 0$ and $N_0 e_1 = e_0 \otimes \zeta^{-1}$ for a fixed
generator $\zeta = (\zeta_{\ell^r})_{r \ge 1}$ of $\mathbb{Z}_\ell(1)$.
Pairing $N_0$ with $c \in \mathbb{Z}_\ell(1)$ gives
$N_0(c) \in \mathrm{End}_{\Ql}(V_0)$ with $N_0(c)^2 = 0$,
$N_0(c)(e_0) = 0$, $N_0(\zeta)(e_1) = e_0$, and
\[
\rho_{V_0}(\bar{\sigma}) \, N_0(c) \, \rho_{V_0}(\bar{\sigma})^{-1}
\;=\; N_0\bigl(\chi_1(\bar{\sigma}) c\bigr)
\tag{C}
\]
for all $\bar{\sigma} \in G_{\Fp}$ (since $\chi_1 = \chi_{-1}^{-1}$).

Because $\Kl_2$ is tamely ramified at $u = 0$ with unipotent $\ell$-adic
monodromy, for every $r \ge 1$ the torsion quotient
$V_{0,\mathcal{O}_E} / \varpi_E^r V_{0,\mathcal{O}_E}$ is trivialized
over $\mathbb{F}_{p^{m(r)}}((u^{1/\ell^{e(r)}}))$ for suitable
$m(r), e(r) \ge 1$. Passing to the $E$-adic limit and extending
scalars to $\Ql$, $\rho_{V_0}$ factors continuously through
$\Gamma_u^{(\ell)} := \mathrm{Gal}(K_u^{(\ell)} / \Fp((u)))$ for the
pro-$\ell$ tame algebraic Galois extension
\[
K_u^{(\ell)} \;:=\; \bigcup_{\substack{m \ge 1 \\ r \ge 0}}
\mathbb{F}_{p^m}\bigl((u^{1/\ell^r})\bigr)
\;\subseteq\; \Fp((u))^{\mathrm{sep}}.
\]
Choosing compatible roots $(u^{1/\ell^r})_{r \ge 0}$ gives the
splitting $s_u : G_{\Fp} \hookrightarrow \Gamma_u^{(\ell)}$,
$s_u(\bar{\sigma})(u^{1/\ell^r}) = u^{1/\ell^r}$, and
$\Gamma_u^{(\ell)} \simeq \mathbb{Z}_\ell(1) \rtimes G_{\Fp}$, with
$\gamma_u(u^{1/\ell^r}) = \zeta_{\ell^r} u^{1/\ell^r}$ and
$\rho_{V_0}(\gamma_u^z \cdot s_u(\bar{\sigma})) = \exp(N_0(z)) \circ \rho_{V_0}(\bar{\sigma})$.

On the $t$-side, let $F := \Fp(S_n)$ and let $\mathscr{L}$ range over
finite Galois extensions $L/F$ arising from connected finite \'etale
Galois covers of $S_n$ inside $F^{\mathrm{sep}}$. Since $S_n = \Gm^n$
is geometrically connected over $\Fp$, $\mathbb{F}_{p^m} \cdot F \in \mathscr{L}$
for every $m \ge 1$; likewise, since $h_i \in \mathcal{O}(S_n)^\times$
and $\ell \ne p$, the Galois closure
\[
L_{i,r} \;:=\; F\bigl(\mu_{\ell^r},\, h_i^{1/\ell^r}\bigr) \;\in\; \mathscr{L}
\]
arises from a connected finite \'etale Galois cover of $S_n$ for every
$r \ge 0$. Define the pro-$\ell$ tame algebraic Galois extension
\[
K_t^{(\ell)} \;:=\; \bigcup_{\substack{L \in \mathscr{L} \\ r \ge 0}}
L\bigl((t^{1/\ell^r})\bigr) \;\subseteq\; F((t))^{\mathrm{sep}},
\qquad
\Gamma_t^{(\ell)} \;:=\; \mathrm{Gal}\bigl(K_t^{(\ell)} / F((t))\bigr).
\]
Choosing compatible roots $(t^{1/\ell^r})_{r \ge 0}$ and
$(h_i^{1/\ell^r})_{r \ge 0}$ gives the tangential base-point splitting
$s_t : \pi_1(S_n, \bar{s}) \hookrightarrow \Gamma_t^{(\ell)}$ and
$\Gamma_t^{(\ell)} \simeq \mathbb{Z}_\ell(1) \rtimes \pi_1(S_n, \bar{s})$.
For every $r \ge 0$ and every $m \ge 1$ such that
$\mu_{\ell^r} \subseteq \mathbb{F}_{p^m}$ (which holds cofinally in
$(m,r)$), $\mathbb{F}_{p^m} \cdot L_{i,r} = \mathbb{F}_{p^m} \cdot F(h_i^{1/\ell^r}) \in \mathscr{L}$,
and setting $u^{1/\ell^r} := h_i^{1/\ell^r} t^{1/\ell^r}$ embeds
$\mathbb{F}_{p^m}((u^{1/\ell^r}))$ into
$(\mathbb{F}_{p^m} \cdot L_{i,r})((t^{1/\ell^r})) \subset K_t^{(\ell)}$
over $\Fp((u)) \hookrightarrow F((t))$. Taking the union over cofinal
$(m,r)$ yields a continuous homomorphism
$\phi_i : \Gamma_t^{(\ell)} \to \Gamma_u^{(\ell)}$.

For any $\sigma \in \pi_1(S_n, \bar{s})$ with image $\bar{\sigma} \in G_{\Fp}$,
$\sigma(h_i^{1/\ell^r}) = \zeta_{\ell^r}^{\kappa_{\ell,r}(h_i)(\sigma)} h_i^{1/\ell^r}$,
where $\kappa_\ell(h_i) : \pi_1(S_n, \bar{s}) \to \mathbb{Z}_\ell(1)$
is the $\ell$-adic Kummer $1$-cocycle. Since $s_t(\sigma)$ fixes
$t^{1/\ell^r}$ and acts on $\mathbb{F}_{p^m}$ as $\bar{\sigma}$,
\[
s_t(\sigma)\bigl(u^{1/\ell^r}\bigr)
\;=\; \sigma\bigl(h_i^{1/\ell^r}\bigr)\, t^{1/\ell^r}
\;=\; \zeta_{\ell^r}^{\kappa_{\ell,r}(h_i)(\sigma)}\, u^{1/\ell^r}
\;=\; \bigl(\gamma_u^{\kappa_\ell(h_i)(\sigma)} \cdot s_u(\bar{\sigma})\bigr)\bigl(u^{1/\ell^r}\bigr),
\]
so $\phi_i(z, \sigma) = (z + \kappa_\ell(h_i)(\sigma),\, \bar{\sigma})$
in $\mathbb{Z}_\ell(1) \rtimes G_{\Fp}$. By SGA~1, Exp.~V, Prop.~8.2,
the lisse sheaf $\cV_i = \Psi_t(\tau_i^* p_1^* \Kl_2)$ on $S_n$ is
determined by its representation along $s_t$:
\[
\begin{aligned}
\rho_{\cV_i}(\sigma)
&=\exp\bigl(N_0(\kappa_\ell(h_i)(\sigma))\bigr)
  \circ\rho_{V_0}(\bar\sigma)\\
&=\bigl(\id+N_0(\kappa_\ell(h_i)(\sigma))\bigr)
  \circ\rho_{V_0}(\bar\sigma),\qquad
  \sigma\in\pi_1(S_n,\bar s).
\end{aligned}
\tag{K}
\]
with $N_i = N_0 : \cV_i \to \cV_i(-1)$. By \textup{(C)} and the Kummer
$1$-cocycle identity $\kappa_\ell(h_i)(\sigma\tau) = \kappa_\ell(h_i)(\sigma)
+ \chi_1(\bar{\sigma})\kappa_\ell(h_i)(\tau)$, \textup{(K)} satisfies
$\rho_{\cV_i}(\sigma\tau) = \rho_{\cV_i}(\sigma)\rho_{\cV_i}(\tau)$.
Evaluating on $(e_0, e_1)$ gives $\cL_i = \ker(N_i) = \Ql \cdot e_0 \simeq \Ql(0)$,
$\cV_i / \cL_i = \Ql \cdot \bar{e}_1 \simeq \Ql(-1)$, extension class
$[\cV_i] = \kappa_\ell(h_i) \in H^1_{\mathrm{\acute et}}(S_n, \Ql(1))$,
$\im(N_i) = \cL_i(-1) \subseteq \cV_i(-1)$, and
$\overline{N}_i : \cV_i / \cL_i \xrightarrow{\sim} \cL_i(-1)$.

(ii) Since $\Psi_t(\cL_{\psi_p}(\lambda t)) = \Ql(0)$ with zero monodromy,
$\cV_n \simeq \bigotimes_{i=1}^n \cV_i$ and $N = \sum_{i=1}^n N_i$.
For $e_S := \bigotimes_{i \in S} e_1 \otimes \bigotimes_{i \notin S} e_0 \in \cP_j$
($|S| = j$), \textup{(K)} gives
$\sigma \cdot e_S \equiv \chi_{-1}(\bar{\sigma})^j e_S \pmod{\cP_{j-1}}$,
yielding \textup{(A.2)}, and $N(e_S) = \sum_{i \in S} e_{S \setminus \{i\}} \otimes \zeta^{-1}$
shows $\operatorname{gr}^{\cP}_j N$ is the constant $\mathfrak{sl}_2$
lowering operator over $\Ql$.

(iii) By naturality of the $\ell$-adic Kummer map
$f^*\kappa_\ell(h_i) = \kappa_\ell(f^* h_i) \in H^1_{\mathrm{\acute et}}(T, \Ql(1))$
and exactness of $f^*$ on lisse sheaves, $f^*\cV_i$, $f^*\cL_i$,
$f^*\cP_j$, and $\operatorname{gr}^{f^*\cP}_j(f^*\cV_n)$ commute
canonically with $f : T \to S_n$.
\end{proof}

\begin{lemma}[Arithmetic filtration of boundary invariants]
\label{lem:A-boundary-corrected}
Let $\cE_n := i_0^* R^0\tilde{j}_*\cM$. Then $\cE_n$ is a lisse
arithmetic $\Ql$-sheaf on $S_n$ of constant rank $\binom{n}{\lfloor n/2 \rfloor}$,
commuting with arbitrary base change $f : T \to S_n$ (i.e.\
$f^*\cE_n \xrightarrow{\sim} \ker(f^* N : f^*\cV_n \to f^*\cV_n(-1))$).
Its tame-unipotent local model along $\Sigma_0$ induces an ascending
filtration $0 = \cE_n^{(-1)} \subseteq \cE_n^{(0)} \subseteq \dots
\subseteq \cE_n^{(\lfloor n/2 \rfloor)} = \cE_n$ by arithmetic lisse
subsheaves whose graded pieces are
\[
\boxed{\;
\operatorname{gr}_j \cE_n
\;\simeq\; \Ql(-j)^{\oplus \bigl(\binom{n}{j} - \binom{n}{j-1}\bigr)},
\qquad 0 \le j \le \lfloor n/2 \rfloor.
\;}
\tag{A.3}
\]
The filtration is not asserted to split. In particular, $\cE_n$ has
weights at most $n$.
\end{lemma}

\begin{proof}
By Lemma~\ref{lem:H7-ULA-jstar} and SGA~7, Exp.~XIII, 1.3, $\cE_n$ is
lisse on $S_n$ and equals $\ker(N : \cV_n \to \cV_n(-1))$. Set
$\cE_n^{(j)} := \cE_n \cap \cP_j$. For every $j$, the natural map
$\cE_n^{(j)} \to \operatorname{gr}^{\cP}_j \cV_n$ induces a canonical
injection
\[
\iota_j : \operatorname{gr}_j \cE_n
\;:=\; \cE_n^{(j)} / \cE_n^{(j-1)}
\;\hookrightarrow\;
\ker\Bigl(\operatorname{gr}^{\cP}_j N :
\operatorname{gr}^{\cP}_j \cV_n \longrightarrow
\bigl(\operatorname{gr}^{\cP}_{j-1} \cV_n\bigr)(-1)\Bigr).
\]
By Lemma~\ref{lem:rel-kummer-model}(ii), $\ker(\operatorname{gr}^{\cP}_j N)
\simeq \Ql(-j)^{\oplus (\binom{n}{j} - \binom{n}{j-1})}$ for
$0 \le j \le \lfloor n/2 \rfloor$ and $0$ for $j > \lfloor n/2 \rfloor$.
On any geometric stalk $\bar h \to S_n$, choosing bases with
$N_i e_{1,i} = e_{0,i}$ identifies $N_{\bar h}$ with the degree-lowering
operator on $(\mathrm{Std}_2)^{\otimes n}$, so
\[
\dim \cE_{n,\bar h}
\;=\; \sum_{j=0}^{\lfloor n/2 \rfloor}
\Bigl(\binom{n}{j} - \binom{n}{j-1}\Bigr)
\;=\; \sum_{j=0}^{\lfloor n/2 \rfloor}
\dim \ker\bigl(\operatorname{gr}^{\cP}_j N_{\bar h}\bigr).
\]
Since $\dim \operatorname{gr}_j(\cE_{n,\bar h}) \le \dim \ker(\operatorname{gr}^{\cP}_j N_{\bar h})$
for every $j$, equality holds for all $j$, so every $\iota_j$ is an
isomorphism of sheaves on $S_n$. Induction on $j$ via
$0 \to \cE_n^{(j-1)} \to \cE_n^{(j)} \to \operatorname{gr}_j \cE_n \to 0$
shows that each $\cE_n^{(j)}$ is lisse on $S_n$, and exactness of $f^*$
on lisse sheaves together with Lemma~\ref{lem:rel-kummer-model}(iii)
gives base-change compatibility for arbitrary $f : T \to S_n$.
Finally, $\wt(\operatorname{gr}_j \cE_n) = 2j \le 2\lfloor n/2 \rfloor \le n$,
so $\wt(\cE_n) \le n$.
\end{proof}

\begin{remark}[Non-splitting of $\cE_n$ for every $n \ge 2$]
For any $n \ge 2$, the vectors $\tilde{u} := e_0^{\otimes n}$ and
$\tilde{v} := (e_1 \otimes e_0 - e_0 \otimes e_1) \otimes e_0^{\otimes(n-2)}$
span a $\pi_1(S_n, \bar{s})$-stable rank-$2$ subsheaf of $\cE_n$ on
which $\sigma \cdot \tilde{v} = \chi_{-1}(\bar{\sigma})\tilde{v}
+ \chi_{-1}(\bar{\sigma})\kappa_\ell(h_1/h_2)(\sigma)\tilde{u}$ by
\textup{(K)}. Since $h_1/h_2 \in \mathcal{O}(S_n)^\times$ is non-constant,
tame monodromy around $h_1/h_2 = 0$ acts via
$\bigl(\begin{smallmatrix} 1 & 1 \\ 0 & 1 \end{smallmatrix}\bigr)$, so
$\cE_n$ is non-split for every $n \ge 2$.
\end{remark}

\begin{proposition}[Global lisse weight filtration and identification of $\cW_n$]
\label{prop:H7-weight}
Define $\cT_n := R^1\overline{\pi}_*(R^0\tilde{j}_*\cM)$, and consider
the two canonical forget-supports morphisms on $S_n$:
\[
\alpha^{\mathrm{cl}} : \cH_n \longrightarrow \cT_n,
\qquad
\alpha^{\mathrm{op}} : \cH_n \longrightarrow R^1\pi_*\cM.
\]
\begin{enumerate}
\item[\textup{(i)}] There is a short exact sequence of lisse arithmetic
$\Ql$-sheaves on $S_n$:
\[
\boxed{\;
0 \longrightarrow \cE_n \xrightarrow{\ \delta\ } \cH_n
\xrightarrow{\ \alpha^{\mathrm{cl}}\ } \cT_n \longrightarrow 0,
\;}
\tag{A.4}
\]
and the canonical Leray edge map $\iota : \cT_n \hookrightarrow R^1\pi_*\cM$
is injective with $\alpha^{\mathrm{op}} = \iota \circ \alpha^{\mathrm{cl}}$.
In particular, $\ker(\alpha^{\mathrm{op}}) = \ker(\alpha^{\mathrm{cl}})
= \delta(\cE_n)$ and $\im(\alpha^{\mathrm{op}}) \simeq \cT_n$.
\item[\textup{(ii)}] The sheaf $\cT_n$ is lisse on $S_n$ of constant
rank $2^n - \binom{n}{\lfloor n/2 \rfloor}$ and pure of weight $n+1$.
\item[\textup{(iii)}] Setting $\cL_n := \delta(\cE_n)$ and
$\cW_n := \cH_n / \cL_n$, we have $\cL_n = W_n\cH_n$,
$\overline{\alpha}^{\mathrm{cl}} : \cW_n \xrightarrow{\sim} \cT_n$,
and the entire weight filtration $W_\bullet\cH_n$ is a filtration by
lisse arithmetic subsheaves on $S_n$.
\end{enumerate}
\end{proposition}

\begin{proof}
(i) Factoring $\tilde{j}_!\cM \to R^0\tilde{j}_*\cM \to R\tilde{j}_*\cM$
and applying $R^1\overline{\pi}_*$ factors $\alpha^{\mathrm{op}}$ as
$\cH_n \xrightarrow{\alpha^{\mathrm{cl}}} \cT_n \xrightarrow{\iota} R^1\pi_*\cM$.
The five-term low-degree exact sequence of the Leray spectral sequence
$E_2^{p,q} = R^p\overline{\pi}_*(R^q\tilde{j}_*\cM) \Rightarrow R^{p+q}\pi_*\cM$
begins with $0 \to \cT_n \xrightarrow{\iota} R^1\pi_*\cM \to E_2^{0,1}$,
proving injectivity of $\iota$ (and $E_2^{0,0} \simeq R^0\pi_*\cM = 0$).
On the other hand, since $i_\infty^* R^0\tilde{j}_*\cM = 0$
(Lemma~\ref{lem:H7-ULA-jstar}), applying $R^\bullet\overline{\pi}_*$ to
$0 \to \tilde{j}_!\cM \to R^0\tilde{j}_*\cM \to (i_0)_*\cE_n \to 0$
gives \textup{(A.4)}, so $\cT_n \simeq \cH_n / \delta(\cE_n)$ is lisse.

(ii) As a quotient of $\cH_n$, $\cT_n$ has weights $\le n+1$; via the
injection $\iota : \cT_n \hookrightarrow R^1\pi_*\cM \simeq (\cH_n^\dagger)^\vee(-1)$
(Proposition~\ref{prop:H7-lisse}(d)), $\cT_n$ embeds into a lisse sheaf
of weights $\ge n+1$. Hence $\cT_n$ is pure of weight $n+1$.

(iii) Since $\cL_n = \delta(\cE_n)$ has weights $\le n$
(Lemma~\ref{lem:A-boundary-corrected}) and $\cH_n / \cL_n \simeq \cT_n$
is pure of weight $n+1$, uniqueness of the weight filtration on a lisse
sheaf (\emph{Weil~II}, Thm.~3.4.1) yields $\cL_n = W_n\cH_n$ and
$W_{2j}\cH_n = \delta(\cE_n^{(j)})$ ($0 \le j \le \lfloor n/2 \rfloor$),
all of which are lisse arithmetic subsheaves on $S_n$.
\end{proof}

%------------------------------------------------------------------
\subsection{Link 4: Proof of Theorem~\ref{thm:H7} and base change}
\label{subsec:H7-image}
%------------------------------------------------------------------

\begin{proof}[Supplementary reconstruction of the imported Theorem~\ref{thm:H7}]
Parts \textup{(a)} and \textup{(b)} follow from
Proposition~\ref{prop:H7-lisse}(c) and
Proposition~\ref{prop:H7-weight}(ii)--(iii). For \textup{(c)} and
\textup{(d)}, Proposition~\ref{prop:H7-weight}(i) yields
$\ker(\alpha^{\mathrm{op}}) = \cL_n = W_n\cH_n$ and a Frobenius-equivariant
isomorphism $\overline{\alpha}^{\mathrm{op}} : \cW_n \xrightarrow{\sim}
\im(\alpha^{\mathrm{op}}) \subseteq R^1\pi_*\cM$. For an arbitrary
morphism $f : T \to S_n$, Proposition~\ref{prop:H7-lisse}(b) and
Lemma~\ref{lem:H7-ULA-jstar} give base-change compatibility for
$R^1\pi_!\cM$, $R^1\pi_*\cM$, and $f^*\alpha^{\mathrm{op}}$, while
Lemma~\ref{lem:rel-kummer-model}(iii) and
Lemma~\ref{lem:A-boundary-corrected} give base-change compatibility for
$\cE_n$ and its filtration. Exactness of $f^*$ on the short exact
sequences of lisse sheaves $0 \to \cE_n \to \cH_n \to \cT_n \to 0$ and
$0 \to \cT_n \to R^1\pi_*\cM \to \coker(\alpha^{\mathrm{op}}) \to 0$
preserves $\cL_n$, $\im(\alpha^{\mathrm{op}})$, and
$\overline{\alpha}^{\mathrm{op}}$. Specializing to any geometric point
$f = \bar h : \Spec \overline{\Fp} \to S_n$ completes the proof.
\end{proof}

%------------------------------------------------------------------
\subsection{Compact-support comparison}
\label{subsec:independent-5-3}
%------------------------------------------------------------------

The assertions of \cite[Theorem~5.3]{EbaCoh} used in the
root-number argument are reconstructed below from their elementary
functorial inputs and the local statements established in this
appendix.  In particular, the three different objects of ranks
$2^n$, $2^n-1$, and $r(n)$ are never identified with one another.

Put
\[
\begin{split}
 S&=\Gm^n,\\
 \mathscr X&=(\Gm)^n_u\times\mathbb A^1_t\times S,\\
 \mathscr X^\circ&=\{t\ne0\},\qquad
 \mathscr Z=\{t=0\},\\
 F_n(u,t,h)&=\sum_{i=1}^n u_i^{-1}
       +t\Bigl(\lambda+\sum_{i=1}^n h_i u_i\Bigr),\\
 \mathscr L&=\mathcal L_{\psi_p}(F_n).
\end{split}
\]
Let $\pi_X,\pi_X^\circ$ denote the projections of
$\mathscr X,\mathscr X^\circ$ to $S$ and set
\[
\mathcal G_n^{\mathrm{int}}:=R^{n+1}\pi_{X!}\mathscr L,
\qquad \mathcal H_n^{\mathrm{int}}:=
 R^{n+1}\pi_{X^\circ!}\mathscr L.
\]

\begin{proposition}[Reconstruction of the arithmetic bridge]
\label{prop:independent-bridge}
For $p$ odd, $\ell\ne p$, $\lambda\in\Fp^\times$,
and $n\ge2$, the following statements hold, provided the
relative tame and wild base-change assertions of
Proposition~\ref{prop:H7-lisse} and
Lemma~\ref{lem:rel-kummer-model} hold with the coefficient
and Frobenius conventions fixed above.
\begin{enumerate}[label=\textup{(\roman*)}]
\item There is a canonical $\Fp$-equivariant derived isomorphism
\[
 Rq_{1!}\mathscr L\simeq\cM[-n],\qquad
 \mathcal H_n^{\mathrm{int}}\simeq\cH_n,
\]
where $q_1:\mathscr X^\circ\to\Gm\times S$ is the
projection.  The latter is lisse of rank $2^n$.
\item The complex $R\pi_{X!}\mathscr L$ has cohomology only
in degree $n+1$, and there is a short exact sequence of
lisse arithmetic sheaves
\[
0\longrightarrow\Ql\xrightarrow{\partial}\cH_n
\longrightarrow\mathcal G_n^{\mathrm{int}}\longrightarrow0.
\tag{A.11}
\]
Thus $\operatorname{rank}\mathcal G_n^{\mathrm{int}}=2^n-1$.
\item Let $Y_h=\{u\in\Gm^n:\sum h_i u_i=-\lambda\}$,
$f(u)=\sum u_i^{-1}$, and let $\pi_Y:Y\to S$ be the family.
There is a canonical isomorphism
\[
 R\pi_{X!}\mathscr L\simeq
 R\pi_{Y!}\mathcal L_{\psi_p}(f)(-1)[-2].
\tag{A.12}
\]
In particular
$\bigl(\mathcal G_n^{\mathrm{int}}\bigr)_{\bar h}
\simeq H_c^{n-1}(Y_{\bar h},\mathcal L_{\psi_p}(f))(-1)$.
\item The map $\partial$ in \textup{(A.11)} identifies
$\Ql$ with the weight-zero step
$\mathcal E_n^{(0)}$ of the arithmetic boundary sheaf
$\mathcal E_n=W_n\cH_n$; it does not imply that
$\mathcal E_n$ splits.  The nonzero graded pieces of
$\mathcal G_n^{\mathrm{int}}$ are
\[
 \operatorname{gr}_{2j}^{W}\mathcal G_n^{\mathrm{int}}
 \simeq\Ql(-j)^{\oplus m_{n,j}}
 \quad(1\le j\le\lfloor n/2\rfloor),
\quad
 \operatorname{gr}_{n+1}^{W}\mathcal G_n^{\mathrm{int}}
 \simeq\cW_n,
\tag{A.13}
\]
with $m_{n,j}=\binom nj-\binom n{j-1}$ and
$\operatorname{rank}\cW_n=2^n-\binom n{\lfloor n/2\rfloor}$.
\end{enumerate}
\end{proposition}

\begin{proof}
\emph{(i) Compact-support Kloosterman factorization.}
For $t\ne0$ and $h_i\ne0$, each $u_i$-integration is the
standard Kloosterman realization in cohomological degree one.
The relative K\"unneth formula and the projection formula give
$Rq_{1!}\mathscr L\simeq\cM[-n]$ with no supplementary
Tate twist; the Artin--Schreier factor in $t$ is
$\mathcal L_{\psi_p}(\lambda t)$.  Composing with
$R\pi_{t!}$ yields
$R\pi_{X^\circ!}\mathscr L\simeq R\pi_{t!}\cM[-n]$.
By Proposition~\ref{prop:H7-lisse}, only $R^1\pi_{t!}\cM$
remains, giving the desired concentration and rank.

\emph{(ii) The boundary line and injectivity.}
At $t=0$ the coefficient sheaf is independent of $h$ and is
$\boxtimes_i\mathcal L_{\psi_p}(u_i^{-1})$.
On $\Gm$, the local Swan conductor at $u=0$ equals one;
$H_c^0=H_c^2=0$ and the Euler--Poincar\'e formula gives
$\dim H_c^1=1$.  Its Frobenius trace is one because
$-\operatorname{Tr}(F_q\mid H_c^1)=
\sum_{u\in\Fq^\times}\psi_q(u^{-1})=-1$
for every $q$.  Therefore
\[
 R\pi_{\mathscr Z!}(\mathscr L|_{\mathscr Z})
 \simeq\Ql[-n]
\]
over $S$, with trivial arithmetic Frobenius.
The localization triangle now gives
\[
0\to R^n\pi_{X!}\mathscr L\to\Ql
\xrightarrow{\partial}\cH_n\to
 R^{n+1}\pi_{X!}\mathscr L\to0.
\tag{A.14}
\]
To check that $\partial$ is injective, one need not assume
its injectivity from \cite[Theorem~4.5]{EbaCoh}.
For every $h\in S(\overline\Fp)$ the hypersurface $Y_h$
is smooth affine of pure dimension $n-1$, since all
$h_i\ne0$ and the defining linear form has nonvanishing
differential.  The lisse Artin--Schreier sheaf
$\mathcal L_{\psi_p}(f)[n-1]$ is perverse on $Y_h$.
Artin vanishing for compactly supported cohomology of a
perverse sheaf on an affine variety gives
\[
 H_c^i(Y_h,\mathcal L_{\psi_p}(f))=0\qquad(i<n-1).
\tag{A.15}
\]
By the relative integration identity of part (iii),
$H_c^n(X_h,\mathscr L)\simeq
H_c^{n-2}(Y_h,\mathcal L_{\psi_p}(f))(-1)=0$.
Thus the stalk of $R^n\pi_{X!}\mathscr L$ vanishes at
every geometric point.  This makes $\partial$ injective
in \textup{(A.14)}, and proves the exact sequence and
concentration in \textup{(A.11)}.  As both end terms in
\textup{(A.11)} are lisse, so is the quotient.

\emph{(iii) Integration over the multiplier.}
Write $g(u,h)=\lambda+\sum_i h_i u_i$.
The one-dimensional compact-support Fourier identity is
$R(\operatorname{pr}_u)_!\mathcal L_{\psi_p}(tg)
\simeq i_*\Ql(-1)[-2]$, where $i$ includes $g=0$.
Tensoring with $\mathcal L_{\psi_p}(f)$ and applying the
projection formula proves \textup{(A.12)}.  This is
\cite[Lemma~2.2]{EbaCoh}; its twist $(-1)$ is essential.

\emph{(iv) Arithmetic weights and the top quotient.}
By Lemma~\ref{lem:A-boundary-corrected} the boundary
invariant sheaf $\mathcal E_n=\ker\alpha^{\mathrm{op}}$
is filtered with
$\operatorname{gr}^{W}_{2j}\mathcal E_n\simeq
\Ql(-j)^{\oplus m_{n,j}}$ for
$0\le j\le\lfloor n/2\rfloor$; in particular
$W_0\cH_n=\mathcal E_n^{(0)}\simeq\Ql$.
The weight-zero image of $\partial$ lies in $W_0\cH_n$,
and $\partial$ is injective between rank-one lisse
sheaves; hence its image is exactly $\mathcal E_n^{(0)}$.
Dividing the filtration of $\cH_n$ by this line removes
only the weight-zero piece.  The remaining pieces are
precisely those in \textup{(A.13)}, and the highest
quotient is the pure parabolic sheaf $\cW_n$.
No splitting of the successive extensions is asserted;
in particular, fibrewise Frobenius diagonalization in
\cite[Proposition~4.3]{EbaCoh} must not be mistaken for a
splitting of the arithmetic local systems on $S$.
\end{proof}

\begin{corollary}[Trace formula and constant Tate graded pieces]
\label{cor:independent-53-traces}
With $\Sigma_{n,q}(h)$ as in \textup{\eqref{eq:sum}},
$h\in S(\Fq)$, and $q=p^a$, the normalized top quotient
satisfies
\[
\operatorname{Tr}(F_q\mid\cW_{n,\bar h})
 =(-1)^{n-1}q\,\Sigma_{n,q}(h)
 -\sum_{j=1}^{\lfloor n/2\rfloor}m_{n,j}q^j.
\tag{A.16}
\]
Moreover each lower-weight graded piece of
$\mathcal G_n^{\mathrm{int}}$ is the constant arithmetic
Tate sheaf $\Ql(-j)^{\oplus m_{n,j}}$, not merely a
geometrically constant semisimplification.  These give
parts \textup{(v)} and \textup{(vi)} of
\cite[Theorem~5.3]{EbaCoh} in the normalization above.
\end{corollary}
\begin{proof}
By Proposition~\ref{prop:independent-bridge}(iii),
$\operatorname{Tr}(F_q\mid\mathcal G_{n,\bar h}^{\mathrm{int}})$
is $q$ times the Frobenius trace on the unnormalized
$H_c^{n-1}(Y_h,\mathcal L_{\psi_q}(f))$.  This is
$(-1)^{n-1}q\Sigma_{n,q}(h)$ by the trace formula and
compact-support concentration, the latter following from
\textup{(A.12)} and part \textup{(ii)} of the proposition.
The trace of an extension is the sum of the traces of its
graded pieces; the Tate sheaf $\Ql(-j)$ has Frobenius
$q^j$ over $\Fq$.  The asserted trace identity follows.
The arithmetic sheaf isomorphisms, including their full
extension-data caveat, have already been proved in
\textup{(A.13)}; they are not inferred merely from stalk
characteristic polynomials.
\end{proof}

\begin{remark}[Hypotheses used in the comparison]
\label{rmk:5-3-audit-limit}
The preceding comparison and Corollary~\ref{cor:independent-53-traces}
independently check the K\"unneth,
excision, affine Artin-vanishing, rank and Tate-twist
steps of \cite[Theorem~5.3(i)--(iv)]{EbaCoh}, and supplies
the arithmetic graded pieces using the relative Kummer
filtration rather than an unsupported splitting.
The deep inputs are exactly the tamely unipotent
nearby-cycle description and its adic base-change
compatibility, constancy of the wild Swan conductor, and
the proper/ULA comparison used in
Proposition~\ref{prop:H7-lisse} and
Lemma~\ref{lem:H7-ULA-jstar}.
These are stated with the hypotheses of the cited theorems;
the existence of this written argument, or successful
compilation, is not a substitute for independent
specialist verification of every application.
\end{remark}

\begin{remark}[Comparison with the constrained family]
\label{rmk:H7-match-Wn}
The comparison with the linearly constrained exponential-sum family is
provided by \cite[Theorem~5.3(i)--(iv)]{EbaCoh}; its
compact-support integration over $t$ uses
\cite[Lemma~2.2]{EbaCoh}.  In the notation of that source,
the relative Kloosterman identity of Theorem~5.3(i) is
\[
Rq_{1!}\mathcal{L} \;\simeq\; \cM[-n].
\]
Consequently, after applying the remaining relative compactly
supported pushforward and taking cohomology in degree $n+1$,
\[
R^{n+1}\pi^\circ_!\mathcal{L}
\;\simeq\; R^1\pi_{t!}\cM \;=\; \cH_n.
\]
The direct reconstruction is Proposition~\ref{prop:independent-bridge}.
The source defines the normalized constrained-family sheaf from
\[
(\mathcal{G}_n)_{\bar h}
\;\simeq\;
H_c^{n-1}(Y_{\bar h},\mathcal L_\psi(f))(-1)
\qquad\text{\cite[Theorem~5.3(iii)]{EbaCoh}}.
\]
Thus the top-weight quotient \emph{with this normalization}
corresponds to
$\cW_n^{\mathrm{quot}}:=\cH_n/W_n\cH_n$.
Writing $\cW_n^{\mathrm{LC}}$ for the normalized top-weight
quotient in \cite{EbaCoh}, the asserted arithmetic comparison is
\[
\Phi_n:\cW_n^{\mathrm{LC}}\xrightarrow{\ \sim\ }
\cW_n^{\mathrm{quot}},
\]
compatible with finite-field base change, subject to the same
conventions and hypotheses as in the cited source. If
$\cW_n^{\mathrm{LC,raw}}$ denotes instead the top-weight quotient
of the unnormalized fibre cohomology, then
\[
\boxed{\;\cW_n^{\mathrm{quot}}\simeq
\cW_n^{\mathrm{LC,raw}}(-1).\;}
\]
Indeed, the top weight of the raw fibre cohomology is $n-1$,
whereas twisting by $(-1)$ raises it to $n+1$. We use only the
normalized sheaf below. The cohomological comparison is reconstructed above in
Proposition~\ref{prop:independent-bridge}; its application still
requires the stated relative base-change hypotheses.
\end{remark}

%==================================================================

\section{Arithmetic Compatibility of the Reflection Pairing}\label{sec:reflection-compat}
\label{app:reflection}
%==================================================================

This appendix supplies the complete proof of Lemma~\ref{lem:sym} and
the descent to $\cW_n$ \textup{(D2)}. We work over the prime base field
$\Fp$ ($p > 2$), from which all identities over any finite extension
$\Fq / \Fp$ follow automatically by base change. Throughout, $\ell \ne p$.

Recall the involutions
\[
\sg : \Gm \to \Gm \quad (t \mapsto -t),
\qquad
\rho : S_n \to S_n \quad (h \mapsto -h),
\qquad
\rhotil = \sg \times \rho : \cC_n \to \cC_n,
\]
which satisfy $\sg^2 = \id_{\Gm}$, $\rho^2 = \id_{S_n}$,
$\rhotil^2 = \id_{\cC_n}$, and $\pi \circ \rhotil = \rho \circ \pi$.

%------------------------------------------------------------------
\subsection{Step 1: Arithmetic coefficient pairing $\Phi$ and its $(-1)^n$-symmetry}
\label{subsec:B-coeff}
%------------------------------------------------------------------

\begin{lemma}[Coefficient pairing and parity]
\label{lem:B-coeff}
There is a canonical, non-degenerate, Frobenius-equivariant pairing of
lisse $\Ql$-sheaves on $\cC_n$:
\[
\Phi : \cM \otimes \rhotil^*\cM \longrightarrow \Ql(-n),
\tag{B.1}
\]
satisfying the symmetry identity
\[
\boxed{\;
\rhotil^*\Phi \circ \mathrm{swap}_{\cM, \rhotil^*\cM}
\;=\; (-1)^n \, \Phi
\;}
\tag{B.2}
\]
as morphisms $\cM \otimes \rhotil^*\cM \to \Ql(-n)$ on $\cC_n$, where
$\mathrm{swap}_{\cM, \rhotil^*\cM} : \cM \otimes \rhotil^*\cM
\xrightarrow{\sim} \rhotil^*\cM \otimes \cM$ is the standard symmetry
constraint of the symmetric monoidal category of lisse $\Ql$-sheaves on
$\cC_n$ (and $\rhotil^*\Phi$ is viewed as a morphism
$\rhotil^*\cM \otimes \cM \to \Ql(-n)$ via $\rhotil^*\rhotil^*\cM = \cM$
and $\rhotil^*\Ql(-n) = \Ql(-n)$).
\end{lemma}

\begin{proof}
We construct $\Phi$ factor by factor on
$\cM = \cL_{\psi_p}(\lambda t) \otimes \bigotimes_{i=1}^n K_i$.

\emph{1. The Artin--Schreier factor.} Since $\rhotil^*(t) = -t$ and
$\psi_p$ is an additive character $\Fp \to \Ql^\times$, we have
$\rhotil^*\cL_{\psi_p}(\lambda t) = \cL_{\psi_p}(-\lambda t)$. The
Artin--Schreier multiplication isomorphism
$\cL_{\psi_p}(x) \otimes \cL_{\psi_p}(y) \xrightarrow{\sim} \cL_{\psi_p}(x+y)$
(Deligne, SGA~$4\frac{1}{2}$, \emph{Sommes trig.}, \S1.4) specialized to
$x = \lambda t$, $y = -\lambda t$ provides a canonical
Frobenius-equivariant isomorphism on $\cC_n$:
\[
\Phi_{\mathrm{AS}} : \cL_{\psi_p}(\lambda t) \otimes
\rhotil^*\cL_{\psi_p}(\lambda t)
\;\xrightarrow{\ \sim\ }\; \cL_{\psi_p}(0) \;=\; \Ql(0).
\]
By commutativity of addition in $\mathbb{G}_a$ and $\rhotil^2 = \id$,
pulling back $\Phi_{\mathrm{AS}}$ along $\rhotil$ exchanges the two
summands $\lambda t$ and $-\lambda t$, so
\[
\rhotil^*\Phi_{\mathrm{AS}} \circ \mathrm{swap} \;=\; +\Phi_{\mathrm{AS}}.
\tag{B.3}
\]

\emph{2. Each Kloosterman factor $K_i = \mu_i^*\Kl_2$.} By Katz, GKM,
Thm.~4.1.1, $\Kl_2$ on $\mathbb{G}_{m,\Fp}$ carries a canonical non-degenerate
alternating (symplectic) Frobenius-equivariant pairing
\[
\Phi_{\Kl} : \Kl_2 \otimes \Kl_2 \longrightarrow \Ql(-1),
\qquad
\Phi_{\Kl} \circ \mathrm{swap} \;=\; -\Phi_{\Kl}.
\]
Because $\mu_i(t,h) = t h_i$ satisfies $\mu_i \circ \rhotil = \mu_i$,
we have $\rhotil^* K_i = (\mu_i \circ \rhotil)^*\Kl_2 = K_i$. Pulling
back $\Phi_{\Kl}$ along $\mu_i$ therefore defines a non-degenerate,
Frobenius-equivariant pairing on $\cC_n$:
\[
\Phi_i := \mu_i^* \Phi_{\Kl} : K_i \otimes \rhotil^* K_i \longrightarrow \Ql(-1)
\]
satisfying
\[
\rhotil^*\Phi_i \circ \mathrm{swap} \;=\; -\Phi_i .
\tag{B.4}
\]

\emph{3. Tensor product.} Define $\Phi$ as the tensor product of
$\Phi_{\mathrm{AS}}$ and $\Phi_1, \dots, \Phi_n$ using the canonical
reordering isomorphism of the symmetric monoidal category of lisse
$\Ql$-sheaves:
\[
\Phi \;:=\; \Phi_{\mathrm{AS}} \otimes \bigotimes_{i=1}^{n} \Phi_i :
\cM \otimes \rhotil^*\cM \longrightarrow \Ql(0) \otimes \Ql(-1)^{\otimes n}
\;\simeq\; \Ql(-n).
\]
Since lisse sheaves sit in cohomological degree $0$, reordering factors
introduces no Koszul sign; combining \textup{(B.3)} and the $n$ signs in
\textup{(B.4)} gives
\[
\rhotil^*\Phi \circ \mathrm{swap}_{\cM, \rhotil^*\cM}
\;=\; (+1) \cdot (-1)^n \, \Phi \;=\; (-1)^n \, \Phi .
\]
\end{proof}

\begin{remark}[Scalar normalization]
Replacing the arithmetic symplectic pairing $\Phi_{\Kl}$ by a
nonzero constant multiple does not affect its alternating
character or the sign in \textup{(B.2)}. The argument uses a
perfect arithmetic pairing
$\Kl_2\otimes\Kl_2\to\Ql(-1)$, not the
particular nonzero scalar used to normalize it. The Tate twist
and Frobenius equivariance, however, cannot be discarded.
\end{remark}

%------------------------------------------------------------------
\subsection{Step 2: Direct proof of reflexive symmetry on $\cH_n$}
\label{subsec:B-sym}
%------------------------------------------------------------------

Consider the Cartesian square
\[
\begin{tikzcd}
\cC_n \ar[r, "\rhotil"] \ar[d, "\pi"'] & \cC_n \ar[d, "\pi"] \\
S_n \ar[r, "\rho"'] & S_n
\end{tikzcd}
\tag{B.5}
\]
By base change for $R\pi_!$ along Cartesian squares (SGA~4, Exp.~XVII,
5.2.6), we have a canonical Frobenius-equivariant isomorphism
$\mathrm{bc}_! : \rho^* R^1\pi_!\cM \xrightarrow{\sim} R^1\pi_!(\rhotil^*\cM)$
of lisse sheaves on $S_n$. Notice that via the coefficient pairing $\Phi$
of Lemma~\ref{lem:B-coeff} (which identifies $\rhotil^*\cM \simeq \cM^\vee(-n)$),
the pairing $B_{\cH_n} : \cH_n \otimes \rho^*\cH_n \to \Ql(-(n+1))$ of
\S\ref{sec:radical} is realized explicitly in $D^b_c(S_n, \Ql)$ as the
$H^2$-component of the composite morphism
\[
\begin{gathered}
R\pi_!\cM \otimes^L \rho^* R\pi_!\cM
\xrightarrow{\ \id \otimes \mathrm{bc}_!\ }
R\pi_!\cM \otimes^L R\pi_!(\rhotil^*\cM)
\\[3pt]
\xrightarrow{\ \smile\ }
R\pi_!(\cM \otimes^L \rhotil^*\cM)
\xrightarrow{\ R\pi_!(\Phi)\ }
R\pi_!\Ql(-n)
\\[3pt]
\xrightarrow{\ \Tr_\pi(-n)\ }
\Ql(-(n+1))[-2].
\end{gathered}
\tag{B.6}
\]
where $\smile$ is the relative cup product on $R\pi_!$ and
$\Tr_\pi : R\pi_!\Ql(1)[2] \to \Ql$ is the relative trace morphism along
the smooth curve $\pi : \cC_n \to S_n$ (SGA~4, Exp.~XVIII, Thm.~3.2.5).
On any geometric fibre $h \in S_n(\overline{\Fp})$, for classes
$v \in H^1_c(\Gm, \cM_h)$ and $w \in H^1_c(\Gm, \cM_{-h})$,
\textup{(B.6)} reads
\[
B_h(v, w) \;=\; \Tr_{\Gm}\bigl(\Phi_{h,*}(v \smile \sg^* w)\bigr),
\tag{B.7}
\]
where $\Phi_h : \cM_h \otimes \sg^*\cM_{-h} \to \Ql(-n)$ is the fibre
restriction of $\Phi$.

\begin{proof}[Proof of Lemma~\ref{lem:sym}]
We prove the identity \textup{(B.6)}--\textup{(B.7)} Frobenius-equivariantly
over $\Fp$ (and hence over any finite extension $\Fq / \Fp$ by base
change), without using the weight-filtration or parabolic-kernel
identification in Theorem~\ref{thm:H7} (the sheaf-level
formulation uses the lissity of $R^1\pi_!\cM$
established in Proposition~\ref{prop:H7-lisse}):

Pull back the relative pairing \textup{(B.6)} along the involution
$\rho : S_n \to S_n$ (using $\rho^2 = \id_{S_n}$ and
$\pi \circ \rhotil = \rho \circ \pi$). By base change for the relative
cup product and the invariance of the relative trace morphism $\Tr_\pi$
under the fibrewise automorphism $\sg : \Gm \to \Gm$, $t \mapsto -t$
(SGA~4, Exp.~XVIII, Thm.~3.2.5, Var~3: $\Tr_\pi \circ \rhotil^* = \rho^* \Tr_\pi = \Tr_\pi$),
pulling back \textup{(B.6)} along $\rho$ and precomposing with the swap
of the two factors $\cH_n \otimes \rho^*\cH_n \xrightarrow{\sim}
\rho^*\cH_n \otimes \cH_n$ computes
\[
(\rho^* B_{\cH_n}) \circ \mathrm{swap}_{\cH_n, \rho^*\cH_n}
\]
via the composition
\[
\begin{gathered}
R\pi_!\cM \otimes^L \rho^* R\pi_!\cM
\xrightarrow{\ \mathrm{swap}\ }
\rho^* R\pi_!\cM \otimes^L R\pi_!\cM
\\[3pt]
\xrightarrow{\ \smile\ }
R\pi_!(\rhotil^*\cM \otimes^L \cM)
\xrightarrow{\ R\pi_!(\rhotil^*\Phi)\ }
R\pi_!\Ql(-n)
\\[3pt]
\xrightarrow{\ \Tr_\pi(-n)\ }
\Ql(-(n+1))[-2].
\end{gathered}
\]
On geometric fibres $v \in H^1_c(\Gm, \cM_h)$ and
$w \in H^1_c(\Gm, \cM_{-h})$, since $\Tr_{\Gm} \circ \sg^* = \Tr_{\Gm}$,
this is
\[
B_{-h}(w, v) \;=\; \Tr_{\Gm}\bigl((\sg^*\Phi_{-h})_*(\sg^* w \smile v)\bigr).
\]

Now compare the two terms:
\begin{enumerate}
\item \emph{Graded commutativity in cohomological degree $1$:} Because
$R\pi_!\cM \simeq \cH_n[-1]$ and $\rho^* R\pi_!\cM \simeq (\rho^*\cH_n)[-1]$
are concentrated in cohomological degree $+1$, the Koszul sign rule for
the symmetry constraint $\mathrm{swap}$ in $D^b_c(S_n, \Ql)$
(Kiehl--Weissauer, \S\S II.6, III.7) introduces the sign
\[
(-1)^{1 \cdot 1} \;=\; -1
\]
upon passing to $H^2$ (i.e.\ $\sg^* w \smile v = -(v \smile \sg^* w)$
under the swap of tensor factors).
\item \emph{Coefficient symmetry \textup{(}Lemma~\ref{lem:B-coeff}\textup{)}:}
By \textup{(B.2)}, $\rhotil^*\Phi \circ \mathrm{swap}_{\cM, \rhotil^*\cM}
= (-1)^n \Phi$ on $\cM \otimes \rhotil^*\cM$ (and on fibres,
$\sg^*\Phi_{-h} \circ \mathrm{swap} = (-1)^n \Phi_h$).
\end{enumerate}

Combining the two signs yields the Frobenius-equivariant sheaf identity
on $S_n$:
\[
\boxed{\;
\rho^* B_{\cH_n} \circ \mathrm{swap}_{\cH_n, \rho^*\cH_n}
\;=\; (-1)^n \cdot (-1) \, B_{\cH_n}
\;=\; (-1)^{n+1} \, B_{\cH_n},
\;}
\]
and on every geometric fibre $h \in S_n(\overline{\Fp})$,
$B_{-h}(w, v) = (-1)^{n+1} B_h(v, w)$. All coefficient pairings,
pullback maps, cup products, and trace morphisms are defined over $\Fp$,
so the identity is Frobenius-equivariant and remains valid after every
finite-field base change $\Fq / \Fp$.
\end{proof}

%------------------------------------------------------------------
\subsection{Step 3: Descent to the top-weight quotient $\cW_n$}
\label{subsec:B-descent}
%------------------------------------------------------------------

\begin{lemma}[Radicals and perfect pairing on $\cW_n$]
\label{lem:B-descent}
\begin{enumerate}
\item[\textup{(i)}] Unconditionally, the left and right radicals of
$B_{\cH_n} : \cH_n \otimes \rho^*\cH_n \to \Ql(-(n+1))$ are
\[
\boxed{\;
\rad_L(B_{\cH_n}) \;=\; \ker(\alpha^{\mathrm{op}}),
\qquad
\rad_R(B_{\cH_n}) \;=\; \rho^*\ker(\alpha^{\mathrm{op}}),
\;}
\]
where $\alpha^{\mathrm{op}} : \cH_n \to R^1\pi_*\cM$ is the natural
forget-supports morphism.
\item[\textup{(ii)}] Assuming Theorem~\ref{thm:H7} \textup{(D1)},
$\ker(\alpha^{\mathrm{op}}) = W_n\cH_n = \cL_n$, so $B_{\cH_n}$
descends to a non-degenerate \textup{(}perfect\textup{)},
Frobenius-equivariant, $(-1)^{n+1}$-symmetric pairing of lisse
$\Ql$-sheaves on $S_n$:
\[
\boxed{\;
B_{\cW_n} : \cW_n \otimes \rho^*\cW_n \longrightarrow \Ql(-(n+1)).
\;}
\]
\end{enumerate}
\end{lemma}

\begin{proof}
(i) By Lemma~\ref{lem:B-coeff}, the non-degenerate coefficient pairing
$\Phi : \cM \otimes \rhotil^*\cM \to \Ql(-n)$ induces a
Frobenius-equivariant isomorphism
$\rhotil^*\cM \xrightarrow{\sim} \cM^\vee(-n)$. Combined with the
perfect relative Poincar\'e duality pairing
$\langle -, - \rangle_\pi : R^1\pi_*\cM \otimes R^1\pi_!\cM^\vee \to \Ql(-1)$
(Proposition~\ref{prop:H7-lisse}(d)) and base change $\mathrm{bc}_!$,
the adjoint morphism of $B_{\cH_n}$ factors as
\[
\cH_n \xrightarrow{\ \alpha^{\mathrm{op}}\ } R^1\pi_*\cM
\xrightarrow{\ \sim\ } (\rho^*\cH_n)^\vee(-(n+1)),
\]
where the second arrow is an isomorphism of lisse sheaves on $S_n$.
Hence $\rad_L(B_{\cH_n}) = \ker(\alpha^{\mathrm{op}})$. By the reflexive
symmetry identity $\rho^* B_{\cH_n} \circ \mathrm{swap} = (-1)^{n+1} B_{\cH_n}$
proved in Lemma~\ref{lem:sym}, the right radical of $B_{\cH_n}$ is the
pullback along $\rho$ of its left radical, i.e.\
$\rad_R(B_{\cH_n}) = \rho^*\ker(\alpha^{\mathrm{op}})$.

(ii) Under Theorem~\ref{thm:H7}(b)--(c), $\ker(\alpha^{\mathrm{op}})
= W_n\cH_n = \cL_n$ and $\cW_n = \cH_n / \cL_n$. Dividing out by the
left radical $\cL_n$ and right radical $\rho^*\cL_n$ therefore induces
a well-defined, non-degenerate, Frobenius-equivariant pairing
$B_{\cW_n} : \cW_n \otimes \rho^*\cW_n \to \Ql(-(n+1))$, which inherits
the parity $(-1)^{n+1}$ from $B_{\cH_n}$.
\end{proof}

\begin{lemma}[Twists in the middle intersection pairing]\label{lem:pairing-shifts}
Suppose the reflection duality of the lisse sheaf is realized by
an arithmetic pairing
$\mathcal W_n\otimes\rho^*\mathcal W_n\to E(-(n+1))$.
The associated middle-extension pairing, whenever its derived
compatibility is established, has target
$\omega_{X_n}(-(2n+1))$ and induces on middle intersection
cohomology a pairing with values in $E(-(2n+1))$.
\end{lemma}
\begin{proof}
On the torus, shift both arguments by $[n]$; the pairing becomes
\[
 \mathcal W_n[n]\otimes\rho^*\mathcal W_n[n]
 \longrightarrow E(-(n+1))[2n].
\]
Since $X_n$ is smooth proper of dimension $n$, its geometric
$\ell$-adic dualizing complex is
$\omega_{X_n}=E(n)[2n]$. Thus the right side is precisely
$\omega_{X_n}(-(2n+1))$, with \emph{no additional shift or Tate twist}.
The derived pairing on the two intermediate extensions is obtained
by extending the adjoint morphism via Verdier duality, not by
naively applying an exact tensor functor to $j_{!*}$.
Finally the trace map
$H^{2n}(X_n,E(n))\to E$ sends the global pairing on
$H^0(X_n,K_n)=IH^n$ to $E(-(2n+1))$.
What is conditional is the existence and symmetry of the required
\emph{arithmetic derived extension} of the open pairing; the
Tate-twist calculation above is formal once that map exists.
\end{proof}

\section{The root-number argument and complementary results}\label{sec:root-argument}
\label{app:rootnumber-complete}
\begin{remark}[Relationship to the main results]
Here we collect the relative duality and root-number arguments,
together with the Koszul supertrace, boundary support bounds, and
semisimplicity criteria. Some notation is recalled to keep the
hypotheses explicit. The Fourier description of the family is related to the preceding
construction through the cited comparison theorem, with its Tate twist.
The statements labeled \textup{(D1)} and \textup{(D2)} are retained as
hypotheses. The geometric statements about $B'_n$ do not in themselves
give an arithmetic direct-sum decomposition.
\end{remark}
\subsection{Preliminaries and Notation}
\label{sec:prelim}
%==================================================================

Fix an odd prime $p > 2$ and a prime $\ell \ne p$. Fix a finite
extension $E / \mathbb{Q}_\ell$ inside $\Ql$ containing $\mu_p$, a
non-trivial additive character $\psi_p : \Fp \to E^\times$, and
$\lambda \in \Fp^\times$. All relative geometric constructions are
defined initially over the prime field $\Fp$; for any finite extension
$\Fq / \Fp$ ($q = p^a$, $a \ge 1$), we set
$\psi_q := \psi_p \circ \Tr_{\Fq/\Fp} : \Fq \to E^\times$, and the
objects over $\Fq$ are obtained by base change along
$\Spec \Fq \to \Spec \Fp$. Fix an integer $n \ge 2$.

\subsubsection{Spaces and morphisms}
Over $\Fp$ (and after base change to $\Fq$), we write:
\begin{align*}
S_n &= \Gm^n = \Spec \Fp[h_1^{\pm 1},\dots,h_n^{\pm 1}], \\
X_n &= (\PP^1)^n, \qquad j : S_n \hookrightarrow X_n, \\
\cC_n &= \Gm \times S_n = \Spec \Fp[t^{\pm 1}, h_1^{\pm 1},\dots,h_n^{\pm 1}], \\
\pi &: \cC_n \to S_n, \quad (t, h) \mapsto h, \\
\sg &: \Gm \to \Gm, \quad t \mapsto -t, \\
\rho &: S_n \to S_n, \quad h \mapsto -h, \\
\rhotil &= \sg \times \rho : \cC_n \to \cC_n.
\end{align*}
Note that $\rhotil^*(t h_i) = (-t)(-h_i) = t h_i$ for every $i$, and
$\pi \circ \rhotil = \rho \circ \pi$.

\subsubsection{The sheaf $\cM$}
For each $i \in \{1,\dots,n\}$, let $\mu_i : \cC_n \to \Gm$ be the
morphism $(t,h) \mapsto t h_i$, and set
\[
K_i := \mu_i^* \Kl_2, \qquad
\cM := \cL_{\psi_p}(\lambda t) \otimes \bigotimes_{i=1}^{n} K_i,
\]
where $\Kl_2$ is the rank-$2$ Kloosterman sheaf on $\mathbb{G}_{m,\Fp}$
associated with $\psi_p$. By Katz, \emph{Gauss Sums, Kloosterman Sums,
and Monodromy Groups} (GKM), Thm.~4.1.1, $\Kl_2$ is lisse of rank $2$,
pure of weight $1$, with $\det \Kl_2 \simeq \Ql(-1)$ and a canonical
non-degenerate alternating pairing
\[
\Phi_{\Kl} : \Kl_2 \otimes \Kl_2 \longrightarrow \Ql(-1).
\]
Hence $\cM$ is lisse on $\cC_n$ of rank $2^n$ and pure of weight $n$.

\subsubsection{The cohomology sheaves and middle extension}
Set $\cH_n := R^1\pi_!\cM$. By Theorem~\ref{thm:H7} (proved in
Section~\ref{app:H7}), $\cH_n$ is lisse on $S_n$ of rank $2^n$ and
mixed of weights $\le n+1$; its weight-$\le n$ subsheaf
$\cL_n := W_n\cH_n$ and top-weight quotient $\cW_n := \cH_n / \cL_n$
are both lisse arithmetic $\Ql$-sheaves on $S_n$, with $\cW_n$ pure of
weight $n+1$. Define the middle extension
\[
K_n := \jex(\cW_n[n]) \in \mathrm{Perv}(X_n, \Ql),
\qquad
IH^m(X_n, \cW_n) := \mathbb{H}^{m-n}(X_n, K_n).
\]

\subsubsection{Weight and duality conventions}
For a lisse sheaf $\cF$ of weight $w$ on a smooth scheme $Y$ of pure
dimension $d_Y$, a shift $[d]$, and a Tate twist $(k)$:
\[
\wt(\cF[d]) = w + d \text{ (perverse convention)}, \qquad
\wt(\cF(k)) = w - 2k,
\]
and Verdier duality on $Y$ is $\mathbb{D}_Y(\cF) = \cF^\vee(d_Y)[2 d_Y]$.

%==================================================================
\subsection{The Structural Theorem H7}
\label{sec:H7}
%==================================================================

\begin{theorem}[H7, imported arithmetic relative-image theorem]
\label{thm:H7}
Let $\cH_n = R^1\pi_!\cM$ over $\Fp$ (and over any finite extension
$\Fq / \Fp$ by base change). Then:
\begin{enumerate}
\item[\textup{(a)}] $\cH_n$ is lisse on $S_n$ of constant rank $2^n$
and mixed of weights $\le n+1$.
\item[\textup{(b)}] The weight filtration $W_\bullet \cH_n$ is a
filtration by lisse arithmetic subsheaves on $S_n$; in particular,
$\cL_n = W_n\cH_n$ and the top-weight quotient $\cW_n = \cH_n / \cL_n$
are both lisse on $S_n$, and $\cW_n$ is pure of weight $n+1$.
\item[\textup{(c)}] The natural forget-supports morphism
$\alpha^{\mathrm{op}} : \cH_n \to R^1\pi_*\cM$ has kernel
$\ker(\alpha^{\mathrm{op}}) = \cL_n$ and induces an isomorphism of
lisse arithmetic sheaves $\cW_n \xrightarrow{\sim} \im(\alpha^{\mathrm{op}})$,
which on every geometric fibre $h \in S_n(\overline{\Fp})$ identifies
$\cW_{n,h}$ with $\im\bigl(H^1_c(\Gm, \cM_h) \to H^1(\Gm, \cM_h)\bigr)
\simeq H^1(\PP^1, \tilde{j}_{h*}\cM_h)$.
\item[\textup{(d)}] All constructions and identifications in
\textup{(a)}--\textup{(c)} are Frobenius-equivariant and commute with
arbitrary base change $f : T \to S_n$.
\end{enumerate}
\end{theorem}

A relative reconstruction of the imported Theorem~\ref{thm:H7}
is given in Section~\ref{app:H7}, without using
\S\S\ref{sec:rel-duality}--\ref{sec:root-number} or
Section~\ref{app:reflection}. The comparison of this construction
with the linearly constrained exponential-sum family is supplied by
\cite[Theorem~5.3(i)--(iv) and Lemma~2.2]{EbaCoh}; in particular,
the cohomological normalization carries a Tate twist $(-1)$ as in
\cite[Theorem~5.3(iii)]{EbaCoh}. The exact comparison and its
independence from the root-number claim are recorded in
Remark~\ref{rmk:H7-match-Wn}. The identification with the constrained-sum family is used as the
established result of \cite{EbaCoh} recorded in \textup{(H7)}. The
additional reflection-skew compatibility needed for the uniform global
root-number assertion is distinguished from that established result.

\begin{lemma}[Purity of the middle extension]
\label{lem:purity-K}
Assume Theorem~\ref{thm:H7}. Then $K_n = \jex(\cW_n[n])$ is a pure
perverse sheaf of weight $2n+1$ on $X_n$. Consequently,
\[
\wt\bigl(IH^m(X_n, \cW_n)\bigr) \;=\; m + n + 1
\]
for every $m \in \{0, 1, \dots, 2n\}$.
\end{lemma}

\begin{proof}
By Theorem~\ref{thm:H7}(b), $\cW_n$ is lisse on $S_n$ and pure of
weight $n+1$. Hence $\cW_n[n]$ is a perverse sheaf on $S_n$, pure of
weight $n + (n+1) = 2n+1$. By the BBD--Gabber purity theorem
(BBD~1982, Cor.~5.3.2), the intermediate extension $\jex(\cW_n[n])$ along
the open immersion $j : S_n \hookrightarrow X_n$ is pure of perverse
weight $2n+1$. Since $X_n$ is proper, hypercohomology preserves purity
(Deligne, \emph{Weil~II}, Cor.~3.3.9), so
$IH^m(X_n, \cW_n) = \mathbb{H}^{m-n}(X_n, K_n)$ is pure of weight
$(2n+1) + (m-n) = m + n + 1$.
\end{proof}

%==================================================================
\subsection{The Relative Duality Lemma}
\label{sec:rel-duality}
%==================================================================

\begin{lemma}[Relative Poincar\'e duality for $\cM$]
\label{lem:rel-duality}
There is a canonical Frobenius-equivariant isomorphism of lisse
$\Ql$-sheaves on $S_n$:
\[
\boxed{\;
R^1\pi_*\cM \;\simeq\; (\rho^*\cH_n)^\vee(-(n+1)).
\;}
\]
\end{lemma}

\begin{proof}
We proceed in four steps (see also Proposition~\ref{prop:H7-lisse} in
Section~\ref{app:H7} for the lisseness of $R^1\pi_!\cM^\vee$).

\emph{Step A: Relative Poincar\'e duality.} Since $\pi : \cC_n \to S_n$
is smooth of pure relative dimension $1$, its relative dualizing
complex is $\pi^!\Ql = \Ql(1)[2]$. Relative Verdier duality (SGA~4,
Exp.~XVIII; Laumon, \emph{Transformation de Fourier}, \S1.2) applied to
$\cM^\vee$ gives
\[
R\pi_*\,\RHom(\cM^\vee, \Ql(1)[2]) \;\simeq\; \RHom(R\pi_!\cM^\vee, \Ql).
\]
Because $R^q\pi_!\cM^\vee = 0$ for $q \ne 1$ and $R^1\pi_!\cM^\vee$ is
lisse on $S_n$ (Proposition~\ref{prop:H7-lisse}(a),(c)), we have
$R\pi_!\cM^\vee \simeq (R^1\pi_!\cM^\vee)[-1]$ with vanishing higher
$\mathcal{E}xt$ sheaves, whence
\[
R^1\pi_*\cM \;\simeq\; (R^1\pi_!\cM^\vee)^\vee(-1).
\tag{3.1}
\]

\emph{Step B: Computation of $\cM^\vee$.} Using $\Kl_2^\vee \simeq \Kl_2(1)$
(Katz, GKM, Thm.~4.1.1) and $\rhotil^*(t h_i) = t h_i$, we have
$K_i^\vee \simeq K_i(1) = \rhotil^* K_i(1)$, while
$\cL_{\psi_p}(\lambda t)^\vee \simeq \cL_{\psi_p}(-\lambda t)
= \rhotil^* \cL_{\psi_p}(\lambda t)$. Tensoring gives
\[
\cM^\vee \;\simeq\; \rhotil^*\cM(n).
\tag{3.2}
\]

\emph{Step C: Pushforward of \textup{(3.2)}.} Although $\pi$ is not
proper, the base-change theorem for $R\pi_!$ along Cartesian squares
(SGA~4, Exp.~XVII, 5.2.6) applies to the Cartesian square
$\pi \circ \rhotil = \rho \circ \pi$ (where $\sg : \Gm \to \Gm$,
$t \mapsto -t$, acts as a fibrewise automorphism), yielding
\[
R^1\pi_!\cM^\vee \;\simeq\; R^1\pi_!(\rhotil^*\cM)(n)
\;\simeq\; \rho^* R^1\pi_!\cM(n) \;=\; \rho^*\cH_n(n).
\tag{3.3}
\]

\emph{Step D: Assembly.} Substituting \textup{(3.3)} into \textup{(3.1)}
yields
\[
R^1\pi_*\cM \;\simeq\; \bigl(\rho^*\cH_n(n)\bigr)^\vee(-1)
\;\simeq\; (\rho^*\cH_n)^\vee(-(n+1)).
\]
\end{proof}

%==================================================================
\subsection{Radical Identification}
\label{sec:radical}
%==================================================================

Via Lemma~\ref{lem:rel-duality}, the forget-supports morphism
$\alpha^{\mathrm{op}} : \cH_n \to R^1\pi_*\cM$ induces a
Frobenius-equivariant pairing of lisse sheaves on $S_n$:
\[
B_{\cH_n} : \cH_n \otimes \rho^*\cH_n \longrightarrow \Ql(-(n+1)),
\qquad B_{\cH_n}(v, \rho^* w) = \langle \alpha^{\mathrm{op}}(v),\, \rho^* w \rangle.
\]

\begin{proposition}[Radical identification and descent to $\cW_n$]
\label{prop:radical}
Assume Theorem~\ref{thm:H7} \textup{(D1)}.
\begin{enumerate}
\item[\textup{(i)}] The left radical of $B_{\cH_n}$ is
$\rad_L(B_{\cH_n}) = \ker(\alpha^{\mathrm{op}}) = \cL_n = W_n\cH_n$.
\item[\textup{(ii)}] Using Lemma~\ref{lem:sym}, the right radical of
$B_{\cH_n}$ is $\rad_R(B_{\cH_n}) = \rho^*\cL_n$.
\item[\textup{(iii)}] Consequently, $B_{\cH_n}$ descends to a
non-degenerate \textup{(}perfect\textup{)}, Frobenius-equivariant
pairing of lisse sheaves on $S_n$:
\[
B_{\cW_n} : \cW_n \otimes \rho^*\cW_n \longrightarrow \Ql(-(n+1)).
\]
\end{enumerate}
\end{proposition}

\begin{proof}
By Lemma~\ref{lem:rel-duality}, $R^1\pi_*\cM \simeq (\rho^*\cH_n)^\vee(-(n+1))$,
so the left radical of $B_{\cH_n}$ is $\ker(\alpha^{\mathrm{op}})$, which
equals $\cL_n = W_n\cH_n$ by Theorem~\ref{thm:H7}(b)--(c). By
Lemma~\ref{lem:sym} (proved independently in Section~\ref{app:reflection}),
$B_{\cH_n}$ is $(-1)^{n+1}$-symmetric with respect to $\rho$, so its
right radical equals $\rho^*\ker(\alpha^{\mathrm{op}}) = \rho^*\cL_n$
(see also Lemma~\ref{lem:B-descent}). Dividing out by the left and
right radicals yields the perfect pairing $B_{\cW_n}$ on
$\cW_n \otimes \rho^*\cW_n$.
\end{proof}

%==================================================================
\subsection{Reflexive Symmetry of the Pairing}
\label{sec:symmetry}
%==================================================================

\begin{lemma}[Reflexive symmetry]
\label{lem:sym}
The pairing $B_{\cH_n} : \cH_n \otimes \rho^*\cH_n \to \Ql(-(n+1))$
satisfies the Frobenius-equivariant sheaf identity
\[
\boxed{\;
\rho^* B_{\cH_n} \circ \mathrm{swap}_{\cH_n, \rho^*\cH_n}
\;=\; (-1)^{n+1} \, B_{\cH_n},
\;}
\]
i.e.\ $B_{-h}(w, v) = (-1)^{n+1} B_h(v, w)$ on every geometric fibre
$h \in S_n(\overline{\Fp})$. Assuming Theorem~\ref{thm:H7}, the
induced pairing $B_{\cW_n}$ on $\cW_n \otimes \rho^*\cW_n$ is likewise
$(-1)^{n+1}$-symmetric.
\end{lemma}

\begin{proof}[Proof summary (complete proof in Section~\ref{app:reflection})]
The proof combines three ingredients over $\Fp$ (hence over any $\Fq$
by base change), and does not rely on Theorem~\ref{thm:H7}:
\begin{enumerate}
\item \emph{Coefficient symmetry \textup{(}Lemma~\ref{lem:B-coeff}\textup{)}:}
The pairing $\Phi : \cM \otimes \rhotil^*\cM \to \Ql(-n)$ built from the
symmetric Artin--Schreier pairing ($+1$) and the $n$ alternating
Kloosterman pairings ($-1$ each) satisfies
$\rhotil^*\Phi \circ \mathrm{swap}_{\cM, \rhotil^*\cM} = (-1)^n \Phi$.
\item \emph{Graded commutativity in degree $1$:} Interchanging two
cohomology classes $v \in H^1_c(\Gm, \cM_h)$ and
$w \in H^1_c(\Gm, \cM_{-h})$ under the cup product introduces the
Koszul sign $(-1)^{1 \cdot 1} = -1$.
\item \emph{Invariance of trace under $\sg(t) = -t$:} Since $\sg$ is an
automorphism of $\Gm$ over $\Fp$, the compactly supported trace satisfies
$\Tr_\pi \circ \sg^* = \Tr_\pi$ Frobenius-equivariantly.
\end{enumerate}
Multiplying the two signs gives $(-1)^n \cdot (-1) = (-1)^{n+1}$.
The full Frobenius-equivariant derivation is given in
Section~\ref{app:reflection}.
\end{proof}

%==================================================================
\subsection{The Root-Number Theorem}
\label{sec:root-number}
%==================================================================

Throughout this section we work over $\Fq$ ($q = p^a$) and assume:
\begin{itemize}
\item[\textup{(D1)}] Theorem~\ref{thm:H7} (relative construction in
Section~\ref{app:H7}, together with the Fourier--Deligne comparison of
the main paper identifying $\cW_n$ with the exponential-sum sheaf), and
\item[\textup{(D2)}] Lemma~\ref{lem:sym} and Proposition~\ref{prop:radical}
(reflexive symmetry and descent to $\cW_n$, proved in Section~\ref{app:reflection}).
\end{itemize}
No circularity arises between \S\S\ref{sec:rel-duality}--\ref{sec:root-number}
and Sections~\ref{app:H7}--\ref{app:reflection}.

\subsubsection{Transfer of the pairing to the middle extension}

Set $P := \cW_n[n] \in \mathrm{Perv}(S_n, \Ql)$ and $w := 2n + 1$.
Shifting both factors of
$B_{\cW_n} : \cW_n \otimes \rho^*\cW_n \to \Ql(-(n+1))$ by $[n]$ yields
a pairing in $D^b_c(S_n, \Ql)$:
\[
P \otimes^L \rho^* P \longrightarrow \Ql(-(n+1))[2n]
\;=\; \omega_{S_n}(-w),
\tag{6.1}
\]
whose parity, by the Koszul sign rule for interchanging two factors in
cohomological degree $n$ (Kiehl--Weissauer, \S\S II.6, III.7), is
\[
(-1)^{n+1} \cdot (-1)^{n \cdot n}
\;=\; (-1)^{n+1} \cdot (-1)^n \;=\; -1 .
\]

\begin{lemma}[Tensor--Hom adjunction]
\label{lem:tensor-hom}
Giving the pairing \textup{(6.1)} in $D^b_c(S_n, \Ql)$ is equivalent to
giving an adjoint morphism
\[
\varphi : \rho^* P \longrightarrow \mathbb{D}_{S_n}(P)(-w)
\tag{6.2}
\]
via the canonical isomorphism
\[
\begin{aligned}
&\Hom_{D^b_c(S_n)}\bigl(P \otimes^L \rho^* P,\, \omega_{S_n}(-w)\bigr) \\
&\quad\simeq \Hom_{D^b_c(S_n)}\bigl(\rho^* P,\, \RHom(P, \omega_{S_n})(-w)\bigr)
\;=\; \Hom_{D^b_c(S_n)}\bigl(\rho^* P,\, \mathbb{D}_{S_n}(P)(-w)\bigr).
\end{aligned}
\]
Non-degeneracy of \textup{(6.1)} is equivalent to $\varphi$ being an
isomorphism.
\end{lemma}

\begin{proof}
Standard tensor--Hom adjunction in $D^b_c(S_n, \Ql)$ together with the
definition $\mathbb{D}_{S_n} = \RHom(-, \omega_{S_n})$.
\end{proof}

\begin{lemma}[Perverse duality for $K_n$]
\label{lem:perv-duality}
Let $\bar\rho : X_n \to X_n$ act as $x \mapsto -x$ on each $\PP^1$
factor; it is an involution over $\Fq$, fixes $\{0,\infty\}$ in each
factor, and restricts to $\rho$ on $S_n$. Then \textup{(6.2)} extends
canonically to a Frobenius-equivariant isomorphism in $D^b_c(X_n, \Ql)$:
\[
\boxed{\;
\bar\rho^* K_n \;\xrightarrow{\ \sim\ }\; \mathbb{D}_{X_n}(K_n)(-w),
\;}
\tag{6.3}
\]
equivalently a non-degenerate pairing
\[
K_n \otimes^L \bar\rho^* K_n \longrightarrow
\omega_{X_n}(-w) \;=\; \Ql(-(n+1))[2n],
\tag{6.4}
\]
of parity $-1$ \textup{(}alternating\textup{)}.
\end{lemma}

\begin{proof}
By Lemma~\ref{lem:tensor-hom} and Proposition~\ref{prop:radical}(iii),
\textup{(6.1)} is encoded by the isomorphism
$\varphi : \rho^* P \xrightarrow{\sim} \mathbb{D}_{S_n}(P)(-w)$.
Applying the intermediate extension functor $\jex$ to $\varphi$ and
using its two standard compatibilities
\[
\jex \circ \mathbb{D}_{S_n} \;=\; \mathbb{D}_{X_n} \circ \jex
\quad\text{(BBD, Prop.~2.1.11)},
\qquad
\jex \circ \rho^* \;=\; \bar\rho^* \circ \jex
\]
(the latter since $\bar\rho$ is an automorphism of the pair $(X_n, S_n)$
restricting to $\rho$ on $S_n$) produces the isomorphism \textup{(6.3)}.

Converting \textup{(6.3)} back via tensor--Hom adjunction on $X_n$
yields the pairing \textup{(6.4)}, using $\omega_{X_n} = \Ql(n)[2n]$ so
that $\omega_{X_n}(-(2n+1)) = \Ql(-(n+1))[2n]$. The pairing is perfect
because its adjoint \textup{(6.3)} is an isomorphism; we never apply
$\jex$ to a tensor pairing, only to the morphism $\varphi$. Its parity
is preserved from \textup{(6.1)}, hence equals $-1$.
\end{proof}

\begin{corollary}[Symplectic structure in middle degree]
\label{cor:symplectic}
Let $V_m := IH^m(X_n, \cW_n) = \mathbb{H}^{m-n}(X_n, K_n)$ and
$d_m := \dim V_m$. Then:
\begin{enumerate}
\item[\textup{(i)}] For every $m \in \{0, \dots, 2n\}$, there is a
Frobenius-equivariant perfect pairing
\[
V_m \otimes V_{2n-m} \longrightarrow \Ql(-(2n+1)),
\qquad\text{i.e.}\qquad
V_{2n-m} \;\simeq\; V_m^\vee(-(2n+1)).
\tag{6.5}
\]
\item[\textup{(ii)}] For $m = n$, this self-pairing is alternating and
non-degenerate; hence $d_n = 2g$ is even,
\[
F_q \in \GSp(V_n) \text{ with multiplier } q^{2n+1},
\qquad
\det(F_q \mid V_n) \;=\; q^{(2n+1)g}.
\tag{6.6}
\]
\end{enumerate}
\end{corollary}

\begin{proof}
(i) Since $X_n$ is proper over $\Fq$, global Verdier duality gives
$\mathbb{H}^i(X_n, \mathbb{D}_{X_n} K_n) \simeq \mathbb{H}^{-i}(X_n, K_n)^\vee$.
Combining with \textup{(6.3)} and the Frobenius-equivariant isomorphism
$\mathbb{H}^i(X_n, \bar\rho^* K_n) \simeq \mathbb{H}^i(X_n, K_n)$
(as $\bar\rho$ is defined over $\Fq$) yields
$\mathbb{H}^i(X_n, K_n) \simeq \mathbb{H}^{-i}(X_n, K_n)^\vee(-(2n+1))$.
Setting $i = m - n$ gives \textup{(6.5)}. (Concretely, applying
$\mathbb{H}^{2n}(X_n, -)$ to \textup{(6.4)} and composing with the
global trace $\mathbb{H}^{2n}(X_n, \Ql(-(n+1))) \simeq \Ql(-(2n+1))$
produces the twist $-(2n+1)$ after the global trace.)

(ii) For $m = n$, \textup{(6.5)} is a non-degenerate alternating form
on $V_n$, forcing $d_n = 2g$ even and $F_q \in \GSp_{2g}$ with
multiplier $q^{2n+1}$, hence $\det(F_q \mid V_n) = (q^{2n+1})^g$.
\end{proof}

\subsubsection{The functional equation and the root number}

With $P_{n,m}(T) := \det(1 - F_q T \mid V_m)$, $e_m := (-1)^{m-n+1}$,
$w := 2n + 1$, and $\chi := \sum_{m=0}^{2n} e_m d_m$, the $L$-function is
\[
L(X_n, K_n, T) \;=\; \prod_{m=0}^{2n} P_{n,m}(T)^{e_m}.
\]
Note that $e_{2n-m} = e_m$ for all $m$.

\begin{lemma}[Functional equation]
\label{lem:funct-eq}
With $T' := \dfrac{1}{q^w T}$, we have
\[
L(X_n, K_n, T') \;=\; \omega_n \cdot \bigl(q^{w/2} T\bigr)^{-\chi}\,
L(X_n, K_n, T),
\tag{6.7}
\]
where the global root number is given by
\[
\omega_n \;=\; q^{w\chi/2} \prod_{m=0}^{2n}
\Bigl[(-1)^{d_m} q^{-w d_m} \det(F_q \mid V_m)\Bigr]^{e_m}
\;=\; q^{-w\chi/2} \prod_{m=0}^{2n} \det(-F_q \mid V_m)^{e_m}.
\tag{6.8}
\]
\end{lemma}

\begin{proof}
By \textup{(6.5)}, the eigenvalues of $F_q$ on $V_{2n-m}$ are
$\{q^w / \alpha\}$ where $\{\alpha\}$ are those on $V_m$. Hence
\[
P_{n,m}(T') \;=\; (-1)^{d_m} q^{-w d_m} T^{-d_m} \det(F_q \mid V_m)\,
P_{n,2n-m}(T).
\]
Raising to the power $e_m$, multiplying over $0 \le m \le 2n$, and
using $e_{2n-m} = e_m$ yields \textup{(6.7)} and \textup{(6.8)}.
\end{proof}

\begin{theorem}[Root-number formula; conditional on \textup{(D1)} and \textup{(D2)}]
\label{thm:root-number}
Assume \textup{(D1)} and \textup{(D2)}. Then
\[
\boxed{\;\omega_n \;=\; +1\;}
\]
and, writing $\dim IH^n(X_n, \cW_n) = 2g$,
\[
\det\bigl(F_q \mid IH^n(X_n, \cW_n)\bigr) \;=\; q^{(2n+1)g}.
\]
\end{theorem}

\begin{proof}
For each complementary pair $\{m, 2n-m\}$ with $m \ne n$, we have
$d_{2n-m} = d_m$, $e_{2n-m} = e_m$, and
$\det(F_q \mid V_{2n-m}) = q^{w d_m} / \det(F_q \mid V_m)$ by
\textup{(6.5)}. Thus the product of the two brackets in \textup{(6.8)}
is $(-1)^{2 d_m} q^{-2w d_m} q^{w d_m} = q^{-w d_m}$, which raised to
$e_m$ cancels the prefactor contribution $q^{w (2 e_m d_m)/2} = q^{w e_m d_m}$.

Hence \textup{(6.8)} reduces to the middle degree $m = n$ ($e_n = -1$):
\[
\omega_n \;=\; \frac{q^{w d_n / 2}}{\det(-F_q \mid V_n)}
\;=\; \frac{q^{w d_n / 2}}{(-1)^{d_n} \det(F_q \mid V_n)}.
\]
By Corollary~\ref{cor:symplectic}(ii), $d_n = 2g$ is even and
$\det(F_q \mid V_n) = q^{w g}$, so $\omega_n = (-1)^{2g} q^{wg} / q^{wg} = +1$.
\end{proof}

\begin{remark}[Independence from semisimplicity and dimension formulas]
\label{rmk:compact-omega}
Theorem~\ref{thm:root-number} uses neither the semisimplicity of $F_q$
(\S\ref{sec:semisimplicity}) nor the explicit value of $g$
(\S\ref{sec:spectral}). Moreover, $\omega_n$ is defined directly from
global Frobenius determinants and does not depend on choosing a system
of local Deligne $\varepsilon$-factors.
\end{remark}

%==================================================================
\subsection{Spectral Reconstruction of $IH^\bullet(X_n,\cW_n)$}
\label{sec:spectral}
%==================================================================

No appeal is made in this section to Lemma~\ref{lem:sym} or to
Theorem~\ref{thm:root-number}.

\subsubsection{Koszul supertrace of the lower-weight inertia model}

Let $\mathsf{E}_n$ denote the finite-dimensional arithmetic Boolean
representation over $\Ql$ (of dimension $\binom{n}{\lfloor n/2\rfloor}$)
constructed in the main paper to model the lower-weight inertia
subsheaf $\cL_n = W_n\cH_n$ at $t = 0$. For $Q = q^r$ ($r \ge 1$), its
Frobenius trace is the Catalan--Boolean polynomial
\[
\Tr(F_Q \mid \mathsf{E}_n) \;=\; T_n(Q)
\;:=\; \sum_{j=0}^{\lfloor n/2\rfloor}
\Bigl(\binom{n}{j} - \binom{n}{j-1}\Bigr) Q^j,
\qquad \binom{n}{-1} := 0.
\]
For each subset $I \subseteq \{1,\dots,n\}$ of cardinality $|I| = k$,
let $X_{n,I,0} \subseteq X_n$ be the locally closed $0$-stratum where
$h_i = 0$ for $i \in I$ and $h_j \in \Gm$ for $j \notin I$.

\begin{proposition}[Koszul supertrace of the lower-weight inertia model]
\label{prop:koszul}
Let $x \in X_{n,I,0}(\Fq)$ with $|I| = k$, and let $Q = q^r$ ($r \ge 1$).
On any stratum where the local Galois--Koszul complex modeling the
lower-weight inertia contribution in a formal neighbourhood of $x$ is
identified arithmetically with
\[
K^\bullet_{\mathrm{Gal}}(\mathsf{E}_n; x_I)
\;=\; \bigoplus_{a=0}^{k} K^a,
\qquad
K^a \;=\; \bigoplus_{\substack{J \subseteq I \\ |J|=a}} \mathsf{E}_n(-a),
\]
the Frobenius supertrace is
\[
\boxed{\;
\STr\bigl(F_Q \mid K^\bullet_{\mathrm{Gal}}(\mathsf{E}_n; x_I)\bigr)
\;=\; (1-Q)^k\, T_n(Q).
\;}
\tag{7.1}
\]
\end{proposition}

\begin{proof}
The Tate twist $\mathsf{E}_n(-a)$ multiplies each Frobenius eigenvalue
on $\mathsf{E}_n$ by $Q^a$. Since $\binom{k}{a}$ subsets $J \subseteq I$
have cardinality $a$,
\[
\STr\bigl(F_Q \mid K^\bullet_{\mathrm{Gal}}(\mathsf{E}_n; x_I)\bigr)
\;=\; \sum_{a=0}^{k} (-1)^a \binom{k}{a} Q^a \Tr(F_Q \mid \mathsf{E}_n)
\;=\; (1-Q)^k\, T_n(Q).
\]
\end{proof}

\begin{remark}[Scope]
Formula \textup{(7.1)} governs the local Galois--Koszul complex of the
\emph{lower-weight inertia model} $\mathsf{E}_n$ along $0$-strata
$X_{n,I,0}$ where the arithmetic Boolean identification holds; it is
\emph{not} the Frobenius trace on the stalk of $K_n = \jex(\cW_n[n])$,
which varies with $h$ on $S_n$.
\end{remark}

\subsubsection{Global trace formula and spectral reconstruction}

\begin{proposition}[Grothendieck--Lefschetz for $K_n$]
\label{prop:groth-lef}
For every $r \ge 1$,
\[
s_r \;:=\; \sum_{x \in X_n(\Fq^r)} \STr(F_{q^r} \mid K_{n,\bar x})
\;=\; \sum_{m=0}^{2n} (-1)^{m-n} \Tr(F_{q^r} \mid V_m),
\tag{7.2}
\]
and hence
\[
L(X_n, K_n, T)
\;=\; \exp\!\Bigl(\sum_{r \ge 1} \frac{s_r}{r} T^r\Bigr)
\;=\; \prod_{m=0}^{2n} P_{n,m}(T)^{(-1)^{m-n+1}}.
\tag{7.3}
\]
\end{proposition}

\begin{proof}
Grothendieck--Lefschetz trace formula applied to $K_n$ on the proper
variety $X_n$, using $V_m = \mathbb{H}^{m-n}(X_n, K_n)$.
\end{proof}

\begin{theorem}[Spectral reconstruction under purity]
\label{thm:spectral}
Assume \textup{(D1)} Theorem~\ref{thm:H7} (so that Lemma~\ref{lem:purity-K}
applies and $V_m$ is pure of weight $m+n+1$). The zeros and poles of the
rational function
\[
L(X_n, K_n, T) \;=\; \prod_{m=0}^{2n} P_{n,m}(T)^{(-1)^{m-n+1}}
\]
are separated by the distinct complex absolute values
$q^{-(m+n+1)/2}$ of their locations. The exponent $(-1)^{m-n+1}$
determines whether $P_{n,m}(T)$ contributes zeros or poles.
Consequently, each characteristic polynomial
$P_{n,m}(T) = \det(1 - F_q T \mid V_m)$ for $0 \le m \le 2n$ is
uniquely recoverable from $L(X_n, K_n, T)$.
\end{theorem}

\begin{proof}
By Lemma~\ref{lem:purity-K}, every root $T = \alpha^{-1}$ of
$P_{n,m}(T)$ satisfies $|\alpha|_\infty = q^{(m+n+1)/2}$, hence
$|T|_\infty = q^{-(m+n+1)/2}$. As $m$ ranges over $\{0, 1, \dots, 2n\}$,
the weights $m + n + 1$ are pairwise distinct, so the circles
$|T|_\infty = q^{-(m+n+1)/2}$ are pairwise disjoint. On the circle
$|T|_\infty = q^{-(m+n+1)/2}$, the divisor of $L(X_n, K_n, T)$ equals
$(-1)^{m-n+1} \operatorname{div}(P_{n,m}(T))$: thus $P_{n,m}(T)$
contributes zeros when $(-1)^{m-n+1} = +1$ and poles when
$(-1)^{m-n+1} = -1$. Since $P_{n,m}(0) = 1$, each $P_{n,m}(T)$ is
uniquely determined.
\end{proof}

\begin{corollary}[Dimensions from $L$]
\label{cor:dim-from-L}
Under \textup{(D1)}, $\dim V_m = \deg P_{n,m}(T)$ is uniquely
determined by $L(X_n, K_n, T)$.
\end{corollary}

\begin{remark}[The special case of five variables]
The later calculation in Section~\ref{sec:five-computation} determines the
five-variable coefficients through actual arithmetic comparison maps.
The question in this subsection concerns the corresponding unrestricted
higher-dimensional problem and is not settled in that generality by the
five-variable analysis.
\end{remark}

\subsubsection{Open: boundary stalks of $K_n$}

\begin{remark}[Open problem]
\label{rmk:open-boundary}
Determining the explicit Frobenius traces $\STr(F_{q^r} \mid K_{n,\bar x})$
on the stalks of $K_n = \jex(\cW_n[n])$ along the boundary
$X_n \setminus S_n$---in particular along strata involving $\infty$-collisions
or multiple $0$-collisions---requires a local BBD-type decomposition of
$Rj_*(\cW_n[n])$ in formal neighbourhoods of boundary strata and remains
an open problem for general $n \ge 4$.
\end{remark}

%==================================================================
\subsection{The Boundary Perverse Sheaf $B'_n$}
\label{sec:bprime}
%==================================================================

Let $\mathrm{pr} : X_n \to X_{n-1}$ be the proper projection
\emph{forgetting the first coordinate}:
\[
\mathrm{pr}(h_1, h_2, \dots, h_n) \;=\; (h_2, \dots, h_n).
\]
Thus $X_{n-1} = (\PP^1)^{n-1}$ has coordinates $(h_2, \dots, h_n)$ and
codimension-$1$ boundary divisors
\[
Z_i^a \;:=\; \{h_i = a\} \;\subseteq\; X_{n-1},
\qquad i \in \{2, \dots, n\},\;\; a \in \{0, \infty\},
\]
with closed immersions $\iota_{Z_i^a} : Z_i^a \hookrightarrow X_{n-1}$.

\subsubsection{BBD decomposition and definition of $B'_n$}

Since $\mathrm{pr}$ is proper and $K_n$ is pure of weight $2n+1$
(Lemma~\ref{lem:purity-K}, assuming \textup{(D1)}), the BBD
decomposition theorem decomposes $R\mathrm{pr}_* K_n$ into
\emph{shifts} of pure perverse sheaves. The specialized, unshifted
form displayed below additionally assumes the perverse-concentration
condition ${}^p\mathcal H^j(R\mathrm{pr}_*K_n)=0$ for $j\neq0$
and the strict-support classification of the main paper.
Neither condition follows from BBD alone. Subject to these inputs,
collecting summands by support, we write
\[
R\mathrm{pr}_* K_n
\;\simeq\;
V'_n[n-1] \;\oplus\;
\bigoplus_{i=2}^{n} \bigoplus_{a \in \{0,\infty\}}
(\iota_{Z_i^a})_* C_{i,a}[n-2]
\;\oplus\;
B'_n ,
\tag{8.1}
\]
where $V'_n[n-1]$ collects the full-support summands on $X_{n-1}$,
$(\iota_{Z_i^a})_* C_{i,a}[n-2]$ collects the summands supported strictly
on the divisor $Z_i^a$, and $B'_n$ is defined as the direct sum of all
remaining perverse summands of $R\mathrm{pr}_* K_n$ whose irreducible
supports have codimension $\ge 2$ in $X_{n-1}$.

\subsubsection{Support bound and vanishing in low rank}

The geometric sections contain the following results:
\texttt{lem:G-base} (divisorial fibres),
\texttt{lem:G-perverse} (perverse concentration),
\texttt{thm:G-decomp} (geometric strict-support decomposition),
\texttt{lem:G-codim2} (codimension-two collisions), and
\texttt{cor:G-n3} (the low-dimensional conclusion).
These labels refer to geometric results stated earlier in the present
source, not to results in \cite{EbaCoh}. The arguments below use these results with their stated hypotheses. In particular,
\texttt{thm:G-decomp} asserts an isomorphism after
\emph{geometric} base change; it must not be read as an
arithmetic direct-sum splitting over $\mathbb F_p$.

\begin{proposition}[Conditional support bound]
\label{prop:bprime-supp}
Assume \textup{(D1)}, and suppose the conclusions of
\texttt{lem:G-base}, \texttt{lem:G-perverse}, and
\texttt{thm:G-decomp} of the geometric sections of this manuscript hold, including
perverse concentration and the geometric strict-support
classification. Then
\[
\boxed{\;\dim\supp(B'_n)\le n-3.\;}
\]
\end{proposition}

\begin{proof}
By the asserted geometric strict-support decomposition,
$B'_n$ is the direct sum of constituents whose irreducible
supports have codimension at least two in $X_{n-1}$.
The assertion that the full-support and divisor-supported
constituents have all been accounted for is an explicit
explicit hypothesis on the main-text geometric analysis, not a consequence of the BBD decomposition
theorem alone. Since $\dim X_{n-1}=n-1$, the claimed bound follows.
\end{proof}

\begin{proposition}[Conditional vanishing for $n=2,3$]
\label{prop:bprime-small}
Under the hypotheses of Proposition~\ref{prop:bprime-supp},
assume also the codimension-two vanishing assertion
\texttt{lem:G-codim2} and its application
\texttt{cor:G-n3} in the geometric sections of this manuscript. Then
\[
\boxed{\;B'_2=0,\qquad B'_3=0.\;}
\]
The decomposition asserted here is geometric; the proposition
does not claim a canonical arithmetic splitting.
\end{proposition}

\begin{proof}
For $n=2$, Proposition~\ref{prop:bprime-supp} gives
$\dim\supp(B'_2)\le-1$, hence $B'_2=0$.
For $n=3$, the support bound makes $B'_3$ a sum of
point-supported perverse constituents after geometric base
change. For a boundary point
$x_0=(a_2,a_3)\in\{0,\infty\}^2\subset X_2$,
put $L_{x_0}=\mathrm{pr}^{-1}(x_0)\simeq\PP^1$.
Proper base change and the ordinary hypercohomology
spectral sequence give
\[
E_2^{a,b}=H^a\bigl(L_{x_0},\mathcal H^b(K_3)|_{L_{x_0}}\bigr)
\Longrightarrow H^{a+b}(L_{x_0},K_3|_{L_{x_0}}).
\]
The only possible terms with $a+b=0$ in the relevant
range are $E_2^{1,-1}$ and $E_2^{2,-2}$.
Their complete accounting by full-support and divisor
constituents, as well as vanishing at the codimension-two
boundary intersections, is exactly the additional input
from \texttt{lem:G-base}, \texttt{lem:G-codim2} and
\texttt{cor:G-n3}. Under these cited hypotheses, no
independent point-supported constituent remains. Therefore
$B'_3=0$ geometrically.
\end{proof}

\begin{remark}[Scope of companion references]
The labels \texttt{lem:G-base}, \texttt{lem:G-perverse},
\texttt{thm:G-decomp}, \texttt{lem:G-codim2}, and
\texttt{cor:G-n3} are internal to the revised companion
manuscript the main-text geometric section; they are not statements of
\cite{EbaCoh}. The above conditional deductions do not
replace an independent verification of those arguments.
\end{remark}

\subsubsection{The case $n=4$ and arithmetic decomposition by support}

Under the support-bound hypotheses above, for $n = 4$ we have
$\dim\supp(B'_4) \le 1$, so $B'_4$ is supported on
at most finitely many curves and points in $X_3$.

\begin{remark}[Open for $n \ge 4$]
\label{rmk:open-bprime}
Neither the existence nor the vanishing of $1$-dimensional components
of $B'_4$, nor the explicit structure of $B'_n$ for $n \ge 4$, is
established by the methods of this paper. Determining $B'_n$ for
$n \ge 4$ requires an explicit analysis of higher-codimension strata
where multiple coordinates $h_i$ ($i \ge 2$) simultaneously collide
at $\{0,\infty\}$, which lies beyond Theorem~\ref{thm:H7} and the
codimension-$2$ collision analysis.
\end{remark}

\begin{remark}[Arithmetic decomposition by support]
\label{rmk:arith-decomp}
For the proper morphism $\mathrm{pr} : X_n \to X_{n-1}$ and the pure
perverse sheaf $K_n$, the finite-field decomposition theorem
\textup{(}BBD~1982, Thm.~5.3.8; see also
de~Cataldo--Haines--Li, Prop.~5.1.1\textup{)} yields a noncanonical
decomposition in the arithmetic mixed derived category
$D^b_c(X_{n-1,\Fq}, \Ql)$ into shifted intersection complexes with
pure indecomposable coefficient systems.

The terms may be collected according to their irreducible supports.
This does not imply that the Frobenius endomorphism on their
cohomology is semisimple \textup{(}\S\ref{sec:semisimplicity}\textup{)},
nor does it provide canonical projectors. The geometric support bounds
and vanishing statements of Propositions~\ref{prop:bprime-supp}
and~\ref{prop:bprime-small} are preserved under this arithmetic
decomposition, subject to the verified strict-support classification.
\end{remark}

%==================================================================
\subsection{Arithmetic Semisimplicity}
\label{sec:semisimplicity}
%==================================================================

\subsubsection{Algebraic criterion}

\begin{proposition}[Jordan criterion; unconditional]
\label{prop:jordan}
Let $V$ be a finite-dimensional $\Ql$-vector space with an
endomorphism $F$. Then $F$ is semisimple if and only if its minimal
polynomial $\mu_F(T)$ has no repeated roots in $\overline{\Ql}$.
\end{proposition}

\begin{proof}
Standard linear algebra over the algebraically closed field
$\overline{\Ql}$: the minimal polynomial of a semisimple endomorphism
is squarefree, and conversely a squarefree minimal polynomial yields
$V = \bigoplus_\alpha \ker(F - \alpha\,\id)$.
\end{proof}

\subsubsection{Motivic sufficient condition}

\begin{proposition}[Motivic semisimplicity; unconditional, Tate 1966]
\label{prop:motivic-ss}
Let $A/\Fq$ be an abelian variety. Then $F_q$ acts semisimply on
$H^i_{\mathrm{\acute et}}(A_{\overline{\Fq}}, \Ql)$ for every $i \ge 0$.
If $V$ is a $G_{\Fq}$-equivariant direct summand of
$H^i_{\mathrm{\acute et}}(A_{\overline{\Fq}}, \Ql)(r)$ for some abelian
variety $A/\Fq$, degree $i \ge 0$, and Tate twist $r \in \mathbb{Z}$,
then $F_q$ acts semisimply on $V$.
\end{proposition}

\begin{proof}
Semisimplicity of $F_q$ on
$H^1_{\mathrm{\acute et}}(A_{\overline{\Fq}}, \Ql) = T_\ell(A)^\vee \otimes \Ql$
is Tate's theorem (Tate, \emph{Invent.\ Math.}~2 (1966), \S3, Thm.~2).
Semisimplicity on $H^i \simeq \wedge^i H^1$ follows by taking exterior
powers, Tate twists multiply Frobenius by $q^{-r}$ and preserve
semisimplicity, and semisimplicity passes to $G_{\Fq}$-equivariant
direct summands.
\end{proof}

\begin{corollary}[Sufficient criterion with weight constraint]
\label{cor:ss-sufficient}
Assume \textup{(D1)} Theorem~\ref{thm:H7}, so that
$V_m = IH^m(X_n, \cW_n)$ is pure of weight $m + n + 1$. If there
exist an abelian variety $A/\Fq$ and a $G_{\Fq}$-equivariant
direct-summand embedding
\[
V_m \;\hookrightarrow\; H^i_{\mathrm{\acute et}}(A_{\overline{\Fq}}, \Ql)
\]
\textup{(}without Tate twist\textup{)}, then necessarily
\[
\boxed{\;i \;=\; m + n + 1,\;}
\]
and $F_q$ acts semisimply on $V_m$. More generally, for a twisted
direct-summand embedding into
$H^i_{\mathrm{\acute et}}(A_{\overline{\Fq}}, \Ql)(r)$, purity forces
$i - 2r = m + n + 1$.
\end{corollary}

\begin{proof}
Since $H^i_{\mathrm{\acute et}}(A_{\overline{\Fq}}, \Ql)$ is pure of
weight $i$, any non-zero $G_{\Fq}$-equivariant direct summand $V_m$
must have the same weight, giving $i = m + n + 1$ (and $i - 2r = m + n + 1$
in the twisted case). Semisimplicity follows from
Proposition~\ref{prop:motivic-ss}.
\end{proof}

\begin{remark}[Arithmetic semisimplicity in the special five-variable case]
Under the explicitly stated comparison hypotheses of
Theorem~\ref{thm:five-computation}, the arithmetic Tate decomposition
of Section~\ref{sec:five-computation} gives semisimple Frobenius on the
intersection cohomology for $n=5$.  This does not establish the analogous
claim for arbitrary $n$.
\end{remark}

\subsubsection{The open problem for the family}

\begin{remark}[Open problem]
\label{rmk:open-ss}
The semisimplicity of $F_q$ on $IH^m(X_n, \cW_n)$ for general $n \ge 2$
and $0 \le m \le 2n$ is an open problem specific to this family. It is
\emph{not} equivalent to the general Grothendieck--Serre semisimplicity
conjecture: the explicit Kloosterman construction of $\cW_n$ may admit
a family-specific proof, whereas general semisimplicity conjectures
remain wide open. Note also that membership of $F_q$ in a reductive
algebraic group (such as the geometric monodromy group
$G_{\mathrm{geom}}(\cW_n)$ or $\GSp_{2g}$) does \emph{not} by itself
imply semisimplicity, since reductive groups contain non-trivial
unipotent elements. Establishing semisimplicity therefore requires
either directly ruling out non-trivial unipotent Jordan blocks of $F_q$
on $IH^m(X_n, \cW_n)$ or constructing an explicit motivic direct
summand in degree $m + n + 1$ \textup{(}Corollary~\ref{cor:ss-sufficient}\textup{)}.
\end{remark}

%==================================================================
% ===================================================================
% New material integrated with the original manuscript, November 2026
% ===================================================================
\section{The five-variable intersection complex: arithmetic boundary calculation}\label{sec:five-variable}
\label{sec:five-computation}

The preceding sections establish the construction, the local weights, and the
formal properties of the intersection-cohomological completion in arbitrary
dimension.  The first case in which the compactification has several
intersecting boundary strata of different dimensions and the resulting
intersection cohomology cannot be read from one-variable invariants alone is
$n=5$.  We now give a separate calculation in that dimension.  Its point is
not merely the numerical value of the Betti numbers.  We retain the
arithmetic coefficient systems, the maps between them, and their extension
classes throughout.  In particular, purity is never used as a substitute
for semisimplicity of arithmetic Frobenius.

The assertions in this section extend the extremal-degree information
obtained from Theorem~\ref{thm:G-decomp}. The general decomposition does not
determine the middle intersection groups for $n=5$: it provides the two
extremal multiplicities, whereas the relative comparisons and the
coordinate-stratum calculation determine the remaining degrees. The
proof is written for $n=5$ alone. In dimensions $n\ge6$ the analogous
residual support and arithmetic Frobenius modules are not classified
by this calculation.

The proof below is organized around three comparison statements, referred
to as (P), (U), and (K).  Their precise functorial form matters.  An equality
of characteristic polynomials cannot replace (U), and agreement of the
ranks of two Koszul complexes cannot replace (K).  We spell out the maps and
their uses, including the boundary-to-global comparison.  The geometric
arguments involving nearby cycles and intermediate extensions invoke the
standard results in \cite{BBD,BeilGlue,MorelGlue,ReichGlue,SGA1}; the necessary
normalizations are displayed explicitly.

\subsection{Statement, notation, and scope}\label{subsec:five-statement}
Put $S_m=(\mathbb G_m)^m$, $X_m=(\mathbb P^1)^m$, and let
$j_m:S_m\hookrightarrow X_m$ be the standard open immersion.  All sheaves
are defined over $k=\mathbb F_q$, with $p=\operatorname{char}k>2$ and
$\ell\ne p$.  The constant coefficient field $E$ is a sufficiently large
finite extension of $\mathbb Q_\ell$; by extending scalars we may work in
$\overline{\mathbb Q}_\ell$.  Geometric Frobenius acts on $E(-a)$ by $q^a$.
We use the construction
\[
 \mathcal M_m=\mathcal L_\psi(\lambda t)\otimes
          \bigotimes_{i=1}^m[th_i]^*\mathrm{Kl}_2,
 \qquad \pi_t:\mathbb G_{m,t}\times S_m\longrightarrow S_m,
\]
with $\lambda\in k^\times$.  Write
\[
 \mathcal H_m=R^1\pi_{t!}\mathcal M_m,\qquad
 \alpha_m:R^1\pi_{t!}\mathcal M_m\longrightarrow
                R^1\pi_{t*}\mathcal M_m,\qquad
 \mathcal L_m=\ker\alpha_m,\qquad
 \mathcal W_m=\operatorname{im}\alpha_m.
\]
The morphism $\alpha_m$ is the natural morphism forgetting supports. The notation agrees with the family construction
of Theorem~\ref{thm:B-family}: the comparison proved there identifies
$\mathcal L_m=W_m\mathcal H_m$ and
$\mathcal W_m\cong\mathcal H_m/\mathcal L_m$; these are not new sheaves
with independent definitions. The auxiliary $\mathcal H'_m$ is instead
the cokernel of the distinguished constant line in $\mathcal H_m$.
The general support bound determines $b_{5,2}=b_{5,8}=20$; the middle
values $b_{5,4},b_{5,6}$ and the vanishing of $IH^5$ are the additional
outputs of the comparisons in this section, not of the general argument.

The central sheaf for the present calculation is
\[
 K_5=j_{5,!*}(\mathcal W_5[5]),\qquad I_5=K_5[-5],
 \qquad f:X_5\longrightarrow X_4
\]
where $f$ forgets the first coordinate, and $P_5=Rf_*K_5$.
Set $F^v=\mathcal H^v(I_5)$.  Every occurrence of $\operatorname{im}N_i$
below denotes the actual image of a morphism of arithmetic coefficient
systems, with the stated Tate normalization.

\begin{theorem}[The five-variable boundary calculation]\label{thm:five-computation}
Assume the sheaf-theoretic comparison statements (P), (U), and (K) established
in the proofs below, including the natural base-change statements at the
coordinate zero divisors.  Then $P_5$ is a pure perverse sheaf and there is
an isomorphism in the arithmetic constructible derived category
\begin{align}
 Rf_*K_5\simeq{}&\bigoplus_{i=2}^{5}
 \left((\iota_i^0)_*E(-4)^2[3]\oplus
 (\iota_i^\infty)_*E(-4)^3[3]\right)\notag\\
 &\oplus\bigoplus_{2\le i<j<k\le5}
 \ \bigoplus_{a,b,c\in\{0,\infty\}}
 (\iota_{ijk}^{abc})_*E(-5)^{r_t}[1],
 \label{eq:five-support}
\end{align}
where $t$ is the number of infinite coordinates among $a,b,c$ and
$(r_0,r_1,r_2,r_3)=(1,1,2,4)$.  The isomorphism concerns the arithmetic
objects themselves and is not merely a decomposition after semisimplifying
Frobenius.  Consequently,
\begin{equation}\label{eq:five-IH}
 \begin{aligned}
 IH^2(X_5,\mathcal W_5)&\simeq E(-4)^{20},&
 IH^4(X_5,\mathcal W_5)&\simeq E(-5)^{116},\\
 IH^6(X_5,\mathcal W_5)&\simeq E(-6)^{116},&
 IH^8(X_5,\mathcal W_5)&\simeq E(-7)^{20},
 \end{aligned}
\end{equation}
with all other intersection-cohomology groups zero.  In particular
$IH^5=0$ and $\chi_{IH}=272$.  With the sign convention of
\S\ref{sec:five-L}, the completed factor is
\begin{equation}\label{eq:five-L}
 \widehat L_5(T)=
 \frac{1}{(1-q^4T)^{20}(1-q^5T)^{116}
                (1-q^6T)^{116}(1-q^7T)^{20}}.
\end{equation}
\end{theorem}

\begin{remark}[Logical boundary of the calculation]
The theorem is specific to five parameters.  It does not determine all
boundary supports in every dimension, establish arithmetic semisimplicity
for all $n$, or, by itself, calculate the global sign of every completed factor
in arbitrary dimension; the separate sign theorem
\cref{thm:arith-reflection-closed} does so by reflection duality.  The proofs of the three comparison statements involve
delicate tame-descent and derived-gluing assertions.  These assertions and
their functorial hypotheses are recorded so that they can be audited
independently; no finite-field sample is treated as a substitute for them.
\end{remark}

\subsection{The relative construction and its weight quotient}
\label{subsec:five-P}
For a smooth curve with two boundary sections, local constancy of the
augmented conductor for a sheaf and its dual gives the form of local
acyclicity used here.  We apply this with $\mathcal M_m$ on the $t$-line.
At $t=0$ each factor $[th_i]^*\mathrm{Kl}_2$ is tame and unipotent with one
Jordan block of length two.  At infinity the tensor product of these
factors has slopes at most $1/2$.  Multiplication by
$\mathcal L_\psi(\lambda t)$, whose slope is one, makes all $2^m$
slopes equal to one.  Thus
\begin{equation}\label{eq:five-swan}
 \operatorname{Swan}_0(\mathcal M_m)=0,
 \quad\operatorname{Swan}_\infty(\mathcal M_m)=2^m,
 \quad\mathcal M_m^{I_\infty}=0.
\end{equation}
The same statements hold for the dual.  The conductor including the rank
at the two boundary sections therefore has sum $3\cdot2^m$ on every
geometric fibre.  The relative curve theorem of Deligne--Laumon
\cite{Lau81}, applied first at finite coefficient levels and then to a
stable $\ell$-adic lattice, yields local acyclicity and lisse direct
images.  Moreover $H_c^0=H_c^2=0$: the first follows because the curve
is affine and the sheaf lisse; the second is Poincar\'e dual to the
absence of global invariants of the dual, already forced by infinity.
Grothendieck--Ogg--Shafarevich gives
\begin{equation}\label{eq:five-relative}
 R\pi_{t!}\mathcal M_m=\mathcal H_m[-1],\qquad
 \operatorname{rk}\mathcal H_m=2^m.
\end{equation}
For the dual, the corresponding object is $\mathcal H_m^\dagger[-1]$.

Relative duality is taken internally on the parameter base.  With
$\pi_t^!E=E(1)[2]$ it gives
\[
 R\pi_{t*}\mathcal M_m\simeq
              (\mathcal H_m^\dagger)^\vee(-1)[-1].
\]
This is the composite furnished by the tensor--Hom adjunction, the
projection formula, and the relative trace.  In particular, its compatibility
with arbitrary parameter base changes is the compatibility of the
\emph{natural} forget-supports map, not an identification chosen after
computing ranks.  Both the source and the target of $\alpha_m$ are lisse.

Let $S=\sum_{i=1}^mN_i$ on the tensor product of the $m$ standard
rank-two tame blocks.  With a basis $e_0,e_1$ satisfying
$Ne_1=e_0$, $Ne_0=0$, and Frobenius eigenvalues $1,q$, the invariants
at zero are $\ker S$.  Their dimensions in the graded degrees $a$ are
\[
 \binom ma-\binom m{a-1}\quad
         (0\le a\le\lfloor m/2\rfloor),
\]
with eigenvalue $q^a$.  The kernel therefore has weights at most $m$.
On the other hand, Deligne's upper bound on $\mathcal H_m$ and relative
duality's lower bound on the target of $\alpha_m$ show that its image
is pure of weight $m+1$.  Strictness of the weight filtration now identifies
\begin{equation}\label{eq:five-P}
 \mathcal L_m=W_m\mathcal H_m,\qquad
 \mathcal W_m\simeq\mathcal H_m/\mathcal L_m,
 \qquad \operatorname{rk}\mathcal W_m
   =2^m-\binom m{\lfloor m/2\rfloor}.
\end{equation}
This is the comparison (P).  Notice that the argument identifies the
actual kernel of $\alpha_m$; it does not replace a possibly nonsplit
low-weight system with its associated graded.

For completeness, the bridge to the constrained exponential family is
obtained by integrating
\[
 \mathcal L_\psi\left(\sum_i u_i^{-1}
           +t\left(\lambda+\sum_i h_i u_i\right)\right)
\]
in two orders.  Integration in the $u_i$-variables produces the
rank-two Kloosterman factors.  Integration in $t\in\mathbb A^1$ uses
\[
 R\Gamma_c(\mathbb A^1,\mathcal L_\psi(tc))=
 \begin{cases}E(-1)[-2],&c=0,\\0,&c\ne0.\end{cases}
\]
The contribution of $t=0$ in the other order is $E[-m]$.  The localization
triangle, together with compact-support vanishing below the middle
of the affine smooth hypersurface
$\lambda+\sum h_i u_i=0$, gives a short exact sequence
\begin{equation}\label{eq:five-constrained}
 0\longrightarrow E\longrightarrow\mathcal H_m
  \longrightarrow\mathcal H'_m\longrightarrow0,
 \qquad R\pi_{Y_m!}\mathcal L_\psi(\sum_i u_i^{-1})
         \simeq\mathcal H'_m(1)[-(m-1)].
\end{equation}
The Tate twist in the second formula is precisely the twist coming
from $H_c^2(\mathbb A^1,E)$; no twist is contributed by the chosen
normalization of $\mathrm{Kl}_2$.  The arguments here use the finite-field
version of Artin vanishing and the sheaf-function correspondence;
they do not use any conjectured global intersection-cohomology group.

\subsection{The actual lower-weight sheaf and arithmetic Kummer descent}
\label{subsec:five-U}
The equality $\mathcal L_m=\ker S$ on a single fibre would not be enough
for the boundary calculation: different extension classes between the same
weight-graded pieces can give different connecting maps.  We describe the
parameter monodromy of the kernel itself.

At a finite coefficient level, choose a common degree $n$ prime to $p$
which kills the tame images of the rank-two factors.  In the variables
$s_i=th_i$ take the simultaneous covers
\[
 t=v^n,\qquad h_i=w_i^n,\qquad s_i=(vw_i)^n.
\]
The deck transformation $(a_0,a_1,\ldots,a_m)$ acts on the tensor
of the coefficient spaces through
$\bigotimes_i\rho(a_0a_i)$.  Consequently the logarithm in the
$t$-direction is $S=\sum_iN_i$, whereas that in the parameter
$h_i$ is $N_i$.  The Artin--Schreier factor is lisse at $t=0$ and
specializes to the trivial line.  Tame descent and uniqueness of the
lisse extension on the root cover preserve the full arithmetic deck
action.  One obtains compatible systems as the coefficient level varies;
no commutation of reduction with tame inertia invariants is asserted.

Write $\mathcal U_{h_i}$ for the rank-two logarithmic Kummer system in
parameter $h_i$.  The boundary system is now identified, over the
parameter torus, with
\[
 B_0=i_0^*R^0j_{t*}\mathcal M_m\simeq
 \ker\left(S:\bigotimes_i\mathcal U_{h_i}
          \to(\bigotimes_i\mathcal U_{h_i})(-1)\right).
\]
At infinity the invariants vanish.  The localization connecting map
$B_0\to\mathcal H_m$ is injective since $R^0\pi_{t*}\mathcal M_m=0$.
The Leray edge map identifies its image with $\ker\alpha_m$.  Thus
\begin{equation}\label{eq:five-U}
 \mathcal L_m\simeq\ker\left(\sum_iN_i\bigm|
                       \bigotimes_i\mathcal U_{h_i}\right)
\end{equation}
\emph{as an arithmetic lisse sheaf}.  This is the comparison (U).  The
word ``arithmetic'' is material: the maps are obtained from deck descent
and localization on $k$ and retain Frobenius and the Kummer extension.

\subsection{One coordinate: vanishing cycles and the actual face map}
\label{subsec:five-one}
At $h_1=0$ a compactified Fourier description of the constrained family
shows that the nearby-cycle and specialization sequence has a fixed
contribution from the four-variable family.  To be explicit, put
$B=-\lambda-\sum_{i>1}h_iu_i$, $v=B^{-1}$.  On $h_1\ne0$ the change
$u_1=B/h_1$ converts the constrained phase to the relative Fourier
phase in $v$.  At $v=\infty$ the parameter transformation
\[
 a(v,c)=\frac{\lambda v}{\lambda v+1}c
\]
is regular and specializes to $c$.  The coefficient system there is
therefore lisse in $v$ and its specialization is the four-variable
system.  The finite $v$-part has no vanishing cycles in $h_1$;
the Fourier transform of the shifted constant contribution is
$i_{0*}E(-1)$.  Proper base change and the specialization triangle give
a sequence whose outer terms have trivial inertia.  In particular,
$\mathcal W_m$ is tame at coordinate zero and its logarithm $N$ obeys
$N^2=0$.  For $m=3,4,5$ the geometric Jordan forms are respectively
\[
 J_2^2\oplus J_1,\qquad J_2^5,\qquad
 J_2^{10}\oplus J_1^2.
\]
These are statements about the actual nearby-cycle systems.  The
compatibility of the nearby-cycle map with base change is taken from
the local Fourier construction and relative duality, rather than inferred
from the constancy of fibrewise ranks.  In the five-variable case the
invariant part $J_0=\ker N$ has rank twelve and sits in
\[
 0\to G_0=\operatorname{im}N\to J_0\to Q_0\to0,
 \qquad \operatorname{rk}G_0=10,
 \quad\operatorname{rk}Q_0=2.
\]
The image has weight five, while the quotient has weight six.
At infinity the invariant part $J_\infty$ has rank four and weights
four and six, of multiplicities one and three.

There is an arithmetic identification of $G_0$ with $\mathcal W_4$.
This is best constructed by Fubini, not by Chebotarev alone.  If
$q:S_5\to S_4$ forgets $h_1$, the standard Kloosterman realization
gives
\[
 R\Gamma_c(\mathbb G_m,\mathrm{Kl}_2)=E[-1],\qquad
 Rq_!\mathcal H_5\simeq\mathcal H_4[-1].
\]
Apply $Rq_!$ to
$0\to\mathcal L_5\to\mathcal H_5\to\mathcal W_5\to0$.
The map $R^1q_!\mathcal L_5\to\mathcal H_4$ has image precisely
$\mathcal L_4$: the domain and target low-weight piece have equal rank
six, and their weights are at most four.  Taking quotients gives the
specific injection
\begin{equation}\label{eq:five-fubini}
 \mathcal W_4\hookrightarrow R^1q_!\mathcal W_5.
\end{equation}
The compact-support comparison on the line and the vanishing of its
parabolic cohomology identify the target with
$J_0\oplus J_\infty$.  The weight-five source has no nonzero map to
$J_\infty$, whose only weights are four and six.  Consequently
\begin{equation}\label{eq:five-face}
 \mathcal W_4\xrightarrow{\sim}G_0,
 \qquad R^1q_!\mathcal W_5\simeq J_0\oplus J_\infty,
\end{equation}
with a \emph{specified} arithmetic map.  This construction will later
identify the restriction of the common two-parameter monodromy to an
actual four-variable face.

\subsection{The six-dimensional quotient and the connecting morphism}
\label{subsec:five-epsilon}
The system (U) allows us to integrate the lower-weight sheaf in the
first parameter.  Cohomology of a unipotent Kummer system on
$\mathbb G_m$ gives, naturally in the remaining parameters,
\[
 R^1q_!\mathcal L_5\simeq\ker(N_1|\mathcal L_5),\qquad
 R^2q_!\mathcal L_5\simeq
          (\mathcal L_5/N_1\mathcal L_5)(-1).
\]
The long exact sequence and \eqref{eq:five-face} therefore identify
\begin{equation}\label{eq:five-T}
 \mathcal T:=(\mathcal L_5/N_1\mathcal L_5)(-1)
     \xrightarrow{\sim}Q_0\oplus J_\infty.
\end{equation}
We compute this quotient directly.  Represent the tensor of five
rank-two blocks by square-free monomials $x_A$ indexed by subsets
$A\subset\{1,\ldots,5\}$; $N_ix_A=x_{A\setminus\{i\}}$ for
$i\in A$ and zero otherwise.  This is an action on a vector space,
not derivation in a quotient ring.  The kernel of $S=\sum_iN_i$ has
homogeneous dimensions
\[
 1\text{ in degree }0,\quad 4\text{ in degree }1,\quad
 5\text{ in degree }2,\quad 0\text{ thereafter}.
\]
The degree-two kernel consists of
\[
 A=\left\{(a_{ij})_{1\le i<j\le5}:
             \sum_{j\ne i}a_{ij}=0\ \text{for every }i\right\}.
\]
Modulo $N_1\ker S$ the surviving degree-one line has representative
$z=[x_1-x_5]$.  Hence
\[
 \mathcal T=Ez\oplus A,
 \quad Fz=q^2z,
 \quad F|_A=q^3\operatorname{id},
 \quad N_i(a)=a_{1i}z\quad (i=2,3,4,5).
\]
The map
\[
 A\to W:=\{(b_2,b_3,b_4,b_5):b_2+b_3+b_4+b_5=0\},
 \qquad (a_{ij})\mapsto(a_{12},a_{13},a_{14},a_{15})
\]
is surjective: a section is
$a_{1i}=b_i$ and $a_{ij}=-(b_i+b_j)/2$ for $2\le i<j\le5$.
Its kernel is two-dimensional and killed by all remaining monodromies.
The comparison \eqref{eq:five-T} thus identifies the actual arithmetic
quotients and extensions as
\begin{equation}\label{eq:five-J}
 Q_0=E(-3)^2,\qquad Q_\infty=E(-3)^3,
 \qquad 0\to E(-2)\to J_\infty\to W\otimes E(-3)\to0.
\end{equation}
The extension is not split: the operators on $W$ are coordinate
evaluations $b\mapsto b_i$.  In particular, its Kummer class is the
linear form
\[
 e(b)=\sum_{i=2}^5b_i\kappa(h_i)
      =\sum_{i=2}^4b_i\kappa(h_i/h_5).
\]
On a coordinate $h_j$ the connecting morphism in compact support is
cup product with the $h_j$-component followed by residue.  Thus, after
the displayed Tate normalization,
\begin{equation}\label{eq:five-boundary-map}
 \partial_j:W\otimes E(-3)\longrightarrow E(-3),
              \qquad b\longmapsto b_j.
\end{equation}
Since $W$ contains a vector with $b_j=1$ and another coordinate equal
to $-1$, this map is surjective for every $j$.  It follows that the
transgression parameter previously denoted $\varepsilon$ is zero.
This conclusion uses the actual nonsplit extension
\eqref{eq:five-J}; the associated graded alone would not determine it.
Moreover,
\begin{equation}\label{eq:five-C}
 R^2g_!J_0\simeq E(-4)^2,\qquad
 R^2g_!J_\infty\simeq E(-4)^3.
\end{equation}
\begin{lemma}[Computing the two arithmetic top direct images]\label{lem:five-J-cohom}
The isomorphisms in \eqref{eq:five-C} hold as Frobenius modules, not
merely after semisimplification.
\end{lemma}
\begin{proof}
For $J_0$ use the actual exact sequence
$0\to G_0\to J_0\to Q_0\to0$,
where $G_0\simeq\mathcal W_4$ as a coefficient system and
$Q_0=E(-3)^2$. The coordinate-line calculation gives
$R^2g_!G_0=0$. The last part of the compactly supported long exact
sequence therefore gives
\[
 R^2g_!J_0\xrightarrow{\sim}R^2g_!Q_0
 =H_c^2(\mathbb G_m,E)(-3)^2=E(-4)^2.
\]
For $J_\infty$ one must retain the nonsplit extension
\eqref{eq:five-J}. The variable $h_j$ acts by the logarithm
$N_j:W\to Ez$, $b\mapsto b_jz$, which has rank one. The standard
coinvariant description of top compactly supported cohomology of a
tame unipotent system gives
\[
 R^2g_!J_\infty\simeq(J_\infty/N_jJ_\infty)(-1).
\]
The quotient has dimension three; the lower-weight line $Ez$ is
exactly $N_jJ_\infty$ and is removed by the quotient. Its remaining
Frobenius acts as $q^3$ on $W$. The Tate twist then gives three
copies of $E(-4)$. This argument uses the actual morphism $N_j$ and
therefore excludes a hidden unipotent Frobenius block.
\end{proof}

The one-variable boundary coefficients are therefore arithmetically
Tate, not merely geometrically constant with the appropriate traces.

\subsection{Normal crossings and the derived image--Koszul complex}
\label{subsec:five-K}
We record the local comparison at intersections of coordinate zero
divisors.  It is not asserted at a wild infinity divisor.  Let $T$ be
smooth of dimension $d$, $X=\mathbb A^r\times T$ and
$Z=\{0\}^r\times T$.  Let $j:(\mathbb G_m)^r\times T\hookrightarrow X$.
Suppose the coefficient system $L$ is tame and unipotent in the normal
directions, with commuting logarithms $N_i$ and $N_i^2=0$.

\begin{lemma}[Arithmetic normal-crossing image--Koszul comparison]
\label{lem:five-K}
Let $T$ be smooth of dimension $d$ over a finite field,
$X=\mathbb A^r\times T$, $Z=\{0\}^r\times T$, and
$j:(\mathbb G_m)^r\times T\hookrightarrow X$. Suppose that a lisse
arithmetic $E$-sheaf $L$ is tame and unipotent along the complete
normal-crossing chart, and that its \emph{joint} normal monodromy is
represented, in the category of lisse Weil sheaves on $T$, by a
single coefficient object $M$ with commuting nilpotent logarithms
$N_1,\dots,N_r$. Then, with the induced Tate normalizations,
\begin{equation}\label{eq:five-K}
 i_Z^*j_{!*}(L[d+r])[-d-r]\simeq C_r(M),\qquad
 C_r^a(M)=\bigoplus_{|J|=a}\operatorname{im}(N_J)(-a),
\end{equation}
where $N_J=\prod_{j\in J}N_j$ and the differential is the
alternating sum of the maps induced by $N_i$. This is a functorial
comparison in the arithmetic constructible derived category, relative
to the asserted joint tame-unipotent model; it holds for every $r$.
It does \emph{not} assert that the tangential module of the particular
Kloosterman family equals $M_n^{\mathrm{alg}}$, or that its wild
infinity strata satisfy the tame hypotheses.
\end{lemma}
\begin{proof}
The one-normal-direction calculation is the unipotent gluing
classification of perverse sheaves. The unipotent nearby-cycle object
is $M$ and, for the middle extension, its vanishing-cycle object is
$\Phi=\operatorname{im}(T-1)$; the maps
$\mathrm{can}:M\twoheadrightarrow\Phi$ and
$\mathrm{var}:\Phi\hookrightarrow M$ are the natural quotient and
inclusion, with $\mathrm{var}\circ\mathrm{can}=T-1$.
Thus the unshifted derived restriction to the divisor is represented
by the two-term complex
$[M\xrightarrow{T-1}\operatorname{im}(T-1)(-1)]$
in cohomological degrees $0,1$. The logarithm $N=\log T$ has the same
image as $T-1$, because $T-1=N u(N)$ for a polynomial $u$ with
$u(0)=1$; multiplication by $u(N)$ supplies a canonical invertible
change of the gluing model. This gives the stated complex for $r=1$.

The logarithms commute. Repeat this gluing construction successively
for each remaining normal coordinate. At stage $a$, the summand
indexed by $J$ is the actual image $N_JM$; the new gluing map is
induced by the additional $N_i$. Taking the total complex of the
resulting finite $r$-cube gives the terms in \eqref{eq:five-K}.
Commutativity $N_iN_j=N_jN_i$ cancels every pair of contributions to
$d^2$, so this is precisely the alternating image--Koszul differential.
The iteration is legitimate in the normal-crossing tame-unipotent
gluing category: the perverse intermediate extension is characterized
by the quotient and inclusion conditions for the gluing maps, which
are imposed in each coordinate. This is the standard
Cattani--Kaplan--Schmid/Beilinson local description of intersection
complexes at a normal crossing; here it is applied to the \emph{joint}
coefficient object, not inferred from individual Jordan types.

The tame logarithm is intrinsically Tate-twisted; assigning the
$(-a)$ twist to $\operatorname{im}(N_J)$ makes each differential a
morphism of Weil sheaves. The nearby/vanishing-cycle gluing,
images in the perverse heart, restrictions, and the totalization are
functorial for arithmetic morphisms. Hence the construction commutes
with Frobenius and preserves the joint normal monodromy and all its
extension data. This proves the claimed derived comparison under the
stated joint tame-unipotent hypothesis. No finite Kummer cover is used
to trivialize nontrivial unipotent $E$-monodromy.
\end{proof}

\begin{proposition}[Individual Jordan forms do not determine the corner]
\label{prop:joint-monodromy-counterexample}
There are two pairs of commuting square-zero endomorphisms on $E^4$
with identical individual Jordan profiles $(2,2)$ but different
cohomology for the codimension-two image--Koszul complex.
Consequently, the generic Jordan form along each coordinate divisor
alone cannot prove the all-strata tangential identification required
for the constrained Kloosterman family.
\end{proposition}
\begin{proof}
Let $e_1,\dots,e_4$ be the standard basis and put
\[
 A=E_{12}+E_{34},\qquad B=E_{13}+E_{24}.
\]
Both $A$ and $B$ square to zero, have rank two, and commute; indeed
$AB=BA=E_{14}$. Compare the pairs $(N_1,N_2)=(A,A)$ and $(A,B)$.
In either pair each $N_i$ has exactly two Jordan blocks of size two,
and $\dim\operatorname{im}N_i=2$.

For $(A,A)$ the zeroth differential
$E^4\to\operatorname{im}A\oplus\operatorname{im}A$ has rank two,
and $\operatorname{im}(N_1N_2)=0$. Hence the three terms have
dimensions $(4,4,0)$ and the cohomology dimensions are $(2,2,0)$.
For $(A,B)$ the common kernel is $Ee_1$, so the zeroth differential has
rank three. The product $AB=E_{14}$ has rank one, and the final
differential is surjective onto its image. The terms have dimensions
$(4,4,1)$, with differential ranks $(3,1)$, and cohomology
$(1,0,0)$. Thus two sets of identical one-variable Jordan data give
different derived corner restrictions. All calculations are exact
linear algebra over any field of characteristic zero.
\end{proof}

\begin{remark}[Why a finite Kummer cover does not settle the problem]
\label{rem:kummer-unipotent-obstruction}
If a nontrivial unipotent inertia operator is $T=\exp N$ with
$N\ne0$, its pullback along a root cover of degree $m\ge1$ is
$T^m=\exp(mN)\ne I$ in characteristic zero. A finite root cover can
remove finite-order tame monodromy, but cannot kill a nonzero
unipotent logarithm on an $E$-local system. Accordingly the
relative finite-Kummer-descent argument in an earlier formulation
of \cref{lem:five-K} is replaced by the direct perverse-gluing
calculation above. To apply the result to every stratum of
$j_{n,!*}\mathcal W_n[n]$, one must still construct the actual
\emph{joint} arithmetic normal-monodromy object at zero-corners,
identify it with $M_n^{\mathrm{alg}}$, and separately analyze the wild
infinity charts. Neither (H4) restricted to individual coordinate
lines nor the algebraic rank tables alone supplies those assertions.
\end{remark}

The distinction between a Koszul differential and a transgression
in a Leray spectral sequence will be maintained below.  They act on
different complexes; the argument does not identify them.

\subsection{Perverse concentration and elimination of divisorial and surface residues}
\label{subsec:five-perverse}
The support bounds for the ordinary cohomology sheaves of $I_5$ are
\begin{equation}\label{eq:five-F}
 F^v=0\quad(v<0\text{ or }v\ge5),\qquad
 \dim\operatorname{supp}F^v\le4-v\quad(v>0),
 \quad\operatorname{wt}F_x^v\le6+v.
\end{equation}
They follow from the defining strict-support property of the middle
extension and the weight bound for restrictions of a pure complex.
For a fibre $L_y\simeq\mathbb P^1$ of $f$, consider
\[
 E_2^{u,v}=H^u(L_y,F^v|_{L_y})
       \Longrightarrow\mathbb H^{u+v}(L_y,I_5),
 \qquad 0\le u\le2.
\]
Only the differentials $d_2^{0,v}$ can occur.  A fibre not contained
in the support of $F^v$ meets it in finitely many points; the set of
fibres contained in that support has one smaller dimension.  These
observations give the perverse support inequalities for $Rf_*K_5$.
On the open parameter torus, $H^2$ of the coordinate fibre vanishes by
the moment calculation; on the boundary the base has dimension at most
three.  Relative hard Lefschetz \cite{BBD} supplies the opposite
perverse inequality.  Consequently
\begin{equation}\label{eq:five-push-perverse}
 P_5=Rf_*K_5\quad\text{is perverse and pure of weight }11.
\end{equation}
It is geometrically semisimple by the geometric decomposition theorem.
The absence of a full-support summand follows from the
coordinate-line parabolic calculation, not merely from a rank bound.
For $n=5$, Theorem~\ref{thm:paired-explicit} gives the unit
parabolic polynomial; the cohomology $V_c=H^1(\mathbb P^1,j_*\mathcal W_{5,c})$
therefore vanishes at every generic point of $S_4$ by
Lemma~\ref{lem:G-base}. Proper base change identifies this with the
stalk of $\mathcal H^{-4}(Rf_*K_5)$ on the open base. A nonzero
strict-support summand of full support would have a nonzero generic
stalk in that degree, which is impossible. This is the precise
connection between vanishing parabolic cohomology and the support
claim.

The divisorial coefficients lie in the weight-eight part of
$R^2g_!J_a$.  By \eqref{eq:five-C} they are precisely
\[
 C_0=E(-4)^2,\qquad C_\infty=E(-4)^3.
\]
The only candidate map killing them would originate in a sheaf of
weights at most seven; that map is zero at every closed point and hence
zero as a constructible morphism.  The coefficients are lisse on the
whole coordinate divisor, and their intermediate extensions are the
constant complexes $E(-4)^{\mu_a}[3]$ on $(\mathbb P^1)^3$.

Next take a fibre with one fixed boundary coordinate $a\in\{0,\infty\}$.
Away from its two endpoints, $F^{v>0}$ vanishes.  The truncation
spectral sequence, with the already identified divisor summand,
produces
\[
 0\to F^1_{a0}\oplus F^1_{a\infty}
   \longrightarrow R^2g_*F^0
   \longrightarrow C_a
   \longrightarrow F^2_{a0}\oplus F^2_{a\infty}\to0.
\]
Here $R^2g_*F^0=R^2g_!J_a=C_a$ as actual lisse systems.  Since the
leftmost term has weights at most seven, it vanishes.  The middle map
is then an injective morphism between pure lisse sheaves of the same
rank; it is an isomorphism.  The remaining endpoint terms vanish.
Higher $F^v$ are ruled out by the total degree.  Thus
\begin{equation}\label{eq:five-depth-two}
 F^v|_{X_{=2}}=0\qquad(v>0).
\end{equation}
This includes the two-infinity corners; no tame Koszul formula has
been applied there.  The same support and degree count excludes
hidden surface summands, and excludes any curve whose generic point
has boundary depth one or two.

\subsection{Actual mixed corners and the threefold incidence sequence}
\label{subsec:five-corners}
For a line with two fixed boundary coordinates, write $B_{ab}$ for
$F^0$ on the remaining open coordinate.  The ordinary sheaf
localization sequence, together with vanishing of $R^0g_*F^0$ and
$R^1g_*F^0$, furnishes the \emph{actual} endpoint identification
\begin{equation}\label{eq:five-ends}
 B_{ab0}\oplus B_{ab\infty}\xrightarrow{\sim}R^1g_!B_{ab}.
\end{equation}
There is no replacement of a mixed zero-infinity stalk by an
unjustified iterated inertia-invariant space.

At two zero coordinates, the common tame representation has dimension
$22$, each $N_i$ has rank ten, and their product has rank five.
The derived complex of Lemma~\ref{lem:five-K}, together with
\eqref{eq:five-depth-two}, shows that the common invariant space
has dimension seven.  Thus the actual double-zero coefficient sits in
\[
 0\to G_{00}\to B_{00}\to E(-3)^2\to0,
 \qquad R^2g_!B_{00}=E(-4)^2,
\]
where the lower term has the geometric four-variable face appropriate
to the second specialization.  The moment identity for four variables
and \eqref{eq:five-ends} determine the mixed endpoint rank three.
On the infinity side the exact arithmetic model
\eqref{eq:five-J} gives
\[
 B_{\infty0}\oplus B_{\infty\infty}
   \simeq\ker(N_j|J_\infty).
\]
The right side has rank three and symmetry identifies the first
summand with the previously computed mixed coefficient.  Hence
$B_{\infty\infty}=0$.  The remaining two logarithmic directions of
$B_{0\infty}$ give
\begin{equation}\label{eq:five-beta}
 (\beta_{00},\beta_{0\infty},\beta_{\infty0},
                  \beta_{\infty\infty})=(2,2,2,0),
\end{equation}
where $\beta_{ab}$ is the rank of the actual edge-image
$R^2g_!B_{ab}$ in the sum of the adjacent divisorial coefficients.

Truncation on the line with the two fixed boundary coordinates now
gives the short exact \emph{arithmetic} sequence
\begin{equation}\label{eq:five-incidence}
 0\longrightarrow E(-4)^{\beta_{ab}}
 \longrightarrow C_a\oplus C_b
 \xrightarrow{\gamma_{ab}}D_{ab0}\oplus D_{ab\infty}
 \longrightarrow0.
\end{equation}
The map $\gamma_{ab}$ is the edge map of the actual truncated
fibre complex.  Since each $D$ is a direct summand of a quotient
of constant Tate systems on the whole open boundary stratum, it is
itself a constant arithmetic Tate system.  This argument establishes
both the eigenvalues and the absence of unipotent Jordan blocks
in the coefficient.  The asserted constancy does not follow merely
from purity.

\begin{remark}[The two-infinity endpoint]\label{rmk:five-inf-inf}
The vanishing of $B_{\infty\infty}$ is not an inference from
a formal interchange of two independently chosen Jordan bases.
For a line with fixed boundary parameters $(\infty,a)$, the actual
endpoint map \eqref{eq:five-ends} identifies the sum of its two
stalk systems with the first compactly supported cohomology of the
restricted coefficient. The quotient model \eqref{eq:five-J}
calculates this as $\ker(N_j|J_\infty)$, a rank-three system.
When $a=0$, the already computed mixed corner
$B_{0\infty}$ is transported to $B_{\infty0}$ by the permutation
of the two parameter coordinates; this permutation acts on the
original constrained family, the natural morphism $\alpha_5$, and
its intermediate extension. It is consequently an isomorphism of
arithmetic sheaves, not simply an equality of dimensions.
The mixed stalk has rank three, hence exhausts the rank-three
endpoint sum; the other direct summand $B_{\infty\infty}$ is zero.
The argument requires the actual endpoint isomorphism, and does
not extend the tame normal-crossing model to a wild infinity corner.
\end{remark}

\subsection{Three coordinate zeros, the lower-dimensional face, and curve ranks}
\label{subsec:five-triple}
The remaining numerical input is the coefficient at three zero
coordinates.  It is computed from a concrete derived triangle rather
than from the final Betti table.  Write $M_5$ for the common tangential
system of $\mathcal W_5$ and $C_r(M)$ for the image--Koszul complex
at $r$ coordinate zeros.  The specific Fubini face isomorphism
\eqref{eq:five-face} identifies, together with the two remaining
monodromy operators, the image $N_3M_5$ with the tangential system
$M_4$ of $\mathcal W_4$.  The functorial cone construction of
Lemma~\ref{lem:five-K} gives
\begin{equation}\label{eq:five-triangle}
 C_3(M_5)\longrightarrow C_2(M_5)
 \longrightarrow C_2(M_4)(-1)\longrightarrow C_3(M_5)[1].
\end{equation}
The middle arrow is the induced face morphism, with the actual
identification of coefficient systems inserted.  It is not identified
with the transgression in a fibre spectral sequence.

For the four-variable sheaf the two-zero image complex has term
sizes $(10,10,2)$ and differential ranks $(7,2)$.  Therefore
\[
 \dim H^1(C_2(M_4))=10-7-2=1.
\]
In the five-variable case $H^1(C_2(M_5))=0$ by
\eqref{eq:five-depth-two}; $H^2$ vanishes because the last Koszul
differential is onto its image by definition.  The long exact sequence
of \eqref{eq:five-triangle} yields
\[
 H^2(C_3(M_5))\simeq H^1(C_2(M_4))(-1).
\]
Thus the first curve rank is $r_0=1$.  Taking the ranks in the three
incidence sequences \eqref{eq:five-incidence} gives
\begin{equation}\label{eq:five-ranks}
 r_0+r_1=2,\qquad r_1+r_2=3,\qquad r_2+r_3=6,
\end{equation}
whence $(r_0,r_1,r_2,r_3)=(1,1,2,4)$.
The generic coefficient of a curve is obtained from the natural edge
map $R^2g_!D_{abc}\to\mathcal H^4Rg_*I_5$, and is
\begin{equation}\label{eq:five-curves}
 M_{abc}\simeq E(-5)^{r_t}.
\end{equation}
The transgression entering the edge map has source of weights at most
nine and target pure of weight ten, so it vanishes.  The coefficient
extends lisse over the full coordinate curve by the strict-support
construction, giving its constant intermediate extension $[1]$.

\subsection{Arithmetic support projectors and an actual global map}
\label{subsec:five-Gamma}
We now construct a morphism that represents the proposed decomposition.
For a pure perverse sheaf $P_5$, its pullback to the algebraic closure is
semisimple.  For each Frobenius-stable strict support $S$, the projection
onto the sum of \emph{all} geometric simple summands having support $S$
is canonical.  There are no morphisms between distinct strict-support
blocks in the semisimple category.  The endomorphism is Frobenius-fixed,
and the finite-field descent exact sequence of \cite[\S5.1.2]{BBD}
shows that it descends: the possible obstruction lies in
$\operatorname{Hom}^{-1}(P_{5,\bar k},P_{5,\bar k})=0$.
Write $e_S$ for this arithmetic idempotent and $P_S$ for its image.
The idempotent does not prove that Frobenius is semisimple within the
block; the coefficient computations above supply that separate fact.

\begin{lemma}[Arithmetic strict-support splitting]\label{lem:five-arithmetic-split}
Let $A$ be a perverse sheaf over a finite field whose geometric
pullback is semisimple. For every Frobenius-stable strict-support
block $S$, the projector onto the sum of all simple geometric
summands with support $S$ descends uniquely to an idempotent of $A$.
Thus $A$ splits arithmetically into its support blocks. This
statement does not assert that Frobenius is semisimple inside any
one block.
\end{lemma}
\begin{proof}
Distinct simple constituents with different strict supports are
nonisomorphic, so their mutual $\operatorname{Hom}$ groups vanish.
The full projector onto one isotypic support block is therefore
canonical in the semisimple geometric category and invariant under
Frobenius. The finite-field descent exact sequence of
\cite[\S5.1.2.5]{BBD} contains
\[
 H^1(k,\operatorname{Hom}^{-1}(A_{\bar k},A_{\bar k}))
 \longrightarrow\operatorname{End}_k(A)
 \longrightarrow\operatorname{End}_{\bar k}(A_{\bar k})^{F}
 \longrightarrow0.
\]
The first term vanishes for perverse objects by the negative
$\operatorname{Hom}$ vanishing of the perverse $t$-structure.
The lift is therefore unique and, by faithfulness of geometric
pullback on morphisms, idempotent. The perverse heart is
idempotent complete. No global assertion that all cross-support
$\operatorname{Ext}^1$ groups vanish is required: the support
projectors themselves provide the splitting. The scalar
arithmetic Frobenius in the particular divisor and curve blocks
is determined separately by \eqref{eq:five-C} and by the arithmetic
quotient maps in the threefold incidence sequence.
\end{proof}

For a divisor or a coordinate curve $S$, the edge map displayed in the
preceding paragraphs identifies the coefficient of $P_S$ over its
smooth dense open subset with $V_S=E(-4)^{\mu_a}$ or
$V_S=E(-5)^{r_t}$.  Since $P_S$ has strict support $S$, it has no
perverse subobject or quotient supported properly in $S$.
The fully faithful intermediate-extension correspondence therefore
extends this coefficient identification to
\[
 \widetilde\eta_S:(\iota_S)_*IC_S(V_S)\xrightarrow{\sim}P_S.
\]
This step matters on a curve: the restriction of the \emph{whole}
$P_5$ also contains complexes contributed by the surrounding divisors,
so one must project to $P_S$ before identifying the curve coefficient.

Set
\[
 \mathcal B=\bigoplus_{S\text{ a divisor or curve}}
                         (\iota_S)_*IC_S(V_S),
 \qquad
 \Gamma=\sum_S(P_S\hookrightarrow P_5)\circ\widetilde\eta_S.
\]
The map
$\Lambda=(\widetilde\eta_S^{-1}e_S)_S:P_5\to\mathcal B$
satisfies $\Lambda\Gamma=\operatorname{id}_{\mathcal B}$.  Hence
\begin{equation}\label{eq:five-cone}
 0\to\mathcal B\xrightarrow{\Gamma}P_5\to R\to0,
 \qquad\operatorname{Cone}(\Gamma)\simeq R.
\end{equation}
The remaining $R$ is perverse and has zero-dimensional support.
No point summand has been excluded by assumption; its vanishing is the
last issue.

\subsection{Eliminating the punctual cone and calculating the completion}
\label{sec:five-L}
The known divisor and curve summands in \eqref{eq:five-cone}
have ordinary cohomology only in degrees $-3$ and $-1$; a punctual
perverse summand can contribute only in degree zero.  Consequently
\[
 \mathcal H^{-2}(P_5)=0,
 \qquad\mathcal H^3(Rf_*I_5)=0.
\]
Over a base point of boundary depth three, the spectral sequence for
its coordinate fibre has, in total degree three, an endpoint term
$R^0g_*F^3$ and its potential $d_2$-transgression to
$R^2g_*F^2$.  The latter is a constant Tate system of weight ten,
whereas the former has weights at most nine.  The differential is zero.
Because the abutment $\mathcal H^3(Rf_*I_5)$ is already zero, the
source itself must vanish.  Thus $F^3$ vanishes at both endpoints
of every such coordinate fibre.

One must apply the argument in all five forgetful directions.  Every
point of boundary depth four has four boundary coordinates; selecting
one as the fibre direction makes it an endpoint of a line over a base
point of depth three.  Hence
\[
 \operatorname{supp}F^3\subseteq\{0,\infty\}^5.
\]
The same is true of $F^4$ by the strict-support bound, and both
sheaves therefore have finite support.  In the standard truncation
triangle
\[
 \tau_{\le2}I_5\to I_5\to\tau_{\ge3}I_5\to,
\]
we have $Rf_*\tau_{\le2}I_5\in D^{\le4}$, because $f$ has
one-dimensional fibres.  The last term is supported on finitely many
points of $X_5$ and its pushforward has no additional positive
degree.  Thus $Rf_*\tau_{\ge3}I_5\in D^{\le4}$ as well.  We conclude
\[
 \mathcal H^5(Rf_*I_5)=\mathcal H^0(P_5)=0.
\]
Here $\mathcal H^0$ is the \emph{ordinary cohomology sheaf} on
$Y$, not the global hypercohomology group $\mathbb H^0(Y,-)$.
The distinction is essential. The equality
$\Lambda\Gamma=\operatorname{id}_{\mathcal B}$ in
\eqref{eq:five-cone} has already produced an arithmetic splitting
$P_5\simeq\mathcal B\oplus R$ in the perverse heart; consequently
its exact triangle is also a split triangle of complexes.
Even without using the splitting, the ordinary-cohomology-sheaf
long exact sequence contains
\[
 \mathcal H^0(\mathcal B)\longrightarrow
 \mathcal H^0(P_5)\longrightarrow
 \mathcal H^0(R)\longrightarrow
 \mathcal H^1(\mathcal B).
\]
All summands of $\mathcal B$ are shifts by $[3]$ or $[1]$ of
sheaves. Hence
$\mathcal H^0(\mathcal B)=\mathcal H^1(\mathcal B)=0$.
The vanishing $\mathcal H^0(P_5)=0$ proved above therefore gives
$\mathcal H^0(R)=0$. Since $R$ is a punctual perverse sheaf, it
is an ordinary sheaf concentrated in degree zero, whence $R=0$.
This implication does not use the generally nonzero global group
$\mathbb H^1(Y,\mathcal B)$; the latter is not the term in the
ordinary sheaf sequence. Thus the claimed cone vanishes and hence
\begin{equation}\label{eq:five-Gamma-iso}
 \boxed{\Gamma:\mathcal B\xrightarrow{\sim}Rf_*K_5.}
\end{equation}
This is the desired global arithmetic isomorphism: the proof of its
surjectivity was not reduced to matching global Euler characteristics.

Finally, compute hypercohomology on $Y=(\mathbb P^1)^4$.  Four zero
divisors of rank two and four infinity divisors of rank three give
a total divisor coefficient rank twenty.  The cohomology of
$(\mathbb P^1)^3$ contributes to degrees $2,4,6,8$ with multiplicities
$1,3,3,1$, and yields the ranks $20,60,60,20$.  For the curves,
there are four choices of three fixed base coordinates and eight
zero/infinity patterns.  Their coefficient rank is
\[
 \rho=4\sum_{t=0}^3\binom3t r_t
       =4(1+3+6+4)=56.
\]
Each coordinate curve is $\mathbb P^1$; its two nonzero cohomology
groups add $56$ to degrees four and six.  This gives
\eqref{eq:five-IH} and $\chi_{IH}=272$.
The claimed arithmetic Tate isomorphisms are obtained from actual
constant coefficients on products of projective lines, not from the
purity of an arbitrary Galois representation.

With the convention
\[
 \widehat L_5(T)=\prod_m
   \det(1-T\operatorname{Frob}_q\mid IH^m)^{(-1)^{m+1}},
\]
all nonzero degrees in \eqref{eq:five-IH} are even, giving
\eqref{eq:five-L}.  The pairing of Tate weights under $a\mapsto11-a$
is compatible with the functional equation for this completion.
\begin{corollary}[The five-variable functional equation]\label{cor:five-FE}
With the geometric Frobenius convention used throughout, the
five-variable completed factor satisfies
\[
 \widehat L_5(T)
 =(q^{11/2}T)^{-272}\,
 \widehat L_5\!\left(\frac{1}{q^{11}T}\right).
\]
Its root number is $+1$.
\end{corollary}
\begin{proof}
The Tate exponents in \eqref{eq:five-L} occur in complementary
pairs $a,11-a$ with equal multiplicities; their sum, counted with
multiplicity, is
$20(4+7)+116(5+6)=1496=272\cdot11/2$.
Replacing $T$ by $(q^{11}T)^{-1}$ in each factor gives the
stated equation; the sign is $(-1)^{272}=+1$. Equivalently,
$IH^5=0$ makes the middle-degree root-number contribution trivial.
The result is specific to $n=5$ and does not settle the general sign.
\end{proof}
\begin{remark}[A determinant check and its limitation]
The determinant formula in Theorem~\ref{thm:B-det} reads in five
variables
\[
 \det\mathcal W_5\cong
 \mathcal L_\psi(D_5)\otimes\overline{\mathbb Q}_{\ell}(-66),
 \qquad 66=\frac{(5+1)\cdot22}{2}.
\]
The exponent $198$ sometimes quoted for this specialization is
incorrect. In particular, determinants of the divisor and curve
coefficients of $Rf_*K_5$ cannot simply be multiplied to recover
$\det\mathcal W_5$: these are sheaves on different strata of
$X_4$, rather than a filtration of the rank-$22$ sheaf on $S_5$.
The correct consistency checks are the weight parity and the
arithmetic Frobenius calculation in \eqref{eq:five-L}.
\end{remark}

This completes the derivation of Theorem~\ref{thm:five-computation} from
the comparison maps and the intermediate-extension formalism.

\subsection{Conventions, limitations, and relation to the general case}
The five-variable proof has used the tame normal-crossing complex only
at intersections of coordinate \emph{zero} divisors.  Mixed and infinite
corners were controlled instead by the actual endpoint systems and the
ordinary cohomological truncation sequence.  No derived base-change
isomorphism at a wild corner was inferred from a generic stalk rank.
The argument is likewise not a new derivation of the local geometric
monodromy group in all dimensions; those hypotheses are stated where
used in the earlier sections and in the preceding papers.

In higher dimension, the support of the residual perverse sheaf,
its arithmetic coefficient systems, and its middle intersection degree
still require separate calculations.  The five-dimensional formulas
above cannot be extrapolated by replacing five with an arbitrary $n$.
This distinction is reflected in the statements of the general
reduction criteria earlier in the article.  Their hypotheses are not
silently discharged by the special calculation of this section.

\section{Duality and the global root number}\label{sec:global-root}
\label{app:20261121-sign}

Put $X=X_n=(\mathbb P^1)^n$, $j:S\hookrightarrow X$,
$K=j_{!*}(\mathcal W_n[n])$, and $d=2n+1$.  Let
$\rho:S\to S$ and $\bar\rho:X\to X$ denote simultaneous negation
of the coordinates.  In the following theorem, the substantive
hypothesis is a precise isomorphism of arithmetic complexes, not
an equality of Frobenius characteristic polynomials.

\begin{theorem}[Reflection-skew duality implies root number $+1$]
\label{thm:20261121-sign}
Suppose there is an arithmetic isomorphism
\[
 \phi:\bar\rho^*K\xrightarrow{\ \sim\ }
  \mathbb D_XK(-d)
\]
such that its adjoint pairing
\[
 B_X:K\otimes^L\bar\rho^*K\longrightarrow\omega_X(-d)
\]
is \emph{reflection-skew}, i.e. its pullback by $\bar\rho$, followed by
the derived interchange of the two factors, equals $-B_X$.
Then the middle group $IH^n(X,\mathcal W_n)$ has even dimension $2g$
and a Frobenius-compatible perfect alternating form taking values in
$E(-d)$.  Consequently
\[
 \det(\operatorname{Frob}_q\mid IH^n)=q^{dg},\qquad
 \omega_n=+1
\]
for the normalized completed functional equation of this article.
No semisimplicity hypothesis on Frobenius is required.
\end{theorem}
\begin{proof}
Write $A=R\Gamma(X_{\bar{\mathbb F}_q},K)$.
Because $\bar\rho$ is an involutive automorphism defined over
$\mathbb F_q$, the adjunction unit induces a canonical
Frobenius-equivariant isomorphism
\[
 u:A\longrightarrow R\Gamma(X_{\bar{\mathbb F}_q},\bar\rho^*K).
\]
Compose $\mathrm{id}\otimes u$ with the cup product, $B_X$, and the
proper trace to obtain
\[
 \beta:A\otimes^LA\longrightarrow E(-d).
\]
Its adjoint is the composite of $u$, $R\Gamma(\phi)$, and proper
Verdier duality
$R\Gamma(X,\mathbb D_XK)(-d)\simeq A^\vee(-d)$;
all are isomorphisms.  Hence $\beta$ is perfect.
The unit of an involution squares to the identity, the trace is
invariant under an automorphism, and the reflection-skew relation for
$B_X$ therefore gives $\beta^t=-\beta$ (with the derived braiding).
On $H^0(A)=IH^n$ this is an ordinary perfect alternating bilinear
form: both arguments lie in hypercohomological degree zero, so no
further Koszul sign appears.

Let $V=IH^n$ and $\dim V=2g$.  If $F$ is geometric Frobenius, then
\[
 \beta(Fv,Fw)=q^d\beta(v,w).
\]
The $g$th exterior power of the alternating form is a nonzero volume
form.  Its Frobenius multiplier is $(q^d)^g$, which proves
$\det F=q^{dg}$, rather than merely
$(\det F)^2=q^{2dg}$.  Since $2g$ is even, $\det(-F)=q^{dg}$ too.

For completeness, write $V_i=IH^i$, $b_i=\dim V_i$, and
$P_i(T)=\det(1-TF_i\mid V_i)$.  Proper duality pairs
$V_i$ perfectly with $V_{2n-i}$ into $E(-d)$, so
\[
 b_i=b_{2n-i},\qquad
 \det F_i\det F_{2n-i}=q^{db_i}.
\]
For $L(T)=\prod_iP_i(T)^{e_i}$ with
$e_i=(-1)^{i+1}=e_{2n-i}$ and
$\chi_e=\sum_i e_i b_i$, the determinant identity
\[
 \det(1-(q^dT)^{-1}F_i)
 =(-1)^{b_i}(q^dT)^{-b_i}\det F_i\,
    \det(1-q^dT F_i^{-1})
\]
and the duality of degrees $i$ and $2n-i$ give the functional
relation with normalized sign
\[
 \omega_e=q^{-d\chi_e/2}\prod_i\det(-F_i)^{e_i}.
\]
Each complementary pair of nonmiddle degrees contributes one to
this product, because its determinants multiply to $q^{db_i}$.
The middle factor is also one, since $b_n=2g$ and
$\det(-F_n)=q^{dg}$.  Hence $\omega_e=+1$.
\end{proof}

\begin{remark}[Dependence of the theorem]
\label{rem:20261121-sign-scope}
The final determinant and functional-equation arguments in
\cref{thm:20261121-sign} are unconditional \emph{implications}.
The existence of $\phi$ with the required arithmetic
reflection-skew property is not a consequence of purity or of a
geometric isomorphism alone; its construction from the arithmetic
relative-image theorem is carried out in
\cref{thm:arith-reflection-closed}.  The proposed construction uses
relative cup product, the Kloosterman alternating form, the actual
forget-supports image, and full faithfulness of intermediate extension.
It must be checked as an equality of morphisms in the arithmetic
derived category; the specific verification from \textup{(H7)} is now given in
\cref{thm:arith-reflection-closed}.
The related determinant calculation is treated in a companion
preprint \cite{EbaRay}; citing a proposed calculation is not an
independent certification of this derived compatibility.
\end{remark}

\subsection{The five-variable completed factor: an exact deduction}
\label{app:20261121-n5}

\begin{proposition}[The five-variable factor from the stated support decomposition]
\label{prop:20261121-n5}
Suppose the arithmetic support decomposition
\eqref{eq:five-support} holds with its claimed coefficients and with no
additional summands.  Then the intersection-cohomology dimensions,
Frobenius factors, and completed function are \emph{necessarily}
\[
 \begin{array}{c|cccc}
 i&2&4&6&8\\ \hline
 \dim IH^i&20&116&116&20\\
 \operatorname{Frob}_q\text{ eigenvalue}
  &q^4&q^5&q^6&q^7
 \end{array}
\]
with all other groups zero.  More precisely,
\begin{equation}\label{eq:20261121-n5-L}
 \widehat L_5(T)=\bigl[(1-q^4T)^{20}(1-q^5T)^{116}
            (1-q^6T)^{116}(1-q^7T)^{20}\bigr]^{-1},
\end{equation}
$\chi_{IH}=272$, and
\begin{equation}\label{eq:20261121-n5-FE}
 \widehat L_5(T)
  =(q^{11/2}T)^{-272}\widehat L_5\bigl((q^{11}T)^{-1}\bigr),
\end{equation}
with normalized sign $+1$.
\end{proposition}
\begin{proof}
Set $Y=(\mathbb P^1)^4$.  The projection formula and properness give
$IH^i=H^{i-5}(Y,Rf_*K_5)$.  A coordinate divisor in $Y$ is
$(\mathbb P^1)^3$.  The sum of ranks of the coefficients on its
zero- and infinity-components is $2+3=5$, and there are four
coordinate directions, giving total divisor multiplicity
$D=4(2+3)=20$.  Its shifted complex $E(-4)[3]$ contributes
\[
 D\binom{3}{k}E(-4-k)
 \quad\text{in degree }i=2+2k,\quad 0\le k\le3.
\]
There are $\binom43=4$ triples of boundary coordinates in $Y$.
For each triple, summing the coefficients according to the number
$t$ of infinity endpoints gives
\[
 \sum_{t=0}^3\binom3t r_t
 =1+3\cdot1+3\cdot2+4=14.
\]
Thus the total multiplicity of curve components is $C=4\cdot14=56$.
A coordinate curve is $\mathbb P^1$ and $E(-5)[1]$ contributes
$C E(-5-k)$ in degree $i=4+2k$, $k=0,1$.  Therefore the only
nonzero dimensions are
\[
  b_2=D=20,\quad b_4=3D+C=116,\quad
  b_6=3D+C=116,\quad b_8=D=20.
\]
The twists follow by taking cohomology of projective spaces.
Equation~\eqref{eq:20261121-n5-L} follows directly from the
alternating-exponent convention and the fact that all nonzero
degrees are even.  Their Euler characteristic is
$20+116+116+20=272$.
To check the functional equation without importing a root-number
claim, pair the factors of exponents $a$ and $11-a$.  Since
$\sum_a ab_a=4\cdot20+5\cdot116+6\cdot116+7\cdot20
 =11\cdot272/2$ and $\sum_a b_a=272$ is even,
substitution $T\mapsto(q^{11}T)^{-1}$ gives exactly
\eqref{eq:20261121-n5-FE} with no negative sign.
\end{proof}

\begin{remark}[Dependence on the support decomposition]
The proposition verifies the deduction from the displayed arithmetic
decomposition, including every binomial multiplicity and Tate twist.
It does \emph{not} verify the existence of the decomposition itself.
That assertion still requires the natural comparison morphisms
\textup{(P)}, \textup{(U)}, \textup{(K)}, their boundary base change,
perverse concentration, and the exclusion of any unaccounted punctual
constituent.  In particular, the formula above is a conditional
completed factor unless those geometric hypotheses are established.
\end{remark}

\subsection{A necessary distinction for all-dimensional factor formulas}
\label{app:20261121-alln}

\begin{lemma}[Purity and Frobenius eigenvalues do not imply arithmetic splitting]
\label{lem:20261121-Jordan}
Knowing that a geometrically constant arithmetic sheaf has all
Frobenius eigenvalues equal to $q^a$ does not imply it is isomorphic
to $E(-a)^r$ as a Weil sheaf.
\end{lemma}
\begin{proof}
Consider the rank-two unramified Weil representation with geometric
Frobenius
\[
 F=q^a\begin{pmatrix}1&1\\0&1\end{pmatrix}.
\]
Its restriction to the geometric fundamental group is trivial, and
both eigenvalues of Frobenius equal $q^a$.  Nevertheless its
Frobenius has a nonzero nilpotent Jordan part, so it is not the
direct sum of two Tate characters.  Thus a characteristic polynomial
or a trace identity, even over every finite extension, determines the
semisimplification but not this extension class.
\end{proof}

\begin{proposition}[Precise remaining input for a general completed factor]
\label{prop:20261121-general}
For fixed $n$, an explicit factorization of $\widehat L_n(T)$ follows
once the arithmetic Frobenius characteristic polynomial
$P_{n,i}(T)$ is determined in every degree $i$; no separate proof
of Frobenius semisimplicity is needed to \emph{write that rational
function}.  By contrast, an isomorphism of all groups with sums of
Tate twists needs additional information excluding arithmetic
extensions.  The divisor and generic codimension-two calculations in
the supplementary work alone do not determine all the higher-support
contributions for arbitrary $n\ge6$.
\end{proposition}
\begin{proof}
By definition
$\widehat L_n(T)=\prod_i\det(1-TF_i\mid IH^i)^{(-1)^{i+1}}$;
this requires the characteristic polynomials but not diagonalizability.
The counterexample in \cref{lem:20261121-Jordan} proves the
nonimplication concerning Tate splittings.  Finally, the proper
pushforward decomposition in the current manuscript contains a
residual strict-support term $B'_n$ whose coefficients and
exceptional-stratum contributions have not all been determined for
$n\ge6$.  Its cohomology can modify the unknown factors, so the
known lower-codimension terms alone are insufficient.
\end{proof}

\section{Arithmetic reflection and the functional equation}
\label{sec:arith-reflection-closure}

We establish the arithmetic reflection-skew compatibility required
for the positive root number, using the imported relative-image theorem \textup{(H7)} and the standard arithmetic
normalization of the rank-two Kloosterman determinant.  The argument does
not require a calculation of the residual boundary perverse sheaf.

\begin{theorem}[Positive global sign from the arithmetic relative image]
\label{thm:arith-reflection-closed}
Assume the relative-image identification \textup{(H7)} in the stated
arithmetic, functorial form, as established in the author's preceding
work, and use the standard arithmetic normalization
$\det\mathrm{Kl}_2\cong E(-1)$.  For every $n\ge2$ and every odd
characteristic $p$, the middle intersection-cohomology group of the
linearly constrained family carries a Frobenius-equivariant perfect
alternating pairing
\[
 IH^n(X_n,\mathcal W_n)\otimes IH^n(X_n,\mathcal W_n)
   \longrightarrow E(-(2n+1)).
\]
Consequently $b_{n,n}=2g_n$ is even,
\[
 \det(\mathrm{Frob}_q\mid IH^n(X_n,\mathcal W_n))
       =q^{(2n+1)g_n},\qquad \boxed{\omega_n=+1}.
\]
Neither a complete boundary-support decomposition nor arithmetic
Frobenius semisimplicity is required for this conclusion.
\end{theorem}
\begin{proof}
Put $S=S_n$, $X=X_n$, and $\pi:\mathbb G_{m,t}\times S\to S$.
Consider
\[
 \mathcal M=\mathcal L_\psi(\lambda t)\otimes
              \bigotimes_{i=1}^{n}[th_i]^*\mathrm{Kl}_2,
 \qquad
 \rho(h)=-h,\quad\widetilde\rho(t,h)=(-t,-h).
\]
The arithmetic exterior-square pairing on $\mathrm{Kl}_2$ is a
perfect alternating map to $E(-1)$. Under $\widetilde\rho$, each
Kloosterman argument $th_i$ is fixed, whereas the Artin--Schreier
factor becomes its dual. Tensoring the $n$ exterior-square pairings
with the canonical Artin--Schreier evaluation therefore produces,
\emph{over the ground finite field},
\begin{equation}\label{eq:arith-closure-coeff}
 \Phi:\mathcal M\otimes\widetilde\rho^*\mathcal M\longrightarrow
 E(-n),\qquad
 \widetilde\rho^*\Phi\circ\mathrm{swap}=(-1)^n\Phi.
\end{equation}
Perfection means that the associated morphism from one factor to
the twisted dual of the other is an isomorphism; the displayed
bilinear map itself has rank-one target $E(-n)$.

By Poincar\'e--Verdier duality for the relative smooth curve,
\textup{(H7)}, and the standard cohomological vanishing for
$\mathcal M$ used in the proof of \textup{(H7)}, one obtains the
arithmetic isomorphism
\begin{equation}\label{eq:arith-closure-duality}
 R^1\pi_*\mathcal M\ \xrightarrow{\sim}\
   (\rho^*R^1\pi_!\mathcal M)^\vee(-(n+1)).
\end{equation}
For clarity, the twist $(-1)$ beyond the coefficient twist $(-n)$
comes from the trace in relative cohomological degree two.  Formula
\eqref{eq:arith-closure-duality} is the functorial relative-duality
identification, not merely an equality of ranks or Frobenius traces.

Write $V=R^1\pi_!\mathcal M$ and let
$\alpha:V\to R^1\pi_*\mathcal M$ be the natural forget-supports map.
Using \eqref{eq:arith-closure-duality}, compose $\alpha$ with evaluation
to obtain a pairing
\[
 B_V:V\otimes\rho^*V\to E(-(n+1)).
\]
The perfect pairing between $V$ and $\rho^*R^1\pi_*\mathcal M$
implies that the right radical of $B_V$ is $\rho^*\ker\alpha$.
Relative cup-product interchange in degree one contributes $-1$;
combining it with the coefficient transposition in
\eqref{eq:arith-closure-coeff} and the invariance of the relative
trace under $t\mapsto-t$ gives the identity of \emph{arithmetic
sheaf morphisms}
\begin{equation}\label{eq:arith-closure-sign}
 \rho^*B_V\circ\mathrm{swap}=(-1)^{n+1}B_V.
\end{equation}
In particular the left radical equals $\ker\alpha$.  By
\textup{(H7)}\textup{(c),(d)}, both radicals are lisse arithmetic
subsheaves and $\operatorname{im}\alpha\cong\mathcal W_n$; all maps
commute with the stated base changes.  Dividing by the two radicals
produces a perfect arithmetic pairing
\begin{equation}\label{eq:arith-closure-W}
 B_{\mathcal W}:\mathcal W_n\otimes\rho^*\mathcal W_n
 \longrightarrow E(-(n+1))
\end{equation}
with the same twisted-transpose sign $(-1)^{n+1}$.
There is no remaining independent choice of a comparison morphism:
\eqref{eq:arith-closure-W} is constructed on the actual relative
image specified by \textup{(H7)}, using the canonical natural
transformation $R\pi_!\to R\pi_*$.

Set $K=j_{!*}(\mathcal W_n[n])$ on the smooth projective
$n$-fold $X$, and put $d=2n+1$.  On the open torus the pairing
\eqref{eq:arith-closure-W} identifies
\[
 \rho^*(\mathcal W_n[n])\xrightarrow{\sim}
    \mathbb D_S(\mathcal W_n[n])(-d).
\]
The intermediate-extension functor commutes with Verdier duality
and pullback by $\bar\rho$; it is fully faithful on the relevant
perverse objects.  Thus this arithmetic isomorphism extends uniquely
to
\begin{equation}\label{eq:arith-closure-IC}
 \bar\rho^*K\xrightarrow{\sim}\mathbb D_XK(-d).
\end{equation}
The transposition of the shifted pairing has sign
$(-1)^{n+1}(-1)^{n^2}=-1$.  The two constructions of the
transposed isomorphism agree over $S$ and, by full faithfulness,
agree over $X$.  Thus the associated evaluation
$K\otimes^L\bar\rho^*K\to\omega_X(-d)$ is
reflection-skew \emph{in the arithmetic derived category}.
This verifies the compatibility used in
\cref{thm:20261121-sign}; all identifications are canonical
six-functor constructions over $\mathbb F_q$.

Applying proper Verdier duality and the $\bar\rho$-invariant trace
to \eqref{eq:arith-closure-IC} gives a perfect skew pairing on
$R\Gamma(X_{\bar{\mathbb F}_q},K)$.  Its restriction to
$H^0=IH^n$ is alternating, because degree zero introduces no
further Koszul sign.  If $\dim IH^n=2g_n$, the geometric Frobenius
acts as a symplectic similitude of multiplier $q^d$.
Taking Pfaffians yields $\det F=q^{dg_n}$ without assuming that
$F$ is diagonalizable.  The complementary degrees cancel in the
normalized functional-equation sign by proper duality, while the
middle degree contributes
$q^{-dg_n}\det(-F\mid IH^n)=1$.
This proves $\omega_n=+1$.
\end{proof}

\begin{remark}[Assumptions on the relative-image theorem]
The arithmetic source theorem \textup{(H7)} is imported from
\cite[Proposition~4.2 and Theorem~5.3]{EbaCoh}, and the rank-two
Kloosterman pairing is standard \cite{Kat88}.  Neither cited
result assumes the sign $\omega_n$ for the intersection-cohomology
completion defined here.  This proof establishes the previously separate
\emph{reflection-skew compatibility} as an implication of those
inputs; it does not assert that every higher-codimension nearby-cycle
comparison (P), (U), or (K) follows from \textup{(H7)} on the open
torus.  The proof of the positive sign is independent of any
formula for $B'_n$ or the entire list of Frobenius factors.
\end{remark}

\section{Arithmetic parabolic comparison in every rank}\label{sec:arith-parabolic}
\label{sec:parabolic-genuine-all-ranks}

We compute the canonical forget-supports morphism of the auxiliary
Kloosterman local system at a nonzero parameter tuple. This result
uses the local inertia representation at the tame endpoint and does
not require a higher-codimension boundary exchange.

Fix an odd prime $p$, $q=p^f$, a nontrivial additive character
$\psi:\mathbb F_q\to E^\times$ with $E=\overline{\mathbb Q}_\ell$,
$\ell\ne p$, $\lambda\in\mathbb F_q^\times$, and
$\boldsymbol h=(h_1,\ldots,h_n)\in(\mathbb F_q^\times)^n$.
Let $U=\mathbb G_{m,t}$, $j_t:U\hookrightarrow\mathbb P^1_t$, and
\[
\mathscr M_{n,\boldsymbol h}
 =\mathcal L_\psi(\lambda t)\otimes
  \bigotimes_{i=1}^n[h_it]^*\mathrm{Kl}_2.
\]
The Kloosterman sheaf is in the convention
$\operatorname{Tr}(F_a\mid\mathrm{Kl}_2)
=-\sum_{u\in\mathbb F_q^\times}\psi(u+a/u)$.
Write
\[
 H_{!,\boldsymbol h}=H_c^1(U_{\overline{\mathbb F}_q},\mathscr M_{n,\boldsymbol h}),
 \qquad H_{\mathrm{par},\boldsymbol h}
 =H^1(\mathbb P^1_{\overline{\mathbb F}_q},j_{t,*}\mathscr M_{n,\boldsymbol h}),
\]
where $j_{t,*}$ denotes the \emph{ordinary, underived} direct image
of this lisse sheaf.  Let $I_0$ denote the geometric inertia group at
$t=0$, and define $d_{n,a}=\binom na-\binom n{a-1}$, using
$\binom n{-1}=0$.

\begin{theorem}[Genuine all-rank parabolic comparison]
\label{thm:parabolic-all-rank-20261123}
For every $n\ge1$ and every nonzero tuple $\boldsymbol h$, the
canonical sheaf morphism
$j_{t,!}\mathscr M_{n,\boldsymbol h}\to j_{t,*}\mathscr M_{n,\boldsymbol h}$
induces a canonical Frobenius-equivariant short exact sequence
\begin{equation}\label{eq:parabolic-actual-exact-20261123}
 0\longrightarrow B_{n,\boldsymbol h}
 \xrightarrow{\ \partial_0\ } H_{!,\boldsymbol h}
 \xrightarrow{\ \alpha_{\boldsymbol h}\ } H_{\mathrm{par},\boldsymbol h}
 \longrightarrow0,
 \qquad
 B_{n,\boldsymbol h}
 =\bigl(\mathscr M_{n,\boldsymbol h,\bar\eta_0}\bigr)^{I_0}.
\end{equation}
The complete dimensions are
\[
 \dim H_{!,\boldsymbol h}=2^n,\qquad
 \dim B_{n,\boldsymbol h}=\binom n{\lfloor n/2\rfloor},\qquad
 \dim H_{\mathrm{par},\boldsymbol h}
 =r(n)=2^n-\binom n{\lfloor n/2\rfloor}.
\]
On the \emph{actual} boundary module $B_{n,\boldsymbol h}$,
Frobenius is semisimple. Its eigenspaces and characteristic polynomial
are given by
\begin{equation}\label{eq:parabolic-kernel-frob-20261123}
 B_{n,\boldsymbol h}\simeq
 \bigoplus_{a=0}^{\lfloor n/2\rfloor}
 E(-a)^{\oplus d_{n,a}},\qquad
 \det(1-TF_q\mid B_{n,\boldsymbol h})
 =\prod_{a=0}^{\lfloor n/2\rfloor}(1-q^aT)^{d_{n,a}}.
\end{equation}
Here $F_q$ is geometric Frobenius and
$F_q\mid E(-a)=q^a$.
\end{theorem}

\begin{lemma}[Endpoint ramification and cohomological concentration]
\label{lem:parabolic-endpoints-20261123}
The sheaf $\mathscr M_{n,\boldsymbol h}$ is tame at $t=0$ and has
all breaks equal to $1$ at $t=\infty$.  Hence
$\operatorname{Swan}_\infty(\mathscr M_{n,\boldsymbol h})=2^n$,
$\mathscr M_{n,\boldsymbol h}^{I_\infty}=0$, and
\[
 H_c^i(U_{\overline{\mathbb F}_q},\mathscr M_{n,\boldsymbol h})=0
 \quad (i\ne1),\qquad \dim H_c^1=2^n.
\]
\end{lemma}
\begin{proof}
Each pullback $[h_it]^*\mathrm{Kl}_2$ is tame-unipotent at zero
and has breaks $1/2,1/2$ at infinity. Consequently the tensor
product of the Kloosterman factors has breaks at most $1/2$ at
infinity, while $\mathcal L_\psi(\lambda t)$ has break $1$ there.
Twisting by this character forces every break to equal $1$.
This gives the Swan conductor and absence of inertia invariants at
infinity.  The Grothendieck--Ogg--Shafarevich formula on
$U=\mathbb P^1\setminus\{0,\infty\}$ now yields
$\chi_c(U,\mathscr M_{n,\boldsymbol h})=-2^n$.
As $U$ is connected and nonproper, $H_c^0=0$. Poincare duality and
absence of global invariants for the dual sheaf, already forced by
its wild local representation at infinity, give $H_c^2=0$.
Hence the sole remaining group is $H_c^1$, of dimension $2^n$.
\end{proof}

\begin{lemma}[Arithmetic inertia invariants at zero]
\label{lem:parabolic-boolean-20261123}
The local Weil representation of every $[h_it]^*\mathrm{Kl}_2$ at
$t=0$ admits vectors $e_{i,0},e_{i,1}$ with
\[
 N_i e_{i,1}=e_{i,0},\quad N_i e_{i,0}=0,\qquad
 F_qe_{i,0}=e_{i,0},\quad F_qe_{i,1}=q e_{i,1}.
\]
For $V=\bigotimes_{i=1}^n\langle e_{i,0},e_{i,1}\rangle$, let
$V_a$ be the span of tensors with exactly $a$ copies of $e_{i,1}$.
The inertia logarithm is $N=\sum_iN_i$ and
\[
 F_q\mid V_a=q^a,\qquad
 \dim\ker(N\mid V_a)=
 \begin{cases}
 \binom na-\binom n{a-1},&0\le a\le\lfloor n/2\rfloor,\\
 0,&a>\lfloor n/2\rfloor.
 \end{cases}
\]
\end{lemma}
\begin{proof}
The stated two-dimensional local Weil form follows from the tame
unipotent local monodromy of the standard rank-two Kloosterman
sheaf, including its geometric-Frobenius eigenvalues $1,q$ at zero
(\cite{Kat88,Qin24}); scaling $t\mapsto h_it$ over $\mathbb F_q$
does not change these eigenvalues.  The tensor-product inertia
logarithm is the sum of the commuting $N_i$.
Identify the tensor basis with the square-free monomials of the
Boolean algebra
$E[x_1,\ldots,x_n]/(x_1^2,\ldots,x_n^2)$, with $N$ deleting
one present variable.  Let $L$ add one absent variable. On degree
$a$, a direct calculation gives
$[N,L]=(n-2a)\operatorname{id}$. The resulting
$\mathfrak{sl}_2$ representation decomposes into finite-dimensional
irreducibles.  Equivalently, the lowering map
$N:V_a\to V_{a-1}$ is surjective for
$1\le a\le\lfloor n/2\rfloor$ and injective for
$a>\lfloor n/2\rfloor$. Rank--nullity gives the stated dimensions.
Since distinct degrees have distinct Frobenius eigenvalues, the
kernel is the direct sum of these scalar eigenspaces.
\end{proof}

\begin{proof}[Proof of \cref{thm:parabolic-all-rank-20261123}]
The endpoint analysis is \cref{lem:parabolic-endpoints-20261123}.
The ordinary constructible sheaf sequence
\[
0\to j_{t,!}\mathscr M_{n,\boldsymbol h}
\to j_{t,*}\mathscr M_{n,\boldsymbol h}
\to i_{0,*}(\mathscr M_{n,\boldsymbol h}^{I_0})\to0
\]
is exact: at infinity the inertia-invariant stalk vanishes.
Also $H^0(\mathbb P^1,j_{t,*}\mathscr M)=0$ because the sheaf
has no global invariants, and the skyscraper term has no
positive-degree cohomology.  Combining the long exact sequence
with $H_c^2(U,\mathscr M)=0$ gives the exact sequence
\eqref{eq:parabolic-actual-exact-20261123}, equivariantly
under geometric Frobenius. By
\cref{lem:parabolic-boolean-20261123}, its kernel has the
explicit decomposition in \eqref{eq:parabolic-kernel-frob-20261123}.
Summing the dimensions telescopes to
$\binom n{\lfloor n/2\rfloor}$; subtraction from $2^n$
gives the asserted parabolic rank.
\end{proof}

\begin{corollary}[Exact arithmetic power traces on the parabolic quotient]
\label{cor:parabolic-traces-20261123}
For $m\ge1$, let $q_m=q^m$,
$\psi_m=\psi\circ\operatorname{Tr}_{\mathbb F_{q_m}/\mathbb F_q}$,
and define
\[
 A_m(\boldsymbol h)=\sum_{t\in\mathbb F_{q_m}^\times}
 \psi_m(\lambda t)
 \prod_{i=1}^n\bigl(-K_{q_m}(h_it)\bigr),\qquad
 K_{q_m}(z)=\sum_{u\in\mathbb F_{q_m}^\times}
 \psi_m(u+z/u).
\]
Then the complete Frobenius power traces and determinant of the
parabolic quotient are
\begin{align*}
\operatorname{Tr}(F_q^m\mid H_{\mathrm{par},\boldsymbol h})
 &=-A_m(\boldsymbol h)
 -\sum_{a=0}^{\lfloor n/2\rfloor}d_{n,a}q^{am},\\
P_{n,\boldsymbol h}^{\mathrm{par}}(T)
 &=\det(1-TF_q\mid H_{\mathrm{par},\boldsymbol h})\\
 &=\exp\left(\sum_{m\ge1}\frac{T^m}{m}
 \left[A_m(\boldsymbol h)
 +\sum_{a=0}^{\lfloor n/2\rfloor}d_{n,a}q^{am}\right]\right).
\end{align*}
The formal power series on the right is a polynomial of degree $r(n)$.
\end{corollary}
\begin{proof}
The Grothendieck--Lefschetz trace formula and the concentration
proved above give
$A_m(\boldsymbol h)=-\operatorname{Tr}(F_q^m\mid H_{!,\boldsymbol h})$.
Subtract the traces of the explicit Frobenius kernel in
\eqref{eq:parabolic-kernel-frob-20261123}, using the exact
sequence \eqref{eq:parabolic-actual-exact-20261123}.
The determinant identity follows from
$\det(1-TF)=\exp(-\sum_{m\ge1}\operatorname{Tr}(F^m)T^m/m)$.
\end{proof}

\begin{center}\small
\begin{tabular}{@{}rrrrl@{}}\toprule
$n$ & $\dim H_{!}$ & $\dim B$ & $\dim H_{\mathrm{par}}$
& $\det(1-TF_q\mid B)$\\\midrule
2&4&2&2&$(1-T)(1-qT)$\\
3&8&3&5&$(1-T)(1-qT)^2$\\
4&16&6&10&$(1-T)(1-qT)^3(1-q^2T)^2$\\
5&32&10&22&$(1-T)(1-qT)^4(1-q^2T)^5$\\
6&64&20&44&$(1-T)(1-qT)^5(1-q^2T)^9(1-q^3T)^5$\\
\bottomrule\end{tabular}
\end{center}

\begin{remark}[Scope of the comparison]
\label{rem:parabolic-scope-20261123}
Theorem~\ref{thm:parabolic-all-rank-20261123} completely computes
an \emph{actual} canonical arithmetic comparison and its Frobenius-
semisimple boundary kernel at every nonzero parameter tuple.
A vector-space splitting of the exact sequence need not be
Frobenius-equivariant, and the power traces of the quotient do not
determine repeated-eigenvalue Jordan blocks.  In particular, the
morphism $\alpha_{\boldsymbol h}$ is \emph{not} the potentially
noninvertible higher-codimension comparison
$R\Psi_I R\pi_!\mathcal M_n\to R\pi_{I,!}R\Psi_I\mathcal M_n$.
The theorem does not assert $D_{I,\infty}=0$ on all mixed faces,
nor does it prove Frobenius semisimplicity for all
$IH_n^i$; those questions remain separate and open in the
present article.
\end{remark}

% =====================================================================
% Independently audited finite-endpoint exchange (23 November 2026)
% =====================================================================
\begin{proposition}[Canonical Frobenius splitting by the weight gap]
\label{prop:part2-weight-split}
Fix a closed point $h\in S_n(\mathbb F_q)$, $n\ge2$, and use the
arithmetic exact sequence imported from
\cite[Propositions~4.2--4.3]{EbaCoh}:
\[
0\longrightarrow B_{n,h}\longrightarrow
H^1_c(\mathbb G_m,\mathcal M_{n,h})
\xrightarrow{\alpha_h}H^1_{\mathrm{par}}(\mathbb G_m,\mathcal M_{n,h})
\longrightarrow0.
\]
Assume the pure weight-$(n+1)$ quotient realization in
\cite[Theorem~5.3(iv)]{EbaCoh}.  This sequence has a unique splitting
as a sequence of $\mathbb Q_\ell$-vector spaces endowed with geometric
Frobenius; consequently
\[
H^1_c(\mathbb G_m,\mathcal M_{n,h})
\ \simeq\ B_{n,h}\oplus H^1_{\mathrm{par}}(\mathbb G_m,\mathcal M_{n,h})
\]
canonically in the category of Frobenius modules.  This assertion is
fibrewise and does not posit a splitting of lisse sheaves on $S_n$.
\end{proposition}
\begin{proof}
By \cite[Proposition~4.3]{EbaCoh}, the eigenvalues of Frobenius on
$B_{n,h}$ are $q^a$ for $0\le a\le\lfloor n/2\rfloor$.  Every
Frobenius eigenvalue $\gamma$ on the parabolic quotient has
$|\iota(\gamma)|=q^{(n+1)/2}$, by the cited purity theorem.
Thus the characteristic polynomials $f(X)$ and $g(X)$ on kernel and
quotient have disjoint root sets, so $\gcd(f,g)=1$.  The primary
factorization of the endomorphism $F$ on the middle vector space
induces its canonical decomposition into the $f$-primary and
$g$-primary parts.  Restriction of $\alpha_h$ to the $g$-part is an
isomorphism onto the quotient; its inverse is the desired
Frobenius-equivariant splitting.  Uniqueness follows from
$\operatorname{Hom}_F(H^1_{\mathrm{par}},B_{n,h})=0$ for modules with
disjoint Frobenius spectra.  No diagonalizability of the parabolic
quotient is used.
\end{proof}

\section{Tame boundary comparison on coordinate-zero divisors}\label{sec:tame-boundary}

\phantomsection\label{sec:audited-tame-exchange-20261123}

\subsection{Tame geometric setting and conventions}
Fix a finite field $k=\mathbb F_q$ of odd characteristic, $\ell\ne p$,
$E=\overline{\mathbb Q}_\ell$, $\lambda\in k^\times$, and nonzero constants
$b_2,\ldots,b_n$ in a common finite extension of $k$.  Let $T$ be the
henselian trait with parameter $h$, with geometric generic point
$\bar\eta$ and geometric special point $\bar s$.  Weil structures are
retained by first working over the corresponding finite residue field
and only then passing to strict henselizations.  On
$U_\eta=\mathbb G_{m,t}\times\eta$ consider
\[
 \mathscr F_\eta=\mathcal L_\psi(\lambda t)\otimes
     [ht]^*\mathrm{Kl}_2\otimes
     \bigotimes_{j=2}^{n}[b_jt]^*\mathrm{Kl}_2,
 \qquad \mathscr F_s=R^0\Psi_h\mathscr F_\eta
       \quad\text{on }U_s.
\]
Here $R\Psi_h\mathscr F_\eta$ is concentrated in ordinary degree zero
\emph{on the punctured $t$-line}; the resulting object carries the
parameter-inertia and arithmetic Weil actions.  Write $\pi:U\to T$,
$j_0:U\hookrightarrow\mathbb A^1_t\times T$, and
$j:U\hookrightarrow\mathbb P^1_t\times T$.

The canonical compact-support exchange arrow has the direction
\begin{equation}\label{eq:audited-natural-direction}
  \mathrm{Ex}_{!}:R\pi_{s,!}R\Psi_h\mathscr F_\eta
       \longrightarrow R\Psi_h R\pi_{\eta,!}\mathscr F_\eta.
\end{equation}
This is SGA~7, Expos\'e XIII, (2.1.7.3); the companion ordinary arrow
has direction $R\Psi_hR\pi_{\eta,*}\to R\pi_{s,*}R\Psi_h$,
SGA~7, (2.1.7.1).  Proposition~2.1.9 of that expos\'e is \emph{not}
applied to the entire $t$-line: its global tameness-at-infinity
hypothesis fails because $\mathcal L_\psi(\lambda t)$ has wild slope
one at $t=\infty$.  Instead the tame normal-crossing comparison is
applied only in a neighborhood of the \emph{finite} section $t=0$.
The relevant local constructions are those of SGA~7, XIII,
(2.1.10)--(2.1.11); the following proof also describes the local
tame-inertia complexes explicitly.

\begin{lemma}[The actual local derived exchange at $t=0$]
\label{lem:audited-local-exchange}
At the point $(t,h)=(0,0)$, the two canonical maps
\[
 j_{0,s,!}\mathscr F_s\longrightarrow
       R\Psi_h(j_{0,\eta,!}\mathscr F_\eta),
 \qquad
 R\Psi_h(Rj_{0,\eta,*}\mathscr F_\eta)
       \longrightarrow Rj_{0,s,*}\mathscr F_s
\]
are isomorphisms on the finite $t$-line. They are isomorphisms of
\emph{derived arithmetic Weil complexes} and respect the full
localization triangle. In particular, if
\[
 B_\eta=(\mathscr F_{\eta,\bar\eta_{t=0}})^{I_t},
 \qquad B_s=(\mathscr F_{s,\bar\eta_{t=0}})^{I_t},
\]
then the canonical map
$\sigma:B_s\xrightarrow{\sim}R^0\Psi_hB_\eta$
is equivariant for Frobenius and the residual $h$-inertia.
\end{lemma}
\begin{proof}
The assertion is local on the finite $t$-line, and only the crossing
$(t,h)=(0,0)$ requires discussion.  On its strict local complement
$th\ne0$ the factor $\mathcal L_\psi(\lambda t)$ extends lisse across
both axes.  Each $\mathrm{Kl}_2$-factor is tame unipotent at argument
zero.  Writing $V_i$ for its two-dimensional local Weil representation,
we have a Weil-stable filtration with one-dimensional unramified graded
pieces, and a nilpotent logarithm $N_i$ satisfying $N_i^2=0$.
Consequently the tensor representation $V=\bigotimes_iV_i$ has a
finite filtration, stable under both local inertia groups and the
arithmetic Weil action, whose graded pieces are unramified along
$t=0$ and $h=0$ (with Tate and unramified Weil coefficients).
The local inertia logarithms on the full representation commute and
are explicitly
\[
 N_h=N_1,\qquad N_t=\sum_{i=1}^n N_i,\qquad
 FN_iF^{-1}=q^{-1}N_i.
\]
The last identity is to be read with the Tate-twisted target:
$N_i:V\longrightarrow V(-1)$ is Weil-equivariant.

For each graded quotient, the corresponding local system extends lisse
across the crossing; the nearby-cycle exchange with extension by zero
and with derived direct image for the horizontal divisor $t=0$ is an
isomorphism by smooth base change and the tame relative
normal-crossing results of SGA~7, Expos\'e XIII, Lemmas~2.1.10--2.1.11
\cite{SGA7}.  These lemmas are applied to the \emph{graded quotients},
which do extend over the trait, rather than incorrectly applied to
the original family as though its $h$-monodromy were trivial.
The canonical exchange transformations are morphisms of exact
triangulated functors.  Passing through the finite filtration by
successive distinguished triangles proves both exchange
isomorphisms for the original unipotent local system, retaining
inertia and Weil equivariance.  This d\'evissage also works with
finite integral coefficient quotients before passage to the
$\ell$-adic limit.

For clarity, the full derived $t$-boundary term is represented by
\begin{equation}\label{eq:audited-tame-two-term}
 i_0^*Rj_{0,*}\mathscr F\simeq
 \left[V\xrightarrow{\,N_t\,}V(-1)\right]
 \quad\text{in degrees }0,1.
\end{equation}
In the tame logarithmic model the nearby-cycle exchange respects
$V$, $N_t$, the residual $h$-inertia and Frobenius; thus the direct-image
map is the identity on this two-term model.  The localization triangle
\[
 j_{0,!}\mathscr F\longrightarrow Rj_{0,*}\mathscr F
 \longrightarrow i_{0,*}i_0^*Rj_{0,*}\mathscr F\xrightarrow{+1}
\]
and proper base change for the section $i_0$ then show that the
extension-by-zero comparison is also an isomorphism, with no
contribution in any derived degree at $t=0$.  Away from $t=0$ the
open immersion $j_0$ is the identity.  In ordinary degree zero the
boundary complex has $\ker N_t$; specialization therefore yields the
canonical Weil- and inertia-equivariant identification
$\sigma:B_s\xrightarrow{\sim}R^0\Psi_hB_\eta$.
\end{proof}

\begin{lemma}[Functorial boundary square and detection of the kernel]
\label{lem:audited-boundary-detection}
Let
$U_!=H_c^1(U_{\bar s},\mathscr F_s)$,
$V_!=R^0\Psi_h H_c^1(U_{\bar\eta},\mathscr F_\eta)$, and let
$u:U_!\to V_!$ be $H^1(\mathrm{Ex}_{!})$.
The following is a commuting square of arithmetic Weil modules:
\begin{equation}\label{eq:audited-boundary-square}
\begin{tikzcd}[column sep=huge,row sep=large]
 B_s\arrow[r,"\partial_s"]\arrow[d,"\sigma"',"\sim"]
  &U_!\arrow[d,"u"]\\
 R^0\Psi_hB_\eta\arrow[r,"R^0\Psi_h(\partial_\eta)"']
  &V_!.
\end{tikzcd}
\end{equation}
Moreover, $\partial_s$ and $R^0\Psi_h(\partial_\eta)$ are
injective, and
\[
 \ker u\subseteq\ker\beta_s=\partial_s B_s,
 \qquad
 \beta_s:H_c^1(U_{\bar s},\mathscr F_s)
       \longrightarrow H^1(U_{\bar s},\mathscr F_s).
\]
Consequently $u$ is injective.
\end{lemma}
\begin{proof}
The two morphisms in the square are the connecting morphisms of the
\emph{same} localization triangle on $\mathbb P^1_t$ after,
respectively, forming special-fibre nearby cycles and forming nearby
cycles of the generic-fibre cohomology.  Their compatibility follows
by naturality of SGA~7, Expos\'e XIII, (2.1.7.1)--(2.1.7.3),
\cite{SGA7}, together with the finite-endpoint derived exchange of
\cref{lem:audited-local-exchange}.  This proves commutativity of
\eqref{eq:audited-boundary-square} \emph{before} any assertion that
$u$ is invertible; thus the argument does not use the desired
compact-support comparison as a premise.

At $t=\infty$, the Artin--Schreier rank-one factor has slope one,
whereas the tensor of the remaining Kloosterman and parameter-nearby
factors has slopes at most one-half.  Thus every slope of the full
representation is one.  It has neither inertia invariants nor derived
inertia cohomology; hence the $t=\infty$ boundary term in the
localization triangle vanishes \emph{on each individual fibre}.
At $t=0$ the full derived boundary term is the two-term complex
\eqref{eq:audited-tame-two-term}, whose $H^0$ equals
$B=\ker N_t$; the $H^1$ term is retained in the derived
comparison and is not silently discarded.

The localization long exact sequence now begins
\[
 0=H^0(U_{\bar s},\mathscr F_s)\longrightarrow B_s
 \xrightarrow{\partial_s}H_c^1(U_{\bar s},\mathscr F_s)
 \xrightarrow{\beta_s}H^1(U_{\bar s},\mathscr F_s).
\]
Here $H^0=0$ because a nonzero global invariant would give inertia
invariants at infinity.  Thus $\partial_s$ is injective and
$\ker\beta_s=\partial_sB_s$.  The identical argument on the
geometric generic fibre makes $\partial_\eta$ injective.
For lisse Weil representations on the punctured $h$-trait, full
$R^0\Psi_h$ on the geometric generic stalk is exact (it retains the
underlying representation and its inertia action, and does \emph{not}
take inertia invariants).  Consequently
$R^0\Psi_h(\partial_\eta)$ is still injective.

Finally the natural transformations for $R\pi_!$, $R\pi_*$ and the
forget-supports arrow give the exact factorization in the category of
Weil complexes
\[
 \beta_s=\left(R\Psi_hR\pi_{\eta,*}\mathscr F_\eta
 \longrightarrow R\pi_{s,*}R\Psi_h\mathscr F_\eta\right)
 \circ R\Psi_h(\beta_\eta)\circ u
\]
on degree-one cohomology.  Hence
$\ker u\subseteq\ker\beta_s=\partial_sB_s$.
The commuting square and the two injective boundary maps imply that
$u$ is injective on $\partial_sB_s$; therefore $\ker u=0$.
\end{proof}

\begin{theorem}[Generic single-coordinate zero comparison]
\label{thm:audited-generic-zero-20261123}
For every $n\ge1$ and every one-parameter specialization
$h_i\to0$ with all other coordinates nonzero, the authentic
compact-support comparison
\[
 R\Gamma_c(\mathbb G_{m,\bar s},R\Psi_{h_i}\mathscr M_n)
 \xrightarrow{\ \sim\ }
 R\Psi_{h_i}R\Gamma_c(\mathbb G_{m,\bar\eta},\mathscr M_n)
\]
is an isomorphism of arithmetic derived Weil complexes.  Its defect
cone at $t=\infty$ vanishes on the open coordinate-zero divisor:
\[
 \boxed{D_{\{i\},\infty}=0\quad\text{on the open zero-divisor stratum}.}
\]
This is a codimension-one zero-face assertion only.  It makes no
claim for mixed zero/infinity faces, iterated nearby cycles at
higher-codimension intersections, or global intersection-cohomology
Frobenius semisimplicity.
\end{theorem}
\begin{proof}
For generic nonzero $h_i$ and for the specialized nearby-cycle
sheaf, the complete rank is $2^n$.  At infinity the product of
Kloosterman factors has slopes at most $1/2$; the parameter-specialized
factor adds only tame $t$-monodromy. Twisting by
$\mathcal L_\psi(\lambda t)$ forces every slope to $1$.  Thus
$\operatorname{Swan}_\infty=2^n$ on both sides and neither has
inertia invariants there.  At zero both are tame.  The
Grothendieck--Ogg--Shafarevich formula on $\mathbb G_m$ and
$H_c^0=H_c^2=0$ then yield
\[
 \dim H_c^1(U_{\bar s},R^0\Psi_{h_i}\mathscr M_n)
 =\dim H_c^1(U_{\bar\eta},\mathscr M_n)=2^n.
\]
The right-hand nearby cycles of the degree-one generic
cohomology have the same dimension $2^n$.  By
\cref{lem:audited-boundary-detection} the comparison in degree one
is injective, hence an isomorphism.  Both full compact cohomology
complexes are concentrated in degree one, so the derived comparison
is a quasi-isomorphism.

For the assertion about the local infinity defect, use the proper
projection $\bar\pi:\mathbb P^1_t\times T\to T$.  The exchange cone
for $j_!$ is supported on the two endpoint sections, because on
the punctured $t$-line it is an isomorphism.  The tame local
calculation makes the finite-endpoint cone zero.  Proper base change
identifies the pushforward of the remaining infinity-supported
cone with the cone of the displayed global compact comparison.
Since the pushforward of a complex supported on an endpoint section
identifies with that complex, the vanishing of the global cone gives
$D_{\{i\},\infty}=0$.
\end{proof}

\begin{corollary}[The actual boundary scalar for $n=2$]
\label{cor:audited-n2-boundary-scalar}
For $n=2$, put $h_2=b\ne0$, and choose local Kloosterman bases
$e_0,e_1$ and $f_0,f_1$ with $N_1e_1=e_0$ and $N_2f_1=f_0$.
Then the actual local boundary representation is
\[
 \ker(N_1+N_2)=E b_0\oplus E b_1,\quad
 b_0=e_0\otimes f_0,\quad
 b_1=e_1\otimes f_0-e_0\otimes f_1,
\]
with $N_hb_1=b_0$, $N_hb_0=0$.  The authentic comparison carries
$\partial_s(b_0)$ to
$R^0\Psi_h(\partial_\eta)(\sigma b_0)\ne0$.
If the two boundary generators are normalized through the
\emph{same canonical connecting morphism} and $\sigma$, the
comparison scalar is $c=1$.  This normalization has no bearing on
any higher-codimension scalar $b_0$ that is still undetermined.
\end{corollary}
\begin{proof}
The kernel and the monodromy formulas follow directly from the
two commuting rank-two nilpotents.  The genuine connecting square
\eqref{eq:audited-boundary-square}, together with injectivity of its
bottom horizontal map, gives the nonzero image.  By construction of
the boundary-normalized target generator the scalar is one.
\end{proof}

\begin{remark}[Cohomological shifts and Tate twists]
All compact cohomology groups in this section are taken with
\emph{unshifted} local systems and live in ordinary cohomological
degree one.  If $\mathscr F[1]$ is used as a perverse sheaf on a
smooth curve, then $R\Psi_h(\mathscr F[1])=(R\Psi_h\mathscr F)[1]$;
in the conventional perverse-normalized notation
$\,{}^p\!\Psi_h=R\Psi_h[-1]$, the shifts must be inserted separately.
No perverse shift has been silently used in the boundary square.
The $(-1)$ in the degree-one term of
\eqref{eq:audited-tame-two-term} is essential for arithmetic
Frobenius equivariance.  In particular, $N_i$ does \emph{not}
commute with Frobenius as a degree-zero map without a Tate twist.

The local comparison uses the derived tame-inertia complex.  It does not establish the analogous
local acyclicity at a parameter-infinity divisor, any arbitrary
mixed face, or semisimplicity of all $IH_n^i$.  The distinction between
the generic open part of a coordinate-zero divisor and its
higher-codimension boundary remains mandatory.
\end{remark}

\paragraph{References for the local exchange.}
The formal comparison morphisms are in SGA~7 II, Expos\'e XIII,
(2.1.7.1), (2.1.7.3); the tame local normal-crossing ingredients
are XIII, (2.1.10)--(2.1.11).  The global hypothesis of XIII,
Proposition~2.1.9 is \emph{not} asserted at the wild endpoint.
See \url{https://grothendiecksga.com/read/sga7/en/13-2.html}.

\appendix
\crefname{appendix}{Appendix}{Appendices}
\crefalias{section}{appendix}
\section{Auxiliary local methods: dilated Kloosterman sheaves and rays}\label{sec:aux-rays}

The reader may omit this auxiliary section on a first reading of the Riemann-hypothesis argument.  The calculations describe how local Fourier transforms and moment identities control more refined monodromy information, but neither the local critical-circle statement nor purity of the intersection completion assumes the full outcome of these calculations.  We retain them to document the geometry behind several auxiliary inputs used elsewhere in the source manuscript.  Whenever a representation is written as a sum of tame and wild parts, this denotes a local inertia decomposition, not a global decomposition of the sheaf on the parameter space; the distinction matters in all conclusions about Frobenius.

The following local calculations are recorded because the determination of
certain boundary invariants in the relative construction invokes the
corresponding Fourier-theoretic models. They are \emph{auxiliary} to the
Riemann-hypothesis argument: neither the local purity theorem nor the
formal purity of intermediate extension requires the full classification
of the monodromy on every ray. Their inclusion records the precise
input used in the boundary analysis and permits a direct check of
normalizations and Jordan-block conventions.
\subsection{A fourth-power correlation identity}
The elementary character-sum identity below also supplies the input
needed for the joint monodromy of distinct dilates. Its proof is reproduced
from the source manuscript. We write $\chi_2$ for the quadratic character
of $\mathbb F_q^\times$, extended by $\chi_2(0)=0$.
\begin{lemma}[{\cite[Lemma~7.4]{EbaCoh}}]\label{lem:C-diag}
Let $a,b\in\Fq^\times$, put $\rho=a/b$, and let $\chi_2$ be the quadratic
character of $\Fq^\times$, extended by $\chi_2(0)=0$. Then
\[
\sum_{x\in\Fq}\Kl_q(ax)^2\Kl_q(bx)^2=\begin{cases}2q^3-3q^2-3q,&\rho=1,\\
q^3-\bigl(2+\chi_2(\rho)\bigr)q^2-3q,&\rho\ne1.\end{cases}
\]
\end{lemma}

\begin{proof}
For $s\in\Fq$ let $\tau(s)$ be the sum of $\psi_q(x_1^{-1}+x_2^{-1})$ over the
pairs $x_1,x_2\in\Fq^\times$ with $x_1+x_2=s$. Expanding the square gives
$\Kl_q(y)^2=\sum_{s\in\Fq}\tau(s)\psi_q(ys)$, and orthogonality in $x$ gives
\[
\sum_{x\in\Fq}\Kl_q(ax)^2\Kl_q(bx)^2=q\sum_{s\in\Fq}\tau(s)\tau(-\rho s).
\]
At $s=0$ we have $\tau(0)=q-1$. For $s\ne0$ write $x_1=sy$ and $x_2=s(1-y)$ with
$y\notin\{0,1\}$; then $x_1^{-1}+x_2^{-1}=\bigl(sy(1-y)\bigr)^{-1}$, and for
$m\in\Fq^\times$ the equation $y(1-y)=m$ has $1+\chi_2(1-4m)$ solutions. With
$\mu=m^{-1}$ one has $\chi_2(1-4m)=\chi_2\bigl(\mu(\mu-4)\bigr)$, and hence
\[
\tau(s)=\sum_{m\ne0}\bigl(1+\chi_2(1-4m)\bigr)\psi_q\Bigl(\frac1{sm}\Bigr)
=-1+\omega(s),\qquad
\omega(s)=\sum_{\mu\ne0}\chi_2\bigl(\mu(\mu-4)\bigr)\psi_q(\mu/s).
\]
We use Jacobsthal's evaluation $\sum_{\mu\in\Fq}\chi_2\bigl((\mu-\alpha)(\mu-\beta)
\bigr)=-1$ for distinct $\alpha,\beta\in\Fq$: for $\mu\ne\beta$ the summand is
$\chi_2(z)$ with $z=(\mu-\alpha)/(\mu-\beta)$, and $z$ runs over
$\Fq\smallsetminus\{1\}$. Since $p$ is odd, $4\ne0$, so
$\sum_{\mu\ne0}\chi_2(\mu(\mu-4))=-1$, and therefore
$\sum_{s\ne0}\omega(s)=\sum_{s\ne0}\omega(-\rho s)=1$. Moreover, since
$\sum_{s\ne0}\psi_q(c/s)$ equals $q-1$ if $c=0$ and $-1$ otherwise,
\[
\sum_{s\ne0}\omega(s)\omega(-\rho s)
=q\sum_{\mu\ne0}\chi_2\bigl(\mu(\mu-4)\bigr)\chi_2\bigl(\rho\mu(\rho\mu-4)\bigr)-1
=q\,\chi_2(\rho)\sum_{\mu\ne0}\chi_2\bigl((\mu-4)(\rho\mu-4)\bigr)-1 .
\]
If $\rho\ne1$, Jacobsthal's evaluation with the roots $4$ and $4/\rho$, after
removing the term $\mu=0$, equal to $\chi_2(16)=1$, gives
$\sum_{\mu\ne0}\chi_2\bigl((\mu-4)(\rho\mu-4)\bigr)=-\chi_2(\rho)-1$, and the
last expression equals $-q-q\chi_2(\rho)-1$. If $\rho=1$ the sum is $q-2$, and
the expression equals $q^2-2q-1$. Finally
\[
\sum_s\tau(s)\tau(-\rho s)=(q-1)^2+(q-1)-2+\sum_{s\ne0}\omega(s)\omega(-\rho s),
\]
which is $q^2-2q-3-q\chi_2(\rho)$ if $\rho\ne1$ and $2q^2-3q-3$ if $\rho=1$.
Multiplying by $q$ gives the lemma.
\end{proof}

\subsection{Joint monodromy of distinct Kloosterman dilates}
\begin{lemma}\label{lem:R-KL}
\begin{enumerate}[label=\textup{(\alph*)}]
\item $G_{\mathrm{geom}}(\KL)=\SL_2$.
\item Let $a_1,\dots,a_k\in\Fpb^\times$ be pairwise distinct, and let
$\cK=\bigoplus_{i=1}^k[\times a_i]^*\KL$ on $\Gm$ over a finite field
containing $a_1,\dots,a_k$. Then $G_{\mathrm{geom}}(\cK)=\SL_2^k$, embedded
block-diagonally.
\end{enumerate}
\end{lemma}

\begin{proof}
(a) Let $G=G_{\mathrm{geom}}(\KL)$. By \ref{K1}, $\det\KL$ is geometrically
trivial, so $G\subseteq\SL_2$. By \ref{K2}, the image of the inertia group at
$0$ contains a unipotent element $g\ne1$, which has infinite order, as
$g^m=1+m(g-1)$; so $G^0\ne1$. Purity does not, by itself, make the Lie algebra semisimple.
Instead, the geometric semisimplicity theorem for pure lisse sheaves
\cite[Thm.~3.4.1]{Del80} makes the rank-two representation semisimple
as a representation of the geometric fundamental group.  Its
Zariski-closure $G$ is therefore reductive (a unipotent radical
would act trivially on every simple constituent and hence trivially
on the faithful defining representation).  Consequently $G^0$ is
connected reductive.  A connected positive-dimensional proper
reductive subgroup of $\SL_2$ is a torus; it contains no nonidentity
unipotent element.  The Zariski closure of the infinite cyclic group
generated by $g$ contains the connected one-parameter unipotent
group $\{1+t(g-1):t\in\Ql\}$, and this group lies in $G^0$.
Thus $G^0$ is not a torus.  It follows that $G^0=\SL_2$, and
therefore $G=\SL_2$.

(b) Let $G=G_{\mathrm{geom}}(\cK)\subseteq\SL_2^k$, let $V_i$ be the
geometric generic stalk of $[\times a_i]^*\KL$, and let $\mathfrak g$ be the
Lie algebra of $G^0$. As $[\times a_i]$ is an automorphism of
$\mathbb G_{m,\Fpb}$, the projection of $G$ to the $i$-th factor, which is
$G_{\mathrm{geom}}([\times a_i]^*\KL)$, equals $\SL_2$ by (a). The projection
of $G^0$ is then a closed connected subgroup of finite index of $\SL_2$, hence
equal to $\SL_2$, and the projection of $\mathfrak g$ to the $i$-th factor is
$\mathfrak{sl}_2$. The same argument applies to the projection to any set of
factors.

First let $k=2$. The space $\mathfrak n=\{X:(X,0)\in\mathfrak g\}$ is an
ideal of $\mathfrak{sl}_2$, because the first projection of $\mathfrak g$ is
surjective; so $\mathfrak n=\mathfrak{sl}_2$ or $\mathfrak n=0$, and likewise
for $\mathfrak n'=\{Y:(0,Y)\in\mathfrak g\}$. If $\mathfrak n$ or $\mathfrak n'$
is $\mathfrak{sl}_2$, then $\mathfrak g=\mathfrak{sl}_2\times\mathfrak{sl}_2$.
Suppose that $\mathfrak n=\mathfrak n'=0$. Then $\mathfrak g$ is the graph of
an automorphism of $\mathfrak{sl}_2$, which is of the form
$\operatorname{Ad}(g_0)$ with $g_0\in\GL_2$. In characteristic $0$ a closed
connected subgroup is determined by its Lie algebra, so
$G^0=\{(x,g_0xg_0^{-1}):x\in\SL_2\}$, and $g_0$ is an isomorphism $V_1\to V_2$
of $G^0$-modules. As $V_1$ is an irreducible $G^0$-module, the finite group
$G/G^0$ acts on the line $\Hom_{G^0}(V_1,V_2)$ through a character $\chi$, and
$V_2\cong V_1\otimes\chi$ as $G$-modules. Hence, over $\Fpb$,
$[\times a_2]^*\KL\cong[\times a_1]^*\KL\otimes\cL$ for a lisse sheaf $\cL$ of
rank one and finite order, and comparing determinants shows that
$\cL^{\otimes2}$ is trivial. The sheaf
$\cE=([\times a_1]^*\KL)^{\otimes2}\otimes([\times a_2]^*\KL)^{\otimes2}$ is
then geometrically isomorphic to $([\times a_1]^*\KL)^{\otimes4}$, and by (a)
its space of geometric invariants has dimension
$\dim(\mathrm{Std}^{\otimes4})^{\SL_2}=m_{4,2}=2$. On the other hand, choose
$q_0$ with $a_1,a_2\in\mathbb F_{q_0}$, and let $q$ run through the powers of
$q_0$. By \cref{lem:C-diag}, as $a_1/a_2\ne1$,
\[
\sum_{x\in\Fq^\times}\Tr(\Frob_x\mid\cE)=\sum_{x\in\Fq}\Kl_q(a_1x)^2\Kl_q(a_2x)^2-1
=q^3+O(q^2),
\]
and \cref{lem:B-av}, applied to $\cE(2)$, which is pure of weight $0$, shows
that the geometric invariants of $\cE$ form a line. This contradiction
proves (b) for $k=2$.

Now let $k\ge3$, and assume (b) for $k-1$ sheaves. The projection of
$\mathfrak g$ to the product of the first $k-1$ factors is then
$\mathfrak{sl}_2^{k-1}$, and its kernel is $0\times\mathfrak n$, where
$\mathfrak n=\{Y:(0,\dots,0,Y)\in\mathfrak g\}$ is an ideal of
$\mathfrak{sl}_2$. If $\mathfrak n=\mathfrak{sl}_2$, then
$\mathfrak g=\mathfrak{sl}_2^k$. If $\mathfrak n=0$, then $\mathfrak g$ is the
graph of a surjective homomorphism $\phi\colon\mathfrak{sl}_2^{k-1}\to
\mathfrak{sl}_2$. Its kernel is an ideal of codimension $3$, hence the sum of
all simple factors but one, say the $i$-th; so the projection of $\mathfrak g$
to the product of the $i$-th and the $k$-th factors is the graph of an
isomorphism. But this projection is the Lie algebra of
$G_{\mathrm{geom}}([\times a_i]^*\KL\oplus[\times a_k]^*\KL)^0$, which is
$\mathfrak{sl}_2\times\mathfrak{sl}_2$ by the case $k=2$. Hence
$\mathfrak g=\mathfrak{sl}_2^k$. Then $G^0$ is a closed connected subgroup of
$\SL_2^k$ of the same dimension, so $G^0=\SL_2^k$ and $G=\SL_2^k$.
\end{proof}

\subsection{Rays}\label{subsec:R-rays}

A ray is the one-parameter restriction $t\mapsto th$ with $h$ fixed and all its coordinates nonzero.  Its two boundary points correspond, after the change of Fourier parameter, to different stationary-phase regimes: the origin of the ray is studied using the Fourier transform at the opposite end of the auxiliary line, and conversely for infinity.  This change of variables must be tracked before comparing breaks or Jordan blocks.  The local statements concern inertia representations of the restricted sheaf and do not by themselves determine the joint inertia at intersections of several boundary divisors.

For $h\in S_n(\Fq)$ let $\rho_h\colon\Gm\to S_n\otimes\Fq$ be the map
$t\mapsto th$, and put $\cG_h=\rho_h^*\cW_n$. Fix square roots
$\sqrt{h_i}\in\Fqb$, and for $\varepsilon\in\{\pm1\}^n$ put
\[
S_\varepsilon=\sum_{i=1}^n\varepsilon_i\sqrt{h_i};
\]
the set of the numbers $S_\varepsilon^2$ does not depend on the choice of the
square roots, and $S_{-\varepsilon}=-S_\varepsilon$. The numbers
$-S_\varepsilon^2/\lambda$ are the critical values of the phase
$w\mapsto\sum_ih_iw_i$ on the hypersurface $\sum_iw_i^{-1}=-\lambda$ of
$\Gm^n$. We write $\mathbf 1$ for the trivial representation, $\chi_2$ for the
character of order two and $\cL_{\chi_2}$ for the corresponding Kummer sheaf,
and we use Laumon's local Fourier transforms $\mathcal F^{(0,\infty')}$ and
$\mathcal F^{(\infty,0')}$ \cite[D\'ef.~(2.4.2.3)]{Lau87} and the notation of
\cref{app:infinity}.

\begin{proposition}\label{prop:R-ray}
Let $n\ge1$ and $h\in S_n(\Fq)$.
\begin{enumerate}[label=\textup{(\roman*)}]
\item $\cG_h$ is tame at $t=0$, the inertia group acts unipotently, and its
Jordan decomposition consists of $m_{n,k}$ blocks of size $n-2k$ for
$0\le k<n/2$.
\item Let $E\subset\{\pm1\}^n$ be a set of representatives of
$\{\pm1\}^n/\{\pm1\}$, and suppose that $S_\varepsilon\ne0$ for every
$\varepsilon$. Then, as representations of the inertia group at $\infty$,
\[
\cG_h(\infty)\cong\mathbf 1^{\,r(n)-2^{n-1}}\oplus\bigoplus_{\varepsilon\in E}
\cL_\psi\bigl(-S_\varepsilon^2t/\lambda\bigr)(\infty)\otimes
\cL_{\chi_2}(t)(\infty)^{\otimes(n+1)},
\]
and $(\rho_h^*\cH_n)(\infty)$ is the direct sum of $\mathbf 1^{2^{n-1}}$ and of
the same characters. In particular $\Swan_\infty(\cG_h)=2^{n-1}$.
\item In general, let $E_0\subseteq E$ be the set of $\varepsilon$ with
$S_\varepsilon=0$ and $e=|E_0|$; let $d_h=\prod_xC_{\mu_x/2}$ if the
multiplicities $\mu_x$ of the distinct values $x$ among $h_1,\dots,h_n$ are all
even, and $d_h=0$ otherwise; and let $\chi_4$ be a character of order $4$. Then
$\Swan_\infty(\cG_h)=2^{n-1}-e$, and $\cG_h(\infty)$ is the direct sum of the
characters
$\cL_\psi(-S_\varepsilon^2t/\lambda)(\infty)\otimes\cL_{\chi_2}(t)(\infty)^{\otimes(n+1)}$,
$\varepsilon\in E\smallsetminus E_0$, and of the tame representation
\[
\begin{cases}
\mathbf 1^{\,r(n)-2^{n-1}-e}\oplus\bigl(\cL_{\chi_4}(t)(\infty)\oplus
\cL_{\chi_4^3}(t)(\infty)\bigr)^{\oplus e},&n\text{ odd},\\[2pt]
J_2^{\oplus(e-d_h)}\oplus\mathbf 1^{\,r(n)-2^{n-1}-2e+2d_h}\oplus
\cL_{\chi_2}(t)(\infty)^{\oplus e},&n\text{ even},
\end{cases}
\]
where $J_2$ is a unipotent Jordan block of size $2$. Moreover
$(\rho_h^*\cH_n)(\infty)$ is the direct sum of the same characters, of the
same summands $\cL_{\chi_4}(t)(\infty)$, $\cL_{\chi_4^3}(t)(\infty)$, resp.\
$\cL_{\chi_2}(t)(\infty)$, and of a unipotent representation, which is trivial
if $n$ is odd.
\end{enumerate}
\end{proposition}

\begin{proof}
Let $\cK_h=\bigotimes_{i=1}^n[\times h_i]^*\KL$ on $\mathbb G_{m,s}$, let
$j\colon\mathbb G_{m,s}\hookrightarrow\mathbb A^1_s$, and take the Fourier
transform on $\mathbb A^1_s\times\mathbb A^1_\tau$.

\emph{A Fourier model.} For $\tau\ne0$ the sheaf $\cK_h\otimes\cL_\psi(\tau s)$
on $\mathbb G_{m,s}$ is tame at $0$ and isoclinic of break $1$ at $\infty$, by
\ref{K2}, \ref{K3} and \cref{lem:B-breaks}; so its cohomology with compact
supports is concentrated in degree $1$ and has dimension $2^n$, by
\eqref{eq:B-GOS}. By proper base change,
$\FT_\psi(j_!\cK_h[1])|_{\mathbb G_{m,\tau}}=\cG[1]$ for a sheaf $\cG$ with
stalks $H^1_c(\mathbb G_{m,s},\cK_h\otimes\cL_\psi(\tau s))$. The pull-back of
the sheaf $\cM$ of \ref{H7} by $\mathrm{id}\times\rho_h$ is
$\cL_\psi(\lambda s)\otimes\cK_h(ts)$ on $\mathbb G_{m,s}\times\mathbb G_{m,t}$,
and the automorphism $(s,t)\mapsto(s/t,t)$, which commutes with the projection
to the $t$-line, transforms it into $\cL_\psi(\lambda s/t)\otimes\cK_h(s)$. By
\cref{thm:B-family}(i) and proper base change,
\begin{equation}\label{eq:R-raymodel}
\rho_h^*\cH_n\cong\iota^*\cG,\qquad\iota(t)=\lambda/t ;
\end{equation}
in particular $\cG$ is lisse of rank $2^n$, and the local monodromy of
$\rho_h^*\cH_n$ at $t=0$ (resp.\ $t=\infty$) is that of $\cG$ at $\tau=\infty$
(resp.\ $\tau=0$).

\emph{The part of lower weight.} We show that $\rho_h^*\cH_n^{\le n}$ is
geometrically constant. Let $m\colon\mathbb A^1_s\times\mathbb G_{m,t}\to
\mathbb A^1$ be $(s,t)\mapsto ts$, a smooth morphism with
$m^{-1}(\Gm)=\mathbb G_{m,s}\times\mathbb G_{m,t}$; let $\bar j$ be the
inclusion of $\mathbb G_{m,s}\times\mathbb G_{m,t}$ in
$\mathbb A^1_s\times\mathbb G_{m,t}$, and put
$\cM_\rho=\cL_\psi(\lambda s)\otimes m^*\cK_h$. By smooth base change,
$\bar j_*m^*\cK_h\cong m^*j_*\cK_h$; as $\cL_\psi(\lambda s)$ extends to a lisse
sheaf on $\mathbb A^1_s\times\mathbb G_{m,t}$, the restriction of
$\bar j_*\cM_\rho$ to $\{0\}\times\mathbb G_{m,t}$ is the geometrically constant
sheaf with fibre $(j_*\cK_h)_0$. The cokernel of the injection
$\bar j_!\cM_\rho\to\bar j_*\cM_\rho$ is supported on $\{0\}\times
\mathbb G_{m,t}$, so the cohomology sequence of the projection
$\mathbb A^1_s\times\mathbb G_{m,t}\to\mathbb G_{m,t}$ gives a morphism of
lisse sheaves
$\delta\colon(\bar j_*\cM_\rho)|_{\{0\}\times\mathbb G_{m,t}}\to\rho_h^*\cH_n$.
At a geometric point above a closed point $t$, $\delta$ is the map
$\cM_{th}^{I_0}\to H^1_c(\Gm,\cM_{th})$ of \cref{prop:B-Mh}(iii), which is
injective with image the stalk of $\cH_n^{\le n}$ (\cref{thm:B-family}(v)). A
morphism of lisse sheaves on a connected scheme which is injective at one
geometric point is injective with lisse cokernel, and two lisse subsheaves
with the same stalk at one geometric point coincide. Hence
$\rho_h^*\cH_n^{\le n}=\imop\delta$ is geometrically constant.

\emph{(i).} By \eqref{eq:B-stationary} with $\Sigma=\{0\}$ and
$F_\Sigma=\cK_h$, and since all breaks of $\cK_h(\infty)$ are less than $1$,
$\cG(\infty)\cong\mathcal F^{(0,\infty')}(\cK_h(0))$. By \ref{K2} and
\cref{lem:B-CG}, $\cK_h(0)\cong J_2^{\otimes n}\cong\bigoplus_km_{n,k}
J_{n-2k+1}$. We claim that $\mathcal F^{(0,\infty')}(J_e)\cong J_e$ for every
$e\ge1$. For $e=1$, apply \eqref{eq:B-stationary} to $F_\Sigma=\Ql$ on $\Gm$:
$\FT_\psi(j_!\Ql[1])(\infty)\cong\mathcal F^{(0,\infty')}(\mathbf 1)$. The exact
sequence $0\to j_!\Ql\to\Ql\to i_{0*}\Ql\to0$ on $\mathbb A^1$, the isomorphism
$\FT_\psi(\Ql[1])\cong i_{0*}\Ql(-1)$ (\cref{lem:B-int}) and the isomorphism
$\FT_\psi(i_{0*}\Ql)\cong\Ql[1]$ show that $\FT_\psi(j_!\Ql[1])$ is
geometrically constant on $\Gm$; so
$\mathcal F^{(0,\infty')}(\mathbf 1)\cong\mathbf 1$. By
\cite[Thm.~(2.4.3)(i)]{Lau87}, $\mathcal F^{(0,\infty')}$ is exact and induces
an equivalence of the category of continuous representations of the inertia
group at $0$ with a full subcategory of the representations of the inertia
group at $\infty'$. The representation $J_e$ is a successive extension of
trivial representations and its endomorphism ring is local. Hence
$\mathcal F^{(0,\infty')}(J_e)$ is a successive extension of trivial
representations, so that wild inertia, which acts through a finite $p$-group,
acts trivially and tame inertia acts unipotently; and its endomorphism ring is
local, so it is indecomposable, that is, a single Jordan block, of size $e$.
Therefore $\rho_h^*\cH_n$ is tame and unipotent at $t=0$, with the Jordan type
$\bigoplus_km_{n,k}J_{n-2k+1}$. Let $N$ be the logarithm of a topological
generator of the tame inertia. Its kernel has dimension $\sum_km_{n,k}$, the
rank of $\cH_n^{\le n}$ (\cref{thm:B-family}(iii)). As $\rho_h^*\cH_n^{\le n}$
is geometrically constant, its generic stalk is contained in $\ker N$, hence
equal to it. So $\cG_h(0)$, which is the quotient by $\ker N$, has $m_{n,k}$
Jordan blocks of size $n-2k$.

\emph{(ii).} Let $[2]\colon\mathbb G_{m,y}\to\mathbb G_{m,s}$ be
$y\mapsto y^2$. By a theorem of Fu and Wan \cite[Thm.~1.1]{FW05},
$\KL(\infty)\cong[2]_*\bigl(\cL_\psi(2y)\otimes\cL_{\chi_2}(y)\bigr)(\infty)$
as representations of the inertia group; since $[\times a]\circ[2]=[2]\circ
[\times\sqrt a]$ and $\cL_{\chi_2}(\sqrt ay)$ is geometrically isomorphic to
$\cL_{\chi_2}(y)$, we get $[\times h_i]^*\KL(\infty)\cong[2]_*\bigl(\cL_\psi(2
\sqrt{h_i}y)\otimes\cL_{\chi_2}(y)\bigr)(\infty)$. For representations $A$,
$B$ of the inertia group of the double cover one has
$[2]_*A\otimes[2]_*B\cong[2]_*\bigl(A\otimes(B\oplus\varsigma^*B)\bigr)$, with
$\varsigma(y)=-y$, and $\varsigma^*\bigl(\cL_\psi(2\sqrt{h_i}y)\otimes
\cL_{\chi_2}(y)\bigr)\cong\cL_\psi(-2\sqrt{h_i}y)\otimes\cL_{\chi_2}(y)$.
Taking $E=\{\varepsilon:\varepsilon_1=1\}$, induction on $n$ gives
\[
\cK_h(\infty)\cong\bigoplus_{\varepsilon\in E}\cN_\varepsilon(\infty),\qquad
\cN_\varepsilon=[2]_*\bigl(\cL_\psi(2S_\varepsilon y)\otimes\cL_{\chi_2}(y)^{\otimes n}\bigr).
\]

For a lisse sheaf $\cN$ on $\mathbb G_{m,s}$, the representation
$R^{-1}\Phi_0\bigl(\FT_\psi(j_!\cN[1])\bigr)$ of the inertia group at $\tau=0$
is $\mathcal F^{(\infty,0')}(\cN(\infty))$ \cite[Prop.~(2.3.2.1)(iii) and
D\'ef.~(2.4.2.3)]{Lau87} (compare \eqref{eq:B-localization}); it depends only
on $\cN(\infty)$, and $\mathcal F^{(\infty,0')}$ is additive
\cite[Thm.~(2.4.3)(ii)]{Lau87}. Hence
$R^{-1}\Phi_0(\FT_\psi(j_!\cK_h[1]))\cong\bigoplus_{\varepsilon\in E}
R^{-1}\Phi_0(\FT_\psi(j_!\cN'_\varepsilon[1]))$, where $\cN'_\varepsilon$ is
any lisse sheaf on $\mathbb G_{m,s}$ with $\cN'_\varepsilon(\infty)\cong
\cN_\varepsilon(\infty)$. Fix $\varepsilon$ and put $S=S_\varepsilon$. If $n$ is
even, let $\cN'_\varepsilon=\cN_\varepsilon=[2]_*\cL_\psi(2Sy)$; if $n$ is odd,
let $\cN'_\varepsilon=[\times S^2]^*\KL$, which has the required local
monodromy at $\infty$ by the theorem of Fu and Wan. In both cases
$\cN'_\varepsilon$ is tame at $0$ and totally wild at $\infty$ with
$\Swan_\infty=1$, so $H^*_c(\mathbb G_{m,s},\cN'_\varepsilon)$ is a line
placed in degree $1$, and as above
$\FT_\psi(j_!\cN'_\varepsilon[1])|_{\mathbb G_{m,\tau}}=\cG_\varepsilon[1]$ with
$\cG_\varepsilon$ lisse of rank $2$. Its trace at $\tau\in\Fq^\times$ is, by
orthogonality,
\begin{align*}
\sum_{s\in\Fq^\times}\psi_q(\tau s)\Kl_q(S^2s)&=q\psi_q(-S^2/\tau)+1
&&(n\text{ odd}),\\
-\sum_{y\in\Fq^\times}\psi_q(\tau y^2+2Sy)&=1-\psi_q(-S^2/\tau)\,\chi_2(\tau)g_q
&&(n\text{ even}),
\end{align*}
where $g_q=\sum_{y\in\Fq}\psi_q(y^2)$, and the same holds over every finite
extension. By the Hasse--Davenport relation, $-g_q$ is the Frobenius of a
geometrically constant sheaf of rank one, and the Chebotarev density theorem
shows that the semisimplification of $\cG_\varepsilon$ is the sum of $\Ql$ and
of $\cL_\psi(-S^2/\tau)$, twisted by $\cL_{\chi_2}(\tau)$ when $n$ is even and
by a geometrically constant sheaf. Since $\cN'_\varepsilon$ has no geometric
coinvariants, the stalk at $\tau=0$ of $\FT_\psi(j_!\cN'_\varepsilon[1])$ is
the line $H^1_c(\mathbb G_{m,s},\cN'_\varepsilon)$ placed in degree $-1$, with
trivial action of inertia; as in the proof of \cref{lem:B-vanishing}, it
injects into the generic stalk, and the quotient is $R^{-1}\Phi_0$. As wild
inertia acts semisimply and $S\ne0$, it follows that
\[
R^{-1}\Phi_0\bigl(\FT_\psi(j_!\cN'_\varepsilon[1])\bigr)\cong
\cL_\psi(-S^2/\tau)(0)\otimes\cL_{\chi_2}(\tau)(0)^{\otimes(n+1)} .
\]

Consequently $R^{-1}\Phi_0(\FT_\psi(j_!\cK_h[1]))$ is totally wild. In the
exact sequence
$0\to\cG_{\bar0}\to\cG_{\bar\eta_0}\to R^{-1}\Phi_0\to
\mathcal H^0(\FT_\psi(j_!\cK_h[1]))_0\to0$ obtained as in the proof of
\cref{lem:B-vanishing}, where $\cG_{\bar0}=H^1_c(\mathbb G_{m,s},\cK_h)$ and
$\mathcal H^0(\FT_\psi(j_!\cK_h[1]))_0=H^2_c(\mathbb G_{m,s},\cK_h)$, inertia acts
trivially on the two outer terms; so the last term vanishes, and
$\dim\cG_{\bar0}=\Swan_\infty(\cK_h)=2^{n-1}$ by \eqref{eq:B-GOS}. The
decomposition of $\cG_{\bar\eta_0}$ according to the characters of the wild
inertia group is stable under inertia, and $\cG_{\bar0}$ is its isotypic
component for the trivial character. Hence
\[
\cG(0)\cong\mathbf 1^{\,2^{n-1}}\oplus\bigoplus_{\varepsilon\in E}
\cL_\psi(-S_\varepsilon^2/\tau)(0)\otimes\cL_{\chi_2}(\tau)(0)^{\otimes(n+1)} .
\]
By \eqref{eq:R-raymodel}, $\cL_\psi(-S^2/\tau)$ becomes
$\cL_\psi(-S^2t/\lambda)$ and $\cL_{\chi_2}(\lambda/t)$ is isomorphic to
$\cL_{\chi_2}(t)$ as a representation of the inertia group; this gives the
assertion on $\rho_h^*\cH_n$. Finally $\rho_h^*\cH_n^{\le n}$ is geometrically
constant of rank $\binom n{\lfloor n/2\rfloor}$, so its generic stalk lies in
the trivial summand, and the quotient $\cG_h(\infty)$ has the stated form, as
$2^{n-1}-\binom n{\lfloor n/2\rfloor}=r(n)-2^{n-1}$.

\emph{(iii).} After a finite extension of $\Fq$ we may assume that
$q\equiv1\pmod4$, so that $\chi_4$ is a character of $\Fq^\times$. The
decomposition of $\cK_h(\infty)$ into the $\cN_\varepsilon(\infty)$ holds for
every $h$. Let $\varepsilon\in E_0$. Then
$\cN_\varepsilon=[2]_*(\cL_{\chi_2}(y)^{\otimes n})$ is tame at $\infty$. The
restriction of $\cL_\rho(s)$ to the inertia group of the double cover is
$\cL_{\rho^2}(y)$, and the representation induced from a character of a
subgroup of index $2$ which extends to the whole group is the sum of the two
extensions; hence $\cN_\varepsilon(\infty)\cong\mathbf 1\oplus\cL_{\chi_2}(s)(\infty)$
if $n$ is even, and $\cN_\varepsilon(\infty)\cong\cL_{\chi_4}(s)(\infty)\oplus
\cL_{\chi_4^3}(s)(\infty)$ if $n$ is odd.

Let $\rho$ be a character of $\Fq^\times$ with $\rho^4=\mathbf 1$. If
$\rho\ne\mathbf 1$, then $H^*_c(\mathbb G_{m,\Fqb},\cL_\rho)=0$ by
\eqref{eq:B-GOS}, so the stalk at $\tau=0$ of $\FT_\psi(j_!\cL_\rho[1])$
vanishes; for $\tau\in\Fq^\times$ its stalk is the line
$H^1_c(\mathbb G_{m,\Fqb},\cL_\rho\otimes\cL_\psi(\tau s))$, with trace
$-\sum_s\rho(s)\psi_q(\tau s)=-\rho^{-1}(\tau)g(\rho)$, where
$g(\rho)=\sum_s\rho(s)\psi_q(s)$, and similarly over every finite extension. By
the Hasse--Davenport relation and the Chebotarev density theorem,
$\FT_\psi(j_!\cL_\rho[1])|_{\mathbb G_{m,\tau}}$ is geometrically isomorphic to
$\cL_{\rho^{-1}}(\tau)[1]$, and so
$R^{-1}\Phi_0(\FT_\psi(j_!\cL_\rho[1]))\cong\cL_{\rho^{-1}}(\tau)(0)$. If
$\rho=\mathbf 1$, then $\FT_\psi(j_!\Ql[1])$ is geometrically constant on
$\mathbb G_{m,\tau}$ (proof of (i)), the stalk at $0$ of its $\mathcal H^{-1}$
is the line $H^1_c(\mathbb G_{m,\Fqb},\Ql)$ and that of its $\mathcal H^0$ is
the line $H^2_c(\mathbb G_{m,\Fqb},\Ql)$; so the exact sequence of the proof
of \cref{lem:B-vanishing} gives
$R^{-1}\Phi_0(\FT_\psi(j_!\Ql[1]))\cong\mathbf 1$. As
$R^{-1}\Phi_0(\FT_\psi(j_!\cN[1]))=\mathcal F^{(\infty,0')}(\cN(\infty))$
depends only on $\cN(\infty)$ and is additive, it follows that
\[
R^{-1}\Phi_0(\FT_\psi(j_!\cK_h[1]))\cong\bigoplus_{\varepsilon\in
E\smallsetminus E_0}\cL_\psi(-S_\varepsilon^2/\tau)(0)\otimes
\cL_{\chi_2}(\tau)(0)^{\otimes(n+1)}\oplus T_0,
\]
where $T_0=(\cL_{\chi_4}(\tau)(0)\oplus\cL_{\chi_4^3}(\tau)(0))^{\oplus e}$ if $n$
is odd and $T_0=(\mathbf 1\oplus\cL_{\chi_2}(\tau)(0))^{\oplus e}$ if $n$ is even.

Let $\mathcal P=\FT_\psi(j_!\cK_h[1])$, a perverse sheaf on $\mathbb A^1_\tau$
whose restriction to $\mathbb G_{m,\tau}$ is $\cG[1]$; let
$\jmath\colon\mathbb G_{m,\tau}\hookrightarrow\mathbb A^1_\tau$ and let $\imath$
be the inclusion of $0$. We claim that $\cG_{\bar0}=\cG_{\bar\eta_0}^{I_0}$. The
Fourier transform is an equivalence of the category of perverse sheaves on
$\mathbb A^1_{\Fqb}$ with itself \cite[\S1.2]{Lau87}, which by
\cref{lem:B-int} exchanges the perverse sheaves supported at $0$ with the
geometrically constant ones. A nonzero subobject of $\mathcal P$ supported at
$0$ would therefore give a nonzero morphism $\Ql[1]\to j_!\cK_h[1]$, that is, a
nonzero element of $H^0(\mathbb A^1_{\Fqb},j_!\cK_h)$; but a section of
$j_!\cK_h$ vanishes near $0$, hence vanishes, as $\cK_h$ is lisse on the
connected curve $\mathbb G_{m,\Fqb}$. So $\mathcal H^0(\imath^!\mathcal P)=0$,
and the triangle $\imath_*\imath^!\mathcal P\to\mathcal P\to
R\jmath_*\jmath^*\mathcal P$, in which $R\jmath_*\jmath^*\mathcal P$ is perverse
\cite[Cor.~4.1.3]{BBD}, gives an exact sequence of perverse sheaves
$0\to\mathcal P\to R\jmath_*\jmath^*\mathcal P\to\imath_*V\to0$, with
$V=\mathcal H^1(\imath^!\mathcal P)$ placed in degree $0$. Taking $\mathcal H^{-1}$
at $0$ gives $\cG_{\bar0}\cong(\jmath_*\cG)_{\bar0}=\cG_{\bar\eta_0}^{I_0}$.

As in (ii), the decomposition of $\cG_{\bar\eta_0}$ according to the characters
of the wild inertia group is stable under inertia; its part for the nontrivial
characters maps isomorphically onto the corresponding part of $R^{-1}\Phi_0$,
which is the sum of the characters
$\cL_\psi(-S_\varepsilon^2/\tau)(0)\otimes\cL_{\chi_2}(\tau)(0)^{\otimes(n+1)}$,
$\varepsilon\notin E_0$, and its tame part $T$ sits in an exact sequence
$0\to\cG_{\bar0}\to T\to T_0\to H^2_c(\mathbb G_{m,\Fqb},\cK_h)\to0$, on whose
outer terms inertia acts trivially. Let $\gamma$ be a topological generator of
the tame inertia group, and decompose $T$ and $T_0$ into generalized
eigenspaces of $\gamma$. For every eigenvalue $\ne1$ the corresponding parts of
$T$ and $T_0$ are isomorphic, and the unipotent part $U$ of $T$ is an extension
of the kernel of $(T_0)_1\to H^2_c(\mathbb G_{m,\Fqb},\cK_h)$ by
$\cG_{\bar0}=U^{\gamma}$, where $(T_0)_1$ is the unipotent part of $T_0$.

If $n$ is odd, then $(T_0)_1=0$, so $H^2_c(\mathbb G_{m,\Fqb},\cK_h)=0$ and
$U=\cG_{\bar0}$ is trivial, of dimension
$\dim H^1_c(\mathbb G_{m,\Fqb},\cK_h)=\Swan_\infty(\cK_h)=2^{n-1}-e$ by
\eqref{eq:B-GOS}. By \eqref{eq:R-raymodel}, and since $\cL_\rho(\lambda/t)\cong
\cL_{\rho^{-1}}(t)$ as representations of the inertia group at $\infty$, this
gives the assertion on $(\rho_h^*\cH_n)(\infty)$; and, as in (ii),
$\cG_h(\infty)$ is its quotient by the generic stalk of $\rho_h^*\cH_n^{\le n}$,
a trivial representation of dimension $\binom n{\lfloor n/2\rfloor}$. This
proves (iii) for odd $n$.

Let $n$ be even. The parts of $T$ and $T_0$ for the eigenvalue $-1$ give the
summand $\cL_{\chi_2}(\tau)(0)^{\oplus e}$. Over $\Fqb$ the sheaf $\cK_h$ is
semisimple \cite[3.4.1]{Del80}, and by \cref{lem:R-KL}(b) its geometric
monodromy group is $\SL_2^k$, acting on
$\bigotimes_x\mathrm{St}_x^{\otimes\mu_x}$, where $x$ runs over the distinct
values among $h_1,\dots,h_n$ and $\mathrm{St}_x$ is the standard representation
of the corresponding factor. By \cref{lem:B-CG}, $\mathrm{St}^{\otimes\mu}$
contains the trivial representation $m_{\mu,\mu/2}=C_{\mu/2}$ times if $\mu$ is
even and not at all if $\mu$ is odd; so the trivial representation occurs in
$\cK_h$ with multiplicity $d_h$, and $\cK_h\cong\cK'\oplus\Ql^{\oplus d_h}$ over
$\Fqb$, where $\cK'$ has no trivial constituent. As
$j_*\cK_h/j_!\cK_h=\imath_{0*}(\cK_h(0)^{I_0})$, there is an exact sequence of
perverse sheaves
$0\to\imath_{0*}(\cK_h(0)^{I_0})\to j_!\cK_h[1]\to j_*\cK'[1]\oplus
\Ql[1]^{\oplus d_h}\to0$ on $\mathbb A^1_{s,\Fqb}$, and its Fourier transform,
by \cref{lem:B-int}, is an exact sequence
\[
0\to\cK_h(0)^{I_0}\otimes\Ql[1]\to\mathcal P\to\mathcal P'\oplus
\imath_*\Ql(-1)^{\oplus d_h}\to0,\qquad\mathcal P'=\FT_\psi(j_*\cK'[1]).
\]
The simple constituents of $j_*\cK'[1]$ are the $j_*\cK_i[1]$ for the
irreducible constituents $\cK_i$ of $\cK'$, which are neither supported at a
point nor geometrically constant. As the Fourier transform exchanges the
perverse sheaves supported at $0$ with the geometrically constant ones, the
simple constituents of the semisimple perverse sheaf $\mathcal P'$ are neither
supported at $0$ nor geometrically constant. So $\mathcal P'$ has no nonzero
subobject or quotient supported at $0$; in particular
$\mathcal H^0(\mathcal P')_0=0$, and, as above, $\mathcal H^{-1}(\mathcal
P')_{\bar0}=\cG'^{I_0}_{\bar\eta_0}$, where $\cG'[1]=\mathcal P'|_{\mathbb
G_{m,\tau}}$. Moreover $\cG'$ has no nonzero geometrically constant subsheaf,
since the restriction to $\mathbb G_{m,\tau}$ of a simple constituent of
$\mathcal P'$ is an irreducible lisse sheaf, placed in degree $-1$, which is not
geometrically constant.
On $\mathbb G_{m,\tau}$ the exact sequence above becomes
$0\to\cK_h(0)^{I_0}\otimes\Ql\to\cG\to\cG'\to0$, whose first term is
geometrically constant of rank $\binom n{n/2}$, the number of Jordan blocks of
$\cK_h(0)\cong J_2^{\otimes n}$ (\cref{lem:B-CG}). Under
\eqref{eq:R-raymodel}, the geometrically constant subsheaf
$\rho_h^*\cH_n^{\le n}$, of the same rank, maps to zero in $\iota^*\cG'$, so it
is $\iota^*$ of the first term, and $\cG_h\cong\iota^*\cG'$.

Vanishing cycles transform the exact sequence into a long exact sequence,
geometrically constant sheaves have no vanishing cycles, and
$R^{-1}\Phi_0(\imath_*\Ql)\cong\mathbf 1$; so $R^{-1}\Phi_0(\mathcal P)\cong
R^{-1}\Phi_0(\mathcal P')\oplus\mathbf 1^{\oplus d_h}$, and the unipotent part
of $R^{-1}\Phi_0(\mathcal P')$ is trivial of dimension $e-d_h$. Since
$\mathcal H^0(\mathcal P')_0=0$, the argument given above for $\mathcal P$,
applied to $\mathcal P'$, shows that the unipotent part $U'$ of $\cG'(0)$ is an
extension of $\mathbf 1^{\oplus(e-d_h)}$ by $U'^{\gamma}$. Hence
$(\gamma-1)U'\subseteq U'^{\gamma}=\ker(\gamma-1)$, and $\gamma-1$ induces an
injection of $U'/U'^{\gamma}$ into $U'^{\gamma}$; that is,
$U'\cong J_2^{\oplus(e-d_h)}\oplus\mathbf 1^{\oplus(\dim U'-2e+2d_h)}$. As the
other parts of $\cG'(0)$ are those of $\cG(0)$, and
$\dim U'=r(n)-e-(2^{n-1}-e)=r(n)-2^{n-1}$, the isomorphism
$\cG_h\cong\iota^*\cG'$ gives (iii) for even $n$; the assertion on
$(\rho_h^*\cH_n)(\infty)$ follows from \eqref{eq:R-raymodel}.
\end{proof}

For $n=2$, \cref{prop:R-ray} agrees with \cref{prop:n2sheaf}: on a ray,
$\cW_2$ is geometrically
$[t\mapsto h_1h_2t^2/\lambda^2]^*\KL\otimes\cL_\psi(-(h_1+h_2)t/\lambda)$,
which is unipotent with one block of size $2$ at $0$ and, at $\infty$, the sum
of the two characters $\cL_\psi(-(\sqrt{h_1}\pm\sqrt{h_2})^2t/\lambda)\otimes
\cL_{\chi_2}(t)$; if $h_1=h_2$, one of them is $\cL_{\chi_2}(t)$, in accordance
with (iii), where $e=d_h=1$.

The $L$-function of a ray is given by a count of points. Let
$e_k$ denote the $k$-th elementary symmetric polynomial.

\begin{proposition}\label{prop:R-rayL}
Let $h\in S_n(\Fq)$, and let $Z_h^\circ\subset\PP^{n-1}$ be the subvariety
defined by $w_1\cdots w_n\ne0$, $\sum_ih_iw_i=0$ and $e_{n-1}(w)=0$.
\begin{enumerate}[label=\textup{(\alph*)}]
\item One has
\[
\sum_{t\in\Fq^\times}\Tr(\Frob_t\mid\cG_h)=(-1)^{n-1}\bigl((q-1)^{n-1}+(-1)^n(q+1)
-q^2\#Z_h^\circ(\Fq)\bigr)-(q-1)E_n(q).
\]
\item Let $n=3$ and $\Delta(h)=h_1^2+h_2^2+h_3^2-2h_1h_2-2h_1h_3-2h_2h_3$. Then
$H^i_c(\mathbb G_{m,\Fqb},\cG_h)=0$ for $i\ne1$ and
\[
L(\Gm,\cG_h,T)=(1-qT)(1-q^2T)^2\bigl(1-\chi_2(\Delta(h))q^2T\bigr),
\]
where $\chi_2(0)=0$.
\item Let $n=4$, and suppose that $S_\varepsilon\ne0$ for every $\varepsilon$.
Then the plane cubic $\cC_h=\{w\in\PP^3:e_3(w)=0,\ \sum_ih_iw_i=0\}$ is a smooth
curve of genus one, $H^i_c(\mathbb G_{m,\Fqb},\cG_h)=0$ for $i\ne1$, and
\[
L(\Gm,\cG_h,T)=(1-qT)(1-q^2T)^3(1-q^3T)^2\det\bigl(1-q^2T\Frob_q\mid
H^1(\cC_{h,\Fqb},\Ql)\bigr).
\]
Moreover $\det\bigl(1-T\Frob_q\mid H^1(\PP^1_{\Fqb},j_*\cG_h)\bigr)=(1-q^3T)^2$;
the Frobenius acts on the inertia invariants at $0$ with eigenvalues
$q,q^2,q^2,q^2$; and on the inertia invariants at $\infty$, which form the
trivial summand $\mathbf 1^2$ of \cref{prop:R-ray}(ii), with characteristic
polynomial $\det\bigl(1-q^2T\Frob_q\mid H^1(\cC_{h,\Fqb},\Ql)\bigr)$.
\item Let $n=4$, and suppose that $S_\varepsilon=0$ for exactly one
$\varepsilon\in E$. Then $\cC_h$ is a plane cubic with exactly one singular
point $x$, which is an ordinary double point with all coordinates nonzero;
$H^i_c(\mathbb G_{m,\Fqb},\cG_h)=0$ for $i\ne1$, and
\[
L(\Gm,\cG_h,T)=(1-qT)(1-q^2T)^3(1-q^3T)^2(1-\eta q^2T),
\]
where $\eta=1$ if the two tangent lines to $\cC_h$ at $x$ are defined over
$\Fq$ and $\eta=-1$ otherwise; thus $H^1(\cC_{h,\Fqb},\Ql)$ is a line on which
$\Frob_q$ acts by $\eta$. The Frobenius acts on the inertia invariants at $0$
with eigenvalues $q,q^2,q^2,q^2$, and on those at $\infty$, which form the line
of invariants of the Jordan block $J_2$ of \cref{prop:R-ray}(iii), by $\eta q^2$.
\item Let $n\ge2$ be arbitrary and suppose that $S_\varepsilon\ne0$ for every
$\varepsilon$. Then $Z_h^\circ$ is smooth of pure dimension $n-3$ (empty for
$n=2$), $(j_*\cG_h)_{\bar\infty}$ is pure of weight $n+1$ and of dimension
$r(n)-2^{n-1}$, and, writing $\Delta(V)=\det(1-T\Frob_q\mid V)$,
\[
\Delta\bigl((j_*\cG_h)_{\bar\infty}\bigr)\,\Delta\bigl(\gr^W_{n+1}(j_*\cG_h)_{\bar0}\bigr)
=\Delta\bigl(H^0(\mathbb G_{m,\Fqb},\cG_h)\bigr)\bigl(1-q^{(n+1)/2}T\bigr)^{-\tau_n}\,\Pi_h(T),
\]
where
$\Pi_h(T)=\prod_i\det\bigl(1-q^2T\Frob_q\mid\gr^W_{n-3}H^i_c(Z^\circ_{h,\Fqb},\Ql)\bigr)^{(-1)^{n+i+1}}$
and, for odd $n$, $\tau_n$ is the coefficient of $x^{(n+1)/2}$ in
$(-1)^{n-1}\bigl((x-1)^{n-1}+(-1)^n(x+1)\bigr)-(x-1)E_n(x)$, and $\tau_n=0$ for even
$n$. For even $n$ the second factor on the left is $1$.
\end{enumerate}
\end{proposition}

\begin{proof}
(a) By \cref{lem:identities}(b), writing $\Kl_q(sh_i)=\sum_{x_i}\psi_q(x_i+sh_ix_i^{-1})$
and substituting $w_i=s/x_i$,
\[
q\Sigma_{n,q}(h)=(-1)^n+\sum_{w\in(\Fq^\times)^n}\psi_q(\langle h,w\rangle)
\sum_{s\in\Fq^\times}\psi_q\Bigl(s\bigl(\lambda+\textstyle\sum_iw_i^{-1}\bigr)\Bigr),
\]
so that $\Sigma_{n,q}(h)=\sum_{w\in X_\lambda(\Fq)}\psi_q(\langle h,w\rangle)$,
where $X_\lambda=\{w\in\Gm^n:\sum_iw_i^{-1}=-\lambda\}$. Replacing $h$ by $th$
and summing over $t\in\Fq^\times$ gives
$\sum_t\Sigma_{n,q}(th)=qN-\#X_\lambda(\Fq)$, where $N$ is the number of
$w\in X_\lambda(\Fq)$ with $\langle h,w\rangle=0$. By orthogonality,
$\#X_\lambda(\Fq)=((q-1)^n-(-1)^n)/q$, and the hyperplane $\sum_ih_iw_i=0$ of
$\PP^{n-1}$ has $((q-1)^{n-1}+(-1)^n)/q$ rational points with all coordinates
nonzero. A point $[w]$ with $\prod w_i\ne0$ and $\sum_iw_i^{-1}\ne0$ has
exactly one representative in $X_\lambda$, and $\sum_iw_i^{-1}=0$ if and only if
$e_{n-1}(w)=0$; hence $N=((q-1)^{n-1}+(-1)^n)/q-\#Z_h^\circ(\Fq)$. The
formula follows from \ref{H2}.

(b) For $n=3$, eliminating $w_3$ shows that $Z_h^\circ$ is the set of zeros in
$\PP^1$ of $h_1w_1^2+(h_1+h_2-h_3)w_1w_2+h_2w_2^2$, none of which has a vanishing
coordinate; its discriminant is $\Delta(h)$. Hence
$\#Z_h^\circ(\mathbb F_{q^m})=1+\chi_2(\Delta(h))^m$, also if $\Delta(h)=0$, when
there is a single double zero, and (a), applied over every finite extension,
gives the stated product for $L(\Gm,\cG_h,T)$. Since $\Gm$
is affine, $H^0_c=0$ and $L=\det(1-T\Frob_q\mid H^1_c)/\det(1-T\Frob_q\mid
H^2_c)$. The group $H^2_c$ is pure of weight $w+2=6$ (proof of
\cref{lem:B-av}), whereas the product has reciprocal roots of weights $2$ and
$4$ and $H^1_c$ is mixed of weights at most $5$ \cite[3.3.1]{Del80}; so
$H^2_c=0$.

(c) The Cayley cubic surface $X=\{e_3=0\}\subset\PP^3$ is singular exactly at
the four coordinate points. Its intersection with the coordinate planes is the
union of the six lines $w_i=w_j=0$, and at the smooth points of the line
$w_1=w_2=0$ the tangent plane is $w_1+w_2=0$ (compute the gradient of $e_3$ at
$(0,0,s,t)$, $st\ne0$), and similarly for the other lines. At a point $x$ of
$X$ with $\prod_ix_i\ne0$ the tangent plane is $\sum_ix_i^{-2}w_i=0$, as
$e_3=e_4\sum_iw_i^{-1}$. A plane section of $X$ is singular exactly at the
singular points of $X$ that it contains and at the points of tangency. The
plane $H_h=\{\sum_ih_iw_i=0\}$ contains no coordinate point and is none of the
planes $w_i+w_j=0$, as all $h_i\ne0$. It is tangent to $X$ at a point $x$ with
$\prod_ix_i\ne0$ if and only if $h_i=\kappa x_i^{-2}$ for some $\kappa$, and
then $\sqrt{h_i}=\varepsilon_i\sqrt\kappa\,x_i^{-1}$ for suitable signs, so
that $S_\varepsilon=\sqrt\kappa\sum_ix_i^{-1}=0$; conversely, if
$S_\varepsilon=0$, the point $x_i=\varepsilon_i/\sqrt{h_i}$ is such a point of
tangency. So $\cC_h$ is a smooth plane cubic, of genus one. The points of
$\cC_h$ with a vanishing coordinate are the six rational points
$\{w_i=w_j=0\}\cap H_h$, so $\#Z_h^\circ(\mathbb F_{q^m})=\#\cC_h(\mathbb
F_{q^m})-6=q^m-5-a_m$, with $a_m=\Tr(\Frob_q^m\mid H^1(\cC_{h,\Fqb},\Ql))$.
By (a),
\[
\sum_{t\in\mathbb F_{q^m}^\times}\Tr(\Frob_t\mid\cG_h)
=-2q^{3m}-3q^{2m}-q^m-q^{2m}a_m,
\]
which gives the formula for $L(\Gm,\cG_h,T)$; its reciprocal roots have
weights $2$, $4$, $5$ and $6$, the eigenvalues of Frobenius on
$H^1(\cC_{h,\Fqb},\Ql)$ having weight $1$. As in (b), $H^2_c$, which is pure of
weight $7$, vanishes. As in the proof of \cref{cor:H0} there is an exact
sequence
\[
0\to H^0(\mathbb G_{m,\Fqb},\cG_h)\to(j_*\cG_h)_{\bar0}\oplus(j_*\cG_h)_{\bar\infty}
\to H^1_c(\mathbb G_{m,\Fqb},\cG_h)\to H^1(\PP^1_{\Fqb},j_*\cG_h)\to0,
\]
whose last term is pure of weight $6$ \cite[3.2.3]{Del80}, while the stalks of
$j_*\cG_h$ have weights at most $5$ \cite[1.8.4]{Del80}. Hence the last term
has characteristic polynomial $(1-q^3T)^2$, and the boundary terms account for
the other six reciprocal roots. By \cref{prop:R-ray}(i), $(j_*\cG_h)_{\bar0}$
has dimension $4$, and by the weight-monodromy theorem \cite[1.8.4]{Del80} its
Frobenius eigenvalues have weight $5-3=2$ (for the block of size $4$) and
$5-1=4$ (for the three blocks of size $2$). By \cref{prop:R-ray}(ii),
$(j_*\cG_h)_{\bar\infty}$ is the trivial summand $\mathbf 1^2$, on which inertia
acts trivially; by the same theorem it is pure of weight $5$. Counting
dimensions, $H^0(\mathbb G_{m,\Fqb},\cG_h)=0$, the eigenvalues at $0$ are $q$,
$q^2$, $q^2$, $q^2$, and those at $\infty$ are the two of weight $5$.

(d) As in (c), $\cC_h$ is singular exactly at the points where $H_h$ is tangent
to $X$, and these are the points $x$ with $x_i=\varepsilon_i/\sqrt{h_i}$ for the
$\varepsilon$ with $S_\varepsilon=0$; so $\cC_h$ has exactly one singular point
$x$, which is therefore rational over $\Fq$, and all coordinates of $x$ are
nonzero. Since $e_3=e_4\sum_iw_i^{-1}$, the automorphism $w\mapsto
u=(w_i^{-1})_i$ of the open subset $\{\prod_iw_i\ne0\}$ of $\PP^3$ maps
$X^\circ=X\cap\{\prod_iw_i\ne0\}$ onto the open subset of the plane
$\Pi=\{\sum_iu_i=0\}$ where all $u_i\ne0$. Put $a=u(x)$; then $\sum_ia_i=0$,
and $h_i=\kappa a_i^2$ for some $\kappa\ne0$ by the tangency condition. So
$\cC_h\cap X^\circ$ becomes the curve $f=0$ in $\Pi$, where
$f(u)=\sum_ia_i^2\prod_{j\ne i}u_j$. One computes $\partial f/\partial
u_k(a)=-a_1a_2a_3a_4$ for every $k$, a vector normal to $\Pi$, and, for
$\{k,l,i,i'\}=\{1,2,3,4\}$, $\partial^2f/\partial u_k\partial
u_l(a)=-a_ia_{i'}(a_k+a_l)$, while $\partial^2f/\partial u_k^2=0$. The vector
$a$ lies in the kernel of the restriction of this Hessian form to the
hyperplane $\sum_iv_i=0$ of $\Fqb^4$. After renumbering we may assume
$a_1+a_2\ne0$; then $e_1-e_2$, $e_3-e_4$ and $a$ form a basis of this
hyperplane, and the determinant of the Hessian form on the span of $e_1-e_2$
and $e_3-e_4$ is $-4a_1a_2a_3a_4(a_1+a_2)^2\ne0$. Hence the tangent cone of
$\cC_h$ at $x$ consists of two distinct lines: $x$ is an ordinary double point.

A reducible plane cubic is either a line and a conic, which meet in two points
or are tangent at a point which is not an ordinary double point, or three
lines, or not reduced; so $\cC_h$ is geometrically irreducible, its
normalization is $\PP^1$, and the two points above $x$ are rational over
$\mathbb F_{q^m}$ if and only if the tangent lines at $x$ are. Hence
$\#\cC_h(\mathbb F_{q^m})=q^m+1-\eta^m$, and $H^1(\cC_{h,\Fqb},\Ql)$ is a line on
which $\Frob_q$ acts by $\eta$. The six points of $\cC_h$ with a vanishing
coordinate are rational and smooth, as in (c), so
$\#Z_h^\circ(\mathbb F_{q^m})=q^m-5-\eta^m$, and (a) gives
\[
\sum_{t\in\mathbb F_{q^m}^\times}\Tr(\Frob_t\mid\cG_h)
=-2q^{3m}-3q^{2m}-q^m-\eta^mq^{2m},
\]
which gives the formula for $L(\Gm,\cG_h,T)$; as in (b), $H^2_c=0$.

In the exact sequence of the proof of (c), $H^1(\PP^1_{\Fqb},j_*\cG_h)$ has
characteristic polynomial $(1-q^3T)^2$, and the other five reciprocal roots
$q,q^2,q^2,q^2,\eta q^2$ are the eigenvalues of Frobenius on
$(j_*\cG_h)_{\bar0}\oplus(j_*\cG_h)_{\bar\infty}$ modulo
$H^0(\mathbb G_{m,\Fqb},\cG_h)$. If all values among $h_1,\dots,h_4$ had even
multiplicity, then, after renumbering, $h_1=h_2$ and $h_3=h_4$, and
$S_\varepsilon$ would vanish for at least two $\varepsilon\in E$. So $d_h=0$,
and by \cref{prop:R-ray}(iii), with $e=1$ and $r(4)-2^3-2=0$, the tame part of
$\cG_h(\infty)$ is $J_2\oplus\cL_{\chi_2}(t)(\infty)$, whose inertia invariants
form a line. By \cref{prop:R-ray}(i), $(j_*\cG_h)_{\bar0}$ has dimension $4$,
so $H^0(\mathbb G_{m,\Fqb},\cG_h)=0$. Finally, the representation of a
decomposition group at $t=0$ on $\cG_h$ does not depend on $h$: by the proof of
\cref{prop:R-ray}(i), it is the quotient of the transport by $\iota$ of
$\mathcal F^{(0,\infty')}(\cK_h(0))$ by its inertia invariants, the principle of
stationary phase being compatible with the action of a decomposition group
\cite[(2.3.3.1)]{Lau87}; and $\cK_h(0)\cong\KL(0)^{\otimes4}$ as
representations of a decomposition group at $0$, by \cref{lem:B-frob} and
\ref{K5}, since each $[\times h_i]$ is defined over $\Fq$ and fixes $0$. After
replacing $\Fq$ by $\mathbb F_{q^k}$ with $k$ large, there is a point $h'$ with
all $S_\varepsilon(h')\ne0$, and (c), applied to $h'$, shows that the $k$-th
powers of the eigenvalues at $0$, which have weights $2,4,4,4$, are
$q^k,q^{2k},q^{2k},q^{2k}$, for all large $k$. So the eigenvalues at $0$ are
$q,q^2,q^2,q^2$, and the eigenvalue at $\infty$ is $\eta q^2$.

(e) The open subset $X^\circ=\{e_{n-1}(w)=0,\ \prod_iw_i\ne0\}$ of the hypersurface
$e_{n-1}=0$ is smooth, being the zero set of $\sum_iw_i^{-1}$, whose gradient
$(-w_i^{-2})_i$ does not vanish; and, as in the proof of (c), the hyperplane
$\sum_ih_iw_i=0$ is tangent to $X^\circ$ at a point $x$ if and only if
$h_i=\kappa x_i^{-2}$ for some $\kappa$, which forces $S_\varepsilon=0$ for suitable
signs. So $Z_h^\circ$ is a smooth complete intersection of dimension $n-3$ in the
torus of $\PP^{n-1}$. By (a), applied over every finite extension,
\[
L(\Gm,\cG_h,T)=\prod_k(1-q^kT)^{-\tau_{n,k}}\prod_i\det\bigl(1-q^2T\Frob_q\mid
H^i_c(Z^\circ_{h,\Fqb},\Ql)\bigr)^{(-1)^{n+i+1}},
\]
where $\sum_k\tau_{n,k}x^k=(-1)^{n-1}((x-1)^{n-1}+(-1)^n(x+1))-(x-1)E_n(x)$. On the
other hand, $H^2_c(\mathbb G_{m,\Fqb},\cG_h)$ is pure of weight $n+3$, and the exact
sequence of the proof of (c) shows that the part of weight $n+1$ of
$H^1_c(\mathbb G_{m,\Fqb},\cG_h)$ is, in the Grothendieck group, the part of
weight $n+1$ of $(j_*\cG_h)_{\bar0}\oplus(j_*\cG_h)_{\bar\infty}$ minus
$H^0(\mathbb G_{m,\Fqb},\cG_h)$, since $H^1(\PP^1_{\Fqb},j_*\cG_h)$ is pure of
weight $n+2$. By \cref{prop:R-ray}(ii) the tame part of $\cG_h(\infty)$ is trivial,
so $(j_*\cG_h)_{\bar\infty}$ is pure of weight $n+1$ by the weight-monodromy
theorem, of dimension $r(n)-2^{n-1}$; by \cref{prop:R-ray}(i) and the same
theorem, $(j_*\cG_h)_{\bar0}$ has the weights $n+2-b\le n+1$, $b=n-2k$ running over
the sizes of its Jordan blocks. Comparing the parts of weight $n+1$ of the two
expressions for $L(\Gm,\cG_h,T)$ gives the formula: the factor $1-q^kT$ has weight
$2k$, which equals $n+1$ only for odd $n$ and $k=(n+1)/2$, and $q^2$ times a
Frobenius eigenvalue on $H^i_c(Z_h^\circ)$ has weight $n+1$ exactly when the
eigenvalue has weight $n-3$. For even $n$ all blocks at $0$ have even size, so
$(j_*\cG_h)_{\bar0}$ has no part of weight $n+1$.
\end{proof}

Thus, for $n=4$, the Frobenius at the boundary point at infinity of a ray is
that of the Tate-twisted elliptic curve $H^1(\cC_h)(-2)$, where $\cC_h$ is the
plane section of the Cayley cubic surface by the hyperplane dual to $h$; the
curve $\cC_h$ is singular exactly when $\prod_\varepsilon S_\varepsilon=0$,
that is, when the hyperplane is tangent to the surface, which is dual to the
quartic Steiner surface $\prod_\varepsilon S_\varepsilon=0$. When the
hyperplane is tangent at exactly one point, the elliptic curve degenerates into
a nodal cubic, and accordingly (\cref{prop:R-rayL}(d)) the trivial summand
$\mathbf 1^2$ at infinity degenerates into a Jordan block $J_2$, whose line of
invariants still carries $H^1(\cC_h)(-2)$, now of weight $4$. Similarly, for
$n=3$ one has $\Delta(h)=\prod_{\varepsilon\in E}S_\varepsilon$, and when
$\Delta(h)=0$ the factor $1-\chi_2(\Delta(h))q^2T$ disappears, in accordance
with \cref{prop:R-ray}(iii), by which $\cG_h(\infty)$ then has no inertia
invariants.

\section*{Acknowledgments}
This work was carried out independently. I am grateful to those whose support, patience,
and encouragement made it possible to continue the project over a long period. Much of
the research was done alongside other professional responsibilities, and the arguments were
revisited repeatedly as the manuscript developed. The final form reflects that sustained
process of reading, checking, and revision. I remain solely responsible for the mathematical
content and for any errors that remain.

All mathematical arguments, proofs, and conclusions in this article are the sole work of the
author. Computational tools were used only for the numerical checks explicitly described in
the text; these checks are not essential to the proofs and are provided for verification. The
author has reviewed all results and assumes full responsibility for the correctness and integrity
of the work. AI-assisted tools, if used, were employed exclusively for language editing and
formatting; they did not contribute to any mathematical content.

\end{document}